%
%
%
%
%
\magnification=\magstep1   
\input amstex
\input pictex
\UseAMSsymbols
\hoffset=0truecm 
\hsize=125mm \vsize=185mm
\vcorrection{1cm}
\NoBlackBoxes
\parindent=8mm
\mathsurround=1pt
\font\gross=cmbx10 scaled\magstep1   \font\abs=cmcsc10
\font\rmk=cmr8  \font\itk=cmti8 \font\bfk=cmbx8   \font\ttk=cmtt8

 \def \bdim{\operatorname{\bold{dim}}}
 \def \bm{\bold{m}}
 \def \bh{\bold{h}}
 \def \br{\bold{r}}
 \def \be{\bold{e}}
 \def \bx{\bold{x}}
 \def \bd{\bold{d}}
 \def \can{\operatorname{can}}
 \def \Cok{\operatorname{Cok}}
 \def \diag{\operatorname{diag}}
 \def \End{\operatorname{End}}
 \def \mod{\operatorname{mod}}
 \def \umod{\operatorname{\underline{mod}}}

 \def \Hom{\operatorname{Hom}}
 \def \I{\operatorname{I}}
 \def \Im{\operatorname{Im}}
 \def \incl{\operatorname{incl}}
 \def \ind{\operatorname{ind}}
 
 \def \Ker{\operatorname{Ker}}
 \def \lfd{\operatorname{lfd}}
 \def \lto#1{\;\mathop{\longrightarrow}\limits^{#1}\;} 
 \def \min{\operatorname{min}}
 \def \Mimo{\operatorname{Mimo}}

 \def \Mono{\operatorname{Mono}}

 \def \rad{\operatorname{rad}}

 \def \Res{\operatorname{Res}}
 \def \sto#1{\;\mathop{\to}\limits^{#1}\;}
 \def \soc{\operatorname{soc}} 
 \def \supp{\operatorname{supp}}

 \def\qed{\phantom{m.} $\!\!\!\!\!\!\!\!$\nolinebreak\hfill\checkmark}

 \def\arr#1#2{\arrow <2mm> [0.25,0.75] from #1 to #2}
 \def\sq{\plot 0 0  1 0  1 1  0 1  0 0 /}
 \def\smallsq#1{\plot 0 0  0.#1 0  0.#1 0.#1  0 0.#1  0 0 /}
\def\hobj#1#2#3{\big( #1 \sto{#2} #3\big)}
\def \tube{\beginpicture\setcoordinatesystem units <0.25cm, 0.25cm> 
        \plot -1 1  -1 -1  1 -1  1 1 / \endpicture}
\def \Utube{\beginpicture\setcoordinatesystem units <0.25cm, 0.25cm> 
        \plot -1 1  -1 -1  -0.6 -1  -0.3 -1.5  
                0 -1  0.2 -1  0.5 -1.5  0.8 -1 1 -1  1 1 / 
        \endpicture}
\def \qtube{\beginpicture\setcoordinatesystem units <0.25cm, 0.25cm> 
        \plot -1 1  -1 -1  -0.6 -1  -0.4 -1.3  
                0 -1.3  0.3 -1.3  0.5 -1.6   1 -1.6  1 1 / 
        \endpicture}
\def \ptube{\beginpicture\setcoordinatesystem units <0.25cm, 0.25cm> 
        \plot -1 1  -1 -1.6  -0.6 -1.6  -0.4 -1.3  
                0 -1.3  0.3 -1.3  0.5 -1  1 -1  1 1 / 
        \endpicture}
\def\Left#1#2#3{\beginpicture
  \setcoordinatesystem units <0.1cm,0.1cm>
  \put{$\sssize #1$} at 1.4 1.6
  \put{$\sssize #2$} at 0 0
  \put{$\sssize #3$} at -.7 -1.6
  \endpicture}
\def\all#1#2#3{\beginpicture
  \setcoordinatesystem units <0.1cm,0.1cm>
  \put{$\sssize #1$} at 1.2 1.6
  \put{$\sssize #2$} at 0 0
  \put{$\sssize #3$} at -1.2 -1.6
  \endpicture}
\def\Eeighttubular{\beginpicture
        \setcoordinatesystem units <0.5cm,0.5cm>
        \put{$\ssize E_8$} at 0 0.3
        \circulararc 300 degrees from 0.4 0.3 center at -.05 0.35
        \endpicture}
 \def\t#1{\;\text{#1}\;}
 
 \def\T#1{\qquad\text{#1}\qquad}
 \def\E#1{{\parindent=1truecm\narrower\narrower\noindent{#1}\smallskip}}  
\headline{\ifnum\pageno=1\hfill %
    \else\ifodd\pageno \hfil\Rechts\hfil \else \hfil\Links\hfil \fi  \fi}
    \def\Links{\abs Ringel, Schmidmeier}
    \def\Rechts{\abs Invariant Subspaces of Nilpotent Linear
 Operators. I.}
%
\phantom{\noindent{\rmk [rs-inv3r, 27.~August 2006]}}
%
%
        \vglue1truecm
\centerline{\gross Invariant Subspaces of Nilpotent Linear Operators. I.}
        \bigskip
\centerline{Claus Michael Ringel and Markus Schmidmeier}
        \bigskip\medskip
\E{{\it Abstract.} 
Let $k$ be a field.
We consider triples $(V,U,T)$, where $V$ is a finite dimensional $k$-space, $U$
a subspace of $V$ and $T \:V \to V$ a linear operator with $T^n = 0$ for some $n$, and
such that $T(U) \subseteq U$. Thus, $T$ is a nilpotent operator on $V$, and $U$ is 
an invariant subspace with respect to $T$. We will discuss the question whether it is
possible to classify these triples. These triples $(V,U,T)$ are the objects of a
category with the Krull-Remak-Schmidt property, thus it will be sufficient to deal
with indecomposable triples. Obviously, the classification problem depends on $n$, and
it will turn out that 
the decisive case is $n=6.$ For $n < 6$, there are only finitely many isomorphism classes of
indecomposables triples, whereas for $n > 6$ we deal with what is called  ``wild''
representation type, so no complete classification can be expected. For $n=6$,
we will exhibit a complete description of all the indecomposable triples.
}

\plainfootnote{}{%
{\rmk 2000 \itk Mathematics Subject Classification. \rmk 
Primary:
47A15, 
15A21. 
Secondary: 
15A04, 
16D70, 
16G30, 
16G60, 
18G25, 
47C99. 

{\itk Keywords:} Nilpotent operators. Invariant subspaces. Auslander-Reiten
        quivers. Tubular categories. Coverings.
}
}

\bigskip
Let $k$ be a field which will be fixed throughout the paper (of interest should be $k = \Bbb C$
or $k = \Bbb R$, but also the cases where $k$ is finite). We consider a finite-dimensional 
$k$-space $V$ with 
a linear operator $T\:V \to V$ which is nilpotent, say satisfying $T^n = 0$ for some
natural number $n$, and we are
interested in the subspaces $U$ of $V$ which are $T$-invariant
(this means that $T(U) \subseteq U$). 
To be more precise, we fix the number $n$, look for 
all such triples $X = (V,U,T)$ and ask for a classification
of the isomorphism classes:  We consider these triples $(V,U,T)$ 
as the objects of a category $\Cal S(n)$ (or $\Cal S_k(n)$ 
in case a reference to the base field
$k$ is needed); for two such triples $(V,U,T)$ and 
$(V',U',T')$, a map $(V,U,T) \to (V',U',T')$ is given by a
$k$-linear map $f\:V \to V'$ with $f(U) \subseteq U'$ and such that $fT = T'f.$  Note that a map 
$f\:(V,U,T) \to (V',U',T')$ in $\Cal S(n)$ is an isomorphism provided $f$ is bijective and
$f(U) = U'$. The category $\Cal S(n)$ is additive, the {\it direct sum}  of two objects 
$X = (V,U,T)$ and $X' = (V',U',T')$ is the triple $X\oplus X' = (V\oplus V',U\oplus U', T\oplus T')$; 
the {\it zero object} is $(0,0,0)$. 
The triple $(V,U,T)$ is said to be {\it indecomposable} provided it is not zero 
and not isomorphic to a direct sum $X\oplus X'$
with non-zero objects $X, X'$. 
The aim of this paper is a study of the categories $\Cal S(n)$, in particular we ask for
a classification of the indecomposable objects.
        
        \bigskip\noindent
{\bf (0.1) Crucial results.}
The difficulty of classifying the indecomposable objects in
$\Cal S(n)$ increases with increasing $n$. Whereas for $n\leq5$ there are only finitely
many isomorphism classes of indecomposable triples, there are infinitely many
such classes for $n\geq 6.$  Our main concern will be the case $n=6$ which
turns out to be a ``tame'' case: here, a complete classification is possible,
as we will show in Part~A of the paper. Part B will provide some 
information on the ``finite type'' cases $n <6$ as well as on the ``wild'' cases $n > 6.$
        
The description of the category $\Cal S_k(6)$ presented in this paper
may look complicated, but it allows to draw a lot
of interesting consequences. Some of these consequences
can be formulated and understood very easily, this we are going to do now ---
but at the moment, no elementary proofs seem to be available.
Let us introduce 
an obvious invariant, namely the {\it dimension pair} $(\dim V, \dim U)$
of a triple $(V,U,T)$ 
in $\Cal S(n)$; it is a pair $(v,u)$ of natural numbers with $u \le v$. 
        \medskip\noindent
{\bf (0.1.1)} {\it A pair $(v,u)$ of natural numbers with $u \le v$ 
is the dimension pair of an
indecomposable object in $\Cal S(6)$ if and only if the following two conditions are
satisfied: first, 
$$
 |v-2u| \le 6,
$$
and second, 
$(v,u)$ is different from $(0,0),\,(7,1),\, (8,1),\,(7,6),\,(8,7).$}
        \medskip
This means that the dimension pairs $(v,u)$ of the indecomposable triples $(V,U,T)$ in
$\Cal S(6)$ lie in the following shaded region (the exceptions are encircled):
$$
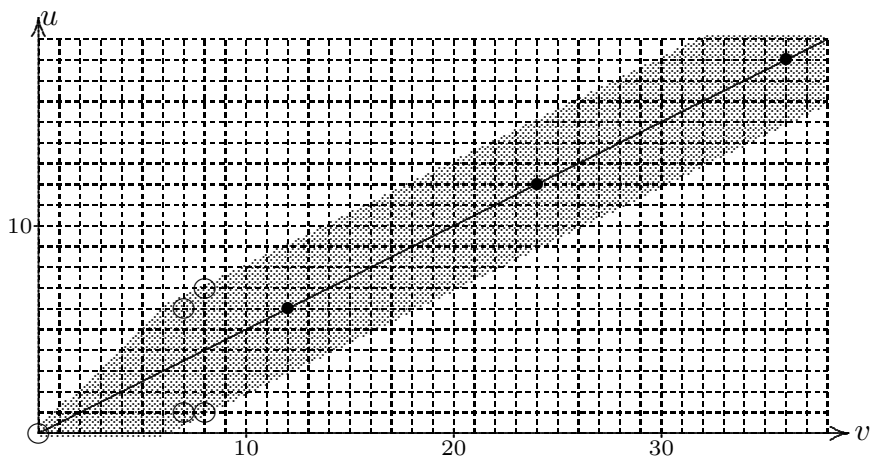
\hbox{\beginpicture
\setcoordinatesystem units <.23cm,.23cm>
\arr{0 0}{0 20}
\arr{0 0}{39 0}
\put{$u$} at  0.5 20
\put{$v$} at 39.7 0
\setdashes <.5mm> 
\setplotarea x from 0 to 38, y from 0 to 19
\grid {38} {19}
\setsolid
\plot 0 0  38  19 /
\setshadegrid span <.3mm>
\vshade 0 -.2 .3 <,z,,> 6 -.2 6.2 <z,z,,> 32 12.8 19.2 <z,,,> 38 15.8 19.2 /  
\plot 10 -0.2  10 0.2 /
\put{$\ssize 10$} at 10 -.7
\plot 20 -0.2  20 0.2 /
\put{$\ssize 20$} at 20 -.7
\plot 30 -0.2  30 0.2 /
\put{$\ssize 30$} at 30 -.7

\plot -0.2 10  0.2 10 /
\put{$\ssize 10$} at -.9 10 
\multiput{$\bullet$} at 12 6  24 12  36  18 /
\multiput{$\ssize \bigcirc$} at 0 0  7 1  8 1  7 6  8 7 /
\endpicture}
$$
One observes that for an indecomposable triple $(V,U,T)$ in $\Cal S(6)$, the
dimension of the subspace 
$U$ is roughly half of the dimension of the total space $V$. Let us stress the symmetry
with respect to the operation $\delta\:(v,u)\mapsto (v,v-u)$ 
(it fixes the points on the line through $(2,1)$ 
and maps a point with vertical distance $\frac12 v-u$ below the line 
to the corresponding point of the same vertical distance above the line).
This operation describes duality:
The categories  $\Cal S(n)$ have the self duality $*=\Hom_k(-,k)$
which maps an object $(V,U,T)$ with dimension pair $(v,u)$ to  the object
$(V,U,T)^*=(V^*,(V/U)^*,T^*)$ with dimension pair $\delta(v,u)$. 

Actually, for $n\leq 6$ we even have $|v-2u|\leq n$, whereas for $n\geq 7$, the
numbers $|v-2u|$ are not bounded.
The black dots in the picture indicate the positions 
which are positive integral multiples of $(12,6)$; the
difference between these dimension pairs and the remaining ones will be discussed next.
        \medskip\noindent
{\bf (0.1.2)}
{\it Assume that $(v,u)$ is not an integral multiple of $(12,6).$ Then the number of isomorphism
classes of indecomposable objects in $\Cal S_k(6)$ with dimension pair 
$(v,u)$ is finite and is independent of $k$.}
        \smallskip
By contrast, the indecomposable objects in $\Cal S_k(6)$ with dimension pair
an integral multiple of $(12,6)$ will depend on the field $k$.
        \medskip\noindent
{\bf (0.1.3)} 
{\it Let $k$ be an infinite field. Let $v,u$ be natural numbers. Then:
There are infinitely many 
isomorphism classes of indecomposable objects in $\Cal S_k(6)$ 
with dimension pair $(v,u)$ if and only if $(v,u)$ is a positive integral multiple of $(12,6).$}
        \smallskip
It may be worthwhile to exhibit here a family of indecomposable triples 
$X_c = (V,U_c,T)$ with dimension pair
$(12,6).$ Let $V$ be a $k$-space with basis $e_1,\dots,e_{12}$, let $T$ be defined by
$T(e_i) = 0$ for $i\in \{1,7,11\}$ and $T(e_i) = e_{i-1}$ otherwise (thus $T$ is given 
by the partition $(6,4,2)$). For $c\in k$, let $U_c$ be the subspace
generated by the three elements $e_3+e_{8}$, $e_{11}-e_8$ and $ce_4+(c-1)e_9+e_{12}$,
and by their images under the powers of $T$.
All these triples $X_c$ are indecomposable and they are pairwise non-isomorphic;
we obtain in this way a {\it one-parameter family indexed by}
$k$.  If we consider indecomposable triples
$(V,U,T)$ with dimension pair $(24,12)$, then there are already two disjoint
one-parameter families, both being  indexed by
$k$; for one family, the space $V$ is given by the partition $(6,6,4,4,2,2)$, for the second
family, by $(6,6,5,3,3,1).$ 
In case $k$ is algebraically closed, there is the following general result:
        \medskip\noindent
{\bf (0.1.4)} {\it Let $k$ be an algebraically closed field, let $t\in \Bbb N_1$. Then there are
$t$ pairwise disjoint one-parameter families of indecomposable objects in $\Cal S_k(6)$ 
with dimension pair $(12t,6t),$ each one being indexed by $k$, and such that there are only 
a finite number of additional isomorphism classes of indecomposable triples with 
dimension pair $(12t,6t).$}
        \smallskip
Assertions (0.1.3) and (0.1.4) show that the classification of the indecomposable objects in
$\Cal S(6)$ is what is called a tame, but non-domestic problem with linear growth (of the
number of one-parameter families). An index set for the indecomposable triples in $\Cal S(6)$
will be provided in this paper. 
For a detailed description and further properties of these indecomposable triples, 
we refer to the forthcoming Part II [12] (in particular, there we
will see that for almost all indecomposable triples $(V,U,T)$ 
with dimension pair $(12t,6t)$,
the dimension of the kernel of $T$ is $3t$, whereas the dimensions 
of the kernels of the endomorphisms induced by $T$ on $U$ as well as on $V/U$ 
both are $2t$).
        \medskip
We now draw the attention to the endomorphism rings $\End(X)$ of indecomposable objects
$X$ in $\Cal S(6).$ These endomorphism rings are usually rather large, however the bulk
of endomorphisms will be nilpotent with nilpotency index at most 8.  
        \medskip\noindent
{\bf (0.1.5)} {\it Let $X$ be an indecomposable object in 
$\Cal S_k(6)$ with dimension pair  $(v,u).$
Then there is an ideal $I$ in $\End(X)$ with $I^8 = 0$ such that
$\End(X)/I$ is a local uniserial ring.}
        \smallskip
The ideal $I$ can be described using the notion of the infinite radical of a category.
Given a Krull-Remak-Schmidt category such that the indecomposable objects have local
endomorphism rings (such as $\Cal S(n)$), its 
{\it radical} $\Cal R$ is by definition the class of all maps
$g\:X\to Y$ such that for any pair of maps $f\:X'\to X$ and $h\:Y\to Y'$ with
$X'$ and $Y'$ indecomposable, the composition $hgf$ is a non-isomorphism.
Then $\Cal R$ forms a categorical ideal.  We are interested in the
{\it infinite radical,} $\Cal I=\bigcap_{i\in \Bbb N} \Cal R^i$.
        \smallskip\noindent
{\bf (0.1.6)} {\it In the category $\Cal S(6)$, the infinite radical $\Cal I$ is
an idempotent ideal, that is, $\Cal I^2=\Cal I$ holds.  If $X$ is an
indecomposable object in $\Cal S(6)$, then the ideal $\Cal I(X,X)$ of $\End(X)$ has
nilpotency index at most $8$.}
        \medskip
The ideal $I$ to be used in Assertion (0.1.5) is $I = \Cal I(X,X).$
        \medskip
An inclusion map $(V',U',T') \to (V,U,T)$ is called an {\it inflation} provided $U' = U\cap V'$, thus
provided the induced map $U/U' \to V/V'$ is injective. If $X' = (V',U',T') \to (V,U,T) = X$ is an
inflation, then we can form $X/X'$ with total space $V/V'$, subspace $U/U'$ and with the 
linear operator on $V/V'$
induced by $T$. Let $X$ be an indecomposable object. Call $X$ {\it $m$-stackable} provided there is 
an indecomposable object $Y$ with a chain of inflations
$$
 0 = Y_0 \to Y_1 \to \cdots \to Y_m = Y
$$
such that all the objects $Y_i/Y_{i-1}$ with $1 \le i \le m$ are isomorphic to $X$. 
        \medskip\noindent
{\bf (0.1.7)} {\it Assume that 
$X$ is an indecomposable object in $\Cal S(6)$ with dimension pair $(v,u).$
If $(v,u)$ is an 
integral multiple of $(12,6)$, then $X$ is $m$-stackable for every $m\in \Bbb N_1$.
If $(v,u)$ is not a
multiple of $(2,1)$, and $X$ is $m$-stackable for some $m\in \Bbb N_1$, then $m \le 6.$}
        \smallskip
Here, we provide no assertion concerning the cases where $(v,u)$ is a multiple of $(2,1)$,
but not an integral multiple of $(12,6).$ It seems to us that in these cases $X$ can be $m$-stackable 
only for $m \le 5.$ Note that in the categories $S(n)$ with $n \le 5$, all the indecomposable objects
can be $m$-stackable only for $m\le n$, and those with dimension pair a multiple of $(2,1)$
can be $m$-stackable only for $m \le n-1.$ 
        \bigskip\noindent
{\bf (0.2) A general frame.} 
We may reformulate (and generalize) the problem to be considered as follows: 
If $\Lambda$ is any ring, let $\Cal S(\Lambda)$ be the category of pairs
$(A,A')$, where $A$ is a finitely generated $\Lambda$-module and 
$A'\subseteq A$ is a submodule of $A$; a map $f\:(A,A') \to (B,B')$ in $\Cal S(\Lambda)$ 
is just a $\Lambda$-linear map $f\:A \to B$ such that $f(A') \subseteq B'$ holds. 
We may call $\Cal S(\Lambda)$ the {\it submodule category} of $\Lambda$-modules.
One will observe that 
$$
 \Cal S_k(n) = \Cal S(k[x]/x^n)
$$ 
where $k[x]$ is the polynomial ring with variable $x$;
namely, a pair $(V,T)$
consisting of a $k$-vector space and a linear operator $T\:V \to V$ with $T^n = 0$ is
in an obvious way a $k[x]/x^n$-module, and a subspace $U$ of $V$ with $T(U) \subseteq U$
leads to a pair $(U,T|U)$ which is just a submodule of $(V,T)$. Note that the ring
$k[x]/x^n$ is a commutative uniserial ring of length $n$ and one may wonder about the
similarities of the various categories $\Cal S(\Lambda)$, with $\Lambda$ 
a commutative uniserial ring of length $n$. In fact, our interest in the problem discussed
in this paper arose from the Birkhoff problem [2] of classifying subgroups
of finite abelian $p$-groups. The Birkhoff problem concerns the categories
$\Cal S(\Bbb Z/p^n)$, thus also here we deal with 
a commutative uniserial ring of length $n$, namely $\Bbb Z/p^n$. The combinatorial 
structure for the categories $\Cal S(n)$ with $n \le 5$ presented in
Section~6 is the same as  for $\Cal S(\Bbb Z/p^n)$.
The structure of the category $\Cal S(\Bbb Z/p^6)$ is not yet known --- it is a 
challenging question whether it looks like the category $\Cal S_k(6)$, for $k = \Bbb Z/p$,
exhibited in this paper. 
        \smallskip
Assume now that $\Lambda$ is an artinian ring (for example a finite-dimensional $k$-algebra, as is the
case when $\Lambda = k[x]/x^n$). We may
consider  $\Cal S(\Lambda)$ as the full 
subcategory of the category $\Cal H(\Lambda)$ of finitely generated modules
over the ring 
$\big[{\Lambda\atop0}\,{\Lambda\atop\Lambda}\big]$ 
of upper triangular matrices as follows: Write the pair 
$(A,A')$ in $\Cal S(\Lambda)$ as column vectors
$\big[{A\atop A'}\big]$ and just use matrix multiplication. 
In this way, one sees immediately that
$\Cal S(\Lambda)$ inherits the Krull-Remak-Schmidt property
from $\Cal H(\Lambda)$: Any object in $\Cal S(\Lambda)$
is the direct sum of a finite number of indecomposable objects and these 
indecomposable direct summands are unique up to isomorphism.
Thus, in order to classify all the
isomorphism classes of triples $(V,U,T)$ in $\Cal S(\Lambda)$, 
it is sufficient to consider the
indecomposable ones. We return now to the case $\Lambda = k[x]/x^n.$
        \bigskip\noindent

{\bf (0.3) Auslander-Reiten theory.} 
The description of the categories $\Cal S(n)$ presented in this paper uses notions
and results from the modern representation theory of finite dimensional $k$-algebras, 
see for example [1] and [9]. 
The basic notions of Auslander-Reiten theory lead to
the possibility of presenting suitable additive categories such as $\Cal S(n)$ at least partly 
by a translation quiver, and the
components of this Auslander-Reiten quiver of $\Cal S(n)$ will be of our concern. 
        \smallskip
For each $n < 6$ we will present the Auslander-Reiten quiver
of $\Cal S(n)$ in Section~6. 
For $n > 6$, we can refer to [10] where it is shown 
that the category $\Cal S_k(n)$ is controlled $k$-wild; in this case the problem of 
classifying all objects is considered infeasible.  Still we can present 
some features of these categories in Section~7. But as we have mentioned already, the
most interesting case seems to be the borderline case $n=6$;
Part A  (as well as the forthcoming
paper [12]) is devoted to exhibit its structure.
        \smallskip
The determination of the structure of the Auslander-Reiten components of $\Cal S_k(6)$
relies on 
methods from the representation theory of quivers which are by now standard: the
knitting of preprojective components, 
the covering theory of the Gabriel school, as well as the use of
one-point extensions, in particular the tubular extension techniques as outlined in [9].
However, all these methods have to be modified since we do not deal with a module
category, but only with an exact category. The necessary modification of these
standard tools will be outlined below. It is amazing that the results obtained
are close to corresponding results concerning module categories.
In particular, the structure of the category $\Cal S(6)$ turns out to be very
similar to the structure of some module categories.
        \smallskip
In order to describe the structure of the category $\Cal S_k(6)$,
one needs the notion of a tube, or better of tubular families. 
Tubular families are sets of 
Auslander-Reiten components all of which are tubes and being indexed
by the projective line $\Bbb P_1(k)$ (in case $k$ is algebraically closed, $\Bbb P_1(k)$
may be considered as the disjoint union of $k$ itself and an additional symbol~$\infty$;
in general, $\Bbb P_1(k)$ 
is the disjoint union of the set of monic irreducible polynomials over $k$ and the
additional symbol~$\infty$). 
Let us try to give at least
an intuitive picture of the global structure 
of the category $\Cal S_k(6)$, it turns out to be a category of tubular type 
$\Eeighttubular$.
Each indecomposable object, up to isomorphism, is uniquely determined
by three parameters 
\item{$\bullet$} the first one is a non-negative rational number $\gamma\in\Bbb Q_0^+$, 
\item{$\bullet$} the second is an element $c\in\Bbb P_1(k)$,
\item{$\bullet$} and as third parameter one can take a natural number $m\in\Bbb N$.
        \par\noindent
Conversely, any combination $(\gamma,c,m)$ occurs in this way, thus we obtain in this
way a complete set of invariants. 
        \smallskip
If we fix the first index $\gamma$, then we obtain the tubular family $\bar{\Cal T}_\gamma$
by taking all the indecomposable triples with first index $\gamma$. Instead of
taking the set $\Bbb Q_0^+$ as first index set, we should better consider the 
topological space 
$\big\{e^{2\pi i\frac{p}{p+q}}\big|\gamma=\frac pq\in\Bbb Q_0^+\big\}$, 
so that the points corresponding to $0$ and $\infty$ 
in $\Bbb Q_0^+\cup\{\infty\}$ get identified:
$$
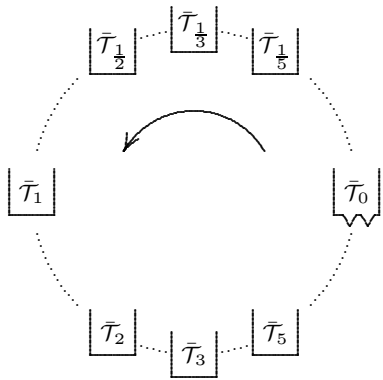
\beginpicture
\setcoordinatesystem units <1.8cm,1.8cm> 
\put{$\ssize {\bar{\Cal T}}_{\frac 13}$} at 0 1
\put{$\ssize {\bar{\Cal T}}_{\frac 15}$} at 0.5 0.866
\put{$\ssize {\bar{\Cal T}}_0$} at 1 0
\put{$\ssize {\bar{\Cal T}}_5$} at 0.5 -0.866
\put{$\ssize {\bar{\Cal T}}_3$} at 0 -1
\put{$\ssize {\bar{\Cal T}}_2$} at -0.5 -0.866
\put{$\ssize {\bar{\Cal T}}_1$} at -1 0
\put{$\ssize {\bar{\Cal T}}_{\frac 12}$} at -0.5 0.866
\put{\tube} at 0 1
\put{\tube} at 0.5 0.866
\put{\Utube} at 1 0
\put{\tube} at 0.5 -0.866
\put{\tube} at 0 -1
\put{\tube} at -0.5 -0.866
\put{\tube} at -1 0
\put{\tube} at -0.5 0.866

\circulararc 120 degrees from 0.433 0.25  center at 0 0
\arr {-0.408 0.2933}  {-0.433 0.25}

\setdots<2pt>
\circulararc -10 degrees from  0.174  0.985 center at 0 0
\circulararc  30 degrees from  0.966  0.259 center at 0 0
\circulararc -30 degrees from  0.966 -0.259 center at 0 0
\circulararc  10 degrees from  0.174 -0.985 center at 0 0
\circulararc -10 degrees from -0.174 -0.985 center at 0 0
\circulararc  30 degrees from -0.966 -0.259 center at 0 0
\circulararc -30 degrees from -0.966  0.259 center at 0 0
\circulararc  10 degrees from -0.174  0.985 center at 0 0

\endpicture
$$
As we have mentioned, for each value of $\gamma$, there is a 
tubular family ${\bar{\Cal T}}_\gamma$, thus a collection of Auslander-Reiten components
which are tubes and which are para\-metrized by the second parameter $c$ from $\Bbb P_1(k)$.
All the tubes in ${\bar{\Cal T}}_\gamma$ but three are homogeneous ones, 
the exceptional ones (say with second parameter $c = 1, 0, \infty$)
are tubes of rank 2, 3, and 6, respectively. 
The third parameter $m$ enumerates the triples in a given tube: there are countable many
isomorphism classes of triples which belong to any  tube. For the homogeneous
tubes, there is indeed a natural bijection with the set of natural numbers, but for the
exceptional tubes the use of natural numbers as the third parameter 
may be considered as rather misleading, it should be considered only as a first 
attempt to get an impression of the diversity of the isomorphism classes in $\Cal S(6)$.
As the proper third index set one should replace $\Bbb N$ by the vertices of a tube
of rank $r(c)$, depending on the second parameter $c$, with $r(\infty) = 6,\ r(0) = 3,\ r(1) = 2$,
and $r(c) = 1$ otherwise. 
        \smallskip
One deviation  has to be mentioned: All but one of the tubular families
$\bar{\Cal T}_\gamma$ have the same shape, the only exception occurs for $\gamma = 0$.
Here, the tube of rank~6 is not stable but contains two projective-injective
vertices; this should explain the curious shape used in the picture above.
The three non-homogeneous tubes in the tubular family $\bar{\Cal T}_0$ (which
we also denote by $\bar{\Cal U}$) will be exhibited in detail in Section~2. 
        \smallskip
We have mentioned already that the distribution 
of the indecomposable triples of $\Cal S_k(6)$
into tubes yields the Auslander-Reiten components of the category $\Cal S_k(6)$,
thus we describe in this way the factor of this category modulo the infinite
radical. But this concerns only the third parameter, whereas the first two parameters
yield information on the infinite radical. The arrangement of the tubular families 
along a circle indicates that the homomorphisms in the infinite radical go 
counterclockwise (see Theorem (1.5.2) and Lemma (2.2.2) for details). The picture presented above
includes a corresponding arrow in order to stress the global direction of 
homomorphisms in the infinite radical.
In particular, the endomorphism of an object $(V,U,T)$ given by multiplication by $T$
always belongs to the infinite radical.
\vfill
        \bigskip\noindent
{\bf (0.4) Covering theory.} 
In order to visualize indecomposable objects in $\Cal S(n)$, we will use two different
ways, namely pairs of sequences of natural numbers, as well as box diagrams (see for
example Section~2). Both ways rely on the use of covering theory which we now want to
explain, since it is the essential tool for our investigation. 
        \smallskip
As we have seen, the objects of the category $\Cal S(n)$ can be described 
as triples $(V,U,T)$ or also as pairs of a $k[x]/x^n$-module and a submodule. But there
is a third description for these objects, namely as representations of the quiver
$$
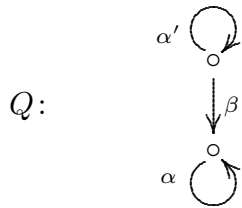
\hbox{\beginpicture
\setcoordinatesystem units <0.5cm,0.5cm>
\put{} at -4 -1
\put{} at 1 3.2
\put{$\circ$} at 0 0
\put{$\circ$} at 0 2
\circulararc 160 degrees from 0 3.2 center at 0 2.7
\circulararc -150 degrees from 0 3.2 center at 0 2.7
\circulararc 150 degrees from 0 -1.2 center at 0 -.7
\circulararc -160 degrees from 0 -1.2 center at 0 -.7
\arr{0 1.6}{0 0.4}
\put{$\ssize \alpha'$} at -1 2.6
\put{$\ssize \alpha$} at -1 -.6
\put{$\ssize \beta$} at 0.4 1
\arr{0.3 2.31} {0.2 2.2}  
\arr{0.3 -.31} {0.2 -.2}
\put{$Q\:$} at -4 1
\endpicture}
$$
which satisfy the relations $\beta\alpha' = \alpha\beta$
and $\alpha^n = 0$, and with the additional requirement that the map $\beta$ is 
an inclusion map.
(We do not have to mention the condition $(\alpha')^n = 0$ since
it follows from the two relations $\beta\alpha' = \alpha\beta$
and $\alpha^n = 0$ in case $\beta$ is assumed to be injective).
We also will consider the following (infinite) quiver (with respect to the
given relations, this is what is called the universal covering for $Q$):
$$
\hbox{\beginpicture
\setcoordinatesystem units <0.5cm,0.5cm>
\put{} at 0 0
\put{} at 5 2
\put{$\circ$} at 0 0
\put{$\circ$} at 2 0
\put{$\circ$} at 4 0
\put{$\circ$} at 0 2
\put{$\circ$} at 2 2
\put{$\circ$} at 4 2
\arr{-0.4 0}{-1.6 0}
\arr{-0.4 2}{-1.6 2}
\arr{1.6 0}{0.4 0}
\arr{1.6 2}{0.4 2}
\arr{3.6 0}{2.4 0}
\arr{3.6 2}{2.4 2}
\arr{5.6 0}{4.4 0}
\arr{5.6 2}{4.4 2}

\arr{0 1.6}{0 0.4}
\arr{2 1.6}{2 0.4}
\arr{4 1.6}{4 0.4}
\put{$\cdots$} at -3 0
\put{$\cdots$} at -3 2

\put{$\cdots$} at 7 0
\put{$\cdots$} at 7 2
\put{$\ssize \alpha'_0$} at -1 2.4 
\put{$\ssize \alpha'_1$} at 1 2.4
\put{$\ssize \alpha'_2$} at 3 2.4
\put{$\ssize \alpha'_3$} at 5 2.4
\put{$\ssize \alpha_0$} at -1 -.4 
\put{$\ssize \alpha_1$} at 1 -.4
\put{$\ssize \alpha_2$} at 3 -.4
\put{$\ssize \alpha_3$} at 5 -.4
\put{$\ssize \beta_0$} at 0.5 1 
\put{$\ssize \beta_1$} at 2.5 1
\put{$\ssize \beta_2$} at 4.5 1
\put{$\ssize 0'$} at  0 2.5
\put{$\ssize 1'$} at  2 2.5
\put{$\ssize 2'$} at  4 2.5

\put{$\ssize 0$} at 0 -.5
\put{$\ssize 1$} at 2 -.5
\put{$\ssize 2$} at 4 -.5
\put{$\widetilde Q\:$} at -6 1
\endpicture}
$$
Usually, we will refrain from mentioning the indices of the arrows 
$\alpha_i,\alpha'_i,\beta_i$ and just write 
$\alpha,\alpha',\beta$.
The representations of $k\widetilde Q$ which satisfy the commutativity relations
$\beta\alpha'=\alpha\beta$ (for all squares), the nilpotency relations
$(\alpha')^n=\alpha^n=0$ (for all compositions of $n$ arrows 
$\alpha'$ or $\alpha$) and for which all the maps $\beta$ are realized
by  monomorphisms form the objects of a category 
which we will denote by $\Cal S(\widetilde n)$, and which will play an
important role throughout the paper.

Of course, we may interprete also $\Cal S(\widetilde n)$ as a submodule
category, namely for an infinite dimensional algebra which we denote by
$kA_\infty^\infty/\alpha^n$,
$$
 \Cal S(\widetilde n)=\Cal S(kA_\infty^\infty/\alpha^n).
$$
Here, $kA_\infty^\infty/\alpha^n$ is the factor algebra of the path 
algebra $kA_\infty^\infty$
of the following infinite quiver
$$
\hbox{\beginpicture
\setcoordinatesystem units <0.5cm,0.5cm>
\put{} at 0 0
\put{} at 5 0
\put{$\circ$} at 0 0
\put{$\circ$} at 2 0
\put{$\circ$} at 4 0
\arr{-0.4 0}{-1.6 0}
\arr{1.6 0}{0.4 0}
\arr{3.6 0}{2.4 0}
\arr{5.6 0}{4.4 0}

\put{$\cdots$} at -3 0

\put{$\cdots$} at 7 0
\put{$A_\infty^\infty:$} at -6 0
\endpicture}
$$
modulo the ideal generated by the relations $\alpha^n = 0$ (all arrows being labelled $\alpha$). 
        \smallskip
The categories $\Cal S(\widetilde n)$ and $\Cal S(n)$ are related by 
an important functor
$$
 \Cal S(\widetilde n) \longrightarrow \Cal S(n),
$$ 
the so-called covering functor.
This functor assigns 
to a representation $(A,A')$ in $\Cal S(\widetilde n)$
the triple $(V,U,T)$ in $\Cal S(n)$ using
for $V$ und $U$ the direct sums $V =\bigoplus_{i\in\Bbb Z}A_i$ and 
$U =\bigoplus_{i\in\Bbb Z} A_i'$, respectively; note that 
the diagonal map 
$\bigoplus_{i\in \Bbb Z}\beta_i\: \bigoplus_{i\in\Bbb Z} A_i'\to
\bigoplus_{i\in\Bbb Z} A_i$ provides an embedding of $U$ into $V$,
and the linear operator $T$ on $V$
is given by the action of the maps corresponding to the arrows $\alpha$.

We will see that for $n \le 6$, the covering functor is dense, thus
any indecomposable triple $(V,U,T)$ in $\Cal S(n)$ is obtained from a 
representation $(A,A')$ of $\widetilde Q$, and this 
representation $(A,A')$ is unique up to the shift, thus its properties
can be used in order to find invariants of $(V,U,T)$. 
        \bigskip\noindent
{\bf (0.5) Contents.}
The paper is organized as follows: As we have mentioned there are two parts.
Part~A with Sections 1, 2, 3 and 4 is devoted to the case $n=6$, whereas the
Sections 6 and 7 are devoted to the cases $n < 6$ and
$n > 6$, respectively. Part B consists of these Sections 6, 7 and an additional Section
5 which provides information concerning some boundary objects in general.  
        \smallskip
In Section~1 we will see that $\Cal S(\widetilde 6)$ has tame  
representation type, more precisely, it is of tubular type $\Eeighttubular$, 
and we exhibit the complete structure of its Auslander-Reiten quiver.
        \smallskip
The topic of Section 2 is the covering functor
$\Cal S(\widetilde n)\to \Cal S(n)$. We will see that in case $n\leq 6$ 
the covering functor is dense, hence it induces a bijection
between the $\Bbb Z$-orbits of indecomposable objects in $\Cal S(\widetilde n)$ 
under the shift $M\mapsto M[1]$, and 
the indecomposable objects in $\Cal S(n)$. In particular, we can use results
from Section~1
in order to obtain not only a full list of the indecomposable objects in
$\Cal S(6)$, but also the Auslander-Reiten quiver for $\Cal S(6)$.
        \smallskip      
In Section 3 we demonstrate the crucial results stated in (0.1).
        \smallskip      
Section 4 provides a slightly different approach for the study of $\Cal S(n)$.
A modification of an essential part of the quiver $\widetilde Q$ leads us to
a finite dimensional algebra $\Omega$ which provides a sort of double
covering of $\Cal S(n)$. This approach has the advantage that we obtain in this way 
an abelian category, namely the category of $\Omega$-modules, which is very similar to
the category $\Cal S(6)$. 
        \smallskip
In Part B we consider the category $\Cal S(n)$ for arbitrary $n$. 
Topic of Section~5 are the objects which form the boundary of the non-stable component of
the Auslander-Reiten quiver of $\Cal S(n)$ (including the projective-injective ones).
They are of similar nature for all the categories
$\Cal S(n)$.
        \smallskip
In Section 6 we assume that $n<6$; then $\Cal S(n)$ has finite representation type.
We specify the Auslander-Reiten quivers and picture each of the
indecomposable objects.  It turns out that the number $s(n)$ of indecomposables 
is given by the formula
$$ 
 s(n) \;=\; 2+2(n-1)\frac 6{6-n}.
$$
We motivate this formula and note that it 
still makes sense for $n=6$ (with both sides being infinite). However for $n = 7$
the right side yields $-70$, and we do not know any interpretation of this number
in terms of representation theory.
        \smallskip
In Section 7 we will consider the case where $n>6$. Then each of the categories
$\Cal S(n)$ has wild representation type;
at least we can determine the shape of the components of the
Auslander-Reiten quiver:  There is one tube containing the boundary objects
on its mouth; all the other components are stable tubes of diameter
a divisor of 6.
Also each corresponding category
$\Cal S(\widetilde n)$ is wild; here the connected components of the Auslander-Reiten
quiver have stable type $\Bbb ZA_\infty$. We can construct in $\Cal S(\widetilde n)$
indecomposable representations of arbitrarily large support.
        \smallskip
There is an Appendix which outlines in which way considerations concerning nilpotent 
operators are of interest when dealing with operators in general.
        \medskip
\noindent{\it Notation.} \/
In this paper we will study representations of the quiver $Q,$ as well as of 
various full subquivers
$Q'$ of the quiver $\widetilde Q$ (taking into account the given relations). 
The projective, injective, and simple representations 
corresponding to
a vertex $z$ of a quiver will be denoted by $P(z)$, $I(z)$, and $S(z)$, respectively.
For $M$ a representation, $M_z$ is the $k$-space at position $z$.
For $M$ a representation of $\widetilde Q$ and $\ell$ an integer, $M[\ell]$ 
is the representation of $\widetilde Q$ obtained by shifting $M$ by $\ell$
steps in direction opposite to $\alpha$. Thus, $M[\ell]_i=M_{i-\ell}$
and $M[\ell]_{i'}=M_{(i-\ell)'}$, for $i\in\Bbb Z$. If $\Cal A$ is a Krull-Remak-Schmidt
category, we denote by $\ind \Cal A$ the class of indecomposable objects in $\Cal A$
(or the choice of one from each isomorphism class). 
For unexplained notions
we refer to [9] and [1].
        \bigskip\noindent
{\bf (0.6) History and related results.} 
As we have mentioned,
it was Birkhoff~[2] who urged to study subcategories of the form 
$\Cal S(\Lambda)$, at least for the rings $\Lambda = \Bbb Z/p^n.$ 
The indecomposable objects in $\Cal S(\Bbb Z/p^n)$
have been determined in [8] in the cases where $n\leq 5$. 
More generally, one may consider not only submodules, but even chains of
submodules, as proposed by Simson in  [14]. For 
each pair $\ell, n$ of natural numbers, he determined the representation type
of the category of chains of length $\ell$ of submodules of $k[x]/x^n$-modules
(in particular, he already has shown that the category $\Cal S(6)$ which we will
study in detail, is tame). 

Also, for $m\leq n$, let us denote by $\Cal S(n;m)$ 
the full subcategory of $\Cal S(n) = \Cal S_k(n)$ of all triples $(V,U,T)$
such that $T^m(U) = 0$. 
In [13] the representation type of $\Cal S(n;m)$ 
has been determined for each pair $(m,n)$ with $m<n$,
and in each finite or tame case a description of the indecomposable
modules has been given.

In a forthcoming paper [12] we are going to expand the study of the category $\Cal S(6)$
considerably. We will outline the relevance of the tubular family $\bar{\Cal T}_1$ and we will show 
in which way the two half circles of tubular families $\{\bar{\Cal T}_\gamma \mid 0 < \gamma < 1\}$
and $\{\bar{\Cal T}_\gamma \mid 1 < \gamma < \infty \}$  are built up from the tubular families
$\bar{\Cal T}_0$ and $\bar{\Cal T}_1.$ Also, we will see that the objects in $\Cal S(6)$ are strongly
related to the representations of a partially ordered set with 9 vertices.
        \bigskip\noindent
{\bf (0.7) Acknowledgement.}
Most of the results presented here have been obtained in the years 2000-2002. They have been
reported in various lectures since that time, and many comments of the audience have 
helped us to improve the results. Unfortunately, it took a long time to complete
a final version. The authors  are grateful to the European Union 
for providing the possibility for the start of this
cooperation (at a rather boring midterm review meeting in 2000 where
the EU representatives insisted that the discussions should only concern 
administrative matters, and not mathematics \dots). 

        \bigskip\bigskip
\vfill\eject
\centerline{\bf Part A\. Operators with Nilpotency Index 6}
        \medskip
We give a description of all indecomposable triples $(V,U,T)$ 
where the linear operator acts with nilpotency index at most 6,
i.e.\ $T^6=0$ holds.  Recall that these triples are the objects of the
category $\Cal S(6).$ But first we study a related category 
$\Cal S(\widetilde 6)$.

        \bigskip\bigskip
{\bf 1. The Category $\Cal S(\widetilde 6)$}
        \medskip
Let $\Cal H(\widetilde 6)$ be the category of all 
(finite dimensional) representations
of $\widetilde Q$ which satisfy the commutativity relations
$\beta\alpha'=\alpha\beta$ and the nilpotence conditions 
$\alpha^6=0=(\alpha')^6$, thus, using notation from (0.2),
$\Cal H(\widetilde6)=\Cal H(kA_\infty^\infty/\alpha^6)$. 
Clearly the category $\Cal S(\widetilde6)$, which we are going to describe
in this section, is the full subcategory of $\Cal H(\widetilde 6)$ of
all representations for which the maps $\beta$ are monomorphisms. The category
$\Cal H(\widetilde 6)$ will be used as a general reference category, its
Grothendieck group $K_0(\Cal H(\widetilde 6))$
(of all objects modulo exact sequences) will be denoted
just by $K_0.$ Note that $K_0$ is free abelian with basis $\be(z)$,
where $z$ is a vertex of $\widetilde Q$. If $M$ is a representation of 
$\widetilde Q$, we will denote by $\bdim M$ the corresponding element in
$K_0$ and call it its {\it dimension vector}. Given a dimension vector $\bd$,
the coefficient with index $z$ (here, $z$ is a vertex of $\widetilde Q$) will be denoted
by $\bd_z.$ 

We consider first the full subcategory $\Cal C\subseteq \Cal
H(\widetilde6)$, which has as objects those representations $M$ 
which satisfy
\item{(a)} $M_z=0$ and $M_{z'}=0$ for all $z\in\{-1,-2,\ldots\}$ and
\item{(b)} the maps $M(\beta_0)$ and $M(\beta_1)$ are the identity maps.

\smallskip
When dealing with $\Cal C$, we do not have to mention
the vertices $0'$ and $1'$, as the maps $\alpha'_1 = \alpha_1$ and 
$\alpha'_2 = \alpha_2\beta_2$ can be recovered.  
In fact, the category $\Cal C$ is equivalent
to the category of representations of the quiver
$$
\hbox{\beginpicture
\setcoordinatesystem units <0.5cm,0.5cm>
\put{} at 0 0
\put{} at 5 2
\put{$\circ$} at 0 0
\put{$\circ$} at 2 0
\put{$\circ$} at 4 0
\put{$\circ$} at 6 0
\put{$\circ$} at 8 0
\put{$\circ$} at 4 2
\put{$\circ$} at 6 2
\put{$\circ$} at 8 2
\arr{1.6 0}{0.4 0}
\arr{3.6 0}{2.4 0}
\arr{5.6 0}{4.4 0}
\arr{7.6 0}{6.4 0}
\arr{9.6 0}{8.4 0}

\arr{5.6 2}{4.4 2}
\arr{7.6 2}{6.4 2}
\arr{9.6 2}{8.4 2}

\arr{4 1.6}{4 0.4}
\arr{6 1.6}{6 0.4}
\arr{8 1.6}{8 0.4}

\put{$\cdots$} at 11 0
\put{$\cdots$} at 11 2
\put{$\ssize \alpha'$} at 5 2.4
\put{$\ssize \alpha'$} at 7 2.4
\put{$\ssize \alpha$} at 1 -.4
\put{$\ssize \alpha$} at 3 -.4
\put{$\ssize \alpha$} at 5 -.4
\put{$\ssize \alpha$} at 7 -.4
\put{$\ssize \beta$} at 4.5 1 
\put{$\ssize \beta$} at 6.5 1
\put{$\ssize \beta$} at 8.5 1
\put{$\ssize 2'$} at  4 2.5
\put{$\ssize 3'$} at  6 2.5
\put{$\ssize 4'$} at  8 2.5

\put{$\ssize 0$} at 0 -.5
\put{$\ssize 1$} at 2 -.5
\put{$\ssize 2$} at 4 -.5
\put{$\ssize 3$} at 6 -.5
\put{$\ssize 4$} at 8 -.5
\endpicture}
$$
with commutativity relations $\beta\alpha' = \alpha\beta$
and nilpotence conditions $(\alpha')^6 = \alpha^6 = 0.$
In order to stress that we consider 
a full subcategory of $\Cal H(\widetilde6)$,
we insert the numbers $\bd_{0'} = \bd_0$ and $\bd_{1'} = \bd_1$ when we display the
dimension vector $\bd$ of a representation of this quiver; for example, the
dimension vector of the indecomposable projective representation corresponding to
the vertex $3$ will be displayed as $\ssize{11000\cdots\atop11110\cdots}$, 
or just as 
$\ssize {1100\atop 1111}$ (we only will
exhibit dimension vectors $\bd$ for representations of $\widetilde Q$ 
which start on the left with the numbers $\bd_{0'}$ and $\bd_0$).

\medskip\noindent {\bf (1.1) The tubular algebra $\Theta$.}
An essential part of the category $\Cal C$, and hence of 
$\Cal S(\widetilde 6)$, consists of modules over the finite dimensional
algebra $\Theta$ which is given by the following subquiver of $\widetilde Q$,
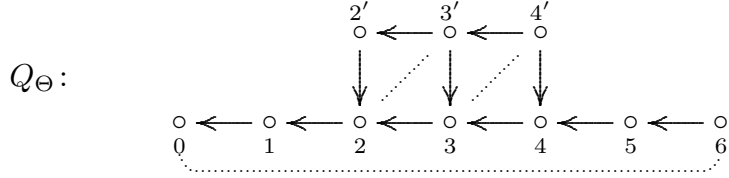
$$
\hbox{\beginpicture
\setcoordinatesystem units <0.5cm,0.5cm>
\put{} at -3 0
\put{} at 5 2
\put{$Q_\Theta\:$} at -3 1
\put{$\circ$} at 0 0
\put{$\circ$} at 2 0
\put{$\circ$} at 4 0
\put{$\circ$} at 6 0
\put{$\circ$} at 8 0
\put{$\circ$} at 10 0
\put{$\circ$} at 12 0
\put{$\circ$} at 4 2
\put{$\circ$} at 6 2
\put{$\circ$} at 8 2
\arr{1.6 0}{0.4 0}
\arr{3.6 0}{2.4 0}
\arr{5.6 0}{4.4 0}
\arr{7.6 0}{6.4 0}
\arr{9.6 0}{8.4 0}
\arr{11.6 0}{10.4 0}

\arr{5.6 2}{4.4 2}
\arr{7.6 2}{6.4 2}

\arr{4 1.6}{4 0.4}
\arr{6 1.6}{6 0.4}
\arr{8 1.6}{8 0.4}


\put{$\ssize 2'$} at  4 2.5
\put{$\ssize 3'$} at  6 2.5
\put{$\ssize 4'$} at  8 2.5

\put{$\ssize 0$} at 0 -.5
\put{$\ssize 1$} at 2 -.5
\put{$\ssize 2$} at 4 -.5
\put{$\ssize 3$} at 6 -.5
\put{$\ssize 4$} at 8 -.5
\put{$\ssize 5$} at 10 -.5
\put{$\ssize 6$} at 12 -.5

\setdots<2pt>
\plot 5.5 1.5  4.5 0.5 /
\plot 7.5 1.5  6.5 0.5 /
\plot 0 -0.7  0.1 -0.9  0.2 -1  0.3 -1.05  11.7 -1.05  11.8 -1  11.9 -0.9  12 -0.7 /
\endpicture}
$$
with relations indicated as usual:  Two commutativity relations and
one zero relation. 
The algebra $\Theta$ is tubular, as we are going to show.

\smallskip
Denote by $\Theta_0$ the algebra obtained from $\Theta$ by deleting the vertices $4'$ und $6$, this
is a tame concealed algebra of type $\widetilde E_7$. In order to see this,
one may use the Happel-Vossieck list (see for example [9], page~365), 
or one may define directly
a tilting functor by forming the pushout of the commutative square, thus
replacing the vertex $3'$ and the four arrows between $3'$ and $2$ 
by a new vertex $3''$ with arrows $2' \to 3''$, $3 \to 3''$, and
$3'' \to 2$. 

\smallskip
We denote by $\Cal P(\Theta_0)$ and $\Cal Q(\Theta_0)$ the preprojective and the 
preinjective component of the Auslander-Reiten quiver $\Gamma(\Theta_0)$ of $\Theta_0$, 
respectively.
The remaining indecomposable $\Theta_0$-modules occur in a $\Bbb P_1(k)$-family 
$\Cal T(\Theta_0)$ of tubes of tubular type $(4,3,2)$.
Almost all of the
tubes are homogeneous, but there are three large tubes of rank $2$, $3$,
and $4$, respectively. Here we picture the mouth for each of these three tubes. 
$$
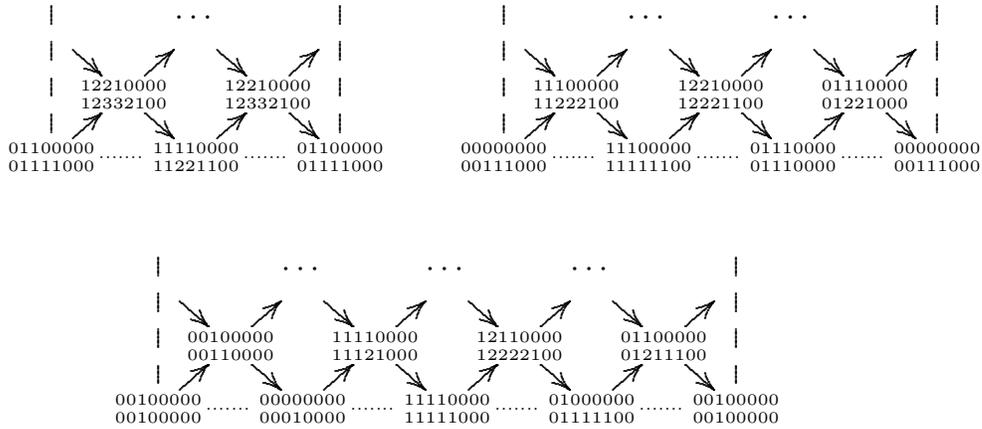
\hbox{\beginpicture
\setcoordinatesystem units <0.8cm,0.7cm>
\put{} at 0 0
\put{} at 4 2.5
\put{$\cdots$} at 2 2.2
\put{$\ssize {01100000    \atop 01111000}$} at 0 0
\put{$\ssize {12210000   \atop  12332100}$} at 1 1
\put{$\ssize {11110000   \atop  11221100}$} at 2 0
\put{$\ssize {12210000   \atop  12332100}$} at 3 1
\put{$\ssize {01100000    \atop 01111000}$} at 4 0

\arr{0.3 0.3} {0.7 0.7}
\arr{1.3 0.7} {1.7 0.3}
\arr{2.3 0.3} {2.7 0.7}
\arr{3.3 0.7} {3.7 0.3}
\arr{0.3 1.7} {0.7 1.3}
\arr{1.3 1.3} {1.7 1.7}
\arr{2.3 1.7} {2.7 1.3}
\arr{3.3 1.3} {3.7 1.7}

\setdots<2pt>
\plot 0.7 0  1.3 0 /
\plot 2.7 0  3.3 0 /

\setdashes <2mm>
\plot 0 0.4  0 2.5 /
\plot 4 0.4  4 2.5 /
\endpicture}
\qquad
\hbox{\beginpicture
\setcoordinatesystem units <0.8cm,0.7cm>
\put{} at 0 0
\put{} at 6 2.5
\put{$\cdots$} at 2 2.2
\put{$\cdots$} at 4 2.2
\put{$\ssize {00000000    \atop 00111000}$} at 0 0
\put{$\ssize {11100000 \atop    11222100}$} at 1 1
\put{$\ssize {11100000   \atop  11111100}$} at 2 0
\put{$\ssize {12210000   \atop  12221100}$} at 3 1
\put{$\ssize {01110000   \atop  01110000}$} at 4 0
\put{$\ssize {01110000   \atop  01221000}$} at 5 1
\put{$\ssize {00000000    \atop 00111000}$} at 6 0

\arr{0.3 0.3} {0.7 0.7}
\arr{1.3 0.7} {1.7 0.3}
\arr{2.3 0.3} {2.7 0.7}
\arr{3.3 0.7} {3.7 0.3}
\arr{4.3 0.3} {4.7 0.7}
\arr{5.3 0.7} {5.7 0.3}
\arr{0.3 1.7} {0.7 1.3}
\arr{1.3 1.3} {1.7 1.7}
\arr{2.3 1.7} {2.7 1.3}
\arr{3.3 1.3} {3.7 1.7}
\arr{4.3 1.7} {4.7 1.3}
\arr{5.3 1.3} {5.7 1.7}

\setdots<2pt>
\plot 0.7 0  1.3 0 /
\plot 2.7 0  3.3 0 /
\plot 4.7 0  5.3 0 /

\setdashes <2mm>
\plot 0 0.4  0 2.5 /
\plot 6 0.4  6 2.5 /
\endpicture}
$$
$$
\hbox{\beginpicture
\setcoordinatesystem units <0.8cm,0.7cm>
\put{} at 0 1
\put{} at 10 3.5
\put{$\cdots$} at 2 3.2
\put{$\cdots$} at 4 3.2
\put{$\cdots$} at 6 3.2
\put{$\ssize {00100000 \atop 00100000 }$} at 0 1
\put{$\ssize {00000000 \atop 00010000}$} at 2 1
\put{$\ssize {11110000 \atop 11111000}$} at 4 1
\put{$\ssize {01000000 \atop 01111100}$} at 6 1
\put{$\ssize {00100000 \atop 00100000}$} at 8 1

\put{$\ssize {00100000 \atop 00110000}$} at 1 2
\put{$\ssize {11110000 \atop 11121000}$} at 3 2
\put{$\ssize {12110000 \atop 12222100}$} at 5 2
\put{$\ssize {01100000 \atop 01211100}$} at 7 2

\arr{0.3 1.3} {0.7 1.7}
\arr{1.3 1.7} {1.7 1.3}
\arr{2.3 1.3} {2.7 1.7}
\arr{3.3 1.7} {3.7 1.3}
\arr{4.3 1.3} {4.7 1.7}
\arr{5.3 1.7} {5.7 1.3}
\arr{6.3 1.3} {6.7 1.7}
\arr{7.3 1.7} {7.7 1.3}

\arr{0.3 2.7} {0.7 2.3}
\arr{1.3 2.3} {1.7 2.7}
\arr{2.3 2.7} {2.7 2.3}
\arr{3.3 2.3} {3.7 2.7}
\arr{4.3 2.7} {4.7 2.3}
\arr{5.3 2.3} {5.7 2.7}
\arr{6.3 2.7} {6.7 2.3}
\arr{7.3 2.3} {7.7 2.7}
\setdots<2pt>
\plot 0.7 1  1.3 1 /
\plot 2.7 1  3.3 1 /
\plot 4.7 1  5.3 1 /
\plot 6.7 1  7.3 1 /

\setdashes <2mm>
\plot 0 1.4  0 3.5 /
\plot 8 1.4  8 3.5 /
\endpicture}
$$

Denote by $R$ the radical of the projective $\Theta$-module $P_\Theta(6)$
and by $R'$ the radical of the projective $\Theta$-module $P_\Theta(4')$, both considered as
$\Theta_0$-modules. These modules $R,R'$ are indecomposable and the dimension vector
of $R$ is $\ssize{ 010000 \atop 011111}$, that of $R'$ is 
$\ssize{ 11110 \atop 11111}$
(here we use already the convention mentioned above).
As both $R$ and $R'$ 
are modules on the mouth of the $4$-tube,
$\Theta$ is obtained from $\Theta_0$ as a two-fold one-point extension 
of a tame concealed algebra and hence is a tubular algebra of
tubular type $(6,3,2)$ ([9], Chapter~5).
The mouth of the big tube for $\Theta_0$ after the twofold ray insertion looks as
follows.
$$
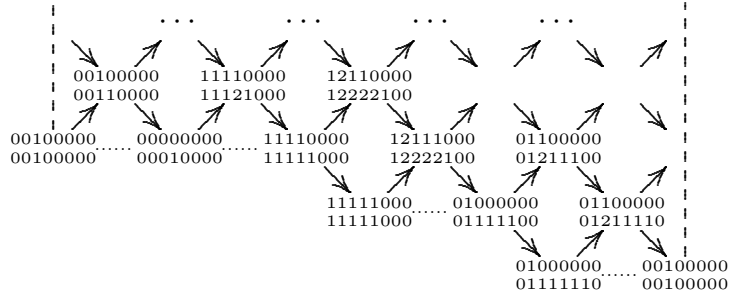
\hbox{\beginpicture
\setcoordinatesystem units <0.7cm,0.7cm>
\put{} at 0 0
\put{} at 10 4.3

\put{$\ssize {00100000 \atop 00100000}$} at 0 2
\put{$\ssize {00100000 \atop 00110000}$} at 1 3
\put{$\ssize {00000000 \atop 00010000}$} at 2 2
\put{$\ssize {11110000 \atop 11121000}$} at 3 3
\put{$\ssize {11110000 \atop 11111000}$} at 4 2
\put{$\ssize {12110000 \atop 12222100}$} at 5 3
\put{$\ssize {11111000 \atop 11111000}$} at 5 1
\put{$\ssize {12111000 \atop 12222100}$} at 6 2
\put{$\ssize {01000000 \atop 01111100}$} at 7 1
\put{$\ssize {01000000 \atop 01111110}$} at 8 0
\put{$\ssize {01100000 \atop 01211100}$} at 8 2
\put{$\ssize {00100000 \atop 00100000}$} at 10 0
\put{$\ssize {01100000 \atop 01211110}$} at 9 1

\arr{7.3 0.7} {7.7 0.3} 
\arr{8.3 0.3} {8.7 0.7} 
\arr{9.3 0.7} {9.7 0.3} 
\arr{4.3 1.7} {4.7 1.3} 
\arr{5.3 1.3} {5.7 1.7} 
\arr{6.3 1.7} {6.7 1.3} 
\arr{7.3 1.3} {7.7 1.7} 
\arr{8.3 1.7} {8.7 1.3} 
\arr{9.3 1.3} {9.7 1.7} 
\arr{0.3 2.3} {0.7 2.7} 
\arr{1.3 2.7} {1.7 2.3} 
\arr{2.3 2.3} {2.7 2.7} 
\arr{3.3 2.7} {3.7 2.3} 
\arr{4.3 2.3} {4.7 2.7} 
\arr{5.3 2.7} {5.7 2.3} 
\arr{6.3 2.3} {6.7 2.7} 
\arr{7.3 2.7} {7.7 2.3} 
\arr{8.3 2.3} {8.7 2.7} 
\arr{9.3 2.7} {9.7 2.3} 
\arr{0.3 3.7} {0.7 3.3} 
\arr{1.3 3.3} {1.7 3.7} 
\arr{2.3 3.7} {2.7 3.3} 
\arr{3.3 3.3} {3.7 3.7} 
\arr{4.3 3.7} {4.7 3.3} 
\arr{5.3 3.3} {5.7 3.7} 
\arr{6.3 3.7} {6.7 3.3} 
\arr{7.3 3.3} {7.7 3.7} 
\arr{8.3 3.7} {8.7 3.3} 
\arr{9.3 3.3} {9.7 3.7} 

\multiput{$\cdots$} at 2 4  4 4  6 4  8 4 /
\setdashes<3pt>
\plot 0 2.3    0 4.3 /
\plot 10 0.3   10 4.3 /
\setdots<2pt>
\plot 0.7 2  1.3 2 /
\plot 2.7 2  3.3 2 /
\plot 5.7 1  6.3 1 /
\plot 8.7 0  9.3 0 /
\endpicture}
$$

Since $\Theta$ is a tubular algebra, there is a second tame concealed algebra involved,
usually denoted by $\Theta_\infty$; 
it is obtained from $\Theta$ by deleting the vertex $0$ and
has type $\widetilde E_8$. The regular $\Theta_\infty$-modules form a 
$\Bbb P_1(k)$-family of tubes of tubular type $(5,3,2)$. 
The tubes of rank $2$ and $3$ are obtained from the corresponding tubes
in $\Cal T(\Theta_0)$ by replacing each representation $M$ by $M[1]$; here is 
the mouth of the tube of rank $5$.
$$
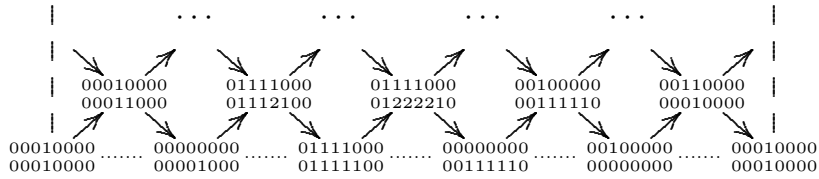
\hbox{\beginpicture
\setcoordinatesystem units <0.8cm,0.7cm>
\put{} at 0 1
\put{} at 10 3.5
\put{$\cdots$} at 2 3.2
\put{$\cdots$} at 4 3.2
\put{$\cdots$} at 6 3.2
\put{$\cdots$} at 8 3.2
\put{$\ssize {00010000 \atop 00010000}$} at 0 1
\put{$\ssize {00000000 \atop 00001000}$} at 2 1
\put{$\ssize {01111000 \atop 01111100}$} at 4 1
\put{$\ssize {00000000 \atop 00111110}$} at 6 1
\put{$\ssize {00100000 \atop 00000000}$} at 8 1
\put{$\ssize {00010000 \atop 00010000}$} at 10 1

\put{$\ssize {00010000  \atop 00011000}$} at 1 2
\put{$\ssize {01111000  \atop 01112100}$} at 3 2
\put{$\ssize {01111000  \atop 01222210}$} at 5 2
\put{$\ssize {00100000  \atop 00111110}$} at 7 2
\put{$\ssize {00110000  \atop 00010000}$} at 9 2

\arr{0.3 1.3} {0.7 1.7}
\arr{1.3 1.7} {1.7 1.3}
\arr{2.3 1.3} {2.7 1.7}
\arr{3.3 1.7} {3.7 1.3}
\arr{4.3 1.3} {4.7 1.7}
\arr{5.3 1.7} {5.7 1.3}
\arr{6.3 1.3} {6.7 1.7}
\arr{7.3 1.7} {7.7 1.3}
\arr{8.3 1.3} {8.7 1.7}
\arr{9.3 1.7} {9.7 1.3}

\arr{0.3 2.7} {0.7 2.3}
\arr{1.3 2.3} {1.7 2.7}
\arr{2.3 2.7} {2.7 2.3}
\arr{3.3 2.3} {3.7 2.7}
\arr{4.3 2.7} {4.7 2.3}
\arr{5.3 2.3} {5.7 2.7}
\arr{6.3 2.7} {6.7 2.3}
\arr{7.3 2.3} {7.7 2.7}
\arr{8.3 2.7} {8.7 2.3}
\arr{9.3 2.3} {9.7 2.7}
\setdots<2pt>
\plot 0.7 1  1.3 1 /
\plot 2.7 1  3.3 1 /
\plot 4.7 1  5.3 1 /
\plot 6.7 1  7.3 1 /
\plot 8.7 1  9.3 1 /

\setdashes <2mm>
\plot 0 1.4  0 3.5 /
\plot 10 1.4  10 3.5 /
\endpicture}
$$

\smallskip
The shape of the category of $\Theta$-modules is as follows:
$$
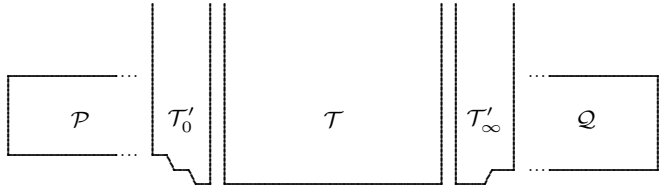
\hbox{\beginpicture
\setcoordinatesystem units <0.8cm,0.8cm>
\put{} at 0 0
\put{} at 9 2.5 
\plot 1.5 0.4  0 0.4  0 1.5  1.5 1.5 /
\setdots<2pt>
\plot 1.5 0.4  1.8 0.4 /
\plot 1.5 1.5  1.8 1.5 /
\setsolid
\plot 2 2.5  2 0.4  2.2 0.4  2.3 .2  2.5 .2  2.6 0  2.8 0  2.8 2.5  /
\plot 3 2.5  3 0  6 0  6 2.5 /
\plot 6.2 2.5  6.2 0  6.6 0  6.7 .2 7 .2  7 2.5 /
\plot 7.5 0.2  9 0.2  9 1.5  7.5 1.5 /
\setdots<2pt>
\plot 7.5 0.2  7.2 0.2 /
\plot 7.5 1.5  7.2 1.5 /

\setsolid

\put{$\ssize \Cal P$} at 1 .9
\put{$\ssize \Cal T_0'$} at 2.4 .9
\put{$\ssize \Cal T$} at 4.5 .9
\put{$\ssize \Cal T'_\infty$} at 6.6 .9
\put{$\ssize \Cal Q$} at 8 .9
\endpicture}
$$

\smallskip
Here, $\Cal P$ is $\Cal P(\Theta_0)$ and $\Cal T_0'$ is obtained from $\Cal T(\Theta_0)$
by the two-fold ray insertion in the $4$-tube at the modules $R$ and $R'$.
Next, $\Cal T$ is the union of countably many stable tubular families
$\Cal T_\gamma$ indexed by $\gamma\in \Bbb Q^+$, each $\Cal T_\gamma$ is a  
$\Bbb P_1(k)$-family of tubes of type $(6,3,2)$.
From the Auslander-Reiten quiver for $\Theta_\infty$ we obtain the 
preinjective component $\Cal Q=\Cal Q(\Theta_\infty)$ for $\Theta$
and the tubular family $\Cal T'_\infty$, which arises 
from the tubular family $\Cal T(\Theta_\infty)$ for $\Theta_\infty$
by the insertion of a coray in the $5$-tube, as pictured below.
Thus, there are two non-stable tubes in the category of $\Theta$-modules,
one in the family of tubes $\Cal T_0'$, the other one in $\Cal T'_\infty$.
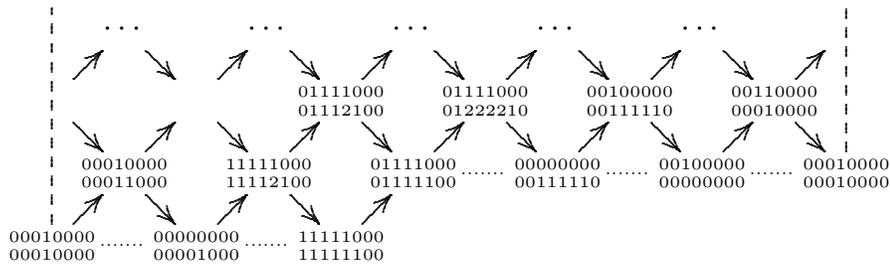
$$
\hbox{\beginpicture
\setcoordinatesystem units <0.8cm,0.8cm>
\put{} at 0 0
\put{} at 0 3

\put{$\ssize {00010000 \atop 00010000}$}   at 0 0
\put{$\ssize {00000000 \atop 00001000}$} at 2 0
\put{$\ssize {11111000 \atop 11111100}$} at 4 0
\put{$\ssize {00010000 \atop 00011000}$} at 1 1
\put{$\ssize {11111000 \atop 11112100}$} at 3 1
\put{$\ssize {01111000 \atop 01111100}$} at 5 1
\put{$\ssize {00000000 \atop 00111110}$} at 7 1
\put{$\ssize {00100000 \atop 00000000}$} at 9 1
\put{$\ssize {00010000 \atop 00010000}$} at 11 1
\put{$\ssize {01111000 \atop 01112100}$} at 4 2
\put{$\ssize {01111000 \atop 01222210}$} at 6 2
\put{$\ssize {00100000 \atop 00111110}$} at 8 2
\put{$\ssize {00110000 \atop 00010000}$} at 10 2

\arr{0.3 0.3} {0.7 0.7} 
\arr{1.3 0.7} {1.7 0.3} 
\arr{2.3 0.3} {2.7 0.7} 
\arr{3.3 0.7} {3.7 0.3} 
\arr{4.3 0.3} {4.7 0.7} 
\arr{0.3 1.7} {0.7 1.3} 
\arr{1.3 1.3} {1.7 1.7} 
\arr{2.3 1.7} {2.7 1.3} 
\arr{3.3 1.3} {3.7 1.7} 
\arr{4.3 1.7} {4.7 1.3} 
\arr{5.3 1.3} {5.7 1.7} 
\arr{6.3 1.7} {6.7 1.3} 
\arr{7.3 1.3} {7.7 1.7} 
\arr{8.3 1.7} {8.7 1.3} 
\arr{9.3 1.3} {9.7 1.7} 
\arr{10.3 1.7} {10.7 1.3}
\arr{0.3 2.3} {0.7 2.7} 
\arr{1.3 2.7} {1.7 2.3} 
\arr{2.3 2.3} {2.7 2.7} 
\arr{3.3 2.7} {3.7 2.3} 
\arr{4.3 2.3} {4.7 2.7} 
\arr{5.3 2.7} {5.7 2.3} 
\arr{6.3 2.3} {6.7 2.7} 
\arr{7.3 2.7} {7.7 2.3} 
\arr{8.3 2.3} {8.7 2.7} 
\arr{9.3 2.7} {9.7 2.3} 
\arr{10.3 2.3} {10.7 2.7} 

\multiput{$\cdots$} at 1 3  3 3  5 3  7 3  9 3 /
\setdashes<3pt>
\plot 0 0.3  0 3.3 /
\plot 11 1.3  11 3.3 /
\setdots<2pt>
\plot 0.7 0  1.3 0 /
\plot 2.7 0  3.3 0 /
\plot 5.7 1  6.3 1 /
\plot 7.7 1  8.3 1 /
\plot 9.7 1  10.3 1 /
\endpicture}
$$

In the sequel we will want to locate indecomposable modules within the
category of $\Theta$-modules, and for this the index functions $\iota_0$ and
$\iota_\infty$ for the tubular algebra $\Theta$ 
can be used, they map the dimension vector of a module
to its dot-product with with 
$\smallmatrix 0&0&1&1&0&0&0&0\\-1&-1&-1&0&1&1&0&0\endsmallmatrix$
and $\smallmatrix 0&0&0&1&1&0&0&0\\0&-1&-1&-1&0&1&1&0\endsmallmatrix$, 
respectively.
We recall from [9], Section 5.2, that an indecomposable $\Theta$-module with dimension
vector $\bd$ is

\smallskip
\centerline{
\vbox{\offinterlineskip
\halign{\strut#&#\hfil&\vrule#& \;\hfil#\hfil&\;\hfil#\hfil&\;\hfil#\hfil&
          \;\hfil#\hfil&\;\hfil#\hfil&\qquad\hfil#\hfil&\;\hfil#\hfil\cr
& in  &height14pt& $\Cal P$ & $\Cal T_0'$ & $\Cal T_\gamma,\gamma\in\Bbb Q^+$ 
    & $\Cal T'_\infty$ & $\Cal Q$ & $\Cal T''_\infty,\Cal U$ & $\Cal C''$ \cr
\noalign{\hrule}
& \omit  &height 14pt& $\iota_0<0$ & $\iota_0=0$ 
    & $\gamma=-\frac{\iota_0}{\iota_\infty}$
    & $\iota_0>0$ & $\iota_0\geq 0$ & $\iota_0>0$ & $\iota_\infty=\iota_0=0$\cr
& if    && and & and & and & and & and & and & or \cr
& \omit && $\iota_\infty\leq 0$ & $\iota_\infty<0$ & $\iota_\infty<0$ 
    & $\iota_\infty=0$ & $\iota_\infty>0$ & $\iota_\infty=0$ 
    & $\iota_\infty>0$ \cr
}}}

\smallskip
\noindent
where we use the abbreviations $\iota_0=\iota_0(\bd)$ and 
$\iota_\infty=\iota_\infty(\bd)$.  (We will come back to the right two 
columns of the table below.)
This finishes our description of the tubular algebra $\Theta$.

\medskip\noindent{\bf (1.2) The category $\Cal C$.}
We obtain $\Cal C$ from the category of $\Theta$-modules
by repeatedly forming one-point extensions.
Note that $\rad P(5')$ with dimension vector $\ssize {111110 \atop 111111}$
and $\rad P(7)$ with dimension vector $\ssize {0000000 \atop 0011111}$
both belong to the big tube in $\Cal T'_\infty,$ thus we can
form the corresponding ray extensions and obtain a component of the following 
form for the two-fold one-point extension algebra $\Theta'$ of the algebra~$\Theta$. 

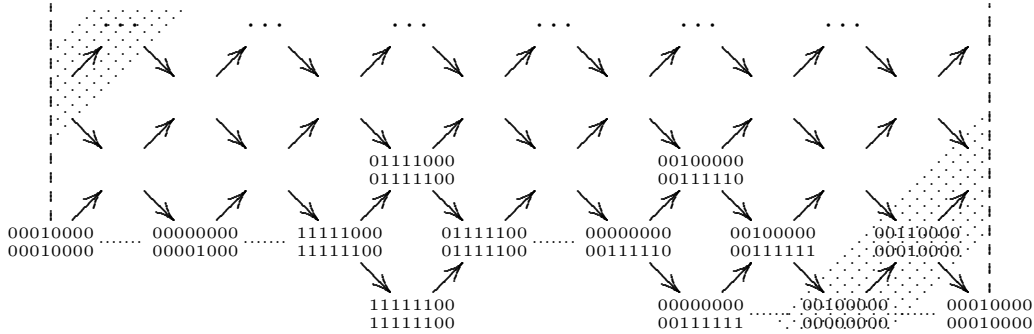
$$
\hbox{\beginpicture
\setcoordinatesystem units <0.8cm,0.8cm>
\put{} at 0 0
\put{} at 0 3

\put{$\ssize {00010000 \atop 00010000}$} at 0 0
\put{$\ssize {00000000 \atop 00001000}$} at 2 0
\put{$\ssize {11111000 \atop 11111100}$} at 4 0
\put{$\ssize {01111000 \atop 01111100}$} at 5 1
\put{$\ssize {11111100 \atop 11111100}$} at 5 -1
\put{$\ssize {01111100 \atop 01111100}$} at 6 0
\put{$\ssize {00000000 \atop 00111110}$} at 8 0
\put{$\ssize {00000000 \atop 00111111}$} at 9 -1

\put{$\ssize {00100000 \atop 00111110}$} at 9 1
\put{$\ssize {00100000 \atop 00111111}$} at 10 0

\put{$\ssize {00100000 \atop 00000000}$} at 11 -1

\put{$\ssize {00110000 \atop 00010000}$} at 12 0
\put{$\ssize {00010000 \atop 00010000}$} at 13 -1
\arr{4.3 -0.3} {4.7 -0.7} 
\arr{5.3 -0.7} {5.7 -0.3} 

\arr{8.3 -0.3} {8.7 -0.7} 
\arr{9.3 -0.7} {9.7 -0.3} 
\arr{10.3 -0.3} {10.7 -0.7}
\arr{11.3 -0.7} {11.7 -0.3} 
\arr{12.3 -0.3} {12.7 -0.7}

\arr{0.3 0.3} {0.7 0.7} 
\arr{1.3 0.7} {1.7 0.3} 
\arr{2.3 0.3} {2.7 0.7} 
\arr{3.3 0.7} {3.7 0.3} 
\arr{4.3 0.3} {4.7 0.7} 
\arr{5.3 0.7} {5.7 0.3} 
\arr{6.3 0.3} {6.7 0.7} 
\arr{7.3 0.7} {7.7 0.3} 
\arr{8.3 0.3} {8.7 0.7} 
\arr{9.3 0.7} {9.7 0.3} 
\arr{10.3 0.3} {10.7 0.7} 
\arr{11.3 0.7} {11.7 0.3} 
\arr{12.3 0.3} {12.7 0.7} 
\arr{0.3 1.7} {0.7 1.3} 
\arr{1.3 1.3} {1.7 1.7} 
\arr{2.3 1.7} {2.7 1.3} 
\arr{3.3 1.3} {3.7 1.7} 
\arr{4.3 1.7} {4.7 1.3} 
\arr{5.3 1.3} {5.7 1.7} 
\arr{6.3 1.7} {6.7 1.3} 
\arr{7.3 1.3} {7.7 1.7} 
\arr{8.3 1.7} {8.7 1.3} 
\arr{9.3 1.3} {9.7 1.7} 
\arr{10.3 1.7} {10.7 1.3}
\arr{11.3 1.3} {11.7 1.7} 
\arr{12.3 1.7} {12.7 1.3}
\arr{0.3 2.3} {0.7 2.7} 
\arr{1.3 2.7} {1.7 2.3} 
\arr{2.3 2.3} {2.7 2.7} 
\arr{3.3 2.7} {3.7 2.3} 
\arr{4.3 2.3} {4.7 2.7} 
\arr{5.3 2.7} {5.7 2.3} 
\arr{6.3 2.3} {6.7 2.7} 
\arr{7.3 2.7} {7.7 2.3} 
\arr{8.3 2.3} {8.7 2.7} 
\arr{9.3 2.7} {9.7 2.3} 
\arr{10.3 2.3} {10.7 2.7} 
\arr{11.3 2.7} {11.7 2.3} 
\arr{12.3 2.3} {12.7 2.7}

\multiput{$\cdots$} at 1 3  3 3  5 3  7 3  9 3  11 3  / 
\setdashes<3pt>
\plot 0 0.3    0 3.3 /
\plot 13 -0.7  13 3.3 /
\setdots<2pt>
\plot 0.7 0  1.3 0 /
\plot 2.7 0  3.3 0 /
\plot 6.7 0  7.3 0 /
\plot 9.7 -1  10.3 -1 /
\plot 11.7 -1  12.3 -1 /
\setshadegrid span <0.8mm>
\vshade 10  -1.3 -1 <,z,,> 11.7 -1.3 0.7 <z,,,> 13  0.2 2 /
\vshade 0 1.2 2.8 <,z,,> 0.5 1.7 3.3 <z,,,> 2.1 3.3 3.3 /
\endpicture}
$$
(The shaded area will be considered later.)
We denote by $\Cal T''_\infty$ the one-parameter family of tubes
obtained from $\Cal T'_\infty$
by forming these two ray extensions.

\smallskip
We claim that the 
category $\Cal C$ has the following shape
$$
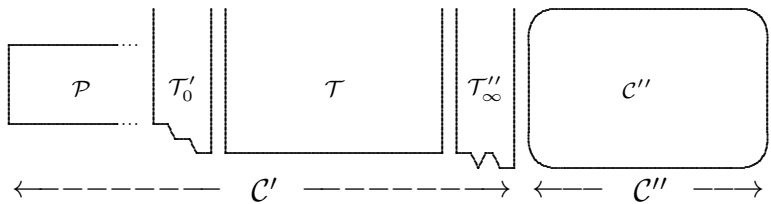
\hbox{\beginpicture
\setcoordinatesystem units <0.8cm,0.8cm>
\put{} at 0 0
\put{} at 9 2.5 
\betweenarrows{$\Cal C'$} from 0 -0.5 to 7 -0.5  
\betweenarrows{$\Cal C''$} from 7.2 -0.5 to 10.5 -0.5  
\plot 1.5 0.4  0 0.4  0 1.5  1.5 1.5 /
\setdots<2pt>
\plot 1.5 0.4  1.8 0.4 /
\plot 1.5 1.5  1.8 1.5 /
\setsolid
\plot 2 2  2 0.4  2.2  0.4  2.3 0.2  2.5 0.2  2.6 0  2.8 0  2.8 2 /
\plot 3 2  3 0  6 0  6 2 /
\plot 6.2 2  6.2 0  6.4 0  6.5 -0.2  6.6 0  6.7 0  6.8 -0.2  7 -0.2  7 2 /
\plot 10 2  7.7 2 /
\plot 10 -.2  7.7 -.2 /
\plot 10.5 0.3  10.5 1.5 /
\plot 7.2 0.3     7.2 1.5 /
\setquadratic
\plot 10 2 10.4 1.9 10.5 1.5 /
\plot 10 -.2 10.4 -.1 10.5 0.3 /
\plot 7.2 1.5  7.3 1.9  7.7 2 /
\plot 7.2 .3  7.3 -.1  7.7 -.2 /

\put{$\ssize \Cal P$} at 1 .9
\put{$\ssize \Cal T_0'$} at 2.4 .9
\put{$\ssize \Cal T$} at 4.5 .9
\put{$\ssize \Cal T''_\infty$} at 6.6 .9
\put{$\ssize \Cal C''$} at 8.7 .9

\endpicture}
$$
where $(\Cal C',\Cal C'')$ is a split torsion pair.
{\it Proof:}\/
For $t \ge 6$, we consider the full subcategory $\Cal C_t$ of $\Cal C$ formed by all
modules with composition factors $S(r),S(q')$ where $0 \le r \le t$ and 
$0 \le q \le t-2.$ Of course, $\Cal C_6$ is just the category of all $\Theta$-modules, 
and $\Cal C_7$ is the category of all $\Theta'$-modules, where $\Theta'$ is obtained
from $\Theta$ by forming the one-point extensions with respect to the two modules
$\rad P(5')$ and $\rad P(7)$, as considered above. In general, $\Cal C_{t+1}$
is obtained from $\Cal C_t$ by forming the one-point extensions with respect to the 
two modules $\rad P((t-1)')$ and $\rad P(t+1)$. Inductively, we claim that
$\Cal C_t$ has a split torsion pair $(\Cal C_t',\Cal C_t'')$ with $\Cal C_t' = \Cal C'$
for $t\ge 7.$ For $t=7$, this follows from the fact that $\Theta'$ is obtained from $\Theta$
using two ray insertions in a tube in $\Cal T'_\infty$. 
It remains to be noted that for $t\geq 7$,
the modules $\rad P((t-1)')$ and $\rad P(t+1)$ belong to $\Cal C_t''.$ 
Indeed, an indecomposable module $M$ in $\Cal C_{t+1}\backslash\Cal C_t$
is given --- as a module over the two-fold one-point extension --- by a 
homomorphism from a non-zero sum of copies of  
$\rad P((t-1)')$ or $\rad P(t+1)$
to a module in $\Cal C_t$, which must be in $\Cal C_t''$ by induction.
So $M$ is in $\Cal C_{t+1}''$. 

\smallskip
As a byproduct of this result we can verify the information in the table
about the extension of the index functions $\iota_0$ and $\iota_\infty$
to the indecomposables in $\Cal C$.
First we observe that for each ray-inserted module $M$ in $\Cal T_\infty''$
both conditions $\iota_0(\bdim M) > 0$ and 
$\iota_\infty(\bdim M)=0$ are satisfied.
Second, for an indecomposable module $M\in\Cal C''$ which is a module over
some iterated one-point extension of $\Theta$, the restriction
$\Res_{Q_\Theta}(M)$ of $M$ to $Q_\Theta$ 
is a (possibly empty) sum of modules in~$\Cal Q$.

\medskip\noindent{\bf (1.3) The category $\Cal C\cap\Cal S(\widetilde 6)$.}
We start with the observation that a module $M\in\Cal C$ which is
not in $\Cal S(\widetilde 6)$ has the property that
$\Hom_{\Cal C}(S(i'),M)\neq 0$ for at least one of the simple modules
$S(i')$ where $i\geq 2$ is a natural number.
Using the index functions we can locate the simple module $S(2')$ 
in the tubular family $\Cal T''_\infty$, and all the 
simple modules $S(t')$ with $t\ge 3$ in the category $\Cal C''$. 
Thus, in $\Cal C'$ there is a single ray
which contains modules not belonging to 
$\Cal S(\widetilde6)$, namely the ray in $\Cal T''_\infty$ starting at 
$S(2')$, as shaded above. If we delete this ray, we obtain the following tube:
$$
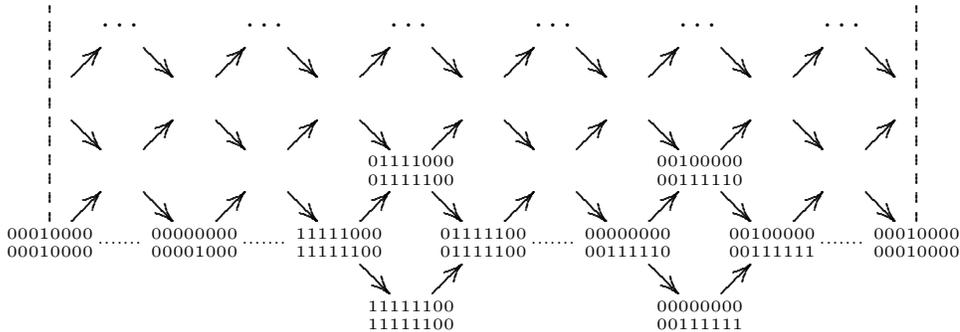
\hbox{\beginpicture
\setcoordinatesystem units <0.8cm,0.8cm>
\put{} at 0 0
\put{} at 0 3

\put{$\ssize {00010000 \atop 00010000}$} at 0 0
\put{$\ssize {00000000 \atop 00001000}$} at 2 0
\put{$\ssize {11111000 \atop 11111100}$} at 4 0
\put{$\ssize {01111000 \atop 01111100}$} at 5 1
\put{$\ssize {11111100 \atop 11111100}$} at 5 -1
\put{$\ssize {01111100 \atop 01111100}$} at 6 0
\put{$\ssize {00000000 \atop 00111110}$} at 8 0
\put{$\ssize {00000000 \atop 00111111}$} at 9 -1
\put{$\ssize {00100000 \atop 00111110}$} at 9 1
\put{$\ssize {00100000 \atop 00111111}$} at 10 0
\put{$\ssize {00010000 \atop 00010000}$} at 12 0
\arr{4.3 -0.3} {4.7 -0.7} 
\arr{5.3 -0.7} {5.7 -0.3} 

\arr{8.3 -0.3} {8.7 -0.7} 
\arr{9.3 -0.7} {9.7 -0.3} 

\arr{0.3 0.3} {0.7 0.7} 
\arr{1.3 0.7} {1.7 0.3} 
\arr{2.3 0.3} {2.7 0.7} 
\arr{3.3 0.7} {3.7 0.3} 
\arr{4.3 0.3} {4.7 0.7} 
\arr{5.3 0.7} {5.7 0.3} 
\arr{6.3 0.3} {6.7 0.7} 
\arr{7.3 0.7} {7.7 0.3} 
\arr{8.3 0.3} {8.7 0.7} 
\arr{9.3 0.7} {9.7 0.3} 
\arr{10.3 0.3} {10.7 0.7} 
\arr{11.3 0.7} {11.7 0.3} 
\arr{0.3 1.7} {0.7 1.3} 
\arr{1.3 1.3} {1.7 1.7} 
\arr{2.3 1.7} {2.7 1.3} 
\arr{3.3 1.3} {3.7 1.7} 
\arr{4.3 1.7} {4.7 1.3} 
\arr{5.3 1.3} {5.7 1.7} 
\arr{6.3 1.7} {6.7 1.3} 
\arr{7.3 1.3} {7.7 1.7} 
\arr{8.3 1.7} {8.7 1.3} 
\arr{9.3 1.3} {9.7 1.7} 
\arr{10.3 1.7} {10.7 1.3}
\arr{11.3 1.3} {11.7 1.7} 
\arr{0.3 2.3} {0.7 2.7} 
\arr{1.3 2.7} {1.7 2.3} 
\arr{2.3 2.3} {2.7 2.7} 
\arr{3.3 2.7} {3.7 2.3} 
\arr{4.3 2.3} {4.7 2.7} 
\arr{5.3 2.7} {5.7 2.3} 
\arr{6.3 2.3} {6.7 2.7} 
\arr{7.3 2.7} {7.7 2.3} 
\arr{8.3 2.3} {8.7 2.7} 
\arr{9.3 2.7} {9.7 2.3} 
\arr{10.3 2.3} {10.7 2.7} 
\arr{11.3 2.7} {11.7 2.3} 

\multiput{$\cdots$} at 1 3  3 3  5 3  7 3  9 3  11 3 /
\setdashes<3pt>
\plot 0 0.3    0 3.3 /
\plot 12 0.3  12 3.3 /
\setdots<2pt>
\plot 0.7 0  1.3 0 /
\plot 2.7 0  3.3 0 /
\plot 6.7 0  7.3 0 /
\plot 10.7 0  11.3 0 /
\endpicture}
$$

\smallskip
After the ray-deletion,
the arrows in this tube represent irreducible maps in the category
$\Cal S(\widetilde6)\cap \Cal C$. Thus the tube becomes an
Auslander-Reiten component for $\Cal S(\widetilde6)\cap\Cal C$:

\smallskip
The functor $\Mono$ takes an object $B=(B'\sto b B'')$ in $\Cal H(\widetilde n)$
to the embedding $\Mono(B)=(\Im b\lto{\incl} B'')$ in $\Cal S(\widetilde n)$.
Then the canonical map $\can\:B\to \Mono(B)$  is a {\it minimal left}
$\Cal S(\widetilde n)$-{\it approximation} for $B$, so $\can$
is left minimal and each test map
$t\:B\to C$ with $C\in\Cal S(\widetilde n)$ factors over it. 
In the above tube in $\Cal C\cap\Cal S(\widetilde 6)$, 
the arrows leaving the ray which is to be deleted represent
maps of type $B\to \Mono(B)$.

Let $0\to A\to B\to C\to 0$ be an Auslander-Reiten sequence in 
$\Cal C$ which starts at an object $A\in\Cal C\cap \Cal S(\widetilde 6)$ that
is not relatively injective. We obtain the
relative Auslander-Reiten sequence in $\Cal C\cap\Cal S(\widetilde 6)$ 
as the lower sequence in the following diagram (cf.~[11], Proposition~3.2).
$$\CD
\ssize 0 @>>> \ssize A  @>>> \ssize B @>>> \ssize C @>>> \ssize 0\cr
   @.              @|              @VV{\can}V       @VV{\can}V        @.\cr
\ssize0 @>>> \ssize A @>>> \ssize \Mono(B) @>>> \ssize \Mono(C) @>>> \ssize 0\cr
\endCD$$
It follows that the ray deletion yields a component in the Auslander-Reiten
quiver for $\Cal C\cap \Cal S(\widetilde 6)$.

\medskip
The two projective modules in this tube,
$P(7)$ and $P(5')$, are projective and injective objects in the category 
$\Cal S(\widetilde 6)$, as we recall from [11], Proposition~1.4.

\medskip\noindent
{\bf (1.3.1) Lemma.} {\it
Let $\Lambda$ be a commutative uniserial ring and $P$ a
projective-injective $\Lambda$-module.
The category $\Cal S(\Lambda)$ has the following projective
or injective objects:

\smallskip\item{\rm (1)} The object $(P,0)$ is projective-injective 
in $\Cal S(\Lambda)$ and has the following source map and sink map
$$(P,0)\lto{\incl} (P,\soc P)\T{and}(\rad P,0)\lto{\incl} (P,0).$$

\item{\rm (2)} The object $(P,P)$ is projective-injective in
$\Cal S(\Lambda)$ and has the following source map and sink map
$$(P,P)\lto{\can} (P/\soc P,P/\soc P)\T{and}(P,\rad P)\lto{\incl} (P,P).$$

}

\smallskip
Denote by $\Cal U$ the intersection of 
$\Cal T''_\infty$ with $\Cal S(\widetilde 6)$, i.e.\ $\Cal U$ is
obtained from $\Cal T''_\infty$ by deleting one ray in the tube of
rank 6, as illustrated above.  
Then the category $\Cal C\cap \Cal S(\widetilde6)$ has the following shape:
$$
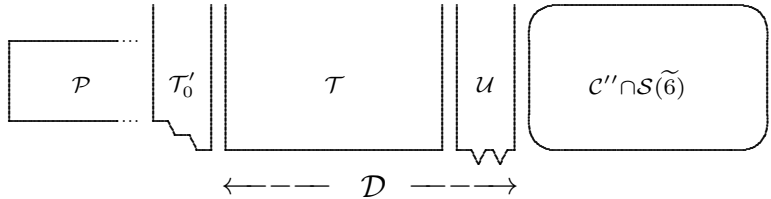
\hbox{\beginpicture
\setcoordinatesystem units <0.8cm,0.8cm>
\put{} at 0 0
\put{} at 9 2.5 
\betweenarrows{$\Cal D$} from 2.9 -0.5 to 7.1 -0.5  
\plot 1.5 0.4  0 0.4  0 1.5  1.5 1.5 /
\setdots<2pt>
\plot 1.5 0.4  1.8 0.4 /
\plot 1.5 1.5  1.8 1.5 /
\setsolid
\plot 2 2  2 0.4  2.2 .4  2.3 .2  2.5 .2  2.6 0  2.8 0  2.8 2 /
\plot 3 2  3 0  6 0  6 2 /
\plot 6.2 2  6.2 0  6.4 0  6.5 -0.2  6.6 0  6.7 0  6.8 -0.2  6.9 0  7 0  7 2 / 

\plot 10 2  7.7 2 /
\plot 10 0  7.7 0 /
\plot 10.5 0.5  10.5 1.5 /
\plot 7.2 0.5     7.2 1.5 /
\setquadratic
\plot 10 2 10.4 1.9 10.5 1.5 /
\plot 10 0 10.4 0.1 10.5 0.5 /
\plot 7.2 1.5  7.3 1.9  7.7 2 /
\plot 7.2 0.5  7.3 0.1  7.7 0 /

\put{$\ssize \Cal P$} at 1 .9
\put{$\ssize \Cal T_0'$} at 2.4 .9
\put{$\ssize \Cal T$} at 4.5 .9
\put{$\ssize \Cal U$} at 6.6 .9
\put{$\ssize \Cal C''\cap \Cal S(\widetilde6)$} at 8.7 .9
\endpicture}
$$

We will see in the next section that 
$\Cal D = \Cal T \sqcup \Cal U$ is a fundamental domain for $\Cal S(\widetilde 6)$
under the shift. For further reference, let us introduce the following notation: Denote by
$\Cal D'$ the objects in $\Cal D$ which are $\Theta$-modules.
Let $\Cal D''$ be the class of the remaining indecomposable objects in $\Cal D$;
they form two rays in the non-stable tube:
one ray consists of all the indecomposable objects $M$ in $\Cal D$ with $M_{5'} \neq 0,$
the other consists of all the indecomposable objects $M$ in $\Cal D$ with $M_7 \neq 0$.
The dimension vectors of the objects in $\Cal D''$ will be discussed in (1.6). 

\bigskip\noindent{\bf (1.4) 
The fundamental domain for $\Cal S(\widetilde 6)$ under the shift.}
It is the part $\Cal D = \Cal T \sqcup \Cal U$ 
which is of interest for us, 
since it will turn out that $\Cal D$ is a fundamental domain 
for the shift on $\widetilde Q$ and 
that the structure of $\Cal S(\widetilde 6)$ is as follows:
$$
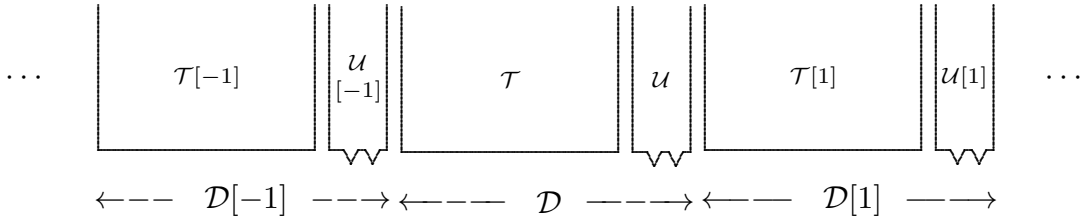
\hbox{\beginpicture
\setcoordinatesystem units <0.8cm,0.8cm>
\put{\beginpicture
\setcoordinatesystem units <0.8cm,0.8cm>
\put{} at 1.5 0
\put{} at 9 2.5 
\put{$\cdots$} at 2 1
\betweenarrows{$\Cal D[-1]$} from 2.9 -0.7 to 7.1 -0.7

\plot 3 2  3 0  6 0  6 2 /
\plot 7 2  7 0 /
\plot 6.2 2  6.2 0 /
\plot 6.2 0  6.4 0 6.5 -0.2 6.6 0 6.7 0 6.8 -0.2 6.9 0  7 0 /

\put{$\ssize \Cal T[-1]$} at 4.5 1
\put{$\ssize \Cal U$} at 6.6 1.2
\put{$\ssize [-1]$} at 6.6 .8
\endpicture} at 0 0 
\put{\beginpicture
\setcoordinatesystem units <0.8cm,0.8cm>
\put{} at 1.5 0
\put{} at 9 2.5 
\betweenarrows{$\Cal D$} from 2.9 -0.7 to 7.1 -0.7  

\plot 3 2  3 0  6 0  6 2 /
\plot 7 2  7 0 /
\plot 6.2 2  6.2 0 /
\plot 6.2 0  6.4 0 6.5 -0.2 6.6 0 6.7 0 6.8 -0.2 6.9 0  7 0 /

\put{$\ssize \Cal T$} at 4.5 1
\put{$\ssize \Cal U$} at 6.6 1
\endpicture} at 4.2 0 
\put{\beginpicture
\setcoordinatesystem units <0.8cm,0.8cm>
\put{} at 1.5 0
\put{} at 9 2.5 
\put{$\cdots$} at 8 1
\betweenarrows{$\Cal D[1]$} from 2.9 -0.7 to 7.1 -0.7  

\plot 3 2  3 0  6 0  6 2 /
\plot 7 2  7 0 /
\plot 6.2 2  6.2 0 /
\plot 6.2 0  6.4 0 6.5 -0.2 6.6 0 6.7 0 6.8 -0.2 6.9 0  7 0 /

\put{$\ssize \Cal T[1]$} at 4.5 1
\put{$\ssize \Cal U[1]$} at 6.6 1
\endpicture} at 8.4 0 
\endpicture}
$$

\medskip\noindent
{\bf (1.4.1) Theorem.} {\it 
For each indecomposable representation $M$ in $\Cal S(\widetilde6)$ there is a 
unique integer $n$ such that $M[n]$ is in $\Cal D$. }

\smallskip\noindent
{\it Proof:}\/ Let $M\in\ind \Cal S(\widetilde 6)$ 
be such that $M_i=0$ for all $i<0$
and that $M_0\neq 0$. The existence of $n$ as in the 
Theorem is a consequence of the following five claims.
In the first three claims we show that either $M$, $M[1]$, or $M[2]$ 
occurs in 
$\Cal C'\cap \Cal S(\widetilde6)
        =\Cal P\sqcup \Cal T_0' \sqcup \Cal T\sqcup \Cal U$; 
in the last two claims we verify that each module in 
$\Cal P\sqcup \Cal T_0'$ has a translate in $\Cal D=\Cal T\sqcup \Cal U$.
We conclude the proof of the Theorem by verifying that the number $n$
is indeed unique.

\medskip\noindent
{\bf Claim 1.} {\it
If $M$ has the property that $\beta_0$ is not surjective, then $M[2]$
is in $\Cal C'$. }

\smallskip\noindent
{\it Proof: }\/
For $M$ such that $\beta_0$ is not surjective, consider the subrepresentation
$N=(\sum M_{i'})$  of $M$ generated by the elements in 
$\sum_{i\in\Bbb Z}M_{i'}$.
The factor $M/N$ has the following properties: 
For each $i\in\Bbb Z$, $(M/N)_{i'}=0$
and $(M/N)_i=M_i/\Im(\beta_i)$.
In particular, $(M/N)_0$ is non-zero.
Viewing $M/N$ as a representation for the full subquiver of $\widetilde Q$
consisting of all points $\{i:i\in\Bbb Z\}$,
we see that there is a non-zero map into the injective representation $P(5)$
for this subquiver.  
The representation $(M/N)[2]$ has the property that $\beta_0$ and $\beta_1$
are identity maps, and hence the composition $M[2]\to(M/N)[2]\to P(7)$
is a non-zero map in $\Cal C$. Since $P(7)$ is in $\Cal C'$, so is $M[2]$.
\qed

\medskip\noindent
{\bf Claim 2.} {\it
If $M$ is such that $\beta_0$ is an isomorphism and $\beta_1$ is not
surjective, then $M[1]$ is in $\Cal C'$. }

\smallskip 
This assertion can be shown in a similar way. \qed

\medskip\noindent
{\bf Claim 3.} {\it
If $M$ is such that both $\beta_0$ and $\beta_1$ are isomorphisms
then $M$ is in $\Cal C'$.  }

\smallskip\noindent
{\it Proof: }\/ If both $\beta_0$ and $\beta_1$ are 
isomorphisms, then $M$ is in $\Cal C$.
Since $M_0\neq 0$ there is a non-zero map $M\to I(0)=P(5')$
to the injective envelope of $S(0)$,
so $M$ is in $\Cal C'$.   \qed

\medskip\noindent
{\bf Claim 4.} {\it
The whole tubular family $\Cal T_0'$ belongs to $\Cal U[-1]$. }

\smallskip\noindent
As a consequence, we can use the notation $\Cal T_0=\Cal U[-1]$.

\smallskip\noindent
{\it Proof:} \/
We denote by $\bh^0=\ssize{122100\atop 123321}$
the radical generator in the Grothendieck group 
$K_0(\mod \Theta_0)$ (we use here an upper index, since one of the vertices of the
quiver is denoted by $0$, thus there could arise some confusion).
The representations in $\Cal T_0'$ with dimension vector a multiple 
of $\bh^0$ are in 1-1 correspondence via the shift $M\mapsto M[1]$
to the representations in $\Cal U$ of dimension type
$\bh^\infty=\ssize{0122100\atop 0123321}$.
(Each module $M$ in $\Cal U$ with dimension vector a multiple of $\bh^\infty$ 
has an isomorphism at $\beta_2$
since this map must be monic because the only modules in 
$\Cal C'$ with  dimension vector a multiple of $\bh^\infty$ 
which do not have a monic map at $\beta_2$ occur in
the ray starting at $S(2')$, which we have deleted.)
Moreover, each module $M$ on the mouth of one of the three large tubes
in $\Cal T_0'$ 
has its translate $M[1]$ in the corresponding big tube in $\Cal U$,
see the diagram above. 
Since each tube is path closed, we conclude that $\Cal T_0'$ is contained
in $\Cal U$.
\qed

\medskip\noindent
{\bf Claim 5.} {\it
Each representation in $\Cal P$ has a translate in $\Cal D$. }

\smallskip\noindent
{\it Proof:} 
First, let us specify the shape of $\Cal P$ including
some of the dimension vectors 
(as we have mentioned before, we always specify the 
numbers $\bd_s$ and $\bd_{s'}$ with $0\le s \le t$ for a suitable $t$). 
We indicate also to which subcategories
$\Cal U[t]$ or $\Cal T[t]$ the modules belong:
$$
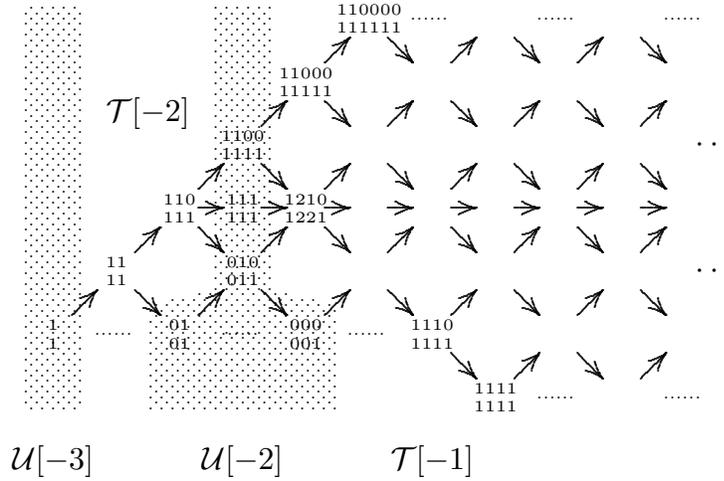
\hbox{\beginpicture
\setcoordinatesystem units <0.7cm,0.7cm>
\put{} at 0 0
\put{} at 7 6
\put{$\ssize {1\atop 1}$}   at 0 1
\put{$\ssize {11\atop 11}$} at 1 2
\put{$\ssize {01\atop 01}$} at 2 1
\put{$\ssize {110\atop 111}$} at 2 3
\put{$\ssize {010\atop 011}$} at 3 2
\put{$\ssize {111\atop 111}$} at 3 3
\put{$\ssize {1100 \atop 1111}$} at 3 4
\put{$\ssize {000 \atop 001}$} at 4 1
\put{$\ssize {1210 \atop 1221}$} at 4 3
\put{$\ssize {11000 \atop 11111}$} at 4 5
\put{$\ssize {110000 \atop 111111}$} at  5 6
\put{$\ssize {1110 \atop 1111}$} at  6 1
\put{$\ssize {1111 \atop 1111}$} at  7 0
\arr{6.3 0.7} {6.7 0.3} 
\arr{7.3 0.3} {7.7 0.7} 
\arr{8.3 0.7} {8.7 0.3} 
\arr{9.3 0.3} {9.7 0.7} 
\arr{0.3 1.3} {0.7 1.7} 
\arr{1.3 1.7} {1.7 1.3} 
\arr{2.3 1.3} {2.7 1.7} 
\arr{3.3 1.7} {3.7 1.3} 
\arr{4.3 1.3} {4.7 1.7} 
\arr{5.3 1.7} {5.7 1.3} 
\arr{6.3 1.3} {6.7 1.7} 
\arr{7.3 1.7} {7.7 1.3} 
\arr{8.3 1.3} {8.7 1.7} 
\arr{9.3 1.7} {9.7 1.3}
\arr{1.3 2.3} {1.7 2.7} 
\arr{2.3 2.7} {2.7 2.3} 
\arr{3.3 2.3} {3.7 2.7} 
\arr{4.3 2.7} {4.7 2.3} 
\arr{5.3 2.3} {5.7 2.7} 
\arr{6.3 2.7} {6.7 2.3} 
\arr{7.3 2.3} {7.7 2.7} 
\arr{8.3 2.7} {8.7 2.3} 
\arr{9.3 2.3} {9.7 2.7} 
\arr{2.3 3.3} {2.7 3.7} 
\arr{3.3 3.7} {3.7 3.3} 
\arr{4.3 3.3} {4.7 3.7} 
\arr{5.3 3.7} {5.7 3.3} 
\arr{6.3 3.3} {6.7 3.7} 
\arr{7.3 3.7} {7.7 3.3} 
\arr{8.3 3.3} {8.7 3.7} 
\arr{9.3 3.7} {9.7 3.3}
\arr{3.3 4.3} {3.7 4.7} 
\arr{4.3 4.7} {4.7 4.3} 
\arr{5.3 4.3} {5.7 4.7} 
\arr{6.3 4.7} {6.7 4.3} 
\arr{7.3 4.3} {7.7 4.7} 
\arr{8.3 4.7} {8.7 4.3} 
\arr{9.3 4.3} {9.7 4.7} 
\arr{4.3 5.3} {4.7 5.7} 
\arr{5.3 5.7} {5.7 5.3} 
\arr{6.3 5.3} {6.7 5.7} 
\arr{7.3 5.7} {7.7 5.3} 
\arr{8.3 5.3} {8.7 5.7} 
\arr{9.3 5.7} {9.7 5.3}
\arr{2.3 3}{2.7 3}
\arr{3.3 3}{3.7 3}
\arr{4.3 3}{4.7 3}
\arr{5.3 3}{5.7 3}
\arr{6.3 3}{6.7 3}
\arr{7.3 3}{7.7 3}
\arr{8.3 3}{8.7 3}
\arr{9.3 3}{9.7 3}
\setdots<2pt>
\plot 0.7 1  1.3 1 /
\plot 2.7 1  3.3 1 /
\plot 4.7 1  5.3 1 /
\plot 7.7 0  8.3 0 /
\plot 9.7 0  10.3 0 /
\plot 5.7 6  6.3 6 /
\plot 7.7 6  8.3 6 /
\plot 9.7 6  10.3 6 /
\setshadegrid span <0.6mm>
\hshade -0.2  1.5 4.5 <,,,z>  1.5 1.5 4.5 <,,z,z> 1.7 2.5 3.5 <,,z,> 6.2 2.5 3.5 /
\hshade -0.2  -0.5 0.5  6.2 -0.5 0.5 /
\put{$\cdots$} at 10.5 2 
\put{$\cdots$} at 10.5 4 
\put{$\Cal U[-3]$} at 0 -1
\put{$\Cal T[-2]$} at 1.5 4.5
\put{$\Cal U[-2]$} at 3 -1
\put{$\Cal T[-1]$} at 6 -1
\endpicture} 
$$
For each module $M\in\Cal P$ the index functions 
$\iota_0$ and $\iota_\infty$ can be used to specify
the tubular family of $M$.
Here we show that most of the modules in $\Cal P$ occur in $\Cal T[-1]$. 
Let $M\in\Cal P$ be such that there is a path 
$X\to\cdots\to M$ where $X\in\Cal P$
is one of the modules with
dimension vector $\ssize {11000\atop 11111}$
or $\ssize {1210\atop 1221}$, and hence 
$X\in\Cal T_{\frac12}[-1]$ or $X\in\Cal T_{\frac14}[-1]$,
respectively, and also a path
from $M$ to one of the modules $Y$ on the mouth of some tube in $\Cal T_0'$.
So we obtain a path (i.e.\ a sequence of non-zero homomorphisms
between indecomposable modules)
$$X[1]\to\cdots\to M[1]\to\cdots\to Y[1]$$
with $X[1]$ in $\Cal T_{\frac 12}$ or in $\Cal T_{\frac 14}$ and
$Y[1]$ in $\Cal U$; since $\Cal D$ is path closed,
$M$ is in $\Cal T[-1]$. \qed

\smallskip
In order to finish the proof of Theorem (1.4.1) we verify that 
at most one translate of an indecomposable module $M$ occurs in $\Cal D$.
Assume that an indecomposable object $M$ is such that both $M$ and
$M[i]$ occur in $\Cal D$ for some $i>0$.  Then there is a path 
$$X\to\cdots \to M\to\cdots\to M[i]\to\cdots\to Y$$ where $X\in \Cal T_0'$
and $Y\in \Cal U$.  Hence there is a path 
$$X[1]\to \cdots \to M[1] \to \cdots \to M[i]\to \cdots \to Y$$ 
where $X[1]\in\Cal U$ by Claim~4. Since $\Cal U$ is path closed, $M[1]$ and
$M[i]$ are in $\Cal U$.  A quick glance back at the index functions
shows that $\iota_0(M)=\iota_\infty (M[1])=0$, contradicting our
assumption that $M$ is in $\Cal D$. 
\qed

\medskip\noindent{\bf (1.5) The Auslander-Reiten quiver for
$\Cal S(\widetilde 6)$.}

\medskip\noindent
{\bf (1.5.1) Theorem.} {\it
Each Auslander-Reiten sequence in the category $\Cal C\cap \Cal S(\widetilde6)$
with modules in $\Cal D$
is an Auslander-Reiten sequence in the category $\Cal S(\widetilde6)$. }

\smallskip\noindent
{\it Proof: }\/
We show that each source map $f:X\to Y$ in $\Cal C$ with $X\in\Cal D$
is left almost split in $\Cal S(\widetilde6)$.
Let $T\in \ind \Cal S(\widetilde6)$
be a module and $t\in\Hom_{\Cal S}(X,T)$ a non-isomorphism.
Then either $T$ is in $\Cal C$ and hence $t$ factors over $f$, or else by Theorem (1.4.1),
there is a positive integer $j$ such that the shifted representation 
$T[j]$ is in $\Cal D$,
so $\Hom_{\Cal S}(X,T)=\Hom_{\Cal C}(X[j],T[j])=0$ 
since $X[j]$ is in $\Cal C''$ and $(\Cal C',\Cal C'')$ 
is a split torsion pair; thus, $t=0$.     \qed

\medskip
In our main result we verify that the overall structure of the category
$\Cal S(\widetilde6)$ is as described in the introduction. 
Consider the following total ordering on the index set $I=\Bbb Z\times
\Bbb Q^+_0$: Let $(m,\gamma)<(m',\gamma')$ if
either $m<m'$ or $m=m'$ and $\gamma<\gamma'$ hold
and define $\Cal T_\mu=\Cal T_\gamma[m]$ for $\mu=(m,\gamma)\in I$.
We have seen that $\ind\Cal S(\widetilde 6)$ is the union 
$\bigsqcup_{\mu\in I}\Cal T_\mu$.

\smallskip\noindent
{\bf (1.5.2) Theorem.} {\it
For each $\mu\in I$, the category $\Cal T_\mu$ is a $\Bbb P_1(k)$-family
of tubes  of type $(6,3,2)$, separating 
$\Cal P_\mu=\bigsqcup_{\lambda<\mu} \Cal T_\lambda$
from $\Cal Q_\mu=\bigsqcup_{\lambda>\mu}\Cal T_\lambda$ in the sense
that there is no non-zero map from $\Cal T_\mu$ to $\Cal P_\mu$, nor 
from $\Cal Q_\mu$ to $\Cal T_\mu$ or $\Cal P_\mu$, 
nor between different tubes in $\Cal T_\mu$,
and that each map from $\Cal P_\mu$ to $\Cal Q_\mu$
factors through each of the tubes in $\Cal T_\mu$. }

\smallskip
In the proof we will use left $\Cal C$-approximations.
For $M\in\Cal S(\widetilde6)$ the minimal left $\Cal C$-approximation 
$M\to{^{\Cal C}M}$ is given by the obvious map:
$$
\hbox{\beginpicture
\setcoordinatesystem units <0.6cm,0.5cm>
\put{} at 0 0
\put{} at 5 2

\arr{-0.5 0}{-1.3 0}
\arr{-0.5 2}{-1.3 2}
\arr{1.5 0}{0.5 0}
\arr{1.5 2}{0.5 2}
\arr{3.5 0}{2.5 0}
\arr{3.5 2}{2.5 2}

\arr{-2 1.6}{-2 0.4}
\arr{0 1.6}{0 0.4}
\arr{2 1.6}{2 0.4}
\arr{4 1.6}{4 0.4}
\put{$\cdots$} at -3 0
\put{$\cdots$} at -3 2

\put{$\cdots$} at 5 0
\put{$\cdots$} at 5 2
\put{$\ssize \alpha'_0$} at -1 2.4 
\put{$\ssize \alpha'_1$} at 1 2.4
\put{$\ssize \alpha'_2$} at 3 2.4
\put{$\ssize \alpha_0$} at -1 -.4 
\put{$\ssize \alpha_1$} at 1 -.4
\put{$\ssize \alpha_2$} at 3 -.4
\put{$\ssize \beta_{-1}$} at -1.5 1 
\put{$\ssize \beta_0$} at 0.5 1 
\put{$\ssize \beta_1$} at 2.5 1
\put{$\ssize \beta_2$} at 4.5 1

\put{$\ssize M_{-1'}$} at -2 2
\put{$\ssize M_{0'}$} at  0 2
\put{$\ssize M_{1'}$} at  2 2
\put{$\ssize M_{2'}$} at  4 2

\put{$\ssize M_{-1}$} at -2 0
\put{$\ssize M_{0}$} at  0 0
\put{$\ssize M_{1}$} at  2 0
\put{$\ssize M_{2}$} at  4 0

\arr{7 1}{9 1}

\arr{12.5 0}{11.5 0}
\arr{12.5 2}{11.5 2}
\arr{14.5 0}{13.5 0}
\arr{14.5 2}{13.5 2}

\plot 10.9 1.6  10.9 0.4 /
\plot 11.0 1.6  11.0 0.4 /
\plot 12.9 1.6  12.9 0.4 /
\plot 13.0 1.6  13.0 0.4 /

\arr{15 1.6}{15 0.4}

\put{$\cdots$} at 16 0
\put{$\cdots$} at 16 2

\put{$\ssize \alpha_1$} at 12 2.4
\put{$\ssize \alpha_2\beta_2$} at 14 2.4 
\put{$\ssize \alpha_1$} at 12 -.4 
\put{$\ssize \alpha_2$} at 14 -.4
\put{$\ssize \beta_2$} at 15.5 1

\put{$\ssize M_{0}$} at 11 2
\put{$\ssize M_{1}$} at 13 2
\put{$\ssize M_{2'}$} at 15 2

\put{$\ssize M_{0}$} at 11 0
\put{$\ssize M_{1}$} at 13 0
\put{$\ssize M_{2}$} at 15 0

\endpicture}
$$
Note that if $M\in\bigsqcup_{m<0}\ind\Cal D[m]$ then the object 
${^{\Cal C}M}$ in $\Cal C$ is a (sum of) modules in $\Cal P\sqcup\Cal T_0'$ 

\smallskip\noindent
{\it Proof:} \/
Let us first consider the tubes in a given category $\Cal T_\mu$.
We claim that each tubular family 
is a $\Bbb P_1(k)$-family of tubes of tubular type $(6,3,2)$. 
Indeed, $\Cal U$ 
being obtained from the category of regular modules over a 
tilted algebra inherits its tubular type and the parameter set $\Bbb P_1(k)$
for the tubes from the corresponding
hereditary algebra;
and the study of wing extensions in [9], Section 3.4, and the 
categorical equivalences induced by the 
shrinking functors in [9], Section 5.4.3, imply that $\Cal T_1$, and hence 
each $\Cal T_\gamma$, 
consists of a $\Bbb P_1(k)$-family of tubes of the same tubular type $(6,3,2)$
([9], Section 5.6).

\smallskip
It remains to check the separation properties.
To simplify the notation we assume that $\mu\in I$ is such that $\Cal T_\mu$ consists
of modules in $\Cal D$.  The separation properties for test modules
in $\Cal D$ follow from the corresponding properties for modules over the
tubular algebra $\Theta$.
Let $M\in\Cal T_\mu$, $M'\in\Cal D[m']$, and 
$M''\in\Cal D[m'']$ where $m'<0<m''$.
Using the fact shown above that $(\Cal C',\Cal C'')$ is a split
torsion pair, it follows that there is no non-zero map $M''\to M$
and --- as in the proof of Theorem (1.5.1) --- that there
is no non-zero map
$M\to M'$.
Moreover, any map $f:M'\to M''$ in $\Cal S(\widetilde6)$ factors over 
the left $\Cal C$-approximation 
$M'\to {^{\Cal C}M'}$ and over the inclusion $\Res_{Q_\Theta}(M'')\to M''$: 
Since ${^{\Cal C}M'}$ is a sum of modules in $\Cal P\sqcup\Cal T_0'$ and 
since $\Res_{Q_\Theta}(M'')$ is  preinjective, $f$ factors through a sum of 
modules in the tube of $M$ in $\Cal T_\mu$ ([9], Section 5.2.2).   \qed
        \medskip\noindent
{\bf (1.6) The dimension vectors of the indecomposable objects.} 
Let $\chi$ be the quadratic form for the quiver $\widetilde Q$ together with the 
commutativity and the nilpotency relations, and $(-,-)$ the corresponding bilinear form.
Thus $\chi$ and $(-,-)$ are defined on $K_0.$ (The bilinear form is defined as follows: 
$2(\be(z),\be(z')) = 1$ for $z = z'$, as well as if $\{z,z'\}$ is equal to $\{i,i+6\}$,
$\{i',(i+6)'\}$ or $\{i,(i+1)'\}$; $2(\be(z),\be(z')) = -1$ in case 
$\{z,z'\}$ is equal to $\{i,i+1\}$,
$\{i',(i+1)'\}$ or $\{i,i'\}$, and $(\be(z),\be(z')) = 0$ otherwise; here, $i$ denotes
an integral number.) 

Given a dimension vector $\bx$, denote by $t(\bx)$ the number of indexes $i\in \Bbb Z$ such
that $\bx_i \neq 0;$ thus in case $\bx$ is the dimension vector of some object in $\Cal S(\widetilde 6)$, 
then $t(\bx)$ measure the size of the support of the total space. Note that for 
indecomposable objects $M$ in $\Cal S(\widetilde 6)$, we always have $t(\bdim M) \le 8,$
and all the indecomposable objects with $t(\bdim M) = 8$ belong to some shift of $\Cal D''.$

        \medskip\noindent
{\bf (1.6.1) Theorem.} {\it Let $M$ be an indecomposable object in $\Cal S(\widetilde 6)$ 
and put $\bx=\bdim M$. Then }
$$\chi(\bx)=\left\{ \matrix 0 & \text{if}\quad \bx\in\langle \bh^0,\bh^\infty\rangle\cr 
        2 & \text{if}\quad t(\bx)=8 \qquad \cr   
        1 & \text{otherwise.}\qquad \endmatrix\right.$$ 
        \medskip\noindent
{\it Proof: } Since $\chi$ is invariant with respect to the shift, we can assume that 
$M$ belongs to $\Cal D$. If $M$ is a $\Theta$-module (thus in $\Cal D'$), then
$\chi(\bx) = \chi_\Theta(\bx)$ and we can use the known information about the 
quadratic form $\chi_\Theta$: we have $\chi_\Theta(\bx) = 0$ in case 
$\bx$ belongs to the radical of $\chi_\Theta$, thus if and only if 
$\bx\in\langle \bh^0,\bh^\infty\rangle$, and $\chi_\Theta(\bdim M) = 1$ for the
remaining indecomposable $\Theta$-modules. 

Thus it remains to consider $M$ in $\Cal D''.$ As we know, $\Cal D''$ consists of 
two rays in the non-stable tube of $\Cal U$: either $M_{5'} \neq 0$ or $M_7 \neq 0$ (and not 
both). We deal with these two rays separately.

First consider the ray starting at $P(5')$, say 
$$
 P(5') = P(5')\{1\} \to P(5')\{2\} \to P(5')\{3\} \to \cdots.
$$
Let $R' = \rad P(5')$ and consider the ray starting in $R'$, say
$$
 R' = R'\{1\} \to R'\{2\} \to R'\{3\} \to \cdots.
$$
Note that $\bdim P(5')\{i\} = \bdim R'\{i\} + \be(5').$
The first six objects $R'\{1\},\dots, R'\{6\}$ have dimension vectors
$$
 {\tsize{11111000\atop 11111100}} \to
 {\tsize{01111000\atop 01111100}} \to
 {\tsize{01111000\atop 01222210}} \to
 {\tsize{01211000\atop 01222210}} \to
 {\tsize{01221000\atop 01232210}} \to
 {\tsize{01221000\atop 01233210}} \to \cdots
$$
the further dimension vectors are obtained by adding multiples of $\bh^\infty = 
\tsize{01221000\atop 01233210}.$ Note that 
$$
 \chi(\bx+\be(5')) = \chi(\bx) + 1 + 2(\bx,\be(5'))
$$
If $\bx = \bdim R'\{i\}$ with $1 \le i \le 5$, then $\chi(\bx) = 1$ and $\bx_{4'} = 1,$
and $\bx_4 = \bx_5$. This shows that $\chi(\bx+\be(5')) = 1$ in these cases.
For $\bx = \bdim R'\{6\} = \bh^\infty$, we have $\chi(\bx) = 0$ and $(\bx,\be(5')) = 0,$ thus
also here we have $\chi(\bx+\be(5')) = 1$. The fact that $(\bh^\infty,\be(5')) = 0$ shows that
$\chi(\bx) = 1$ for all dimension vectors $\bx$ on the ray starting at $P(5').$

Second, consider the ray starting at $P(7),$
$$
 P(7) = P(7)\{1\} \to P(7)\{2\} \to P(7)\{3\} \to \cdots.
$$
Let $R = \rad P(7)$ and consider the ray starting in $R$, say
$$
 R = R\{1\} \to R\{2\} \to R\{3\} \to \cdots.
$$
We have $\bdim P(7)\{i\} = \bdim R\{i\} + \be(7).$
The first six objects $R\{1\},\dots, R\{6\}$ have dimension vectors
$$
 {\tsize{00000000\atop 00111110}} \to
 {\tsize{00100000\atop 00111110}} \to
 {\tsize{00110000\atop 00121110}} \to
 {\tsize{00110000\atop 00122110}} \to
 {\tsize{11221000\atop 11233210}} \to
 {\tsize{01221000\atop 01233210}} \to \cdots
$$
the further dimension vectors are again obtained by adding multiples of $\bh^\infty = 
\tsize{01221000\atop 01233210}.$ As in the previous case, we use the formula
$$
 \chi(\bx+\be(7)) = \chi(\bx) + 1 + 2(\bx,\be(7)).
$$
For the dimension vectors $\bx = \bdim R\{i\}$ with $1 \le i \le 5$, we have $\chi(\bx) =1$,
whereas $\chi(\bdim R\{6\}) = 0.$ This time, we have $2(\bx,\be(7)) = -1$ for $1 \le i \le 4$,
but $2(\bx,\be(7)) = 0$ for $5 \le i \le 6$. This implies that 
$$
 \chi(\bdim P(7)\{5\}) = 2, \t{and} \chi(\bdim P(7)\{i\}) = 1 \t{for} i = 1,2,3,4,6.
$$
Again, we have $(\bh^\infty,\be(7)) = 0$, thus 
$$
 \chi(\bdim P(7)\{5+6n\}) = 2, \t{and} \chi(\bdim P(7)\{i+6n\}) = 1 \t{for} i = 1,2,3,4,6,
$$
for all $n\in \Bbb N_1.$

Altogether we see: The two rays which form $\Cal D''$ do not contain objects with
dimension vector $\bx$ in $\langle \bh^0,\bh^\infty\rangle$. But they contain all the objects
with $t(\bx) = 8.$ The objects $M$ with $t(\bx) = 8$ are just the objects of the form
$P(7)\{5+6n\}$, and then $\chi(\bx) =2$. Otherwise $\chi(\bx) = 1.$ 
        \qed

        \bigskip\bigskip

{\bf 2. Coverings for Submodule Categories}
        \bigskip

In this section we show that the covering functor 
$$
 \pi\:\Cal S(\widetilde n)\to \Cal S(n)
$$
is dense in case $n\leq 6$.  
As a consequence, we can describe the overall structure of the submodule 
category $\Cal S(6)$ and 
obtain a detailed picture for the indecomposable modules in the 
non-stable family of tubes.

        \medskip\noindent
{\bf (2.1) Coverings.} %
The categories $\Cal S(\widetilde n)$ and $\Cal S(n)$ 
are contained in the categories
$\Cal H(\widetilde n)$ and $\Cal H(n)$ of modules over the 
triangular matrix algebras with coefficients in 
$kA_\infty^\infty/\alpha^n$ and $k[x]/x^n$, respectively.
These two categories are related by the covering functor
$$\pi\:\Cal H(\widetilde n)\to \Cal H(n)$$
which maps an object $\big(A'\sto\beta A\big)\in\Cal H(\widetilde n)$
where $A=(A_i)$ and $A'=(A_i')$ are $kA_\infty^\infty$-modules
to the map $\big(\bigoplus A_i'\lto{\diag\beta}\bigoplus A_i\big)\in\Cal H(n)$
where the $k[x]$-module structure on $\bigoplus A_i$ and $\bigoplus A_i'$
is given by the operation of $\alpha$ on $A$ and $A'$.
The group $\Bbb Z$ acts freely via the shift on the isoclasses of modules 
in $\Cal H(\widetilde n)$. 
According to [6], Section~3, there is an injection from the $\Bbb Z$-orbits
of isoclasses of indecomposable objects in $\Cal H(\widetilde n)$ into the
isoclasses of indecomposable objects in $\Cal H(n)$. 

\smallskip Also the corresponding statement for the subobject categories
$\Cal S(\widetilde n)$ and $\Cal S(n)$ holds:

\medskip\noindent{\bf (2.1.1) Proposition.} {\it
The subcategory $\Cal S(\widetilde n)$ of 
$\Cal H(\widetilde n)$ 
is closed under the operation of the group $G=\Bbb Z$. 
The restriction of the covering functor $\pi$ to 
$\Cal S(\widetilde n)$
maps an object into the corresponding subcategory $\Cal S(n)$
of $\Cal H(n)$.
Hence $\pi|_{\Cal S}$ induces an injection from the $G$-orbits of
isoclasses of indecomposable objects in $\Cal S(\widetilde n)$
into the isoclasses of indecomposable objects in
$\Cal S(n)$. \qed
}

\smallskip
$$\CD \Cal S(\widetilde n) @>\incl>> 
               \Cal H(\widetilde n)\\
        @V \pi|_{\Cal S} VV  @VV \pi V  \\
    \Cal S(n) @>\incl>> 
               \Cal H(n)\\
  \endCD$$

\smallskip 
We are interested in the case when the injection given by $\pi|_{\Cal S}$ 
in Proposition (2.1.1)
is bijective.  For this we adapt the definition of {\it locally
support finite} to our situation.

\smallskip \noindent
{\it Definition.} The {\it support}\/ $\supp M$ of a representation 
$M\in\Cal H(\widetilde n)$ is the set of all points $z$ 
of the quiver $\widetilde Q$ for which $M_z\neq 0$ holds.
The category $\Cal H(\widetilde n)$ is said to be 
{\it locally support finite with respect to}
$\Cal S(\widetilde n)$ provided the union
$$S_z=\bigcup\{\supp M:M\in\ind\Cal S(\widetilde n)\t{and}M_z\neq0\}$$
is a finite set for each point $z$ in the quiver.  
For a set $S$ of points in $Q$, we denote by $S^+$ the set of 
all points $y\in Q$ for which there exists a non-zero path to or from
an element in $S$ (thus, $S^+$ is the union of the supports of the
projective and the injective indecomposable objects in $\Cal
H(\widetilde n)$ which correspond to a point in $S$). 
Note that if $S$ is finite then so is $S^+$.

\medskip\noindent
{\bf (2.1.2) Theorem.} {\it Assume that $\Cal H(\widetilde n)$ is
locally support finite with respect to $\Cal S(\widetilde n)$.
Then the restriction $\pi|_{\Cal S}$ of the covering functor 
$\pi\:\Cal H(\widetilde n)\to \Cal H(n)$ induces
a bijection between the $G$-orbits of isoclasses of indecomposables in
$\Cal S(\widetilde n)$ and the isoclasses of indecomposables in 
$\Cal S(n)$.  
}

\medskip
For the proof of Theorem (2.1.2) we use a lemma.

\medskip\noindent
{\bf (2.1.3) Lemma.} {\it
Let $S$ be a finite set of points in $\widetilde Q$.
Suppose that a locally finite dimensional module 
$A\in\Cal S(\widetilde n)^{\lfd}$ has the property that 
its restriction to $S^+$ decomposes in $\Cal H(\widetilde n)$
as $\Res_{S^+}A=B\oplus C$ where $B$ has support in $S$. 
Then $A$ itself has a decomposition $A=B\oplus D$ where
$B\in\Cal S(\widetilde n)$
and $D\in\Cal S(\widetilde n)^{\lfd}$.
}

\smallskip Clearly,  $\Cal S(\widetilde n)^{\lfd}$ 
denotes the category of all locally finite dimensional 
representations of $\widetilde Q$
which satisfy the commutativity and nilpotency relations and the
extra condition that all maps $\beta_i$ are monic. 

\smallskip\noindent
{\it Proof of Lemma (2.1.3):}\/
If $\Res_{S^+}A=B\oplus C$ is a decomposition as in the lemma, 
then the required module 
$D\in\Cal S(\widetilde n)^{\lfd}$ is constructed as in [5], Lemma 2:
For each point $z\in Q$ put $D_z=C_z$ if $z\not\in S^+$ and 
$D_z=A_z$ otherwise.
Then $D$ is a module in  $\Cal H(\widetilde n)$ 
satisfying $A=B\oplus D$.
Clearly, $B\in\Cal S(\widetilde n)$ and 
$D\in\Cal S(\widetilde n)^{\lfd}$.  \qed

\smallskip\noindent
{\it Proof of Theorem (2.1.2):}\/
We need to show that  the restriction of the covering functor
$\pi\:\Cal S(\widetilde n)\to\Cal S(n)$
is dense.
Let $M=(V,U,T)$ be an indecomposable object in $\Cal S(n)$. 
The pull up $(A,A')$ in $\Cal S(\widetilde n)^{\lfd}$ 
is the locally finite dimensional representation 
such that $A_i=V$ and $A'_i=U$ for each $i$ and such that the operation
of each arrow $\alpha$ is given by $T$. 
Let $z$ be a point of the quiver $\widetilde Q$ which is 
in the support of $(A,A')$.
Since $\Cal H(\widetilde n)$ is locally support finite with respect to
$\Cal S(\widetilde n)$, the set $S_z$ and its extension $S_z^+$
are finite.  Hence the restriction $\Res_{S_z^+}(A,A')$ 
is a direct sum of indecomposable objects in 
$\Cal S(\widetilde n)$.
Let $B$ be an indecomposable direct summand so that 
$B_z\neq 0$, then the support of $B$ is contained in $S_z$.
Hence $B$ is an indecomposable finite dimensional direct summand of 
$(A,A')$, according to
Lemma (2.1.3). It follows from [5], Lemma~1, that $M\cong\pi B$.
\qed

\bigskip
Let $\Cal S$ be one of the categories $\Cal S(n)$ or 
$\Cal S(\widetilde n)$ and let $\Cal H$ be the corresponding
module category $\Cal H(n)$ or $\Cal H(\widetilde n)$, respectively. 
We recall from [11], Proposition~3.2, how Auslander-Reiten sequences
in the submodule category $\Cal S$ are constructed from the corresponding
Auslander-Reiten sequences in $\Cal H$.
Consider the $\Mono$-functor
$$\Mono\:\quad \Cal H\to \Cal S,
\quad \hobj {B'}bB \;\mapsto
        \hobj {\Im b}\incl B.$$ 
If $\hobj {A'}aA$ is an object in $\Cal S$ which is not 
relatively injective in $\Cal S$ and if 
$$0\;\to \;\hobj{A'}aA\;\to\;\hobj{B'}bB\;\to\;\hobj{C'}cC\;\to\; 0$$ 
is an Auslander-Reiten sequence 
in the module category $\Cal H$, then the corresponding 
Auslander-Reiten sequence in the category $\Cal S$ is
obtained
by applying the functor $\Mono$ to the middle and end terms:
$$0\;\to \;\hobj{A'}aA\;\to\;
        \hobj{\Im b}\incl B\;\to\;\hobj{\Im c}\incl C\;\to\; 0$$ 
Since the covering functor $\pi\:\Cal H(\widetilde n)\to\Cal H(n)$
preserves Auslander-Reiten sequences, and since it commutes 
with the $\Mono$-functor
in the sense that the following diagram gives rise to an isomorphism for
each object in $\Cal H(\widetilde n)$,
$$\CD \ssize \Cal H(\widetilde n) @>\pi>> \ssize\Cal H(n)\cr
      @V\Mono VV                      @VV\Mono V \cr
      \ssize \Cal S(\widetilde n) @>>\pi|_{\Cal S}> \ssize\Cal S(n)\cr
\endCD
$$
we obtain the following consequence.

\medskip\noindent
{\bf (2.1.4) Proposition.} {\it
The covering functor $\Cal S(\widetilde n)\to 
\Cal S(n)$ preserves Auslander-Reiten sequences. \qed
}

\medskip
In Lemma (1.3.1) we have obtained a 
description of the projective and the injective objects in 
$\Cal S$ and their sink and source maps.
It follows that the covering functor 
$\Cal S(\widetilde n)\to \Cal S(n)$ 
maps the projective and injective objects in $\Cal S(\widetilde n)$
to the corresponding objects in $\Cal S(n)$ and preserves their
sink and source maps, respectively.

\medskip
As a consequence of Theorem (2.1.2), of Proposition (2.1.4), and of the above
remark we obtain the following relative version of [6], Theorem 3.6.

\smallskip\noindent
{\bf (2.1.5) Theorem.} {\it
The restriction of the covering functor $\pi|_{\Cal S}\: 
\Cal S(\widetilde n)\to \Cal S(n)$
induces an isomorphism of the quotient 
$\Gamma(\widetilde n)/\Bbb Z$ of the Auslander-Reiten quiver $\Gamma(\widetilde n)$
of $\Cal S(\widetilde n)$ modulo the shift
onto the union of some connected components of the 
Auslander-Reiten quiver
$\Gamma(n)$ of $\Cal S(n)$.
This functor $\pi|_{\Cal S}$ is dense if
$\Cal H(\widetilde n)$ is locally support finite with respect to 
$\Cal S(\widetilde n)$. \qed
}

\bigskip\noindent{\bf (2.2) The category $\Cal S(6)$.}
First we observe that $\Cal H(\widetilde 6)$ is locally support
finite with respect to $\Cal S(\widetilde 6)$.

\smallskip\noindent
{\bf (2.2.1) Lemma.} {\it
\item{1.} The support of an indecomposable module in the
non-stable tube in $\Cal U$ is contained in
$$[0,\ldots,7]\cup[0',\ldots,5'].$$
\item{2.} The support of an indecomposable object in a
stable tube in $\Cal U$ is contained in
$$[1,\ldots,6]\cup[1',\ldots,4'].$$
\item{3.} The support of an indecomposable object in 
$\Cal T=\bigsqcup_{\gamma\in\Bbb Q^+}\Cal T_\gamma$
is contained in 
$$[0,\ldots,6]\cup[0',\ldots,4'].$$
}

\noindent
{\it Proof:} The possible support of a module 
in the category $\Cal S(\widetilde n)$ is determined by 
its position relatively to the projective and injective objects
in the Auslander-Reiten quiver.
We have located all the projective and injective modules in Section~1.
\qed

\medskip
The picture above Theorem (1.4.1) describes the overall structure of the
category $\Cal S(\widetilde 6)$.  
According to Lemma (2.2.1) and Theorem (2.1.5), the covering functor
$\Cal S(\widetilde 6 )\to \Cal S(6)$ is dense.  
As a consequence, the indecomposable objects 
in $\Cal S(6)$ are in one-to-one correspondence
to the indecomposables in the fundamental domain $\Cal D$ and 
the Auslander-Reiten quiver for $\Cal S(6)$ 
is obtained from the Auslander-Reiten quiver for $\Cal S(\widetilde 6)$
by identifying the shifted copies of $\Cal D$. 
For each tubular family $\Cal T_\gamma$, $\gamma\in\Bbb Q_0^+$, 
in $\Cal D$, denote the image under the covering functor by 
$\bar{\Cal T}_\gamma$.  
We arrange the tubular families along the unit circle 
by placing $\bar{\Cal T}_\gamma$ for $\gamma=\frac pq$ at the point
$e^{2\pi i\frac p{p+q}}$.

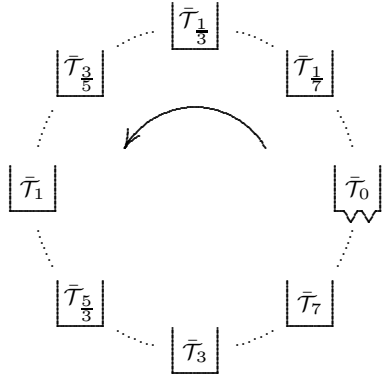
$$\beginpicture
\setcoordinatesystem units <1.8cm,1.8cm> 
\put{} at -1.5 -1.5
\put{} at 1.5 1.5
\put{$\ssize \bar{\Cal T}_{\frac 13}$} at 0 1
\put{$\ssize \bar{\Cal T}_{\frac 17}$} at 0.707 0.707
\put{$\ssize \bar{\Cal T}_0$} at 1 0
\put{$\ssize \bar{\Cal T}_7$} at 0.707 -0.707
\put{$\ssize \bar{\Cal T}_3$} at 0 -1
\put{$\ssize \bar{\Cal T}_{\frac 53}$} at -0.707 -0.707
\put{$\ssize \bar{\Cal T}_1$} at -1 0
\put{$\ssize \bar{\Cal T}_{\frac 35}$} at -0.707 0.707
\put{\tube} at 0 1
\put{\tube} at 0.707 0.707
\put{\Utube} at 1 0
\put{\tube} at 0.707 -0.707
\put{\tube} at 0 -1
\put{\tube} at -0.707 -0.707
\put{\tube} at -1 0
\put{\tube} at -0.707 0.707

\circulararc 120 degrees from 0.433 0.25  center at 0 0
\arr {-0.408 0.2933}  {-0.433 0.25}

\setdots<2pt>
\circulararc -15 degrees from  0.259  0.966 center at 0 0
\circulararc  15 degrees from  0.966  0.259 center at 0 0
\circulararc -15 degrees from  0.966 -0.259 center at 0 0
\circulararc  15 degrees from  0.259 -0.966 center at 0 0
\circulararc -15 degrees from -0.259 -0.966 center at 0 0
\circulararc  15 degrees from -0.966 -0.259 center at 0 0
\circulararc -15 degrees from -0.966  0.259 center at 0 0
\circulararc  15 degrees from -0.259  0.966 center at 0 0

\endpicture$$

Thus, homomorphisms in the infinite radical of $\Cal S(6)$ 
go counterclockwise, as indicated by the arrow. 
More precisely, the homomorphisms in the category $\Cal S(n)$
can be computed from those in $\Cal S(\widetilde n)$ by using 
the following formula.

\smallskip\noindent
{\bf (2.2.2) Lemma.} {\it
For $M,N\in \Cal H(\widetilde n)$ there is an isomorphism of $k$-spaces
$$ \Hom_{\Cal H(n)}(\pi M,\pi N)\;\cong\;
   \bigoplus_{g\in \Bbb Z}\Hom_{\Cal H(\widetilde n)}(M,N[g]).$$
}

\smallskip\noindent
{\it Proof:} \/ The assertion holds since the covering functor
$\pi$ is left adjoint to the pull up functor $F_\bullet$ ([3], Section~3.2) and since
there is a canonical isomorphism $F_\bullet\pi N\cong \bigoplus_{g\in \Bbb Z}N[g]$
([6], Lemma 3.2.) \qed

\medskip\noindent{\bf (2.3) The non-stable tubular family.}
We conclude this section with a detailed picture of the
tubes in $\bar{\Cal T}_0$. We present an object $(A,A')$
in $\Cal S(n)$ as follows.

\smallskip 
The $k[x]$-module $A$ has the form 
$A\cong\bigoplus_{i=1}^tk[x]/x^{\lambda_i}$
where the $\lambda_i$ form a partition with all parts bounded by $n$.
We picture $A$ by rotating the Young diagram for the partition
$(\lambda_i)$ by $90^{\text o}$ and by suitably adjusting the columns
in height and order.  
The submodule $A'$ is represented by the image of a set of  generators:
If one generator has a non-zero entry in the $i$-th summand of $A$,
say the entry is in the $j$-th power of the radical, but not the $(j+1)$-st
power, then we put a symbol $\bullet$ in the $j$-th box of the $i$-th 
column to indicate the entry $x^j$.  Different entries, say $cx^j$ where
$c$ is a unit, are indicated by writing the symbol $c$ into the $j$-th box.
Symbols belonging to the same generator of $A'$ are connected
by a line.  In fact, we will arrange the columns of $A$ in height and
order in such a way that these connecting lines are horizontal and
short, as good as possible.

\smallskip\noindent
{\it Example.}
The only objects in $\Cal D$ where the total space has maximal
support $[0,\ldots,7]$ occur in the non-stable tube in $\Cal U$
at each of the intersections of the ray starting at $P(7)$ with
the coray ending at $P(5')$.
The first such object $M$ has dimension vector $\sssize{11221000\atop 11233211}$.
Let $(A,A')=\pi(M)$ be the corresponding
object in $\Cal S(k[x]/x^6)$.  It follows that 
$A= k[x]/x^6\oplus k[x]/x^2\oplus k[x]/x^6$ so the partition 
representing $A$ is $(6,6,2)$. The submodule $A'$ is isomorphic to
$k[x]/x^5\oplus k[x]/x^2$, the two canonical generators
are mapped to $(0,1,x)$ and $(x^4,x,0)$ in $A$. 
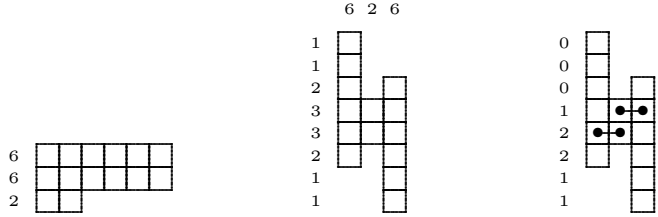
$$      \beginpicture 
        \setcoordinatesystem units <0.25cm,0.25cm>
        \multiput{\phantom{x}} at 0 -2  6 2 /
        \multiput{\sq} at 0 0  0 1  0 2  1 0  1 1  1 2  2 1  2 2  3 1  
                          3 2  4 1  4 2  5 1  5 2 /
        \put{$\sssize 6$} at -1 2
        \put{$\sssize 6$} at -1 1
        \put{$\sssize 2$} at -1 0
        \endpicture
   \qquad\qquad
        \beginpicture 
        \setcoordinatesystem units <0.25cm,0.25cm>
        \multiput{} at -1 0  4 8 /
        \multiput{\sq} at 0 2  0 3  0 4  0 5  0 6  0 7  1 3  1 4  
                2 0  2 1  2 2  2 3  2 4  2 5 /
        \put{$\sssize 1$} at -1 0
        \put{$\sssize 1$} at -1 1
        \put{$\sssize 2$} at -1 2
        \put{$\sssize 3$} at -1 3
        \put{$\sssize 3$} at -1 4
        \put{$\sssize 2$} at -1 5
        \put{$\sssize 1$} at -1 6
        \put{$\sssize 1$} at -1 7
        \put{$\sssize 6$} at 0.5 8.5
        \put{$\sssize 2$} at 1.5 8.5
        \put{$\sssize 6$} at 2.5 8.5
        \endpicture
   \qquad\qquad
        \beginpicture 
        \setcoordinatesystem units <0.25cm,0.25cm>
        \multiput{} at -1 0  4 7 /
        \multiput{\sq} at 0 2  0 3  0 4  0 5  0 6  0 7  1 3  1 4  
                2 0  2 1  2 2  2 3  2 4  2 5 /
        \multiput{$\ssize \bullet$} at .5 3  1.5 3  1.5 4  2.5 4 /
        \plot 0.5 3  1.5 3 /
        \plot 1.5 4  2.5 4 /
        \put{$\sssize 1$} at -1 0
        \put{$\sssize 1$} at -1 1
        \put{$\sssize 2$} at -1 2
        \put{$\sssize 2$} at -1 3
        \put{$\sssize 1$} at -1 4
        \put{$\sssize 0$} at -1 5
        \put{$\sssize 0$} at -1 6
        \put{$\sssize 0$} at -1 7
        \endpicture
$$
The pictures represent the Young diagram for $A$, the rotated Young diagram
for $A$ with the columns suitably adjusted in height and order, 
and the same diagram with the generators
of the subgroup $A'$ indicated.  On the left hand side of the middle and
the last diagram we list the entries of the dimension vector of $M$.

%
%
%
%
%
%
%
%
\def\scale#1{\beginpicture\setcoordinatesystem units <1.5mm,1.5mm>#1
        \endpicture}

\newpage
\centerline{\bf The non-stable tube}

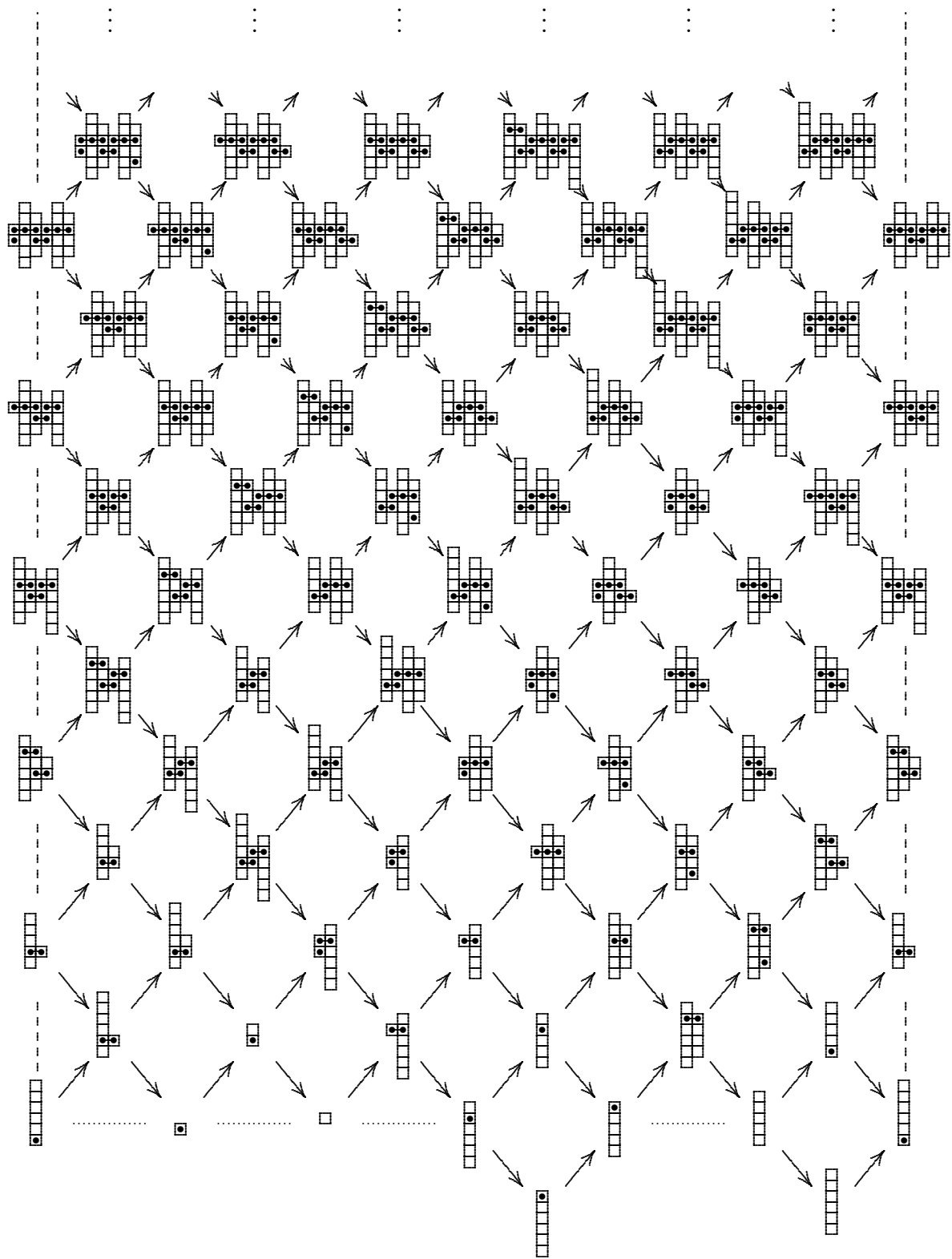
\def\Xooooo{\multiput{} at 0 -2  0 7 /
        \multiput{\sq} at 0 -1 0 0  0 1  0 2  0 3  0 4 /
        \put{$\sssize\bullet$} at .5 4 }
\def\Oooooo{\multiput{} at 0 -2  0 7 /
        \multiput{\sq} at 0 1  0 2  0 3  0 4  0 5  0 6 / }
\def\Oxoooo{\multiput{} at 0 -2  0 7 /
        \multiput{\sq} at 0 -1 0 0  0 1  0 2  0 3  0 4 /
        \put{$\sssize\bullet$} at .5 3 }
\def\cOxoooo{\multiput{} at 0 -2  0 7 /
        \multiput{\sq} at 0 -1 0 0  0 1  0 2  0 3  0 4 /
        \multiput{$\sssize\bullet$} at .5 3  -.5 3 /
        \plot -.5 3  .5 3 / }
\def\Xoooo{\multiput{} at 0 -2  0 7 /
        \multiput{\sq} at 0 0  0 1  0 2  0 3  0 4 /
        \put{$\sssize\bullet$} at .5 4 }
\def\Ooooo{\multiput{} at 0 -2  0 7 /
        \multiput{\sq} at 0 1  0 2  0 3  0 4  0 5 /  }
\def\x{\multiput{} at 0 -2  0 7 /
        \put{$\sssize\bullet$} at .5 1 }
\def\X{\multiput{} at 0 -2  0 7 /
        \put{\sq} at 0 2
        \put{$\sssize \bullet$} at .5 2 }
\def\cX{\multiput{} at 0 -2  0 7 /
        \put{\sq} at 0 2
        \multiput{$\sssize \bullet$} at -.5 2  .5 2 /
        \plot -.5 2  .5 2 / }
\def\O{\multiput{} at 0 -2  0 7 /
        \put{\sq} at 0 3 }
\def\Oxooo{\multiput{} at 0 -2  0 7 /
        \multiput{\sq} at 0 0  0 1  0 2  0 3  0 4 /
        \put{$\sssize\bullet$} at .5 3 }
\def\cOxooo{\multiput{} at 0 -2  0 7 /
        \multiput{\sq} at 0 0  0 1  0 2  0 3  0 4 /
        \multiput{$\sssize\bullet$} at -.5 3  .5 3 /
        \plot -.5 3  .5 3 / }
\def\OcooooXooo{\multiput{} at 0 -2  0 7 /
        \multiput{\sq} at 0 0  0 1  0 2  0 3  0 4  0 5  1 1  1 2  1 3  1 4 /
        \multiput{$\sssize\bullet$} at .5 4  1.5 4 /
        \plot .5 4  1.5 4 / }
\def\Ox{\multiput{} at 0 -2  0 7 /
        \multiput{\sq} at 0 2  0 3 /
        \put{$\sssize\bullet$} at .5 2 }
\def\cOx{\multiput{} at 0 -2  0 7 /
        \multiput{\sq} at 0 2  0 3 /
        \multiput{$\sssize\bullet$} at -.5 2  .5 2 / 
        \plot -.5 2  .5 2 / }
\def\OocoooOxoo{\multiput{} at 0 -2  0 7 /
        \multiput{\sq} at 0 0  0 1  0 2  0 3  0 4  0 5  1 1  1 2  1 3  1 4 /
        \multiput{$\sssize\bullet$} at .5 3  1.5 3 /
        \plot .5 3  1.5 3 / }
\def\cOocoooOxoo{\multiput{} at 0 -2  0 7 /
        \multiput{\sq} at 0 0  0 1  0 2  0 3  0 4  0 5  1 1  1 2  1 3  1 4 /
        \multiput{$\sssize\bullet$} at -.5 3  .5 3  1.5 3 /
        \plot -.5 3  1.5 3 / }
\def\hga{\scale{\put{\OocoooOxoo} at 0 0 
                        \put{\cOx} at 2 0 
                        \put{\cOxooo} at 3 0 }}
\def\hgb{\scale{\put{\OcooooXooo} at 0 0 
                        \put{\cOx} at 2 0 
                        \put{\cOocoooOxoo} at 3 0 }}
\def\hgc{\scale{\put{\Ooooo} at 0 0 
                        \put{\cOx} at 1 0 
                        \put{\cOocoooOxoo} at 2 0 
                        \put{\x} at 3 0 }}
\def\hgd{\scale{\put{\Oooooo} at 0 0 
                        \put{\cOx} at 1 0 
                        \put{\cOocoooOxoo} at 2 0 
                        \put{\cX} at 4 0 }}
\def\hge{\scale{\put{\Ox} at 0 0 
                        \put{\cOocoooOxoo} at 1 0 
                        \put{\cOx} at 3 0 }}
\def\hgf{\scale{\put{\O} at 0 0
                        \put{\cOocoooOxoo} at 1 0 
                        \put{\cOx} at 3 0
                        \put{\cOxoooo} at 4 0 }}

\def\hha{\scale{\put{\O} at 0 0
                        \put{\cOocoooOxoo} at 1 0 
                        \put{\cOx} at 3 0
                        \put{\cOxooo} at 4 0 }}
\def\hhb{\scale{\put{\OocoooOxoo} at 0 0 
                        \put{\cOx} at 2 0 
                        \put{\cOocoooOxoo} at 3 0 }}
\def\hhc{\scale{\put{\OcooooXooo} at 0 0 
                        \put{\cOx} at 2 0 
                        \put{\cOocoooOxoo} at 3 0 
                        \put{\x} at 4 0 }}
\def\hhd{\scale{\put{\Ooooo} at 0 0 
                        \put{\cOx} at 1 0 
                        \put{\cOocoooOxoo} at 2 0 
                        \put{\cX} at 4 0 }}
\def\hhe{\scale{\put{\Oooooo} at 0 0 
                        \put{\cOx} at 1 0 
                        \put{\cOocoooOxoo} at 2 0 
                        \put{\cOx} at 4 0 }}
\def\hhf{\scale{\put{\Ox} at 0 0 
                        \put{\cOocoooOxoo} at 1 0 
                        \put{\cOx} at 3 0
                        \put{\cOxoooo} at 4 0 }}
\def\hia{\scale{\put{\O} at 0 0
                        \put{\cOocoooOxoo} at 1 0 
                        \put{\cOx} at 3 0
                        \put{\cOocoooOxoo} at 4 0 }}
\def\hib{\scale{\put{\OocoooOxoo} at 0 0 
                        \put{\cOx} at 2 0 
                        \put{\cOocoooOxoo} at 3 0 
                        \put{\x} at 4 0 }}
\def\hic{\scale{\put{\OcooooXooo} at 0 0 
                        \put{\cOx} at 2 0 
                        \put{\cOocoooOxoo} at 3 0 
                        \put{\cX} at 5 0 }}
\def\hid{\scale{\put{\Ooooo} at 0 0 
                        \put{\cOx} at 1 0 
                        \put{\cOocoooOxoo} at 2 0 
                        \put{\cOx} at 4 0 }}
\def\hie{\scale{\put{\Oooooo} at 0 0 
                        \put{\cOx} at 1 0 
                        \put{\cOocoooOxoo} at 2 0 
                        \put{\cOx} at 4 0
                        \put{\cOxoooo} at 5 0 }}
\def\hif{\scale{\put{\Ox} at 0 0 
                        \put{\cOocoooOxoo} at 1 0 
                        \put{\cOx} at 3 0
                        \put{\cOxooo} at 4 0 }}
\def\hja{\scale{\put{\Ox} at 0 0 
                        \put{\cOocoooOxoo} at 1 0 
                        \put{\cOx} at 3 0
                        \put{\cOocoooOxoo} at 4 0 }}
\def\hjb{\scale{\put{\O} at 0 0
                        \put{\cOocoooOxoo} at 1 0 
                        \put{\cOx} at 3 0
                        \put{\cOocoooOxoo} at 4 0
                        \put{\x} at 5 0 }}
\def\hjc{\scale{\put{\OocoooOxoo} at 0 0 
                        \put{\cOx} at 2 0 
                        \put{\cOocoooOxoo} at 3 0 
                        \put{\cX} at 5 0 }}
\def\hjd{\scale{\put{\OcooooXooo} at 0 0 
                        \put{\cOx} at 2 0 
                        \put{\cOocoooOxoo} at 3 0 
                        \put{\cOx} at 5 0 }}
\def\hje{\scale{\put{\Ooooo} at 0 0 
                        \put{\cOx} at 1 0 
                        \put{\cOocoooOxoo} at 2 0 
                        \put{\cOx} at 4 0
                        \put{\cOxoooo} at 5 0 }}
\def\hjf{\scale{\put{\Oooooo} at 0 0 
                        \put{\cOx} at 1 0 
                        \put{\cOocoooOxoo} at 2 0 
                        \put{\cOx} at 4 0
                        \put{\cOxooo} at 5 0 }}
\def\hka{\scale{\put{\Ox} at 0 0 
                        \put{\cOocoooOxoo} at 1 0 
                        \put{\cOx} at 3 0
                        \put{\cOocoooOxoo} at 4 0
                        \put{\x} at 5 0 }}
\def\hkb{\scale{\put{\O} at 0 0
                        \put{\cOocoooOxoo} at 1 0 
                        \put{\cOx} at 3 0
                        \put{\cOocoooOxoo} at 4 0
                        \put{\cX} at 6 0 }}
\def\hkc{\scale{\put{\OocoooOxoo} at 0 0 
                        \put{\cOx} at 2 0 
                        \put{\cOocoooOxoo} at 3 0 
                        \put{\cOx} at 5 0 }}
\def\hkd{\scale{\put{\OcooooXooo} at 0 0 
                        \put{\cOx} at 2 0 
                        \put{\cOocoooOxoo} at 3 0 
                        \put{\cOx} at 5 0
                        \put{\cOxoooo} at 6 0 }}
\def\hke{\scale{\put{\Ooooo} at 0 0 
                        \put{\cOx} at 1 0 
                        \put{\cOocoooOxoo} at 2 0 
                        \put{\cOx} at 4 0
                        \put{\cOxooo} at 5 0 }}
\def\hkf{\scale{\put{\Oooooo} at 0 0 
                        \put{\cOx} at 1 0 
                        \put{\cOocoooOxoo} at 2 0 
                        \put{\cOx} at 4 0
                        \put{\cOocoooOxoo} at 5 0 }}

$$\beginpicture\setcoordinatesystem units <1cm,1.23cm>
        \put{\scale{\put{\Xooooo} at 0 0 }} at 6.9 -1
        \put{\scale{\put{\Oooooo} at 0 0 }} at 10.9 -1
        \multiput{\scale{\put{\Oooooo} at 0 0 
                        \put{\x} at 0 0 }} at -.1 0  11.9 0 /
        \put{\scale{\put{\X} at 0 0 }} at 1.9 0
        \put{\scale{\put{\O} at 0 0 }} at 3.9 0
        \put{\scale{\put{\Oxoooo} at 0 0 }} at 5.9 0
        \put{\scale{\put{\Xoooo} at 0 0 }} at 7.9 0
        \put{\scale{\put{\Ooooo} at 0 0 }} at 9.9 0
        \put{\scale{\put{\Oooooo} at 0 0 
                        \put{\cX} at 1 0 }} at .9 1
        \put{\scale{\put{\Ox} at 0 0 }} at 2.9 1
        \put{\scale{\put{\O} at 0 0
                        \put{\cOxoooo} at 1 0 }} at 4.9 1
        \put{\scale{\put{\Oxooo} at 0 0 }} at 6.9 1
        \put{\scale{\put{\OcooooXooo} at 0 0 }} at 8.9 1
        \put{\scale{\put{\Ooooo} at 0 0 
                        \put{\x} at 0 0 }} at 10.9 1
        \multiput{\scale{\put{\Ooooo} at 0 0 
                        \put{\cX} at 1 0 }} at -.1 2  11.9 2 /
        \put{\scale{\put{\Oooooo} at 0 0 
                        \put{\cOx} at 1 0 }} at 1.9 2
        \put{\scale{\put{\Ox} at 0 0 
                        \put{\cOxoooo} at 1 0 }} at 3.9 2
        \put{\scale{\put{\O} at 0 0
                        \put{\cOxooo} at 1 0 }} at 5.9 2
        \put{\scale{\put{\OocoooOxoo} at 0 0 }} at 7.9 2
        \put{\scale{\put{\OcooooXooo} at 0 0 
                        \put{\x} at 1 0 }} at 9.9 2
        \put{\scale{\put{\Ooooo} at 0 0 
                        \put{\cOx} at 1 0 }} at .9 3
        \put{\scale{\put{\Oooooo} at 0 0 
                        \put{\cOx} at 1 0 
                        \put{\cOxoooo} at 2 0 }} at 2.9 3
        \put{\scale{\put{\Ox} at 0 0 
                        \put{\cOxooo} at 1 0 }} at 4.9 3
        \put{\scale{\put{\O} at 0 0
                        \put{\cOocoooOxoo} at 1 0 }} at 6.9 3
        \put{\scale{\put{\OocoooOxoo} at 0 0 
                        \put{\x} at 1 0 }} at 8.9 3
        \put{\scale{\put{\OcooooXooo} at 0 0 
                        \put{\cX} at 2 0 }} at 10.9 3
        \multiput{\scale{\put{\OcooooXooo} at 0 0 
                        \put{\cOx} at 2 0 }} at -.1 4  11.9 4 /
        \put{\scale{\put{\Ooooo} at 0 0 
                        \put{\cOx} at 1 0 
                        \put{\cOxoooo} at 2 0 }} at 1.9 4
        \put{\scale{\put{\Oooooo} at 0 0 
                        \put{\cOx} at 1 0 
                        \put{\cOxooo} at 2 0 }} at 3.9 4
        \put{\scale{\put{\Ox} at 0 0 
                        \put{\cOocoooOxoo} at 1 0 }} at 5.9 4
        \put{\scale{\put{\O} at 0 0
                        \put{\cOocoooOxoo} at 1 0 
                        \put{\x} at 2 0 }} at 7.9 4
        \put{\scale{\put{\OocoooOxoo} at 0 0 
                        \put{\cX} at 2 0 }} at 9.9 4
        \put{\scale{\put{\OcooooXooo} at 0 0 
                        \put{\cOx} at 2 0 
                        \put{\cOxoooo} at 3 0 }} at .9 5 
        \put{\scale{\put{\Ooooo} at 0 0 
                        \put{\cOx} at 1 0 
                        \put{\cOxooo} at 2 0 }} at 2.9 5
        \put{\scale{\put{\Oooooo} at 0 0 
                        \put{\cOx} at 1 0 
                        \put{\cOocoooOxoo} at 2 0 }} at 4.9 5
        \put{\scale{\put{\Ox} at 0 0 
                        \put{\cOocoooOxoo} at 1 0 
                        \put{\x} at 2 0 }} at 6.9 5
        \put{\scale{\put{\O} at 0 0
                        \put{\cOocoooOxoo} at 1 0 
                        \put{\cX} at 3 0 }} at 8.9 5
        \put{\scale{\put{\OocoooOxoo} at 0 0 
                        \put{\cOx} at 2 0 }} at 10.9 5

        \multiput{\scale{\put{\OocoooOxoo} at 0 0 
                        \put{\cOx} at 2 0 
                        \put{\cOxoooo} at 3 0 }} at -.1 6  11.9 6 /
        \put{\scale{\put{\OcooooXooo} at 0 0 
                        \put{\cOx} at 2 0 
                        \put{\cOxooo} at 3 0 }} at 1.9 6 
        \put{\scale{\put{\Ooooo} at 0 0 
                        \put{\cOx} at 1 0 
                        \put{\cOocoooOxoo} at 2 0 }} at 3.9 6
        \put{\scale{\put{\Oooooo} at 0 0 
                        \put{\cOx} at 1 0 
                        \put{\cOocoooOxoo} at 2 0 
                        \put{\x} at 3 0 }} at 5.9 6
        \put{\scale{\put{\Ox} at 0 0 
                        \put{\cOocoooOxoo} at 1 0 
                        \put{\cX} at 3 0 }} at 7.9 6
        \put{\scale{\put{\O} at 0 0
                        \put{\cOocoooOxoo} at 1 0 
                        \put{\cOx} at 3 0 }} at 9.9 6
        \put{\hga} at .9 7
        \put{\hgb} at 2.9 7
        \put{\hgc} at 4.9 7
        \put{\hgd} at 6.9 7 
        \put{\hge} at 8.9 7
        \put{\hgf} at 10.9 7

        \multiput{\hha} at -.1 8  12 8 /
        \put{\hhb} at 1.9 8
        \put{\hhc} at 3.9 8
        \put{\hhd} at 5.9 8
        \put{\hhe} at 7.9 8
        \put{\hhf} at 9.9 8

        \put{\hia} at .9 9
        \put{\hib} at 2.9 9
        \put{\hic} at 4.9 9
        \put{\hid} at 6.9 9 
        \put{\hie} at 8.9 9
        \put{\hif} at 10.9 9

        \multiput{\hja} at -.1 10  12 10 /
        \put{\hjb} at 1.9 10
        \put{\hjc} at 3.9 10
        \put{\hjd} at 5.9 10
        \put{\hje} at 7.9 10
        \put{\hjf} at 9.9 10

        \put{\hka} at .9 11
        \put{\hkb} at 2.9 11
        \put{\hkc} at 4.9 11
        \put{\hkd} at 6.9 11
        \put{\hke} at 8.9 11
        \put{\hkf} at 10.9 11
\arr{6.3 -.3} {6.7 -.7}
\arr{7.3 -.7} {7.7 -.3}
\arr{10.3 -.3} {10.7 -.7}
\arr{11.3 -.7} {11.7 -.3}
\arr{0.3 0.3} {0.7 0.7}
\arr{2.3 0.3} {2.7 0.7}
\arr{4.3 0.3} {4.7 0.7}
\arr{1.3 0.7} {1.7 0.3}
\arr{3.3 0.7} {3.7 0.3}
\arr{5.3 0.7} {5.7 0.3}
\arr{6.3 0.3} {6.7 0.7}
\arr{8.3 0.3} {8.7 0.7}
\arr{10.3 0.3} {10.7 0.7}
\arr{7.3 0.7} {7.7 0.3}
\arr{9.3 0.7} {9.7 0.3}
\arr{11.3 0.7} {11.7 0.3}
\arr{0.3 1.7} {0.7 1.3}
\arr{2.3 1.7} {2.7 1.3}
\arr{4.3 1.7} {4.7 1.3}
\arr{1.3 1.3} {1.7 1.7}
\arr{3.3 1.3} {3.7 1.7}
\arr{5.3 1.3} {5.7 1.7}
\arr{6.3 1.7} {6.7 1.3}
\arr{8.3 1.7} {8.7 1.3}
\arr{10.3 1.7} {10.7 1.3}
\arr{7.3 1.3} {7.7 1.7}
\arr{9.3 1.3} {9.7 1.7}
\arr{11.3 1.3} {11.7 1.7}
\arr{0.3 2.3} {0.7 2.7}
\arr{1.3 2.7} {1.7 2.3}
\arr{2.3 2.3} {2.7 2.7}
\arr{3.3 2.7} {3.7 2.3}
\arr{4.3 2.3} {4.7 2.7}
\arr{5.3 2.7} {5.7 2.3}
\arr{6.3 2.3} {6.7 2.7}
\arr{7.3 2.7} {7.7 2.3}
\arr{8.3 2.3} {8.7 2.7}
\arr{9.3 2.7} {9.7 2.3}
\arr{10.3 2.3} {10.6 2.6}
\arr{11.3 2.7} {11.7 2.3}
\arr{0.3 3.7} {0.7 3.3}
\arr{1.3 3.3} {1.7 3.7}
\arr{2.35 3.65} {2.65 3.35}
\arr{3.3 3.3} {3.7 3.7}
\arr{4.3 3.7} {4.7 3.3}
\arr{5.3 3.3} {5.7 3.7}
\arr{6.3 3.7} {6.7 3.3}
\arr{7.4 3.4} {7.7 3.7}
\arr{8.3 3.7} {8.7 3.3}
\arr{9.3 3.3} {9.6 3.6}
\arr{10.3 3.7} {10.6 3.4}
\arr{11.3 3.3} {11.6 3.6}
\arr{0.3 4.3} {0.6 4.6}
\arr{1.4 4.6} {1.65 4.35}
\arr{2.3 4.3} {2.7 4.7}
\arr{3.3 4.7} {3.6 4.4}
\arr{4.3 4.3} {4.7 4.7}
\arr{5.3 4.7} {5.7 4.3}
\arr{6.3 4.3} {6.7 4.7}
\arr{7.3 4.7} {7.7 4.3}
\arr{8.3 4.3} {8.7 4.7}
\arr{9.3 4.7} {9.6 4.4}
\arr{10.3 4.3} {10.6 4.6}
\arr{11.3 4.7} {11.6 4.4}
\arr{0.4 5.6} {0.6 5.4}
\arr{1.35 5.35} {1.6 5.6}
\arr{2.4 5.6} {2.6 5.4}
\arr{3.35 5.35} {3.65 5.65}
\arr{4.4 5.6} {4.6 5.4}
\arr{5.4 5.4} {5.65 5.65}
\arr{6.4 5.6} {6.65 5.35}
\arr{7.35 5.35} {7.65 5.65}
\arr{8.35 5.65} {8.65 5.35}
\arr{9.35 5.35} {9.65 5.65}
\arr{10.35 5.65} {10.6 5.4}
\arr{11.35 5.35} {11.6 5.6}
\arr{0.35 6.35} {0.6 6.6}
\arr{1.4 6.6} {1.6 6.4}
\arr{2.35 6.35} {2.6 6.6}
\arr{3.4 6.6} {3.6 6.4}
\arr{4.4 6.4} {4.6 6.6}
\arr{5.4 6.6} {5.6 6.4}
\arr{6.4 6.4} {6.6 6.6}
\arr{7.4 6.6} {7.6 6.4}
\arr{8.35 6.35} {8.65 6.65}
\arr{9.35 6.65} {9.65 6.35}
\arr{10.35 6.35} {10.6 6.6}
\arr{11.4 6.6} {11.6 6.4}
\arr{0.4 7.6} {0.6 7.4}
\arr{1.4 7.4} {1.6 7.6}
\arr{2.4 7.6} {2.6 7.4}
\arr{3.4 7.4} {3.6 7.6}
\arr{4.4 7.6} {4.6 7.4}
\arr{5.4 7.4} {5.6 7.6}
\arr{6.35 7.65} {6.55 7.45}
\arr{7.35 7.35} {7.65 7.65}
\arr{8.35 7.65} {8.65 7.35}
\arr{9.35 7.35} {9.6 7.6}
\arr{10.45 7.55} {10.6 7.4}
\arr{11.4 7.4} {11.65 7.65}
\arr{0.4 8.4} {0.6 8.6}
\arr{1.4 8.6} {1.6 8.4}
\arr{2.4 8.4} {2.55 8.55}
\arr{3.4 8.6} {3.55 8.45}
\arr{4.4 8.4} {4.55 8.55}
\arr{5.35 8.65} {5.55 8.45}
\arr{6.35 8.35} {6.6 8.6}
\arr{7.35 8.65} {7.55 8.45}
\arr{8.35 8.35} {8.65 8.65}
\arr{9.5 8.5} {9.6 8.4}
\arr{10.35 8.35} {10.6 8.6}
\arr{11.4 8.6} {11.65 8.35}
\arr{0.4 9.6} {0.6 9.4}
\arr{1.4 9.4} {1.6 9.6}
\arr{2.4 9.6} {2.6 9.4}
\arr{3.4 9.4} {3.6 9.6}
\arr{4.4 9.6} {4.6 9.4}
\arr{5.4 9.4} {5.6 9.6}
\arr{6.4 9.6} {6.6 9.4}
\arr{7.4 9.4} {7.6 9.6}
\arr{8.4 9.6} {8.55 9.45}
\arr{9.4 9.4} {9.6 9.6}
\arr{10.4 9.6} {10.6 9.4}
\arr{11.4 9.4} {11.6 9.6}
\arr{0.4 10.4} {0.6 10.6}
\arr{1.4 10.6} {1.6 10.4}
\arr{2.4 10.4} {2.6 10.6}
\arr{3.4 10.6} {3.6 10.4}
\arr{4.4 10.4} {4.6 10.6}
\arr{5.4 10.6} {5.6 10.4}
\arr{6.4 10.4} {6.6 10.6}
\arr{7.5 10.5} {7.6 10.4}
\arr{8.4 10.4} {8.6 10.6}
\arr{9.4 10.6} {9.5 10.5}
\arr{10.4 10.4} {10.6 10.6}
\arr{11.4 10.6} {11.6 10.4}
\arr{0.4 11.6} {0.6 11.4}
\arr{1.4 11.4} {1.6 11.6}
\arr{2.4 11.6} {2.55 11.45}
\arr{3.4 11.4} {3.6 11.6}
\arr{4.4 11.6} {4.55 11.45}
\arr{5.4 11.4} {5.6 11.6}
\arr{6.4 11.6} {6.55 11.45}
\arr{7.4 11.4} {7.6 11.6}
\arr{8.4 11.6} {8.55 11.45}
\arr{9.4 11.4} {9.6 11.6}
\arr{10.35 11.65} {10.45 11.55}
\arr{11.45 11.45} {11.6 11.6}
\setdashes<3pt>
\plot 0 0.6  0 1.4 /
\plot 0 2.6  0 3.4 /
\plot 0 4.6  0 5.4 /
\plot 0 6.6  0 7.4 /
\plot 0 8.6  0 9.4 /
\plot 0 10.6  0 12.5 /
\plot 12 0.6  12 1.4 /
\plot 12 2.6  12 3.4 /
\plot 12 4.6  12 5.4 /
\plot 12 6.6  12 7.4 /
\plot 12 8.6  12 9.4 /
\plot 12 10.6  12 12.5 /
\setdots<2pt>
\plot 0.5 0  1.5 0 /
\plot 2.5 0  3.5 0 /
\plot 4.5 0  5.5 0 /
\plot 8.5 0  9.5 0 /
\multiput{$\vdots$} at 1 12.5  3 12.5  5 12.5  7 12.5  9 12.5  11 12.5 /
\endpicture$$
%
%

\newpage
\centerline{\bf The stable non-homogeneous tubes in $\bar \Cal T_0$}

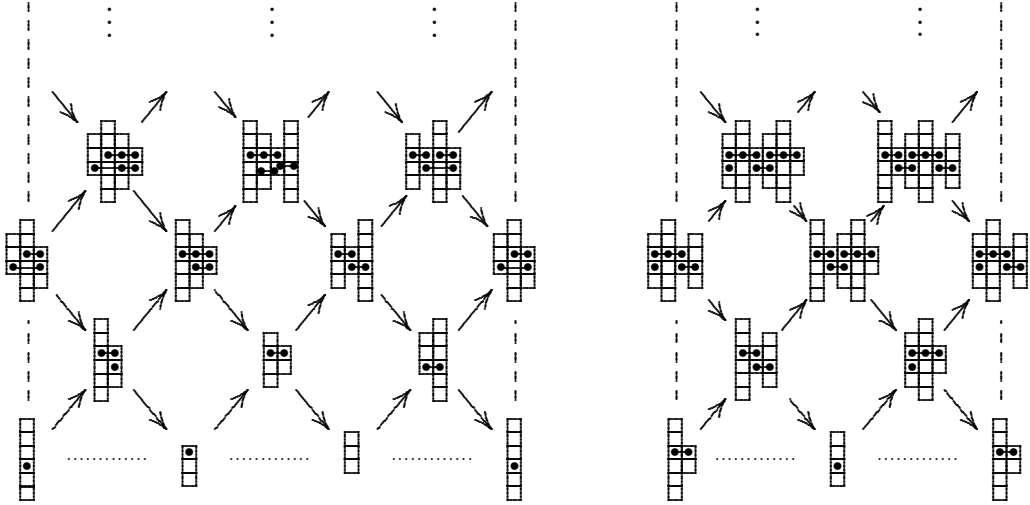
\def\Oooxoo{\multiput{} at 0 -1  0 6 /
        \multiput{\sq} at 0 0  0 1  0 2  0 3  0 4  0 5 /
        \put{$\sssize\bullet$} at .5 2 }
\def\OooxooR{\multiput{} at 0 -1  0 6 /
        \multiput{\sq} at 0 0  0 1  0 2  0 3  0 4  0 5 /
        \put{$\sssize\bullet$} at 1.5 2 }
\def\OooxooRR{\multiput{} at 0 -1  0 6 /
        \multiput{\sq} at 0 0  0 1  0 2  0 3  0 4  0 5 /
        \multiput{$\sssize\bullet$} at 1.5 2  2.5 2 /
        \plot 1.5 2  2.5 2 / }
\def\OooxooRS{\multiput{} at 0 -1  0 6 /
        \multiput{\sq} at 0 0  0 1  0 2  0 3  0 4  0 5 /
        \multiput{$\sssize\bullet$} at 1.3 1.8  2.3 1.8 /
        \plot 1.3 1.8  2.3 1.8 / }
\def\cOooxoo{\multiput{} at 0 -1  0 6 /
        \multiput{\sq} at 0 0  0 1  0 2  0 3  0 4  0 5 /
        \multiput{$\sssize\bullet$} at -.5 2  1.5 2 /
        \plot -.5 2  1.5 2 / }
\def\cOooxooT{\multiput{} at 0 -1  0 6 /
        \multiput{\sq} at 0 0  0 1  0 2  0 3  0 4  0 5 /
        \multiput{$\sssize\bullet$} at -.5 2  1.5 2  2.5 2 /
        \plot -.5 2  2.5 2 / }
\def\cOooxooL{\multiput{} at 0 -1  0 6 /
        \multiput{\sq} at 0 0  0 1  0 2  0 3  0 4  0 5 /
        \multiput{$\sssize\bullet$} at -.5 2  .5 2 /
        \plot -.5 2  .5 2 / }
\def\cOooxooLS{\multiput{} at 0 -1  0 6 /
        \multiput{\sq} at 0 0  0 1  0 2  0 3  0 4  0 5 /
        \multiput{$\sssize\bullet$} at -.3 2.2  .7 2.2 /
        \plot -.3 2.2  .7 2.2 / }
\def\Ooo{\multiput{} at 0 -1  0 6 /
        \multiput{\sq} at 0 2  0 3  0 4  / }
\def\Xoo{\multiput{} at 0 -1  0 6 /
        \multiput{\sq} at 0 1  0 2  0 3 /
        \put{$\sssize\bullet$} at .5 3 }
\def\cXoo{\multiput{} at 0 -1  0 6 /
        \multiput{\sq} at 0 1  0 2  0 3 /
        \multiput{$\sssize\bullet$} at -.5 3  .5 3 /
        \plot -.5 3  .5 3 / }
\def\OcooXo{\multiput{} at 0 -1  0 6 /
        \multiput{\sq} at 0 1  0 2  0 3  0 4  1 2  1 3 /
        \multiput{$\sssize\bullet$} at .5 3  1.5 3 /
        \plot .5 3  1.5 3 / }
\def\cOcooXo{\multiput{} at 0 -1  0 6 /
        \multiput{\sq} at 0 1  0 2  0 3  0 4  1 2  1 3 /
        \multiput{$\sssize\bullet$} at -.5 3  .5 3  1.5 3 /
        \plot -.5 3  1.5 3 / }

%
%
\def\OocoooXo{\multiput{} at 0 -1  0 6 /
        \multiput{\sq} at 0 0  0 1  0 2  0 3  0 4  0 5  1 2  1 3 /
        \multiput{$\sssize\bullet$} at .5 3  1.5 3 /
        \plot .5 3  1.5 3 / }
\def\Ooxo{\multiput{} at 0 -1  0 6 /
        \multiput{\sq} at 0 1  0 2  0 3  0 4 /
        \put{$\sssize\bullet$} at .5 2  }
\def\cOocoooXo{\multiput{} at 0 -1  0 6 /
        \multiput{\sq} at 0 0  0 1  0 2  0 3  0 4  0 5  1 2  1 3 /
        \multiput{$\sssize\bullet$} at -.5 3  .5 3  1.5 3 /
        \plot -.5 3  1.5 3 / }
\def\cOoxo{\multiput{} at 0 -1  0 6 /
        \multiput{\sq} at 0 1  0 2  0 3  0 4 /
        \multiput{$\sssize\bullet$} at -.5 2  .5 2 / 
        \plot -.5 2  .5 2 / }
%
%
$$\hbox{\beginpicture\setcoordinatesystem units <.9cm,1.1cm>
        \multiput{\scale{\put{\Oooxoo} at 0 0  }} at -.1 0  5.9 0 /
        \put{\scale{\put{\Xoo} at 0 0 }} at 1.9 0
        \put{\scale{\put{\Ooo} at 0 0 }} at 3.9 0
        \put{\scale{\put{\OooxooR} at 0 0  
                        \put{\cXoo} at 1 0 }} at .9 1
        \put{\scale{\put{\OcooXo} at 0 0 }} at 2.9 1
        \put{\scale{\put{\Ooo} at 0 0 
                        \put{\cOooxooL} at 1 0 }} at 4.9 1
        \multiput{\scale{\put{\Ooo} at 0 0 
                        \put{\cOooxoo} at 1 0 
                        \put{\cXoo} at 2 0 }} at -.1 2  5.9 2 /
        \put{\scale{\put{\OooxooRR} at 0 0  
                        \put{\cOcooXo} at 1 0 }} at 1.9 2
        \put{\scale{\put{\OcooXo} at 0 0 
                        \put{\cOooxooL} at 2 0 }} at 3.9 2
        \put{\scale{\put{\Ooo} at 0 0 
                        \put{\cOooxooT} at 1 0 
                        \put{\cOcooXo} at 2 0 }} at .9 3
        \put{\scale{\put{\OooxooRS} at 0 0  
                        \put{\cOcooXo} at 1 0 
                        \put{\cOooxooLS} at 3 0 }} at 2.9 3
        \put{\scale{\put{\OcooXo} at 0 0 
                        \put{\cOooxoo} at 2 0 
                        \put{\cXoo} at 3 0 }} at 4.9 3
\arr{0.3 0.3} {0.7 0.7}
\arr{2.3 0.3} {2.7 0.7}
\arr{4.3 0.3} {4.7 0.7}
\arr{1.3 0.7} {1.7 0.3}
\arr{3.3 0.7} {3.7 0.3}
\arr{5.3 0.7} {5.7 0.3}
\arr{0.35 1.65} {0.7 1.3}
\arr{2.3 1.7} {2.7 1.3}
\arr{4.35 1.65} {4.7 1.3}
\arr{1.3 1.3} {1.7 1.7}
\arr{3.3 1.3} {3.65 1.65}
\arr{5.3 1.3} {5.7 1.7}
\arr{0.3 2.3} {0.7 2.7}
\arr{1.3 2.7} {1.7 2.3}
\arr{2.3 2.3} {2.55 2.55}
\arr{3.4 2.6} {3.65 2.35}
\arr{4.35 2.35} {4.65 2.65}
\arr{5.35 2.65} {5.65 2.35}
\arr{0.3 3.7} {0.6 3.4}
\arr{1.4 3.4} {1.7 3.7}
\arr{2.3 3.7} {2.55 3.45}
\arr{3.45 3.45} {3.7 3.7}
\arr{4.3 3.7} {4.6 3.4}
\arr{5.3 3.3} {5.7 3.7}
\setdashes<3pt>
\plot 0 0.6  0 1.4 /
\plot 0 2.6  0 4.6 /
\plot 6 0.6  6 1.4 /
\plot 6 2.6  6 4.6 /
\setdots<2pt>
\plot 0.5 0  1.5 0 /
\plot 2.5 0  3.5 0 /
\plot 4.5 0  5.5 0 /
\multiput{$\vdots$} at 1 4.5  3 4.5  5 4.5 /
\endpicture} \qquad\qquad
%
%
\hbox{\beginpicture\setcoordinatesystem units <.9cm,1.1cm>
        \multiput{\scale{\put{\OocoooXo} at 0 0  }} at -.1 0  3.9 0 /
        \put{\scale{\put{\Ooxo} at 0 0 }} at 1.9 0
        \put{\scale{\put{\OocoooXo} at 0 0  
                        \put{\cOoxo} at 2 0  }} at .9 1
        \put{\scale{\put{\Ooxo} at 0 0 
                        \put{\cOocoooXo} at 1 0 }} at 2.9 1
        \multiput{\scale{\put{\Ooxo} at 0 0 
                        \put{\cOocoooXo} at 1 0 
                        \put{\cOoxo} at 3 0 }} at -.1 2  3.9 2 /
        \put{\scale{\put{\OocoooXo} at 0 0  
                        \put{\cOoxo} at 2 0 
                        \put{\cOocoooXo} at 3 0 }} at 1.9 2 
        \put{\scale{\put{\Ooxo} at 0 0 
                        \put{\cOocoooXo} at 1 0 
                        \put{\cOoxo} at 3 0
                        \put{\cOocoooXo} at 4 0 }} at .9 3
        \put{\scale{\put{\OocoooXo} at 0 0  
                        \put{\cOoxo} at 2 0 
                        \put{\cOocoooXo} at 3 0 
                        \put{\cOoxo} at 5 0 }} at 2.9 3 
\arr{0.3 0.3} {0.6 0.6}
\arr{2.3 0.3} {2.7 0.7}
\arr{1.4 0.6} {1.7 0.3}
\arr{3.3 0.7} {3.7 0.3}
\arr{0.4 1.6} {0.6 1.4}
\arr{2.4 1.6} {2.7 1.3}
\arr{1.3 1.3} {1.6 1.6}
\arr{3.3 1.3} {3.6 1.6}
\arr{0.4 2.4} {0.6 2.6}
\arr{1.45 2.55} {1.6 2.4}
\arr{2.4 2.4} {2.55 2.55}
\arr{3.4 2.6} {3.6 2.4}
\arr{0.3 3.7} {0.6 3.4}
\arr{1.5 3.5} {1.7 3.7}
\arr{2.3 3.7} {2.5 3.5}
\arr{3.4 3.4} {3.7 3.7}
\setdashes<3pt>
\plot 0 0.6  0 1.4 /
\plot 0 2.6  0 4.6 /
\plot 4 0.6  4 1.4 /
\plot 4 2.6  4 4.6 /
\setdots<2pt>
\plot 0.5 0  1.5 0 /
\plot 2.5 0  3.5 0 /
\multiput{$\vdots$} at 1 4.5  3 4.5 /
\endpicture}$$

The remaining tubes in the tubular family $\bar{\Cal T}_0$
are in one-to-one correspondence with the monic irreducible
polynomials in $k[x]$ different from $x$ and $x-1$. 
In the case where $k$ is an algebraically closed field, the modules on the
mouth of those tubes are all as specified in Assertion (0.1.3).
Here we picture the first three modules in the tube corresponding to
$c\in k\backslash\{0,1\}$. We put $c'=1-c$.
{
\def\sixfourtwo{\multiput{} at 0 -1  0 6 /
        \multiput{\sq} at 0 0  0 1  0 2  0 3  0 4  0 5 
                1 1  1 2  1 3  1 4    2 2  2 3 /
        \multiput{$\bullet$} at .5 2  1.5 2 / 
        \plot .5 2  1.5 2 / }
\def\lambdaone{\multiput{} at 0 -1  0 6 /
        \multiput{$\bullet$} at 2.5 3 /
        \put{$\sssize c$} at .5 3
        \put{$\sssize c'$} at 1.5 3.1 
        \plot .8 3  1.1 3 /
        \plot 1.8 3  2.5 3 / }
\def\lambdaonelow{\multiput{} at 0 -1  0 6 /
        \multiput{$\bullet$} at 2.3 2.8 /
        \put{$\sssize c$} at .5 2.8
        \put{$\sssize c'$} at 1.5 2.9
        \plot .8 2.8  1.1 2.8 /
        \plot 1.8 2.8  2.3 2.8 / }
\def\lambdatwohigh{\multiput{} at 0 -1  0 6 /
        \multiput{$\bullet$} at -.3 3.15  2.3 3.15 /
        \put{$\sssize c$} at .5 3.15
        \put{$\sssize c'$} at 1.5 3.15 
        \plot -.3 3.15  .2 3.15 /
        \plot .8 3.15  1.1 3.15 /
        \plot 1.8 3.15  2.3 3.15 / }
\def\lambdatwolow{\multiput{} at 0 -1  0 6 /
        \multiput{$\bullet$} at -.3 2.8  2.3 2.8 /
        \put{$\sssize c$} at .5 2.8
        \put{$\sssize c'$} at 1.5 2.9
        \plot -.3 2.8  .2 2.8 /
        \plot .8 2.8  1.1 2.8 /
        \plot 1.8 2.8  2.3 2.8 / }

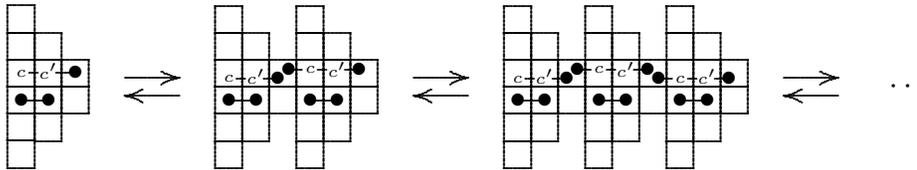
\def\scale#1{\beginpicture\setcoordinatesystem units <3mm,3mm>#1
        \endpicture}
%
%
$$\hbox{\beginpicture\setcoordinatesystem units <1cm,1cm>
        \put{\scale{\put{\sixfourtwo} at 0 0 
                    \put{\lambdaone} at 0 0 }} at 0 0 
        \put{\scale{\multiput{\sixfourtwo} at 0 0  3 0 /
                    \put{\lambdaonelow} at 0 0 
                    \put{\lambdatwohigh} at 3 0 }} at 2.75 0 
        \put{\scale{\multiput{\sixfourtwo} at 0 0  3 0  6 0 / 
                    \put{\lambdaonelow} at 0 0 
                    \put{\lambdatwohigh} at 3 0 
                    \put{\lambdatwolow} at 6 0 }} at 6.4 0 
        \arr{1.3 .1}{1.9 .1}
        \arr{1.9 -.1}{1.3 -.1}
        \arr{4.5 .1}{5.1 .1}
        \arr{5.1 -.1}{4.5 -.1}
        \arr{8.6 .1}{9.2 .1}
        \arr{9.2 -.1}{8.6 -.1}
        \put{$\cdots$} at 10 0
        \endpicture}$$
}

       \bigskip\noindent
{\bf (2.4) Dimension vectors under the covering functor.} In order to analyze the
covering functor $\pi\:\Cal S(\widetilde 6) \to \Cal S(6)$, we are going to look at 
dimension vectors. The covering functor induces a linear map from $K_0 = K_0(\mod \widetilde Q)$
to $K_0(\mod Q)$, which we also denote by $\pi$; here $\pi(\bx) = (\sum \bx_i,\sum \bx_{i'}).$
        \medskip\noindent
{\bf (2.4.1) Theorem.} {\it Let $X$ be an indecomposable object in $\Cal S(\widetilde 6)$
and put $\bx=\bdim X$.  Then }
$$\chi(\bx)=0 \qquad\text{if and only if} \qquad \pi(\bx) = s\cdot(12,6)\quad
        \text{for some $s\in\Bbb N$}.$$

        \medskip
We first verify the following 
        \smallskip
\noindent{\bf Claim.} {\it Let $(V,U,T)$ be an indecomposable triple in $\Cal S(6)$.
Then the dimension pair $(v,u)=\bdim(V,U,T)$ is in the shaded region pictured in (0.1). }

        \smallskip
According to Theorem (2.1.5), there is an  object $M\in\Cal D$ in the fundamental
domain for $\Cal S(\tilde6)$, unique up to isomorphism, which corresponds to
the triple $\pi(M)=(V,U,T)$ under the covering functor. 

Most of the objects in $\Cal D$ are modules over the algebra $\Theta$, the remaining
ones are on one of the two rays which we have inserted to obtain the tube
$\Cal T_\infty''$.  For the modules on the rays starting at
$\ssize{11111100 \atop 11111100}$ and at $\ssize{00000000 \atop 00111111}$,
the covering functor yields objects which have the following dimension pairs:
$$(6,6)\quad (5,5)\quad (10,5)\quad (10,6) \quad(11,7)\quad (12,7)\quad (18,12)
        \quad (17,11)\quad (22,11)\cdots$$
and
$$
 (6,0)\quad (6,1)\quad (7,2)\quad (8,2)\quad (14,7)\quad (13,6)\quad (18,6)
        \quad (18,7)\quad (19,8)\cdots
$$

The remaining indecomposable 
objects in $\Cal D$ are $\Theta$-modules.
The category of modules over the tubular algebra $\Theta$ is controlled by an
integral quadratic form~$\chi_\Theta$ ([9], 5.2 Theorem 6).
Thus, the dimension vector of an indecomposable $\Theta$-module is either a positive
root or a positive radical vector for $\chi_\Theta$. 
The radical for $\chi_\Theta$ is the free $\Bbb Z$-module generated by 
$$\bh^0={\ssize{12210\phantom{00}\atop 1233210}}\qquad \text{and}\qquad
        \bh^\infty={\ssize{01221\phantom{00}\atop 0123321}}.$$
Thus, if $M$ is a $\Theta$-module with dimension vector 
$\bdim M=a_0 \bh^0+a_\infty \bh^\infty$
in the radical of $\chi_\Theta$, then the corresponding object $\pi(M)\in\Cal S(6)$
has dimension pair $(a_0+a_\infty)\cdot (12,6)$.  It remains to deal with positive roots
of $\chi_\Theta$.  For $\bdim M$  such a root, put 
$\bm= \bdim M-(\dim M_0) \bh^0-(\dim M_{4'}) \bh^\infty$.
Since $\bh^0$, $\bh^\infty$ are in the radical, we have the equation
$$1=\chi_\Theta(\bdim M) = \chi_\Theta(\bm).$$
Thus, $\bm$ is a (positive or negative) root of $\chi_\Theta$ and,
since $\bm_0=\bm_{4'}=0$ also a root for the quadratic form of
the algebra $\Xi$ given by the following quiver with relations.

$$
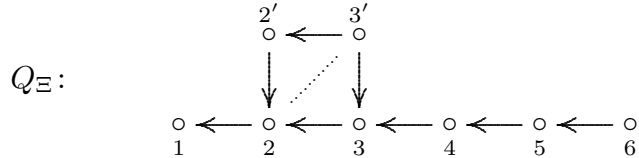
\hbox{\beginpicture
\setcoordinatesystem units <0.5cm,0.5cm>
\put{} at -1 0
\put{} at 5 2
\put{$Q_\Xi\:$} at -1 1
\put{$\circ$} at 2 0
\put{$\circ$} at 4 0
\put{$\circ$} at 6 0
\put{$\circ$} at 8 0
\put{$\circ$} at 10 0
\put{$\circ$} at 12 0
\put{$\circ$} at 4 2
\put{$\circ$} at 6 2
\arr{3.6 0}{2.4 0}
\arr{5.6 0}{4.4 0}
\arr{7.6 0}{6.4 0}
\arr{9.6 0}{8.4 0}
\arr{11.6 0}{10.4 0}

\arr{5.6 2}{4.4 2}

\arr{4 1.6}{4 0.4}
\arr{6 1.6}{6 0.4}


\put{$\ssize 2'$} at  4 2.5
\put{$\ssize 3'$} at  6 2.5

\put{$\ssize 1$} at 2 -.5
\put{$\ssize 2$} at 4 -.5
\put{$\ssize 3$} at 6 -.5
\put{$\ssize 4$} at 8 -.5
\put{$\ssize 5$} at 10 -.5
\put{$\ssize 6$} at 12 -.5

\setdots<2pt>
\plot 5.5 1.5  4.5 0.5 /
\endpicture}
$$
This algebra is tilted of type $E_8$ so the roots of $\chi_\Theta$ are obtained 
from the roots of $E_8$.
Namely, if $\ssize{\phantom{a\,b}\,b'\phantom{c\,d\,e\,f}\atop a\,b\,g\,c\,d\,e\,f}$ 
is a root for
$E_8$ then the corresponding root for $\chi_\Theta$ is
$\ssize{\phantom{a}\,b'\!c'\phantom{d\,e\,f}\atop a\,b\,c\,d\,e\,f}$ where
$c'=b'+c-g$, and gives rise to the dimension pair $(v,u)=(a+b+c+d+e+f,a+b'+c')$.
The roots for $E_8$ can be easily calculated, or else they can be found in [4], Planche VII.

In the following diagram we take each of the 120 positive roots of $E_8$ 
and mark the corresponding dimension pair by putting a dot at position $(v,u)$.
        For each root corresponding to this position
        we label the dot by the largest entry $m$ of the root.

        \smallskip
        Most roots (e.g.~all roots where $m\geq 2$) correspond
        to indecomposable objects in $\Cal S(\widetilde 6)$
        and hence via $\pi$ to indecomposable triples
        in $\Cal S(6)$.
        If $(V,U,T)$ is such a triple then $m$ is just the
        Krull-Remak-Schmidt multiplicity of the total space,
        i.e.~the number of indecomposable summands of the
        $k[x]$-module $(V,T)$.

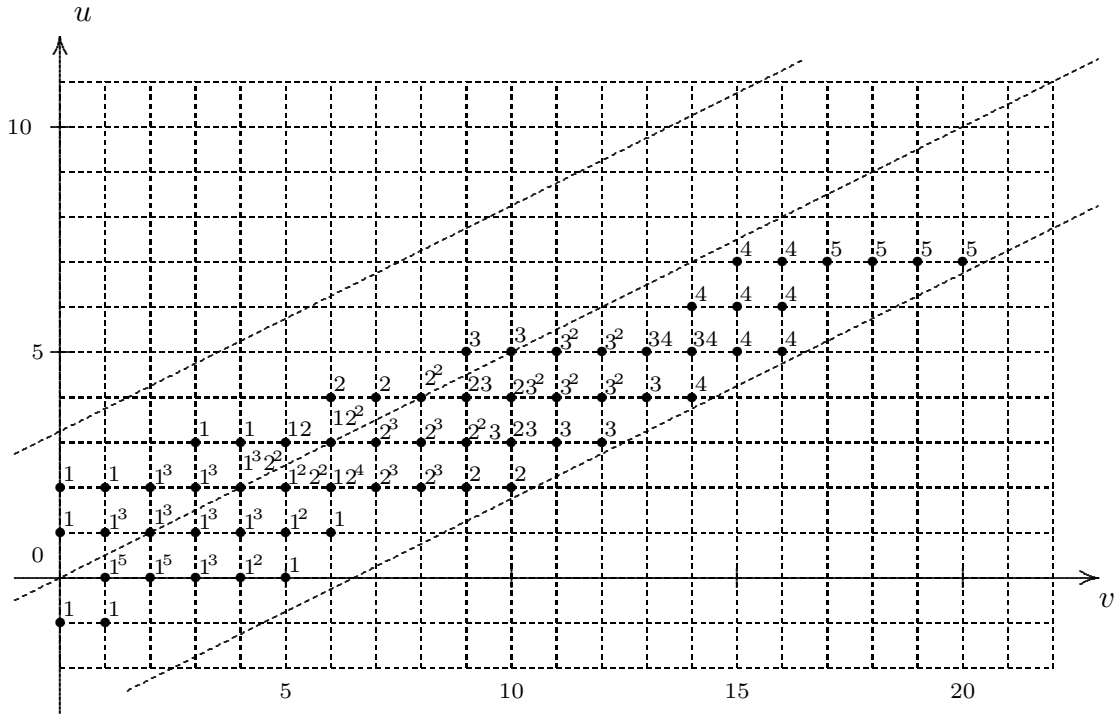
$$\hbox{\beginpicture
\setcoordinatesystem units <.5cm,.5cm>
\arr{0 -3}{0 12}
\arr{-1 0}{23 0}
\put{$u$} at  0.5 12.5
\put{$v$} at 23.2 -.5
\setdashes <.5mm> 
\setplotarea x from 0 to 22, y from -2 to 11
\grid {22} {13}
\plot -1 -.5  23 11.5 /
\plot 1.5 -2.5   23 8.25 /
\plot -1 2.75  16.5 11.5 /
\setsolid
\plot 5 -0.2  5 0.2 /
\put{$\ssize 5$} at 5 -2.5
\plot 10 -0.2  10 0.2 /
\put{$\ssize 10$} at 10 -2.5
\plot 15 -0.2  15 0.2 /
\put{$\ssize 15$} at 15 -2.5
\plot 20 -0.2  20 0.2 /
\put{$\ssize 20$} at 20 -2.5
\plot -0.2 5  0.2 5 /
\put{$\ssize 5$} at -.5 5 
\put{$\ssize 0$} at -.5 .5 
\plot -0.2 10  0.2 10 /
\put{$\ssize 10$} at -.9 10 
%
\multiput{$\ssize \bullet$} at 0 -1  0 1  0 2  1 -1  1 0  1 1  1 2  
                               2 0  2 1  2 2  3 0  3 1  3 2  3 3
                               4 0  4 1  4 2  4 3  5 0  5 1  5 2  5 3  
                               6 1  6 2  6 3  6 4  7 2  7 3  7 4  
                               8 2  8 3  8 4  9 2  9 3  9 4  9 5
                               10 2  10 3  10 4  10 5  11 3  11 4  11 5  
                               12 3  12 4  12 5  13 4  13 5
                               14 4  14 5  14 6  15 5  15 6  15 7
                               16 5  16 6  16 7  17 7  18 7  19 7  20 7 /
\multiput{$\ssize 5$} at 17.2 7.3  18.2 7.3  19.2 7.3  20.2 7.3 /
\multiput{$\ssize 4$} at 14.2 4.3  14.2 6.3  15.2 5.3  15.2 6.3  15.2 7.3
                         16.2 5.3  16.2 6.3  16.2 7.3 /
\multiput{$\ssize 34$} at 13.3 5.3  14.3 5.3 /
\multiput{$\ssize 3^{\!2}$} at 11.3 4.3  11.3 5.3  12.3 4.3  12.3 5.3 /
\multiput{$\ssize 3$} at 9.2 5.3  10.2 5.4  11.2 3.3  12.2 3.3  13.2 4.3 /
\multiput{$\ssize 23$} at 9.3 4.3  10.3 3.3 /
\multiput{$\ssize 23^{\!2}$} at 10.4 4.3 /
\multiput{$\ssize 2^{\!2}3$} at 9.4 3.3 /
\multiput{$\ssize 2$} at 6.2 4.3  7.2 4.3  9.2 2.3  10.2 2.3 /
\multiput{$\ssize 2^{\!2}$} at 8.3 4.5 /
\multiput{$\ssize 2^{\!3}$} at 7.3 2.3  7.3 3.3  8.3 2.3  8.3 3.3  /
\multiput{$\ssize 12$} at 5.3 3.3  /
\multiput{$\ssize 12^{\!2}$} at 6.4 3.6 /
\multiput{$\ssize 1^{\!3}2^{\!2}$} at 4.5 2.6 /
\multiput{$\ssize 1^{\!2}2^{\!2}$} at 5.5 2.3 / 
\multiput{$\ssize 12^{\!4}$} at 6.4  2.3 /
\multiput{$\ssize 1$} at .2 -.7   .2 1.3   .2 2.3  
                         1.2 -.7  1.2 2.3  3.2 3.3  4.2 3.3  5.2 0.3  6.2 1.3 /
\multiput{$\ssize 1^{\!2}$} at 4.3 .3  5.3 1.3  /
\multiput{$\ssize 1^{\!3}$} at 1.3 1.3  2.3 1.4  2.3 2.3  3.3 0.3  
                               3.3 1.3  3.3 2.3  4.3 1.3  /
\multiput{$\ssize 1^{\!4}$} at   /
\multiput{$\ssize 1^{\!5}$} at 1.3 0.3  2.3 0.3 /

\endpicture}
$$

Returning to the triple $(V,U,T)=\pi(M)$ in $\Cal D$,
we need to show that
$$\big|\dim V-2\dim U\big|\leq 6,$$
or, equivalently, 
$$\big|\sum_i \dim M_i - 2\sum_i\dim M_{i'}\big|\leq 6.$$
Let us first consider the modules on the two rays. Given such a module, 
the next module on that ray 
is either obtained as an extension
with a module on the mouth of the tube, 
or by applying a canonical map with kernel with dimension pair
$(1,1)$.  After six such subsequent operations, the dimension pair increases by $(12,6)$. 
It follows that the dimensions of all modules on both rays satisfy the above inequality.

We now deal with the remaining indecomposable 
objects in $\Cal D$, which are $\Theta$-modules.

Clearly, all roots of $\chi_\Theta$
lie in the region
given by $|v-2u|\leq 6$.  This finishes the proof of the claim.

        \bigskip\bigskip
{\bf 3. Proof of the Crucial Results.}
        \smallskip
In this chapter we verify the assertions stated in (0.1) in the introduction.
        \medskip\noindent
{\bf (3.1) Proof of Assertion (0.1.1).} 
We have seen in (2.4) that every dimension pair $(v,u)$ of an indecomposable
triple $(V,U,T)$ in $\Cal S(6)$ is contained in the shaded region pictured in the
introduction.
It remains to verify that every point $(v,u)$ in the shaded region can be
realized as the dimension pair $\bdim(V,U,T)=(v,u)$ of such an object.
Since there is the duality $*$, we only have to consider pairs
$(v,u)$ with $2u \le v.$ 
Moreover, if $u = 0$, then $1 \le v \le 6$ and the triple $(V,0,T)$ where
$V$ is given by the partition $(v)$ is indecomposable.
Similarly, if $u = 1,$ then $2 \le v \le 6$ and again the triple $(V,\ker T,T)$ 
where $V$ has partition $(v)$ is indecomposable. 
Thus we can assume that $u \ge 2.$
        \medskip
Consider again the algebra $\Xi$ given by the following quiver with relations.
$$
\hbox{\beginpicture
\setcoordinatesystem units <0.5cm,0.5cm>
\put{} at -1 0
\put{} at 5 2
\put{$Q_\Xi\:$} at -1 1
\put{$\circ$} at 2 0
\put{$\circ$} at 4 0
\put{$\circ$} at 6 0
\put{$\circ$} at 8 0
\put{$\circ$} at 10 0
\put{$\circ$} at 12 0
\put{$\circ$} at 4 2
\put{$\circ$} at 6 2
\arr{3.6 0}{2.4 0}
\arr{5.6 0}{4.4 0}
\arr{7.6 0}{6.4 0}
\arr{9.6 0}{8.4 0}
\arr{11.6 0}{10.4 0}

\arr{5.6 2}{4.4 2}

\arr{4 1.6}{4 0.4}
\arr{6 1.6}{6 0.4}


\put{$\ssize 2'$} at  4 2.5
\put{$\ssize 3'$} at  6 2.5

\put{$\ssize 1$} at 2 -.5
\put{$\ssize 2$} at 4 -.5
\put{$\ssize 3$} at 6 -.5
\put{$\ssize 4$} at 8 -.5
\put{$\ssize 5$} at 10 -.5
\put{$\ssize 6$} at 12 -.5

\setdots<2pt>
\plot 5.5 1.5  4.5 0.5 /
\endpicture}
$$
Note that there are indecomposable $\Xi$-modules which give rise under $\pi$ to objects
in $\Cal S(6)$ with dimension pairs at the following positions. 
$$\hbox{\beginpicture
\setcoordinatesystem units <.4cm,.4cm>
\arr{0 -1}{0 12}
\arr{-1 0}{23 0}
\put{$u$} at  0.5 12.5
\put{$v$} at 23.2 -.5
\setdashes <.5mm> 
\setplotarea x from 0 to 22, y from 0 to 11
\grid {22} {11}
\plot -1 -.5  23 11.5 /
\plot 5.5 -.5   23 8.25 /
\setsolid
\plot 5 -0.2  5 0.2 /
\put{$\ssize 5$} at 5 -.5
\plot 10 -0.2  10 0.2 /
\put{$\ssize 10$} at 10 -.5
\plot 15 -0.2  15 0.2 /
\put{$\ssize 15$} at 15 -.5
\plot 20 -0.2  20 0.2 /
\put{$\ssize 20$} at 20 -.5
\plot -0.2 5  0.2 5 /
\put{$\ssize 5$} at -.5 5 
\put{$\ssize 0$} at -.5 .5 
\plot -0.2 10  0.2 10 /
\put{$\ssize 10$} at -.9 10 
%
\multiput{$\ssize \bullet$} at 1 0  5 0 
                               2 1  
                               4 2   5 2 
                               6 2  6 3   7 2  7 3 
                               8 2  8 3  8 4  9 2  9 3  9 4 
                               10 2  10 3  10 4  10 5  11 3  11 4  11 5  
                               12 3  12 4  12 5  13 4  13 5
                               14 4  14 5  14 6  15 5  15 6  15 7
                               16 5  16 6  16 7  17 7  18 7  19 7  20 7 /
\multiput{$\ssize \bigcirc$} at 13 6  14 7  17 6 /
\endpicture}
$$
(In addition to $\Xi$-modules which give rise to triples in $\Cal S(6)$ with dimension pairs
$(v,u)$ where $u \ge 2$, we need the three dimension pairs $(1,0),\ (2,1),\ (5,0)$ in
order to fill the holes $(13,6),\ (14,7),\ (17,6)$ indicated by circles.)

The algebra $\Xi$ has the following algebra $\Xi'$ as one-point extension:
$$
\hbox{\beginpicture
\setcoordinatesystem units <0.5cm,0.5cm>
\put{} at -1 0
\put{} at 5 2
\put{$Q_{\Xi'}\:$} at -1 1
\put{$\circ$} at 2 0
\put{$\circ$} at 4 0
\put{$\circ$} at 6 0
\put{$\circ$} at 8 0
\put{$\circ$} at 10 0
\put{$\circ$} at 12 0
\put{$\circ$} at 4 2
\put{$\circ$} at 6 2
\put{$\circ$} at 8 2
\arr{3.6 0}{2.4 0}
\arr{5.6 0}{4.4 0}
\arr{7.6 0}{6.4 0}
\arr{9.6 0}{8.4 0}
\arr{11.6 0}{10.4 0}

\arr{5.6 2}{4.4 2}
\arr{7.6 2}{6.4 2}

\arr{4 1.6}{4 0.4}
\arr{6 1.6}{6 0.4}
\arr{8 1.6}{8 0.4}

\put{$\ssize 2'$} at  4 2.5
\put{$\ssize 3'$} at  6 2.5
\put{$\ssize 4'$} at  8 2.5

\put{$\ssize 1$} at 2 -.5
\put{$\ssize 2$} at 4 -.5
\put{$\ssize 3$} at 6 -.5
\put{$\ssize 4$} at 8 -.5
\put{$\ssize 5$} at 10 -.5
\put{$\ssize 6$} at 12 -.5

\setdots<2pt>
\plot 5.5 1.5  4.5 0.5 /
\plot 7.5 1.5  6.5 0.5 /
\endpicture}
$$
This is a tame concealed algebra with radical generator 
$$\bh^0=\ssize{122100\atop 123321}.$$
        \medskip
Let $M$ be an indecomposable $\Xi$-module such that the corresponding dimension pair $(v,u)$ is labelled
by the symbol $\bullet$ in the above diagram.  If $M$ is preprojective or preinjective as a $\Xi'$-module
then the lemma below yields a realization for each dimension pair
in the orbit of $(v,u)$ under the shift by $(12,6)$.  In case $M$ is a regular $\Xi'$-module,
then the tube of $M$ contains indecomposable modules which realize each dimension pair in the 
orbit of $(v,u)$.  So it remains to deal with the 
orbits of the pairs
$$
 (12,6), (18,6).
$$
        \bigskip
For the orbit of $(12,6)$, take the indecomposable $\Xi'$-modules with dimension vector in the radical.
        \bigskip
Finally we consider  $(18,6)$. For any $t \ge 0$, one finds an indecomposable object 
with dimension pair $(18,6)+(12t,6t)$ in the non-stable tube in $\Cal U$. 
Note that these objects do not come from $\Xi'$-modules.
For example, the unique object in the non-stable tube with dimension pair $(18,6)$ is obtained via the
covering functor from a representation with dimension vector
$$\ssize{1221000\atop 1344321}.$$
        \medskip
Recall that a tame concealed algebra $\Lambda$ is the endomorphism ring of a preprojective 
tilting module over a tame hereditary algebra. An indecomposable $\Lambda$-module is said
to be {\it homogeneous} provided it is fixed under the Auslander-Reiten translation.

        \medskip\noindent
{\bf (3.1.1) Lemma.} {\it Let $\Lambda$ be a tame concealed algebra and
$H$ an indecomposable homogeneous $\Lambda$-module. 
For $M$ an indecomposable $\Lambda$-module in the preprojective or the preinjective
component of the Auslander-Reiten quiver there is an exact sequence
$$
 0 \to M \to N \to H \to 0 \qquad\text{or} \qquad 0 \to H \to N \to M \to 0,
$$
respectively, with $N$ indecomposable.}
        \medskip\noindent
{\it Proof:} 
This is well-known for $\Lambda$ being hereditary and we will use tilting theory for the 
general case. Consider first the case where $M$ is in the preprojective component 
of the Auslander-Reiten quiver of $\Lambda$. Write $\Lambda$ as the
endomorphism ring of a preprojective tilting $\Lambda'$-module $T$, where $\Lambda'$ is
hereditary. Consider the torsion pair $(\Cal F,\Cal G)$ in $\mod \Lambda'$, where $\Cal G$
are the modules generated by $T$. Both $M$ and $H$ are in the image of the tilting functor
$\Hom(T,-)\:\Cal G \to \mod\Lambda$, thus there are indecomposable $\Lambda'$-modules
$M', H'$ with $M = \Hom(T,M')$ and $H = \Hom(T,H')$, and $M'$ is a preprojective or regular
$\Lambda'$-module, whereas $H'$ is again homogeneous. In $\mod \Lambda'$, we find an exact
sequence $0 \to M' \to N' \to H' \to 0$ with $N'$ indecomposable. Take $N = \Hom(T,N');$
this is an indecomposable $\Lambda$-module and under $\Hom(T,-)$ we obtain an exact sequence
of the form $0 \to M \to N \to H \to 0$. 
This completes the proof in the first case. In case $M$ belongs to the preinjective component 
of $\mod\Lambda$, we just use duality and obtain an exact sequence of the second kind.
\qed

        \bigskip\noindent
{\bf (3.2) Proof of Assertions (0.1.2) and (0.1.3).} 
In essence, both statements are a consequence of [9], Theorem~5.2.6, which states that 
the category of modules over a tubular algebra $\Theta$
is controlled by the associated integral quadratic form~$\chi_\Theta$.

\smallskip Let $M$ be an indecomposable $\Theta$-module where 
$\Theta$ is the tubular algebra from Section~1. By [9], Theorem~5.2.6, 
the dimension vector $\bd$ of $M$ is either a root or a radical
vector for $\chi_\Theta$.  In the first case $M$ is determined 
uniquely, up to isomorphism, by $\bd$ and conversely,
every vector $\bd$ which is a root of $\chi_\Theta$
can be realized by an indecomposable module.   
If $\bd$ is a radical vector for $\chi_\Theta$ then 
there are infinitely many isomorphism classes
of indecomposable modules of dimension vector $\bd$ provided
only that the underlying base field $k$ has infinite cardinality.

\smallskip Let us deal with radical vectors first.
According to [9], Theorem 5.1.3(c), the radical vectors
corresponding to index $\gamma$ form a free abelian group
of rank 1 containing 
$$\br_\gamma\quad =\quad q\cdot \bh^0
        + p\cdot \bh^\infty$$
where $\gamma=\frac pq$ with $p$ and $q$ non-negative integers 
which are relatively prime. 
It follows that $\br_\gamma$ is a generator of this group.
Suppose that $k$ is a field of infinite cardinality,
$\gamma=\frac pq\in \Bbb Q_0^+$ as above, and $m\in\Bbb N$.
Then there are infinitely
many indecomposable objects, up to isomorphism, 
of dimension vector $m \br_\gamma$ in $\Cal T_0'\sqcup\Cal T$.
Under the covering functor, they correspond to an infinite
family of pairwise non-isomorphic objects with dimension pair
$m (p+q)\cdot (12,6)$ in $\bar{\Cal T}_\gamma$. 
In particular if $\gamma=0$, then for each $m\in\Bbb N$ 
and $c\in k\backslash\{0,1\}$ we have pictured at the end of (2.3)
an indecomposable module with dimension pair $m\cdot (12,6)$
which occurs in the homogeneous tube in $\bar{\Cal T}_0$ 
corresponding to the monic irreducible polynomial $x-c$.

\smallskip 
Now suppose that a dimension pair $(v,u)$ is given which is not
an integer multiple of $(12,6)$.  We show that there are 
only finitely many indecomposable subspace configurations
with this dimension pair, up to isomorphism, and that this number 
$n(v,u)$ does not depend on the choice of the base field $k$.

\smallskip Let $n'(v,u)$ be the number of roots of $\chi_\Theta$ 
for which the corresponding index $\gamma$ is in $\Bbb Q^+$.
(Clearly, this number is finite since up to the shift, there are only finitely
many connected dimension vectors which are mapped under $\pi$ to $(v,u)$.)
This is the number of indecomposable objects in 
$\Cal T=\bigsqcup_{\gamma\in\Bbb Q^+}\Cal T_\gamma$ which under the
covering functor correspond to an indecomposable 
subspace configuration with dimension pair $(v,u)$. 
There may be further subspace configurations of this dimension, 
but they must correspond to objects in the tubular family
$\bar{\Cal T}_0$, and hence occur in one of the three
non-homogeneous tubes which we have pictured in (2.3). 
Say there are $n_0(v,u)$ such objects; again, this number
does not depend on the choice of the base field $k$ 
because the operations which we performed on the tubular
family $\Cal T'_\infty$ work for every base field, 
and the number is finite because
within each tube the dimensions tend to infinity with 
increasing distance to the mouth.  Then we can compute
$n(v,u)$ as the sum $n'(v,u)+n_0(v,u)$. We have seen that
this number is finite and independent of $k$.

        \medskip \noindent
{\bf (3.3) Proof of Assertion (0.1.4).}
We start with the following observation.  Let $\gamma=\frac pq\in \Bbb Q_0^+$
be written in lowest terms.  The modules on the mouth of the homogeneous
tubes in $\Cal T_\gamma$, and some modules in the non-homogeneous tubes,
form a one-parameter family of indecomposables of dimension vector
$q\bh^0 + p\bh^\infty$.
Moreover, the remaining modules with this dimension vector in $\Cal T_\gamma$
are in the non-homogeneous tubes; there are $6+3+2$ of them if $\gamma\neq 0$
and $4+3+2$ if $\gamma=0$. 
Under the covering functor, the family in $\Cal T_\gamma$ corresponds 
to a family in $\bar \Cal T_\gamma$
with dimension pair $(p+q)\cdot (12,6)$. 
Similarly for $s>0$, the $s$-fold selfextension of each of these modules is an object
with dimension vector 
$$
 s(q\bh^0+p\bh^\infty) 
$$
in $\Cal T_\gamma$.
Those modules give rise to a family of indecomposables with dimension pair
$s(p+q)\cdot (12,6)$.

        \smallskip
Thus given $t>0$, any $p\in\{0,1,\ldots,t-1\}$ gives rise to a family of indecomposables
with dimension pair $t\cdot (12,6)$: Put $\gamma=\frac p{t-p}$ and let $s$ be the
greatest common divisor of $p$ and $t-p$.
Then there is a one-parameter family in $\bar\Cal T_\gamma$ containing 
the $s$-fold selfextensions of the modules on the mouths of the homogeneous tubes.

        \smallskip
It remains to show that there are only finitely many indecomposables with
dimension pair $t\cdot (12,6)$ which are not in one of the families.
In fact, we can describe the possible exceptions:  If $\gamma=\frac p{t-p}$ as above,
then for each family in $\Cal T_\gamma$, 
there are at most $6+3+2$ (or $4+3+2$ if $\gamma=0$) indecomposables in the
non-homogeneous tubes in $\Cal T_\gamma$ which are not in the family, 
their distance from the mouth is the product of
the rank with the greatest common divisor of $p$ and $t-p$.
        \smallskip
Suppose a module $\bar M$ in one of the non-homogeneous tubes has distance
from the mouth which is not a multiple of the rank of the tube.  Then there
is a corresponding module $M\in\Cal D$ such that $\pi(M)=\bar M$; this module $M$
occurs in a non-homogeneous tube of the same rank 
and has the same distance from the mouth of its tube.
Then either $M$ is a $\Theta$-module and $\bdim M$ is a root of the quadratic form 
for $\Theta$ or $M$ is in one of the inserted rays.  We have seen in the proof of
Assertion (0.1.1) that in neither of the two cases the dimension pair of $\bar M$ 
is a multiple of $(12,6)$. 

        \medskip \noindent
{\bf (3.4) Proof of Assertions (0.1.5) and (0.1.6).} As in the case of tubular algebras one
observes that the infinite radical $\widetilde {\Cal I}$
of $\Cal S(\widetilde 6)$ is the following ideal: Let $M$ belong to $\Cal T_\gamma$
and $M'$ to $\Cal T_{\gamma'}$, then
$$
 \widetilde{\Cal I}(M,M') = \left\{\matrix 0 &\text{if} & \gamma \ge \gamma' \cr
                       \Hom(M,M') &\text{if} & \gamma < \gamma'\endmatrix \right. .
$$
Also, the ideal $\widetilde{\Cal I}$ is idempotent. Namely, consider a map
$f\:M \to M'$ where $\gamma < \gamma'$, and choose some rational number $\gamma''$
with $\gamma < \gamma'' < \gamma.$ One knows that $\Cal T_{\gamma''}$ 
is a separating tubular family, thus there is a module $M''$ in $\Cal T_{\gamma''}$ 
such that $f$ factors through $M''$, say $f = f_2f_1$ with $f_1\:M \to M''$
and $f_2\:M'' \to M'.$ Since both $f_1$ and $f_2$ belong to $\widetilde{\Cal I}$,
we see that $f$ belongs to ${\widetilde{\Cal I}}^{\,2}$.

Under the covering functor $\pi\:\Cal S(\widetilde 6) \to \Cal S(6)$, the infinite radical
$\widetilde {\Cal I}$ of $\Cal S(\widetilde 6)$ is mapped onto the infinite radical
$\Cal I$ of $\Cal S(6)$. Since $\widetilde {\Cal I}$ is idempotent, also $\Cal I$ is
idempotent. 

Let $X$ be an indecomposable object of $\Cal S(6)$. Choose an indecomposable object $M$ of 
$\Cal S(\widetilde 6)$ with $\pi(M) = X.$ The covering property of $\pi$ asserts that
one can identify $\End(X) = \Hom(X,X)$ under $\pi$ with the graded ring
$$
 \Pi = \bigoplus_{i\in \Bbb Z} \Pi_i, \t{where} \Pi_i = \Hom(M,M[i]).
$$
Here, the product of $f\in \Hom(M,M[i])$ and $g\in \Hom(M,M[j])$ is 
$f[j]\circ g \in \Hom(M,M[i+j]).$ As we know, $\Hom(M,M[i]) = 0$ for $i < 0$, thus
$\Pi = \bigoplus_{i\in \Bbb N_0} \Pi_i.$ Also, $\Hom(M,M[i]) = 0$ for 
$i \ge 8$, since we can assume
that the support of the total space of $M$ is contained in $[0,7].$
This shows that the ideal $\Pi_+ = \bigoplus_{i>0} \Pi_i$  of $\Pi$ is nilpotent with
nilpotency index at most $8$. However under the identification of $\Pi$ and $\End(X)$, 
the ideal $\Pi_+$ corresponds to $I = \Cal I(X,X)$. Also, under this identification, 
$\Pi/\Pi_+ = \End(M)$ is isomorphic to $\End(X)/I$. But $\End(M)$ is a local uniserial
ring for any indecomposable object in $\Cal S(\widetilde n).$ This completes the proof.

Note that for most of the indecomposable objects $X$ in $\Cal S(6)$ we even must have
$\Cal I(X,X)^7 = 0,$ since usually the support of the total space is $M$ is 
of length at most 7, so that $\Pi_+$  has nilpotency
index at most 7. The only indecomposable objects $M$ in $\Cal S(\widetilde 6)$ such that
the total space of $M$ has support of length 8 have been exhibited in Section (2.3).

%
%
Let us present here an object $X\in\Cal S(6)$ and a non-zero endomorphism
$$\varepsilon\in\Cal I(X,X)^7.$$
Let $M$ in $\Cal S(\widetilde 6)$ be the object in the non-stable tube in $\Cal D$ 
at the second intersection of the ray starting at $P(7)$ and the coray ending at $P(5')$. 
Then $X=(A,A')=\pi(M)$ can be pictured as follows. (The shaded box will be explained below.)
$$      
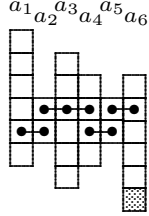
\beginpicture 
        \setcoordinatesystem units <0.25cm,0.25cm>
        \multiput{} at -1 0  7 8.5 /
        \multiput{\sq} at 0 2  0 3  0 4  0 5  0 6  0 7  1 3  1 4  
                2 1  2 2  2 3  2 4  2 5  2 6  3 2  3 3  3 4  3 5
                4 3  4 4  5 0  5 1  5 2  5 3  5 4  5 5  / 
        \multiput{$\ssize \bullet$} at .5 3  1.5 3  1.5 4  2.5 4  3.5 4
                3.5 3  4.5 3  4.5 4  5.5 4 /
        \plot 0.5 3  1.5 3 /
        \plot 1.5 4  3.5 4 /
        \plot 3.5 3  4.5 3 /
        \plot 4.5 4  5.5 4 /
        \put{$\ssize a_1$} at .5 8.5 
        \put{$\ssize a_2$} at 1.6 8
        \put{$\ssize a_3$} at 2.5 8.5
        \put{$\ssize a_4$} at 3.6 8
        \put{$\ssize a_5$} at 4.5 8.5
        \put{$\ssize a_6$} at 5.6 8
        \setshadegrid span <.3mm>
        \vshade  5  -.5 .5 <,z,,> 
                 6  -.5 .5 /
        \endpicture
$$
Write $A=k[x]/x^6\oplus k[x]/x^2\oplus k[x]/x^6\oplus k[x]/x^4\oplus k[x]/x^2\oplus k[x]/x^6$
and denote by $a_i$ a generator of the $i$-th summand.  Then $A'$ is generated by
$a_1x^4+a_2x$, $a_2+a_3x^2+a_4x$, $a_4x^2+a_5x$, and $a_5+a_6x$.
Define $\varphi:A\to A$ by putting $\varphi(a_1)=a_3$, $\varphi(a_2)=a_4x^2$, $\varphi(a_3)=a_6$
and $\varphi(a_i)=0$ for $i\geq 4$, this defines an endomorphism of $(A,A')$. 
Note that $\varphi^2$ is the map given by $\varphi^2(a_1)=a_6$ and $\varphi^2(a_i)=0$ for $i\geq2$.
In fact, $\varphi=\pi(f)$ for a map $f:M\to M[1]$, so $\varphi$ is in the infinite radical.
Let $T\in\End(A,A')$ be given by multiplication by $x$, so $T$ is in the infinite radical, too. 
Then the composition
$$\varepsilon=T^5\circ\varphi^2\in\Cal I(X,X)^7$$
is non-zero, more precisely, the image of $\varepsilon$ is just $k\,a_6x^5$,
pictured as the shaded box in the above diagram.

        \medskip \noindent
{\bf (3.5) Proof of Assertion (0.1.7).} Let $X$ be an indecomposable object in $\Cal S(6)$
with dimension pair $(v,u).$  
First, assume that $(v,u)$ is an integral multiple of $(12,6).$ Let $M$ be an indecomposable
object in $\Cal S(\widetilde 6)$ with $\pi(M) = X.$ Theorem (2.4.1) asserts $\chi(\bdim M) = 0.$
But then there is an increasing chain of inclusions
$$
 0 = N_0 \subset N_1 \subset \cdots \subset N_m \subset \cdots
$$
of indecomposable objects in $\Cal S(\widetilde 6)$ such that $N_i/N_{i-1}$ is isomorphic to $M$
for all $i \ge 1.$ Using the covering functor $\pi$ we obtain 
indecomposable objects $Y_i = \pi(N_i)$ in $\Cal S(6)$. Since $\pi$
is exact, it follows that $Y_i/Y_{i-1}$ is isomorphic to $X$ for all $i\ge 1$.
        \smallskip
Next we deal with the case where $(v,u)$ is not a multiple of $(2,1)$. Assume 
there is a chain of inflations 
$$
 0 = Y_0 \subset Y_1 \subset \cdots \subset Y_m
$$
with $Y_i/Y_{i-1}$ isomorphic to $X$ for all $1\le i \le m.$ The dimension pair of $Y_m$ is
$(mv,mu)$ and 
$$
 |mv - 2mvu| = m|v-2u| \ge m,
$$
sind $v-2u$ is non-zero. According to Assertion (0.1.1), we see that $m \le 6$.
        \medskip
We have mentioned in (0.1) that one may expect that an indecomposable object in $\Cal S(6)$
with dimension pair not an integral multiple of $(12,6)$ can be $m$-stackable only for $m\le 6.$
If one tries to prove this assertion, one has to be aware that when dealing with a covering functor,
usually new exact sequences are created: In general, monomorphisms cannot be lifted to monomorphisms.
        \bigskip\bigskip

{\bf 4. A Sort of Double Covering for $\Cal S(6)$}
        \medskip

In Section~1, all stable tubes in the fundamental domain
of the category $\Cal S(\widetilde6)$ 
have been obtained as tubes of regular modules over the
tubular algebra $\Theta$.  In this section we exhibit a slightly more complicated
finite dimensional algebra $\Omega$
which provides a corresponding model for two adjacent fundamental domains of 
$\Cal S(\widetilde 6)$. Thus the regular $\Omega$-modules form a sort of double covering for 
$\Cal S(6)$.
\smallskip 
\noindent {\bf (4.1) The algebra $\Omega$.}
Here is the algebra $\Omega$, given by the following quiver with relations
$$
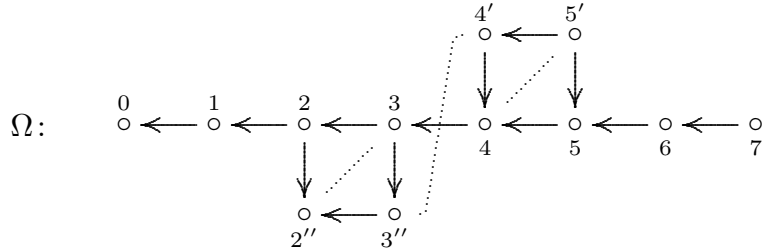
\hbox{\beginpicture
\setcoordinatesystem units <0.5cm,0.5cm>
\put{} at 0  -2
\put{} at 14  2
\put{$\Omega\:$} at -2 0
\put{$\circ$} at 0 0
\put{$\circ$} at 2 0
\put{$\circ$} at 4 0
\put{$\circ$} at 6 0
\put{$\circ$} at 8 0
\put{$\circ$} at 10 0
\put{$\circ$} at 12 0
\put{$\circ$} at 14 0

\put{$\circ$} at 8 2
\put{$\circ$} at 10 2

\put{$\circ$} at 4 -2
\put{$\circ$} at 6 -2

\arr{1.6 0}{0.4 0}
\arr{3.6 0}{2.4 0}
\arr{5.6 0}{4.4 0}
\arr{7.6 0}{6.4 0}
\arr{9.6 0}{8.4 0}
\arr{11.6 0}{10.4 0}
\arr{13.6 0}{12.4 0}

\arr{5.6 -2}{4.4 -2}
\arr{9.6 2}{8.4 2}

\arr{8 1.6}{8 0.4}
\arr{10 1.6}{10 0.4}

\arr{4 -.4}{4 -1.6}
\arr{6 -.4}{6 -1.6}

\put{$\ssize 4'$} at  8 2.5
\put{$\ssize 5'$} at  10 2.5

\put{$\ssize 2''$} at  4 -2.5
\put{$\ssize 3''$} at  6 -2.5

\put{$\ssize 0$} at 0 .5
\put{$\ssize 1$} at 2 .5
\put{$\ssize 2$} at 4 .5
\put{$\ssize 3$} at 6 .5
\put{$\ssize 4$} at 8 -.5
\put{$\ssize 5$} at 10 -.5
\put{$\ssize 6$} at 12 -.5
\put{$\ssize 7$} at 14 -.5

\setdots<2pt>
\plot 9.5 1.5  8.5 0.5 /
\plot 5.5 -0.5  4.5 -1.5 /
\plot 7.5 2  7.3 1.9  6.7 -1.9  6.5 -2 /
\endpicture}
$$
such that in addition any composition of six of the horizontal maps
in the middle is zero. In case we need a notation for the arrows, we will
use the following notation:
\item{$\bullet$} the horizontal arrows $i\!-\!1 \leftarrow i$ will be
denoted by $p_{i}$ or just $p$,
\item{$\bullet$} the horizontal arrow $4'\leftarrow 5'$ by $p'_5$,
\item{$\bullet$} the horizontal arrow $2'' \leftarrow 3''$ 
by $p''_3$.
\item{$\bullet$} the vertical arrows $i \leftarrow i'$ will be
denoted by $\mu_{i}$.
\item{$\bullet$} the vertical arrows $i'' \leftarrow i$ will be
denoted by $\varepsilon_{i}$.

\medskip
In order to see the relationship to $\Cal S(\widetilde6)$
one has to deal with the kernels 
$$
 V'_2 = \Ker \varepsilon_2 \quad\text{and}\quad 
 V'_3 = \Ker \varepsilon_3.
$$
of the two lower vertical maps $\varepsilon_i$; 
thus a representation $V$ of $\Omega$ with both $\mu_4$ and $\mu_5$ injective 
yields the corresponding object $E(V)$ in $\Cal S(\widetilde6).$
$$
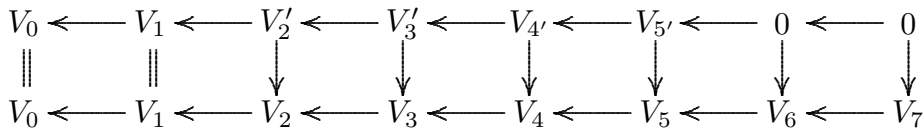
\hbox{\beginpicture
\setcoordinatesystem units <0.7cm,0.5cm>
\put{} at 0   0
\put{} at 14  2
\put{$V_0$} at 0 0
\put{$V_1$} at 2 0
\put{$V_2$} at 4 0
\put{$V_3$} at 6 0
\put{$V_4$} at 8 0
\put{$V_5$} at 10 0
\put{$V_6$} at 12 0
\put{$V_7$} at 14 0

\put{$V_0$} at 0 2
\put{$V_1$} at 2 2
\put{$V'_2$} at 4 2
\put{$V'_3$} at 6 2
\put{$V_{4'}$} at 8 2
\put{$V_{5'}$} at 10 2
\put{$0$} at 12 2
\put{$0$} at 14 2

\arr{1.6 0}{0.4 0}
\arr{3.6 0}{2.4 0}
\arr{5.6 0}{4.4 0}
\arr{7.6 0}{6.4 0}
\arr{9.6 0}{8.4 0}
\arr{11.6 0}{10.4 0}
\arr{13.6 0}{12.4 0}

\arr{1.6 2}{0.4 2}
\arr{3.6 2}{2.4 2}
\arr{5.6 2}{4.4 2}
\arr{7.6 2}{6.4 2}
\arr{9.6 2}{8.4 2}
\arr{11.6 2}{10.4 2}
\arr{13.6 2}{12.4 2}

\plot 0 1.4  0 .6 /
\plot .1 1.4 .1 .6 /
\plot 2 1.4  2 .6 /
\plot 2.1 1.4  2.1 .6 /
\arr{4 1.6}{4 0.4}
\arr{6 1.6}{6 0.4}
\arr{8 1.6}{8 0.4}
\arr{10 1.6}{10 0.4}
\arr{12 1.6}{12 0.4}
\arr{14 1.6}{14 0.4}

\endpicture}
$$

        \medskip\noindent{\bf (4.2) 
The structure of the category of $\Omega$-modules.} 
Consider the following three subquivers which yield algebras
$\Omega_0,\Omega_1,\Omega_2$ that are tame concealed of type~$\widetilde E_7$:
$$
\hbox{\beginpicture
\setcoordinatesystem units <0.5cm,0.5cm>
\put{\beginpicture
\put{} at 0  -2
\put{} at 14  2
\put{$\Omega_0\:$} at -2 -1
\put{$\circ$} at 0 0
\put{$\circ$} at 2 0
\put{$\circ$} at 4 0
\put{$\circ$} at 6 0
\put{$\circ$} at 8 0
\put{$\circ$} at 10 0

\put{$\circ$} at 4 -2
\put{$\circ$} at 6 -2

\arr{1.6 0}{0.4 0}
\arr{3.6 0}{2.4 0}
\arr{5.6 0}{4.4 0}
\arr{7.6 0}{6.4 0}
\arr{9.6 0}{8.4 0}

\arr{5.6 -2}{4.4 -2}

\arr{4 -.4}{4 -1.6}
\arr{6 -.4}{6 -1.6}

\put{$\ssize 2''$} at  4 -2.5
\put{$\ssize 3''$} at  6 -2.5

\put{$\ssize 0$} at 0 -.5
\put{$\ssize 1$} at 2 -.5
\put{$\ssize 2$} at 3.7 -.5
\put{$\ssize 3$} at 5.7 -.5
\put{$\ssize 4$} at 8 -.5
\put{$\ssize 5$} at 10 -.5

\setdots<2pt>
\plot 5.5 -0.5  4.5 -1.5 /

\endpicture} at 0 10
\put{\beginpicture
\put{} at 0  -2
\put{} at 14  2
\put{$\Omega_1\:$} at -2 0

\put{$\circ$} at 2 0
\put{$\circ$} at 4 0
\put{$\circ$} at 6 0
\put{$\circ$} at 8 0
\put{$\circ$} at 10 0
\put{$\circ$} at 12 0

\put{$\circ$} at 8 2

\put{$\circ$} at 6 -2

\arr{3.6 0}{2.4 0}
\arr{5.6 0}{4.4 0}
\arr{7.6 0}{6.4 0}
\arr{9.6 0}{8.4 0}
\arr{11.6 0}{10.4 0}

\arr{6 -.4}{6 -1.6}

\arr{8 1.6}{8 0.4}

\put{$\ssize 4'$} at  8 2.5

\put{$\ssize 3''$} at  6 -2.5

\put{$\ssize 1$} at 2 -.5
\put{$\ssize 2$} at 4 -.5
\put{$\ssize 3$} at 5.7 -.5
\put{$\ssize 4$} at 8 -.5
\put{$\ssize 5$} at 10 -.5
\put{$\ssize 6$} at 12 -.5

\setdots<2pt>
\plot 7.5  2  7.3 1.9  6.7 -1.9   
      6.5 -2 /
\endpicture} at 0 5
\put{\beginpicture
\put{} at 0  -2
\put{} at 14  2

\put{$\Omega_2\:$} at -2 1
\put{$\circ$} at 4 0
\put{$\circ$} at 6 0
\put{$\circ$} at 8 0
\put{$\circ$} at 10 0
\put{$\circ$} at 12 0
\put{$\circ$} at 14 0

\put{$\circ$} at 8 2
\put{$\circ$} at 10 2

\arr{5.6 0}{4.4 0}
\arr{7.6 0}{6.4 0}
\arr{9.6 0}{8.4 0}
\arr{11.6 0}{10.4 0}
\arr{13.6 0}{12.4 0}

\arr{9.6 2}{8.4 2}

\arr{8 1.6}{8 0.4}
\arr{10 1.6}{10 0.4}

\put{$\ssize 4'$} at  8 2.5
\put{$\ssize 5'$} at  10 2.5

\put{$\ssize 2$} at 4 -.5
\put{$\ssize 3$} at 6 -.5
\put{$\ssize 4$} at 8 -.5
\put{$\ssize 5$} at 10 -.5
\put{$\ssize 6$} at 12 -.5
\put{$\ssize 7$} at 14 -.5

\setdots<2pt>
\plot 9.5 1.5  8.5 0.5 /

\endpicture} at 0 0
\endpicture}
$$
        \medskip
The preprojective component of $\Omega_0$ is a component of $\Omega$:
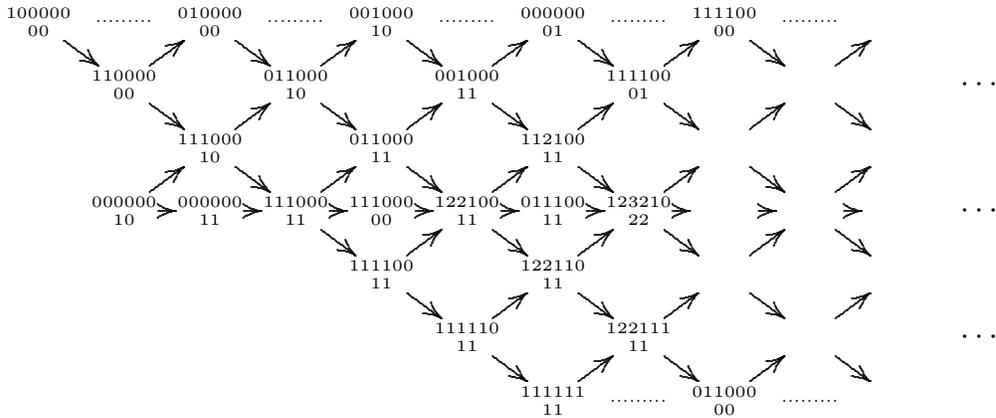
$$
\hbox{\beginpicture
\setcoordinatesystem units <0.95cm,0.7cm>
\put{} at 0 0
\put{} at 12 7
\put{$\ssize {100000 \atop 00}$}   at 0 6
\put{$\ssize {010000 \atop 00}$}   at 2 6
\put{$\ssize {001000 \atop 10}$}   at 4 6
\put{$\ssize {000000 \atop 01}$}   at 6 6
\put{$\ssize {111100 \atop 00}$}   at 8 6

\put{$\ssize {110000 \atop 00}$}   at 1 5
\put{$\ssize {011000 \atop 10}$}   at 3 5
\put{$\ssize {001000 \atop 11}$}   at 5 5
\put{$\ssize {111100 \atop 01}$}   at 7 5

\put{$\ssize {111000 \atop 10}$}   at 2 4
\put{$\ssize {011000 \atop 11}$}   at 4 4
\put{$\ssize {112100 \atop 11}$}   at 6 4

\put{$\ssize {000000 \atop 10}$}   at 1 3
\put{$\ssize {000000 \atop 11}$}   at 2 3
\put{$\ssize {111000 \atop 11}$}   at 3 3
\put{$\ssize {111000 \atop 00}$}   at 4 3
\put{$\ssize {122100 \atop 11}$}   at 5 3
\put{$\ssize {011100 \atop 11}$}   at 6 3
\put{$\ssize {123210 \atop 22}$}   at 7 3

\put{$\ssize {111100 \atop 11}$}   at 4 2
\put{$\ssize {122110 \atop 11}$}   at 6 2

\put{$\ssize {111110 \atop 11}$}   at 5 1
\put{$\ssize {122111 \atop 11}$}   at 7 1

\put{$\ssize {111111 \atop 11}$}   at 6 0 
\put{$\ssize {011000 \atop 00}$}   at 8 0

\arr{6.3 0.3} {6.7 0.7}
\arr{8.3 0.3} {8.7 0.7}

\arr{5.3 0.7} {5.7 0.3}
\arr{7.3 0.7} {7.7 0.3}
\arr{9.3 0.7} {9.7 0.3}

\arr{4.3 1.7} {4.7 1.3} 
\arr{6.3 1.7} {6.7 1.3} 
\arr{8.3 1.7} {8.7 1.3} 

\arr{5.3 1.3} {5.7 1.7} 
\arr{7.3 1.3} {7.7 1.7} 
\arr{9.3 1.3} {9.7 1.7} 

\arr{1.4 3.0} {1.6 3.0}
\arr{2.4 3.0} {2.6 3.0}
\arr{3.4 3.0} {3.6 3.0}
\arr{4.4 3.0} {4.6 3.0}
\arr{5.4 3.0} {5.6 3.0}
\arr{6.4 3.0} {6.6 3.0}
\arr{7.4 3.0} {7.6 3.0}
\arr{8.4 3.0} {8.6 3.0}
\arr{9.4 3.0} {9.6 3.0}

\arr{4.3 2.3} {4.7 2.7} 
\arr{6.3 2.3} {6.7 2.7} 
\arr{8.3 2.3} {8.7 2.7} 

\arr{3.3 2.7} {3.7 2.3} 
\arr{5.3 2.7} {5.7 2.3} 
\arr{7.3 2.7} {7.7 2.3} 
\arr{9.3 2.7} {9.7 2.3} 

\arr{2.3 3.7} {2.7 3.3} 
\arr{4.3 3.7} {4.7 3.3} 
\arr{6.3 3.7} {6.7 3.3} 
\arr{8.3 3.7} {8.7 3.3} 

\arr{1.3 3.3} {1.7 3.7} 
\arr{3.3 3.3} {3.7 3.7} 
\arr{5.3 3.3} {5.7 3.7} 
\arr{7.3 3.3} {7.7 3.7} 
\arr{9.3 3.3} {9.7 3.7} 

\arr{2.3 4.3} {2.7 4.7} 
\arr{4.3 4.3} {4.7 4.7} 
\arr{6.3 4.3} {6.7 4.7} 
\arr{8.3 4.3} {8.7 4.7} 

\arr{1.3 4.7} {1.7 4.3} 
\arr{3.3 4.7} {3.7 4.3} 
\arr{5.3 4.7} {5.7 4.3} 
\arr{7.3 4.7} {7.7 4.3} 
\arr{9.3 4.7} {9.7 4.3} 

\arr{0.3 5.7} {0.7 5.3} 
\arr{2.3 5.7} {2.7 5.3} 
\arr{4.3 5.7} {4.7 5.3} 
\arr{6.3 5.7} {6.7 5.3} 
\arr{8.3 5.7} {8.7 5.3} 

\arr{1.3 5.3} {1.7 5.7} 
\arr{3.3 5.3} {3.7 5.7} 
\arr{5.3 5.3} {5.7 5.7} 
\arr{7.3 5.3} {7.7 5.7} 
\arr{9.3 5.3} {9.7 5.7} 

\setdots<2pt>
\plot 6.7 0  7.3 0 /
\plot 8.7 0  9.3 0 /

\plot 0.7 6  1.3 6 /
\plot 2.7 6  3.3 6 /
\plot 4.7 6  5.3 6 /
\plot 6.7 6  7.3 6 /
\plot 8.7 6  9.3 6 /
\multiput{$\cdots$} at 11 1  11 3  11 5 /
\endpicture} 
$$
The radicals of the projective modules $P(6)$ and $P(4')$ belong to the
mouth of the 4-tube of $\Omega_0$, here is the corresponding component:
        \medskip

$$
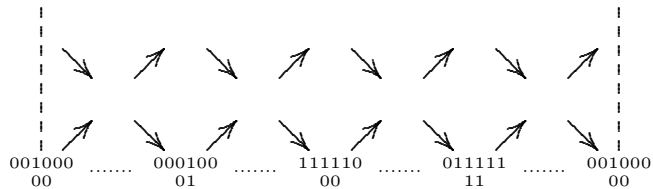
\hbox{\beginpicture
\setcoordinatesystem units <0.8cm,0.8cm>
\put{} at 0 2
\put{} at 0 4
\put{$\ssize {001000 \atop 00}$} at 0 2
\put{$\ssize {000100 \atop 01}$} at 2 2
\put{$\ssize {111110 \atop 00}$} at 4 2
\put{$\ssize {011111 \atop 11}$} at 6 2
\put{$\ssize {001000 \atop 00}$} at 8 2
\arr{0.3 2.3} {0.7 2.7} 
\arr{1.3 2.7} {1.7 2.3} 
\arr{2.3 2.3} {2.7 2.7} 
\arr{3.3 2.7} {3.7 2.3} 
\arr{4.3 2.3} {4.7 2.7} 
\arr{5.3 2.7} {5.7 2.3} 
\arr{6.3 2.3} {6.7 2.7} 
\arr{7.3 2.7} {7.7 2.3} 
\arr{0.3 3.7} {0.7 3.3} 
\arr{1.3 3.3} {1.7 3.7} 
\arr{2.3 3.7} {2.7 3.3} 
\arr{3.3 3.3} {3.7 3.7} 
\arr{4.3 3.7} {4.7 3.3} 
\arr{5.3 3.3} {5.7 3.7} 
\arr{6.3 3.7} {6.7 3.3} 
\arr{7.3 3.3} {7.7 3.7} 

\setdashes<3pt>
\plot 0 2.3  0 4.3 /
\plot 8 2.3  8 4.3 /

\setdots<2pt>
\plot 0.7 2  1.3 2 /
\plot 2.7 2  3.3 2 /
\plot 4.7 2  5.3 2 /
\plot 6.7 2  7.3 2 /
\endpicture}
$$
If we insert the rays starting at $P(4')$ and $P(6)$, we obtain the 
following tube:
$$
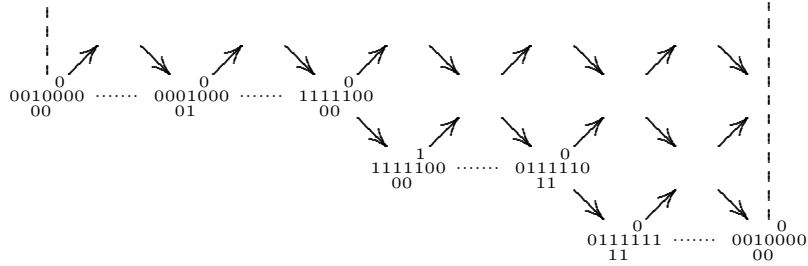
\hbox{\beginpicture
\setcoordinatesystem units <0.8cm,0.8cm>
\put{} at 0 0
\put{} at 0 3
\put{$\Left{0}{0010000}{00}$} at 0 2
\put{$\Left{0}{0001000}{01}$} at 2 2
\put{$\Left{0}{1111100}{00}$} at 4 2
\put{$\Left{1}{1111100}{00}$} at 5 1
\put{$\Left{0}{0111110}{11}$} at 7 1
\put{$\Left{0}{0111111}{11}$} at 8 0
\put{$\Left{0}{0010000}{00}$} at 10 0

\arr{7.3 0.7} {7.7 0.3} 
\arr{8.3 0.3} {8.7 0.7} 
\arr{9.3 0.7} {9.7 0.3} 
\arr{4.3 1.7} {4.7 1.3} 
\arr{5.3 1.3} {5.7 1.7} 
\arr{6.3 1.7} {6.7 1.3} 
\arr{7.3 1.3} {7.7 1.7} 
\arr{8.3 1.7} {8.7 1.3} 
\arr{9.3 1.3} {9.7 1.7} 
\arr{0.3 2.3} {0.7 2.7} 
\arr{1.3 2.7} {1.7 2.3} 
\arr{2.3 2.3} {2.7 2.7} 
\arr{3.3 2.7} {3.7 2.3} 
\arr{4.3 2.3} {4.7 2.7} 
\arr{5.3 2.7} {5.7 2.3} 
\arr{6.3 2.3} {6.7 2.7} 
\arr{7.3 2.7} {7.7 2.3} 
\arr{8.3 2.3} {8.7 2.7} 
\arr{9.3 2.7} {9.7 2.3} 

\setdashes<3pt>
\plot 0 2.3    0 3.3 /
\plot 10 0.3  10 3.3 /

\setdots<2pt>
\plot 0.7 2  1.3 2 /
\plot 2.7 2  3.3 2 /
\plot 5.7 1  6.3 1 /
\plot 8.7 0  9.3 0 /
\endpicture}
$$
        \medskip
It follows that the union of $\Omega_0$ and $\Omega_1$ is a tubular algebra of tubular
type $(6,3,2)$. By duality, the union of $\Omega_1$ and $\Omega_2$ is also a
tubular algebra of tubular type $(6,3,2)$. 
        \medskip
The 4-tube of $\Omega_1$ is
contained in the following 6-tube of $\Omega$, the coray-insertiones for 
adding the vertices $0$ and $2''$ and the ray insertions for adding the 
vertices $7$ and $5'$ yield two projective-injective modules:
        \medskip
$$
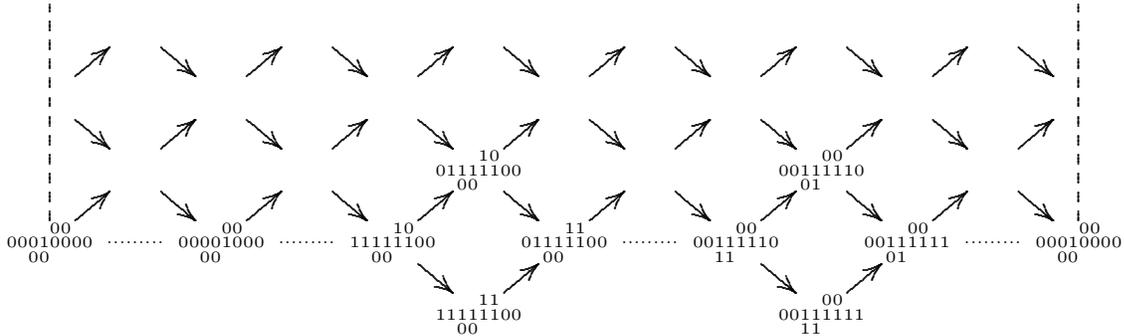
\hbox{\beginpicture
\setcoordinatesystem units <.95cm,0.8cm>
\put{} at 0 0
\put{} at 0 3
\put{$\all{00}{00010000}{00}$} at 0 0
\put{$\all{00}{00001000}{00}$} at 2 0
\put{$\all{10}{11111100}{00}$} at 4 0
\put{$\all{10}{01111100}{00}$} at 5 1
\put{$\all{11}{11111100}{00}$} at 5 -1
\put{$\all{11}{01111100}{00}$} at 6 0
\put{$\all{00}{00111110}{11}$} at 8 0
\put{$\all{00}{00111111}{11}$} at 9 -1
\put{$\all{00}{00111110}{01}$} at 9 1
\put{$\all{00}{00111111}{01}$} at 10 0
\put{$\all{00}{00010000}{00}$} at 12 0
\arr{4.3 -0.3} {4.7 -0.7} 
\arr{5.3 -0.7} {5.7 -0.3} 

\arr{8.3 -0.3} {8.7 -0.7} 
\arr{9.3 -0.7} {9.7 -0.3} 

\arr{0.3 0.3} {0.7 0.7} 
\arr{1.3 0.7} {1.7 0.3} 
\arr{2.3 0.3} {2.7 0.7} 
\arr{3.3 0.7} {3.7 0.3} 
\arr{4.3 0.3} {4.7 0.7} 
\arr{5.3 0.7} {5.7 0.3} 
\arr{6.3 0.3} {6.7 0.7} 
\arr{7.3 0.7} {7.7 0.3} 
\arr{8.3 0.3} {8.7 0.7} 
\arr{9.3 0.7} {9.7 0.3} 
\arr{10.3 0.3} {10.7 0.7} 
\arr{11.3 0.7} {11.7 0.3} 
\arr{0.3 1.7} {0.7 1.3} 
\arr{1.3 1.3} {1.7 1.7} 
\arr{2.3 1.7} {2.7 1.3} 
\arr{3.3 1.3} {3.7 1.7} 
\arr{4.3 1.7} {4.7 1.3} 
\arr{5.3 1.3} {5.7 1.7} 
\arr{6.3 1.7} {6.7 1.3} 
\arr{7.3 1.3} {7.7 1.7} 
\arr{8.3 1.7} {8.7 1.3} 
\arr{9.3 1.3} {9.7 1.7} 
\arr{10.3 1.7} {10.7 1.3}
\arr{11.3 1.3} {11.7 1.7} 
\arr{0.3 2.3} {0.7 2.7} 
\arr{1.3 2.7} {1.7 2.3} 
\arr{2.3 2.3} {2.7 2.7} 
\arr{3.3 2.7} {3.7 2.3} 
\arr{4.3 2.3} {4.7 2.7} 
\arr{5.3 2.7} {5.7 2.3} 
\arr{6.3 2.3} {6.7 2.7} 
\arr{7.3 2.7} {7.7 2.3} 
\arr{8.3 2.3} {8.7 2.7} 
\arr{9.3 2.7} {9.7 2.3} 
\arr{10.3 2.3} {10.7 2.7} 
\arr{11.3 2.7} {11.7 2.3} 
\setdashes<3pt>
\plot 0 0.3    0 3.3 /
\plot 12 0.3  12 3.3 /
\setdots<2pt>
\plot 0.7 0  1.3 0 /
\plot 2.7 0  3.3 0 /
\plot 6.7 0  7.3 0 /
\plot 10.7 0  11.3 0 /
\endpicture}
$$

        \medskip
In conclusion, the category $\mod \Omega$ has the following shape.  Note that the indecomposable
regular $\Omega_i$-modules are contained in the tubular families labelled $\Cal T_i$ for
$i=0,1,2$. 
$$
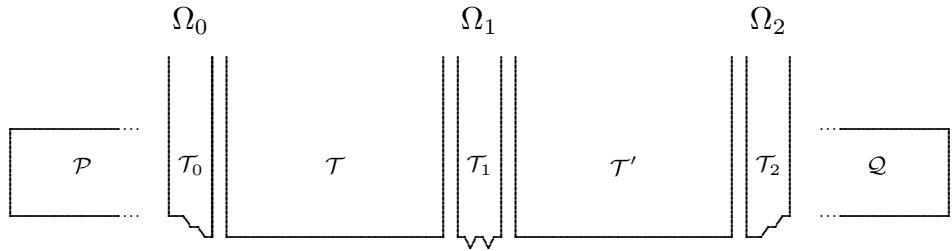
\hbox{\beginpicture
\setcoordinatesystem units <0.8cm,0.8cm>
\put{} at 0 0
\put{} at 9 2.5 
\plot 1.5 0.3  0 0.3  0 1.5  1.5 1.5 /
\setdots<2pt>
\plot 1.5 0.3  1.8 0.3 /
\plot 1.5 1.5  1.8 1.5 /
\setsolid
\plot 2.2 2.5  2.2 0.3  2.4 0.3  2.5 0.15  2.6 0.15  2.7 0  2.8 0
         2.8 2.5  2.8 0 /
\plot 3 2.5  3 0  6 0  6 2.5 /

\plot 6.2 2.5  6.2 0.0  6.3 0.0  6.375 -0.15  6.45 0  6.55 0
 6.625 -0.15  6.7 0  6.8 0  6.8 2.5 /

\plot 7 2.5  7 0  10 0  10 2.5 /

\plot 10.2 2.5  10.2 0  10.4 0  10.5 .15  10.6 .15  10.7 .3  10.8 .3
        10.8 2.5 /

\plot 11.5 0.3  13 0.3  13 1.5  11.5 1.5 /
\setdots<2pt>
\plot 11.5 0.3  11.2 0.3 /
\plot 11.5 1.5  11.2 1.5 /

\setsolid

\put{$\ssize \Cal P$} at 1 1
\put{$\ssize \Cal T_0$} at 2.5 1
\put{$\ssize \Cal T$} at 4.5 1
\put{$\ssize \Cal T_1$} at 6.5 1
\put{$\ssize \Cal T'$} at 8.5 1
\put{$\ssize \Cal T_2$} at 10.5 1
\put{$\ssize \Cal Q$} at 12 1

\put{$\Omega_0$} at 2.5 3
\put{$\Omega_1$} at 6.5 3
\put{$\Omega_2$} at 10.5 3
\endpicture}
$$
        \medskip\noindent{\bf (4.3) The representation embedding into
$\Cal S(\widetilde6)$.}
Let us call an indecomposable $\Omega$-module {\it regular} provided none of its
summands occurs in one of the components $\Cal P$ or $\Cal Q$.
Then the correspondence $V\mapsto E(V)$ from the introduction to this
section gives rise to a representation embedding from the regular
$\Omega$-modules into $\Cal S(\widetilde6)$. 
For this we show that every regular representation
has the property that the maps
$\varepsilon_2$ and $\varepsilon_3$ are epimorphisms and the maps
$\mu_4$ and $\mu_5$ are monomorphisms.  
Clearly, for a representation $M\in\mod \Omega$ the following conditions
are equivalent.
$$\hbox{\beginpicture
\setcoordinatesystem units <1cm,.5cm>
\put{$\varepsilon_2$ surjective} [l] at 1 6  \put{iff} [l] at 4 6
   \put{$\Hom(M,\all{00}{00000000}{11}) = 0$} [l] at 5 6
\put{$\varepsilon_3$ surjective} [l] at 1 5  
   \put{$\Hom(M,S(3'')) = 0$} [l] at 5 5
\put{$\mu_4$ injective} [l] at 1 4  
   \put{$\Hom(S(4'),M) = 0$} [l] at 5 4
\put{$\mu_5$ injective} [l] at 1 3  
   \put{$\Hom(\all{11}{00000000}{00},M) = 0$} [l] at 5 3
\endpicture}
$$
We have already seen that the modules $S(3'')$ and $\all{00}{00000000}{11}$ occur
in the component $\Cal P$, and dually, the modules $S(4')$ and $\all{11}{00000000}{00}$
are both in $\Cal Q$.  

\medskip
Thus, the regular $\Omega$-modules are embedded in $\Cal S(\widetilde6)$ 
as follows:
The tubular families $\Cal T_0$ and $\Cal T_2$ in $\Omega$-mod
are contained in the subcategories $\Cal U[-1]$ and $\Cal U[1]$ 
of $\Cal S(\widetilde6)$, and there are equivalences 
between each of the three full subcategories
$\Cal T$, $\Cal T_1$, and $\Cal T'$ of $\Omega$-mod and the 
corresponding subcategories $\Cal T$, $\Cal U$, and 
$\Cal T[1]$ of $\Cal S(\widetilde6)$:  

\smallskip
The regular $\Omega_0$-modules correspond to the regular $\Theta_0$-modules 
(see Section~1) by taking cokernels; and hence,
via the representation embedding $E$, they give rise to all
the indecomposable modules in the tubes in $\Cal U[-1]$ 
in $\Cal S(\widetilde6)$,
with the exception of several modules on two rays and two corays 
in the tube of rank~6.  
Thus $\Cal T_0$, being obtained from the tubular family in $\Omega_0$-mod
by a twofold ray insertion, is a full subcategory of $\Cal U[-1]$ in
$\Cal S(\widetilde6)$,
and dually, the tubular family $\Cal T_2$ of $\Omega$-mod becomes
a full subcategory of $\Cal U[1]$. 
On the categories $\Cal T$, $\Cal T_1$ and $\Cal T'$, the 
representation embedding $E$ induces an equivalence.
To verify density it
suffices to note that
the support of objects in $\Cal T\sqcup\Cal U\sqcup\Cal T[1]$
coincides with the support of $\Omega$ since the injective
$I(-1)=P(4')$ lies in $\Cal U[-1]$ in $\Cal S(\widetilde 6)$ 
and the projectives $P(6')$ and $P(8)$ are in $\Cal U[1]$; 
moreover,  each object $M$ in  $\Cal T\sqcup\Cal U\sqcup\Cal T[1]$
has the property that the maps $\mu_0$ and $\mu_1$ are isomorphisms,
this follows from the injective modules $I(0)'=P(5)$ and $I(1)'=P(6)$
being contained in $\Cal U[-2]$ and $\Cal U[-1]$ in 
$\Cal S(\widetilde6)$, respectively.
\medskip
        \vfill\eject

Combining the embedding of the regular $\Omega$-modules
in $\Cal S(\widetilde6)$ with the covering functor $\pi\:S(\widetilde6)\to
\Cal S(6)$ we obtain the double covering mentioned in the
introduction.
$$
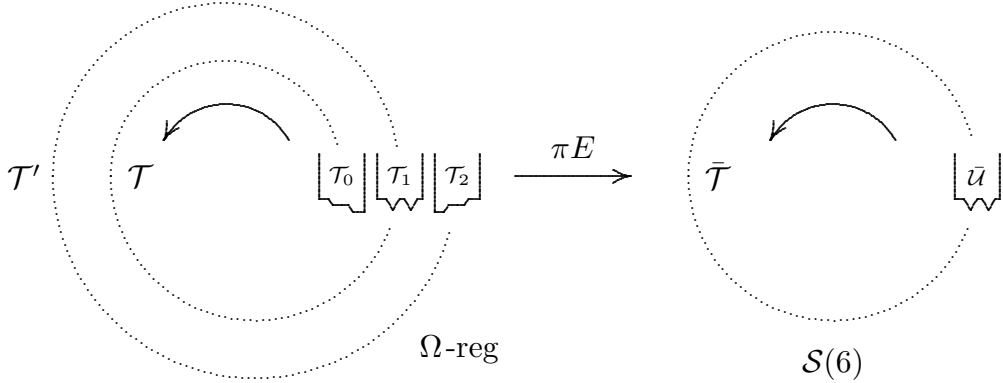
\beginpicture
\setcoordinatesystem units <1.6cm,1.6cm> 
\put{} at -.4 -1.4
\put{} at 6.2 1.2

\circulararc 120 degrees from 1.233 0.25  center at .8 0
\arr {.392 0.2933}  {.367 0.25}
\put{\qtube} at 1.6 0
\put{$\ssize \Cal T_0$} at 1.6 0
\put{\Utube} at 2 0
\put{$\ssize \Cal T_1$} at 2 0
\put{\ptube} at 2.4 0 
\put{$\ssize \Cal T_2$} at 2.4 0

\put{$\Cal T$} at .2 0
\put{$\Cal T'$} at -.6 0
\put{$\Omega$-reg} at 2.4 -1.2

\arr{2.8 0} {3.6 0}
\put{$\pi E$} at 3.2 .2

\setdots<2pt>
\circulararc -165 degrees from 0 0 center at .8 0
\circulararc 160 degrees from 0 0 center at 1 0
\circulararc -170 degrees from -.4 0 center at .8 0
\circulararc 165 degrees from -.4 0 center at 1 0

\setsolid
\put{$\ssize \bar{\Cal U}$} at 6 0
\put{\Utube} at 6 0
\circulararc 120 degrees from 5.433 0.25  center at 5 0
\arr {4.592 0.2933}  {4.567 0.25}
\put{$\bar{\Cal T}$} at 4.2 0
\put{$\Cal S(6)$} at 5 -1.3
\setdots<2pt>
\circulararc -165 degrees from  4 0 center at 5 0
\circulararc  160 degrees from  4 0 center at 5 0
\endpicture$$


        \bigskip\bigskip
\centerline{\bf Part B. Operators with nilpotency index
        different from 6.}

\bigskip

{\bf 5. General considerations.}

\medskip\noindent{\bf (5.1) The boundary modules.}
We will extract a lot of information about the shape of the submodule category
$\Cal S(n)$ from our study of the {\it boundary modules:}
They form the $\tau_{\Cal S}$-orbit to which
the two projective-injective objects are attached. 
In this section, let $\Lambda$ be a commutative uniserial ring.
We consider the two projective objects in $\Cal S(\Lambda)$, $P$ and
$P'$, their radicals $R=\rad P$ 
and $R'=\rad P'$,  the end terms of their source maps $J$ and $P'/K$,
and the simple and the quasi-simple representations $S$ and $K$.

\smallskip
\centerline{\vbox{\halign{\hfil#&$\;=\;$#\hfil&\qquad\hfil#&$\;=\;$#\hfil\cr
$P$ & $(\Lambda,0)$ & $P'$ & $(\Lambda,\Lambda)$\cr
$R$ & $(\rad\Lambda,0)$ & $R'$ 
        &$(\Lambda,\rad\Lambda)$\cr
$J$ & $(\Lambda,\soc\Lambda)$ & 
$P'/K$ & $(\Lambda/\soc\Lambda,\Lambda/\soc\Lambda)$\cr
$S$ & $(\soc\Lambda,0)$ & $K$ & $(\soc\Lambda,\soc\Lambda)$\cr }}}

\smallskip
Note that in case $n=2$, there are three pairs of isomorphic
representations in $\Cal S(\Lambda)$;
otherwise for $n > 2$, 
all the above representations are pairwise non-isomorphic.
Using arguments from [11] 
we show that these modules form a part of a component of
the Auslander-Reiten quiver for $\Cal S(\Lambda)$ of the following shape.

$$
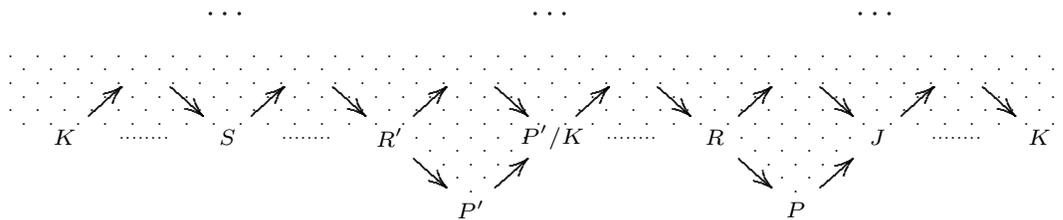
\hbox{\beginpicture
\setcoordinatesystem units <0.9cm,0.8cm>
\put{} at -0.7 -1
\put{} at 12.7 1.7
\put{$\ssize P'/K$} at 6 0
\put{$\ssize R$} at 8 0
\put{$\ssize P$} at 9 -1
\put{$\ssize J$} at 10 0
\put{$\ssize K$} at 12 0
\put{$\ssize K$} at 0 0
\put{$\ssize S$} at 2 0
\put{$\ssize R'$} at 4 0
\put{$\ssize P'$} at 5 -1
\arr{8.3 -0.3} {8.7 -0.7} 
\arr{9.3 -0.7} {9.7 -0.3} 
\arr{4.3 -0.3} {4.7 -0.7} 
\arr{5.3 -0.7} {5.7 -0.3} 
\arr{0.3 0.3} {0.7 0.7} 
\arr{1.3 0.7} {1.7 0.3} 
\arr{2.3 0.3} {2.7 0.7} 
\arr{3.3 0.7} {3.7 0.3} 
\arr{4.3 0.3} {4.7 0.7} 
\arr{5.3 0.7} {5.7 0.3} 
\arr{6.3 0.3} {6.7 0.7} 
\arr{7.3 0.7} {7.7 0.3} 
\arr{8.3 0.3} {8.7 0.7} 
\arr{9.3 0.7} {9.7 0.3} 
\arr{10.3 0.3} {10.7 0.7} 
\arr{11.3 0.7} {11.7 0.3} 
\setdots<2pt>
\plot 6.7 0  7.3 0 /
\plot 10.7 0  11.3 0 /
\plot 2.7 0  3.3 0 /
\plot 0.7 0  1.3 0 /
\put{$\cdots$} at 2 1.7
\put{$\cdots$} at 6 1.7
\put{$\cdots$} at 10 1.7 

\setshadegrid span <1.5mm>
\vshade  -0.7 0 1.3 <,z,,> 
          4 0 1.3  <z,z,,> 
          5 -1 1.3 <z,z,,>
          6 0 1.3 <z,z,,>
          8 0 1.3  <z,z,,>
          9 -1 1.3 <z,z,,>
         10 0 1.3 <z,,,>
         12.5 0 1.3 /

\endpicture}
$$
Clearly, the inclusions of the radicals $R'\to P'$ and $R\to P$ 
are sink maps in
$\Cal S(\Lambda)$; and according to Lemma (1.3.1), the maps
$\incl\: P\to J$ and $\can\: P'\to P'/K$ are source maps.
In order to verify that the non-projective modules form an orbit under
the Auslander-Reiten translation $\tau_{\Cal S}$ we recall from [11],
Theorem 5.1, the formula for the translate of an indecomposable
non-projective object $c\in \Cal S(\Lambda)$:
$$\tau_{\Cal S}(c)=\Mimo\tau_\Lambda\Cok (c)$$
We explain the three steps used for the computation of $\tau_{\Cal S}(c)$.
For this, we display objects in $\Cal S(\Lambda)$ 
as maps between $\Lambda$-modules,
so $c$ is given by an inclusion map $(c\:C'\to C)$.  
Let $(c''\:C\to C'')$ be the cokernel of $c$, then the Auslander-Reiten
translate in the category $\mod \Lambda$ gives rise to a morphism
$\tau_\Lambda(c'')$ in the category $\umod \Lambda$
which we represent by a map $(\tilde a\:A'\to \widetilde A)$ in the 
full subcategory $\mod'\Lambda$ of all $\Lambda$-modules which have no
non-zero injective direct summand; this map $\tilde a$ is determined uniquely,
up to a morphism which factors through an injective module.
From this we obtain $\tau_{\Cal S}(c)$ by making $\tilde a$
a ``minimal monomorphism'':
$$\Mimo(\tilde a)\;=\;[\tilde a\;e]\:A'\to \widetilde A\oplus \I(\Ker\tilde a)$$
where $e$ is an extension of the inclusion map 
$\Ker \tilde a\to \I(\Ker\tilde a)$ to $A'$. 
According to [11], Proposition 4.1, the isomorphism class of $\Mimo(\tilde a)$
in $\Cal S(n)$ does not depend on the choice for $\tilde a$
nor on the choice for $e$. 

\smallskip For example, consider the case that $c=K=(1_k\:k\to k)$
where $k=\soc\Lambda$. 
Then $\tilde a$  is just the cokernel map $(0\: k\to 0)$, and making this
map a monomorphism yields 
$\Mimo(k\to 0)=(\incl\:k\to \Lambda)$, so 
the translate is given by:  $\tau_{\Cal S}(K)=J$.
Similarly,
$$
\tau_{\Cal S}(J)=\Mimo(\tau_\Lambda(\Lambda\to \Lambda/k))
  =\Mimo(0\to \Lambda/k)
  = (0\to \Lambda/k)=R\qquad\text{etc.}$$

\medskip  In the situation where $\Lambda=k[x]/x^n$,
the boundary modules in $\Cal S(n)$ occur in the image
of the covering functor $\pi|_{\Cal S}\:\Cal S(\widetilde n)\to \Cal S(n)$.
Thus, for each boundary module $\bar B\in\Cal S(n)$ there is an orbit
under the shift of boundary modules $B[i]$, $i\in\Bbb Z$, in
the category $\Cal S(\widetilde n)$.  We index the boundary
modules such that 
$i=\max\{j\in\Bbb Z:B[i]_j\neq 0\}$ holds.
We obtain the following (part of an) orbit of boundary modules 
in $\Cal S(\widetilde n)$.
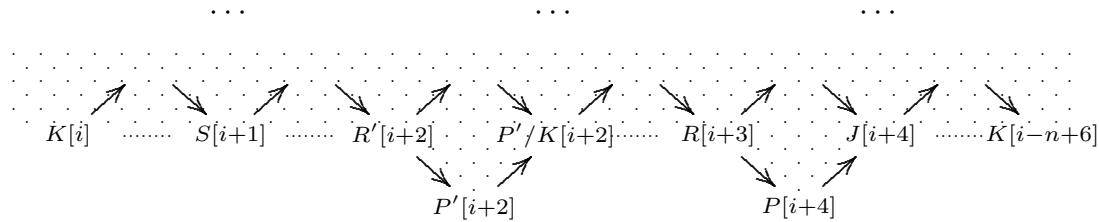
$$
\hbox{\beginpicture
\setcoordinatesystem units <0.9cm,0.8cm>
\put{} at 0 0
\put{} at 0 1.7
\put{$\ssize R[i+3]$} at 8 0
\put{$\ssize P[i+4]$} at 9 -1
\put{$\ssize J[i+4]$} at 10 0
\put{$\ssize K[i-n+6]$} at 12 0
\put{$\ssize K[i]$} at 0 0
\put{$\ssize S[i+1]$} at 2 0
\put{$\ssize R'[i+2]$} at 4 0
\put{$\ssize P'[i+2]$} at 5 -1
\put{$\ssize P'/K[i+2]$} at 6 0
\arr{8.3 -0.3} {8.7 -0.7} 
\arr{9.3 -0.7} {9.7 -0.3} 
\arr{4.3 -0.3} {4.7 -0.7} 
\arr{5.3 -0.7} {5.7 -0.3} 
\arr{0.3 0.3} {0.7 0.7} 
\arr{1.3 0.7} {1.7 0.3} 
\arr{2.3 0.3} {2.7 0.7} 
\arr{3.3 0.7} {3.7 0.3} 
\arr{4.3 0.3} {4.7 0.7} 
\arr{5.3 0.7} {5.7 0.3} 
\arr{6.3 0.3} {6.7 0.7} 
\arr{7.3 0.7} {7.7 0.3} 
\arr{8.3 0.3} {8.7 0.7} 
\arr{9.3 0.7} {9.7 0.3} 
\arr{10.3 0.3} {10.7 0.7} 
\arr{11.3 0.7} {11.7 0.3} 
\setdots<2pt>
\plot 6.7 0  7.3 0 /
\plot 10.7 0  11.3 0 /
\plot 2.7 0  3.3 0 /
\plot 0.7 0  1.3 0 /
\put{$\cdots$} at 2 1.7 
\put{$\cdots$} at 6 1.7
\put{$\cdots$} at 10 1.7

\setshadegrid span <1.5mm>
\vshade  -0.7 0 1.3 <,z,,> 
          4 0 1.3  <z,z,,> 
          5 -1 1.3 <z,z,,>
          6 0 1.3  <z,z,,>
          8 0 1.3  <z,z,,>
          9 -1 1.3 <z,z,,>
         10 0 1.3  <z,,,>
         12.5 0 1.3 /

\endpicture}
$$  
Thus, the boundary modules in $\Cal S(\widetilde n)$ exhibit
a shift $M\mapsto M[n-6]$ given by the sixfold
application of $\tau_{\Cal S}$; we will see how this shift characterizes
the shape of the Auslander-Reiten quiver of the category
$\Cal S(n)$, and how it contributes to the number of indecomposable
objects.

\medskip\noindent{\bf (5.2) The stable part.}
Frequently we will use [7], Theorem 2:
A connected component $\Cal C$ of the 
Auslander-Reiten quiver $\Gamma$ which contains 
a $\tau$-periodic module has stable part of the form
$\Bbb Z\Delta/\varphi$ where $\Delta$ is a Dynkin diagram or $A_\infty$ and $\varphi$
is an admissible automorphism of $\Bbb Z\Delta$.
This result applies in particular to each of the categories
$\Cal S(n)$ since here every module has $\tau_{\Cal S}$-period
a divisor of 6 ([11], Corollary~6.5).
        \bigskip
\bigskip
{\bf 6. Submodule Categories of Finite Representation Type}
        \medskip
For each $n\leq 5$, the category $\Cal S(n)$ is of finite
representation type as we will see.  
In each case, we present the Auslander-Reiten
quiver and display every one of the indecomposable modules.
It turns out that the number $s(n)$ of indecomposable
objects in $\Cal S(n)$, up to isomorphism, is given by the formula
$$ s(n)\;=\;2+2(n-1)\frac 6{6-n}$$
which we will discuss in (6.6).

\medskip
\medskip\noindent{\bf (6.1) The case $n = 1$.}
The category $\Cal S(1)$ is semisimple, the only indecomposable triples are $(k,0,0)$
and $(k,k,1_k)$.

\medskip\noindent{\bf (6.2) The case $n = 2$.}
In the category $\Cal S(2)$, among the eight boundary modules 
there are the two projective-injective representations and
in addition three pairs of isomorphic representations. 
If $X$ is a non-projective boundary module in $\Cal S(\widetilde 2)$
then the isomorphism $\tau_{\Cal S}^3X\cong X[-2]$ implies that
there are two $\tau_{\Cal S}$-orbits of stable boundary modules
in $\Cal S(\widetilde 2)$, one containing $X$, the other
one containing $X[1]$. 
These two orbits are connected as for example the 
Auslander-Reiten sequence ending at $K$ demonstrates,
and together with the projective-injective modules form
the following component of the Auslander-Reiten quiver
of $\Cal S(\widetilde 2)$:
$$
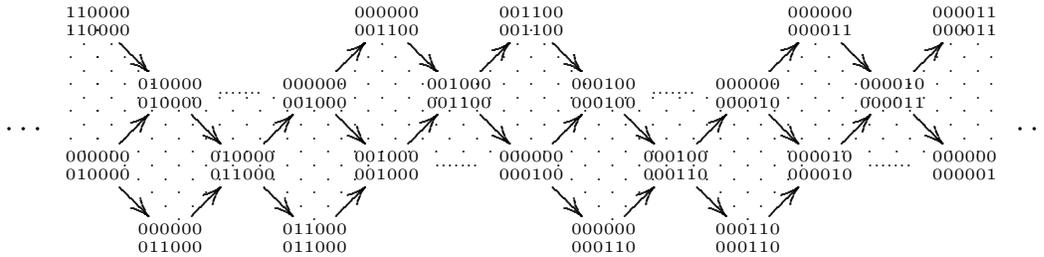
\hbox{\beginpicture
\setcoordinatesystem units <0.8cm,0.8cm>
\put{} at -1.2 -1.2
\put{} at 13.2 2.2
\put{$\cdots$} at -1 .5
\put{$\cdots$} at 13 .5
\put{$\ssize {000000 \atop 010000}$} at 0 0
\put{$\ssize {000000 \atop 011000}$} at 1 -1
\put{$\ssize {010000 \atop 011000}$} at 2 0
\put{$\ssize {011000 \atop 011000}$} at 3 -1
\put{$\ssize {001000 \atop 001000}$} at 4 0
\put{$\ssize {000000 \atop 000100}$} at 6 0
\put{$\ssize {000000 \atop 000110}$} at 7 -1
\put{$\ssize {000100 \atop 000110}$} at 8 0
\put{$\ssize {000110 \atop 000110}$} at 9 -1
\put{$\ssize {000010 \atop 000010}$} at 10 0
\put{$\ssize {000000 \atop 000001}$} at 12 0

\put{$\ssize {110000 \atop 110000}$} at 0 2
\put{$\ssize {010000 \atop 010000}$} at 1 1
\put{$\ssize {000000 \atop 001000}$} at 3 1
\put{$\ssize {000000 \atop 001100}$} at 4 2
\put{$\ssize {001000 \atop 001100}$} at 5 1
\put{$\ssize {001100 \atop 001100}$} at 6 2
\put{$\ssize {000100 \atop 000100}$} at 7 1
\put{$\ssize {000000 \atop 000010}$} at 9 1
\put{$\ssize {000000 \atop 000011}$} at 10 2
\put{$\ssize {000010 \atop 000011}$} at 11 1
\put{$\ssize {000011 \atop 000011}$} at 12 2
\arr{0.3 -0.3} {0.7 -0.7} 
\arr{1.3 -0.7} {1.7 -0.3} 
\arr{2.3 -0.3} {2.7 -0.7} 
\arr{3.3 -0.7} {3.7 -0.3} 
\arr{6.3 -0.3} {6.7 -0.7} 
\arr{7.3 -0.7} {7.7 -0.3} 
\arr{8.3 -0.3} {8.7 -0.7} 
\arr{9.3 -0.7} {9.7 -0.3} 
\arr{0.3 0.3} {0.7 0.7} 
\arr{1.3 0.7} {1.7 0.3} 
\arr{2.3 0.3} {2.7 0.7} 
\arr{3.3 0.7} {3.7 0.3} 
\arr{4.3 0.3} {4.7 0.7} 
\arr{5.3 0.7} {5.7 0.3} 
\arr{6.3 0.3} {6.7 0.7} 
\arr{7.3 0.7} {7.7 0.3} 
\arr{8.3 0.3} {8.7 0.7} 
\arr{9.3 0.7} {9.7 0.3} 
\arr{10.3 0.3} {10.7 0.7} 
\arr{11.3 0.7} {11.7 0.3} 
\arr{0.3 1.7} {0.7 1.3}
\arr{3.3 1.3} {3.7 1.7}
\arr{4.3 1.7} {4.7 1.3}
\arr{5.3 1.3} {5.7 1.7}
\arr{6.3 1.7} {6.7 1.3}
\arr{9.3 1.3} {9.7 1.7}
\arr{10.3 1.7} {10.7 1.3}
\arr{11.3 1.3} {11.7 1.7}
\setdots<2pt>
\plot 1.7 1  2.3 1 /
\plot 4.7 0  5.3 0 /
\plot 7.7 1  8.3 1 /
\plot 10.7 0  11.3 0 /

\setshadegrid span <1.5mm>
\vshade  -0.5 0 2 <,z,,> 
          0   0 2  <z,z,,> 
          1  -1 1  <z,z,,>
          2   0 1  <z,z,,> 
          3  -1 1  <z,z,,>
          4   0 2  <z,z,,> 
          5   0 1  <z,z,,>
          6   0 2  <z,z,,> 
          7  -1 1  <z,z,,>
          8   0 1  <z,z,,> 
          9  -1 1  <z,z,,>
         10   0 2  <z,z,,> 
         11   0 1  <z,z,,>
         12   0 2  <z,z,,>
         12.5 0 2 /

\endpicture}
$$
Thus, all indecomposables in $\Cal S(\widetilde2)$ are boundary modules.
The Auslander-Reiten quiver for $\Cal S(2)$ is as follows.
Here and for the  following representation finite categories we use
two presentations for the objects in the Auslander-Reiten quiver:
First, we represent each module by the dimension vector of a corresponding
object in $\Cal S(\widetilde n)$ and second, 
we display the embedding of the submodule by using 
box notation as introduced in (2.3).
$$
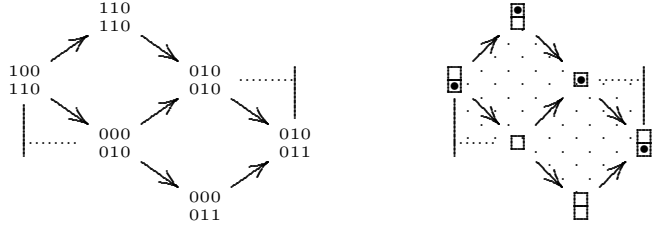
\hbox{\beginpicture
\setcoordinatesystem units <1cm,0.7cm>
\put{} at 0 0
\put{} at 3 3
\put{$\ssize {100\atop 110}$}   at 0 2
\put{$\ssize {000\atop 010}$}   at 1 1
\put{$\ssize {110\atop 110}$}   at 1 3
\put{$\ssize {000\atop 011}$} at 2 0
\put{$\ssize {010\atop 010}$}  at 2 2
\put{$\ssize {010\atop 011}$}    at 3 1

\arr{0.3 1.7} {0.7 1.3} 
\arr{1.3 1.3} {1.7 1.7} 
\arr{2.3 1.7} {2.7 1.3} 
\arr{1.3 0.7} {1.7 0.3} 
\arr{2.3 0.3} {2.7 0.7} 
\arr{0.3 2.3} {0.7 2.7} 
\arr{1.3 2.7} {1.7 2.3} 
\setdots<2pt>
\plot 0 1   0.6 1 /
\plot 2.4 2  3 2 /
\setsolid
\plot 0 .8  0 1.6 /
\plot 3 1.4  3 2.2 /
\endpicture} 
\qquad\qquad
\hbox{\beginpicture
\setcoordinatesystem units <.7cm,0.7cm>
\put{} at 0 0
\put{} at 3 3
\multiput{$\smallsq2$} at -.1 1.9  -.1 2.1 /
\put{$\sssize \bullet$} at  0 1.9
\put{$\smallsq2$} at .9 1
\multiput{$\smallsq2$} at .9 2.9  .9 3.1 /
\put{$\sssize \bullet$} at 1 3.1
\multiput{$\smallsq2$} at 1.9 -.1  1.9 .1 /
\put{$\smallsq2$} at 1.9 2
\put{$\sssize \bullet$} at 2 2
\multiput{$\smallsq2$} at 2.9 0.9  2.9 1.1 /
\put{$\sssize \bullet$} at 3 .9

\arr{0.3 1.7} {0.7 1.3} 
\arr{1.3 1.3} {1.7 1.7} 
\arr{2.3 1.7} {2.7 1.3} 
\arr{1.3 0.7} {1.7 0.3} 
\arr{2.3 0.3} {2.7 0.7} 
\arr{0.3 2.3} {0.7 2.7} 
\arr{1.3 2.7} {1.7 2.3} 
\setdots<2pt>
\plot 0 1   0.7 1 /
\plot 2.3 2  3 2 /
\setsolid
\plot 0 .8  0 1.7 /
\plot 3 1.3  3 2.2 /
\setshadegrid span <1.5mm>
\vshade   0   1 2 <,z,,> 
          1   1 3  <z,z,,> 
          2   0 2  <z,z,,>
          3   1 2 /
\endpicture} 
$$

For each of the categories $\Cal S(n)$ where $2\leq n\leq 5$ we picture
the distribution of the possible dimension pairs $(v,u)=\bdim(V,U,T)$ 
for the indecomposable triples $(V,U,T)$.
We label each point $(v,u)$ with the Krull-Remak-Schmidt multiplicity 
(which is the number of indecomposable direct summands) of the 
$k[x]$-module given by $(V,T)$. 
Note that in each case the dimension pairs lie in the stripe bounded by the
lines $u=\frac12 v\pm \frac12n$. Here is the distribution for $n=2$:
$$
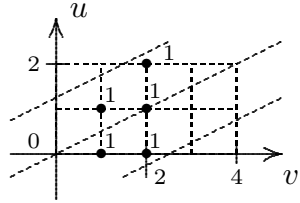
\hbox{\beginpicture
\setcoordinatesystem units <.5cm,.5cm>
\arr{0 -1}{0 3}
\arr{-1 0}{5 0}
\put{$u$} at  0.5 3.2
\put{$v$} at 5.2 -.5
\setdashes <.5mm> 
\setplotarea x from 0 to 4, y from 0 to 2 
\grid {4} {2}
\plot -1 -.5  5 2.5 /
\plot -1 .75  2.5 2.5 /
\plot 1.5 -.5  5 1.25 /
\setsolid
\plot 2 -0.5  2 0.2 /
\put{$\ssize 2$} at 2.3 -.5
\plot 4 -0.2  4 0.2 /
\put{$\ssize 4$} at 4 -.5
\plot -0.2 2  0.2 2 /
\put{$\ssize 2$} at -.5 2 
\put{$\ssize 0$} at -.5 .3 
\multiput{$\ssize \bullet$} at 1 0  2 0   1 1  2 1  2 2 /
\multiput{$\ssize 1$} at 1.2 .3   2.2 .3  1.2 1.3  2.2 1.4  2.5 2.25 /
\endpicture}
$$

\medskip\noindent{\bf (6.3) The case $n = 3$.} 
We derive the shape of the Auslander-Reiten quiver for $\Cal S(3)$
from our knowledge about the boundary modules.
Since the boundary modules are pairwise non-isomorphic,
the formula $\tau_{\Cal S}^6X\cong X[-3]$ 
for a non-projective boundary module $X$ in the covering category
$\Cal S(\widetilde3)$ implies that 
there are $3$ orbits of boundary modules in the
stable part $\Bbb Z\Delta$ of $\Cal S(\widetilde3)$,
more precisely, the orbits containing $X$, $X[1]$ and $X[2]$, 
are pairwise disjoint.
The automorphism of $\Bbb Z\Delta$ given by $M\mapsto M[1]$
permutes these three orbits and 
hence gives rise to an automorphism $\rho$ of the Dynkin diagram
$\Delta$ of order three.  The only Dynkin diagram $\Delta$ admitting such an automorphism
is $\Delta=D_4$.  Thus, the shape of the stable part of $\Gamma(3)$
is determined as $\Bbb ZD_4/\tau_{\Cal S}^2\rho$.  
Here is the Auslander-Reiten quiver:
$$
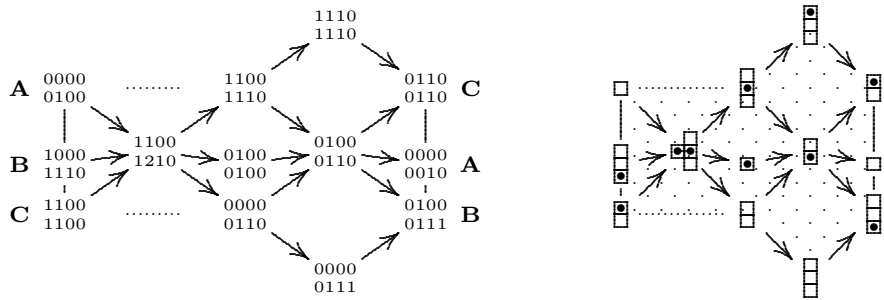
\hbox{\beginpicture
\setcoordinatesystem units <1cm,0.7cm> 
\put{} at 0 0
\put{} at 4 4
\multiput{$\ssize \bold A$} at -.5 3  4.5 1.8 /
\multiput{$\ssize \bold B$} at -.5 1.8  4.5 1 /
\multiput{$\ssize \bold C$} at -.5 1  4.5 3 /

\put{$\ssize {1100\atop 1100}$}   at 0 1
\put{$\ssize {1000\atop 1110}$} at 0 1.8
\put{$\ssize {0000\atop 0100}$}   at 0 3

\put{$\ssize {1100\atop 1210}$}   at 1 2

\put{$\ssize {0000 \atop 0110}$} at 2 1
\put{$\ssize {0100\atop 0100}$}    at 2 1.8
\put{$\ssize {1100\atop 1110}$}  at 2 3

\put{$\ssize {0000\atop 0111}$} at 3 0
\put{$\ssize {0100\atop 0110}$}   at 3 2
\put{$\ssize {1110\atop 1110}$} at 3 4

\put{$\ssize {0100\atop 0111}$} at 4 1
\put{$\ssize {0000\atop 0010}$} at 4 1.8
\put{$\ssize {0110\atop 0110}$} at 4 3

\arr{2.3 0.7} {2.7 0.3} 
\arr{3.3 0.3} {3.7 0.7} 
\arr{0.3 1.3} {0.7 1.7} 
\arr{1.3 1.7} {1.7 1.3} 
\arr{2.3 1.3} {2.7 1.7} 
\arr{3.3 1.7} {3.7 1.3} 
\arr{0.3 2.7} {0.7 2.3} 
\arr{1.3 2.3} {1.7 2.7} 
\arr{2.3 2.7} {2.7 2.3} 
\arr{3.3 2.3} {3.7 2.7} 
\arr{2.3 3.3} {2.7 3.7} 
\arr{3.3 3.7} {3.7 3.3} 
\arr{0.3 1.86}{0.7 1.94}
\arr{1.3 1.94}{1.7 1.86}
\arr{2.3 1.86}{2.7 1.94}
\arr{3.3 1.94}{3.7 1.86}
\setdots<2pt>
\plot 0.7 1  1.3 1 /
\plot 0.7 3  1.3 3 /
\setsolid
\plot 0 1.35  0 1.45 /
\plot 0 2.15  0 2.65 /
\plot 4 1.35  4 1.45 /
\plot 4 2.15  4 2.65 /
\endpicture} 
\qquad\qquad
\hbox{\beginpicture
\setcoordinatesystem units <.7cm,0.7cm> 
\put{} at 0 0
\put{} at 4 4
\multiput{$\smallsq2$} at -.1 .9  -.1 1.1 /
\put{$\sssize \bullet$} at  0 1.1 
\multiput{$\smallsq2$} at -.1 1.6  -.1 1.8  -.1 2 /
\put{$\sssize \bullet$} at  0 1.6
\put{$\smallsq2$} at -.1 3
\multiput{$\smallsq2$} at .8 2  1 1.8  1 2  1 2.2 /
\multiput{$\sssize \bullet$} at .9 2  1.1 2 /
\plot .9 2  1.1 2 /
\multiput{$\smallsq2$} at 1.9 0.9  1.9 1.1 /
\put{$\smallsq2$} at 1.9 1.8
\put{$\sssize \bullet$} at 2 1.8
\multiput{$\smallsq2$} at 1.9 2.8  1.9 3  1.9 3.2 /
\put{$\sssize\bullet$} at 2 3
\multiput{$\smallsq2$} at 2.9 -.2  2.9 0  2.9 .2 /
\multiput{$\smallsq2$} at 2.9 1.9  2.9 2.1 /
\put{$\sssize\bullet$} at 3 1.9
\multiput{$\smallsq2$} at 2.9 3.8  2.9 4  2.9 4.2 /
\put{$\sssize\bullet$} at 3 4.2
\multiput{$\smallsq2$} at 3.9 0.8  3.9 1  3.9 1.2 /
\put{$\sssize \bullet$} at 4 0.8
\put{$\smallsq2$} at 3.9 1.8
\multiput{$\smallsq2$} at 3.9 2.9  3.9 3.1 /
\put{$\sssize \bullet$} at 4 3.1
\arr{2.3 0.7} {2.7 0.3} 
\arr{3.3 0.3} {3.7 0.7} 
\arr{0.3 1.3} {0.7 1.7} 
\arr{1.3 1.7} {1.7 1.3} 
\arr{2.3 1.3} {2.7 1.7} 
\arr{3.3 1.7} {3.7 1.3} 
\arr{0.3 2.7} {0.7 2.3} 
\arr{1.3 2.3} {1.7 2.7} 
\arr{2.3 2.7} {2.7 2.3} 
\arr{3.3 2.3} {3.7 2.7} 
\arr{2.3 3.3} {2.7 3.7} 
\arr{3.3 3.7} {3.7 3.3} 
\arr{0.3 1.86}{0.7 1.94}
\arr{1.3 1.94}{1.7 1.86}
\arr{2.3 1.86}{2.7 1.94}
\arr{3.3 1.94}{3.7 1.86}
\setdots<2pt>
\plot 0.3 1  1.7 1 /
\plot 0.3 3  1.7 3 /
\setsolid
\plot 0 1.3  0 1.4 /
\plot 0 2.2  0 2.8 /
\plot 4 1.4  4 1.6 /
\plot 4 2  4 2.7 /
\setshadegrid span <1.5mm>
\vshade   0   1 3 <,z,,> 
          2   1 3  <z,z,,> 
          3   0 4  <z,z,,>
          4   1 3 /

\endpicture} 
$$

We conclude this section with the distribution of the dimension pairs
of the indecomposable triples in the category $\Cal S(3)$.
$$
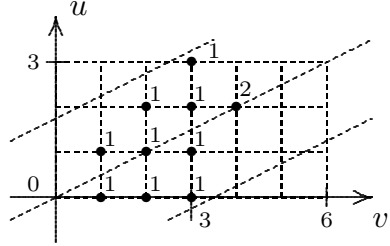
\hbox{\beginpicture
\setcoordinatesystem units <.5cm,.5cm>
\arr{0 -1}{0 4}
\arr{-1 0}{7 0}
\put{$u$} at  0.5 4.2
\put{$v$} at 7.2 -.5
\setdashes <.5mm> 
\setplotarea x from 0 to 6, y from 0 to 3 
\grid {6} {3}
\plot -1 -.5  7 3.5 /
\plot -1 1.25   3.5 3.5 /
\plot 2.5 -.5  7 1.75 /
\setsolid
\plot 3 -0.5  3 0.2 /
\put{$\ssize 3$} at 3.3 -.5
\plot 6 -0.2  6 0.2 /
\put{$\ssize 6$} at 6 -.5
\plot -0.2 3  0.2 3 /
\put{$\ssize 3$} at -.5 3 
\put{$\ssize 0$} at -.5 .3
\multiput{$\ssize \bullet$} at 1 0  2 0  3 0  1 1  2 1  3 1  2 2  3 2  3 3  4 2 /
\multiput{$\ssize 1$} at 1.2 .3   2.2 .3  3.2 .3  1.2 1.3  2.2 1.4  3.2 1.3  
                                      2.2 2.3  3.2 2.3  3.5 3.25  /
\multiput{$\ssize 2$} at 4.2 2.4 /
\endpicture}
$$

\medskip\noindent{\bf (6.4) The case $n = 4$.}
Here again the boundary modules are pairwise non-isomorphic,
and the formula $\tau_{\Cal S}^6X\cong X[-2]$ for a non-projective
module in the covering category $\Cal S(\widetilde 4)$ 
yields two orbits of boundary modules in the Auslander-Reiten
quiver for $\Cal S(\widetilde4)$.
The automorphism on the stable part $\Bbb Z\Delta$ of this Auslander-Reiten quiver
given by the shift $M\mapsto M[1]$ permutes these two orbits and 
induces a non-trivial action $\sigma$ of order two on the Dynkin diagram $\Delta$.
Hence, $\Delta$ is of type $A_n$, $D_n$, or $E_6$.
First we compute the dimension vectors of the indecomposables in
the Auslander-Reiten sequences containing a boundary module.

$$
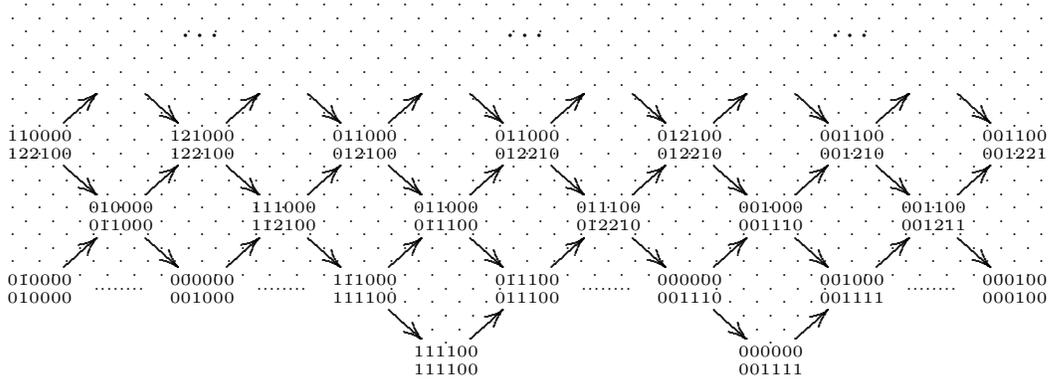
\hbox{\beginpicture
\setcoordinatesystem units <0.9cm,0.8cm>
\put{} at 0 0
\put{} at 0 1.7
\put{$\ssize {010000 \atop 010000}$} at 0 0
\put{$\ssize {000000 \atop 001000}$} at 2 0
\put{$\ssize {111000 \atop 111100}$} at 4 0
\put{$\ssize {111100 \atop 111100}$} at 5 -1
\put{$\ssize {011100 \atop 011100}$} at 6 0
\put{$\ssize {000000 \atop 001110}$} at 8 0
\put{$\ssize {000000 \atop 001111}$} at 9 -1
\put{$\ssize {001000 \atop 001111}$} at 10 0
\put{$\ssize {000100 \atop 000100}$} at 12 0

\put{$\ssize {110000 \atop 122100}$} at 0 2
\put{$\ssize {010000 \atop 011000}$} at 1 1
\put{$\ssize {121000 \atop 122100}$} at 2 2
\put{$\ssize {111000 \atop 112100}$} at 3 1
\put{$\ssize {011000 \atop 012100}$} at 4 2
\put{$\ssize {011000 \atop 011100}$} at 5 1
\put{$\ssize {011000 \atop 012210}$} at 6 2
\put{$\ssize {011100 \atop 012210}$} at 7 1
\put{$\ssize {012100 \atop 012210}$} at 8 2
\put{$\ssize {001000 \atop 001110}$} at 9 1
\put{$\ssize {001100 \atop 001210}$} at 10 2
\put{$\ssize {001100 \atop 001211}$} at 11 1
\put{$\ssize {001100 \atop 001221}$} at 12 2
\arr{4.3 -0.3} {4.7 -0.7} 
\arr{5.3 -0.7} {5.7 -0.3} 
\arr{8.3 -0.3} {8.7 -0.7} 
\arr{9.3 -0.7} {9.7 -0.3} 
\arr{0.3 0.3} {0.7 0.7} 
\arr{1.3 0.7} {1.7 0.3} 
\arr{2.3 0.3} {2.7 0.7} 
\arr{3.3 0.7} {3.7 0.3} 
\arr{4.3 0.3} {4.7 0.7} 
\arr{5.3 0.7} {5.7 0.3} 
\arr{6.3 0.3} {6.7 0.7} 
\arr{7.3 0.7} {7.7 0.3} 
\arr{8.3 0.3} {8.7 0.7} 
\arr{9.3 0.7} {9.7 0.3} 
\arr{10.3 0.3} {10.7 0.7} 
\arr{11.3 0.7} {11.7 0.3} 
\arr{0.3 1.7} {0.7 1.3}
\arr{1.3 1.3} {1.7 1.7}
\arr{2.3 1.7} {2.7 1.3}
\arr{3.3 1.3} {3.7 1.7}
\arr{4.3 1.7} {4.7 1.3}
\arr{5.3 1.3} {5.7 1.7}
\arr{6.3 1.7} {6.7 1.3}
\arr{7.3 1.3} {7.7 1.7}
\arr{8.3 1.7} {8.7 1.3}
\arr{9.3 1.3} {9.7 1.7}
\arr{10.3 1.7} {10.7 1.3}
\arr{11.3 1.3} {11.7 1.7}
\arr{0.3 2.3} {0.7 2.7} 
\arr{1.3 2.7} {1.7 2.3} 
\arr{2.3 2.3} {2.7 2.7} 
\arr{3.3 2.7} {3.7 2.3} 
\arr{4.3 2.3} {4.7 2.7} 
\arr{5.3 2.7} {5.7 2.3} 
\arr{6.3 2.3} {6.7 2.7} 
\arr{7.3 2.7} {7.7 2.3} 
\arr{8.3 2.3} {8.7 2.7} 
\arr{9.3 2.7} {9.7 2.3} 
\arr{10.3 2.3} {10.7 2.7} 
\arr{11.3 2.7} {11.7 2.3} 
\setdots<2pt>
\plot 0.7 0  1.3 0 /
\plot 2.7 0  3.3 0 /
\plot 6.7 0  7.3 0 /
\plot 10.7 0  11.3 0 /

\setshadegrid span <1.5mm>
\vshade  -0.5 0 4 <,z,,> 
          4   0 4  <z,z,,> 
          5  -1 4  <z,z,,>
          6   0 4  <z,z,,> 
          8   0 4  <z,z,,> 
          9  -1 4  <z,z,,>
         10   0 4  <z,z,,> 
         12.5 0 4 /
\put{$\cdots$} at 2 3.5
\put{$\cdots$} at 6 3.5
\put{$\cdots$} at 10 3.5
\endpicture}
$$
The dimension vectors in the $\tau_{\Cal S}$-orbit next to the boundary modules
correspond to indecomposables (since at least the first one does); they
have period six, up to the shift. This excludes the case that $\Delta$ 
has type $D_n$ for some $n$, as the shift $M\mapsto M[1]$ maps the central
axis of the diagram $\Bbb ZD_n$ onto itself.  
Hence, the modules in the second orbit next to the 
boundary modules must be indecomposable.
They occur with period three, up to the shift, so either $\Delta=A_5$ or $\Delta=E_6$ holds.
For dimension reasons, 
the case $\Delta=A_5$ is not possible:
The Auslander-Reiten sequence starting at $\ssize {1100 \atop 1221}$
would have to have a middle term which is the direct sum of the extension of 
$\ssize {01 \atop 01}$ and $\ssize {000 \atop 001}$
and the extension of $\ssize {111\atop 111}$ and $\ssize {0000\atop 0111}$;
but clearly,  there is no exact sequence 
$ 0\to {\ssize{1100\atop 1221}} 
   \to {\ssize{010\atop 011}}\oplus{\ssize{1110\atop 1221}}
   \to {\ssize{1210\atop 1221}} \to 0$, a contradiction.
Hence, $\Delta=E_6$ and 
the Auslander-Reiten quiver for $\Cal S(4)$ has 
stable part $\Bbb ZE_6/\tau_{\Cal S}^3\sigma$:
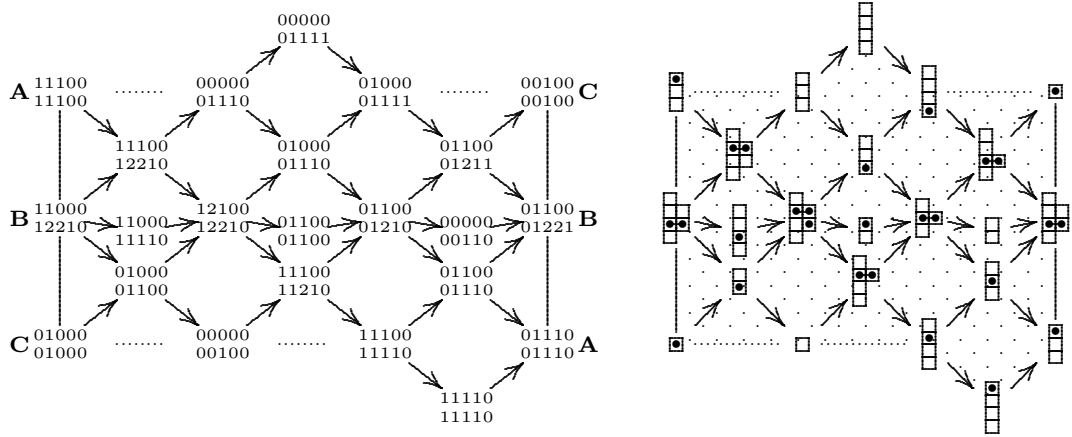
$$
\hbox{\beginpicture
\setcoordinatesystem units <.9cm,0.7cm>
\put{} at 0 0
\put{} at 7 6
\multiput{$\ssize\bold A$} at -.5 5  6.5 1 /
\multiput{$\ssize\bold B$} at -.5 3  6.5 3 /
\multiput{$\ssize\bold C$} at -.5 1  6.5 5 /
\put{$\ssize {11100\atop 11100}$}   at 0 5
\put{$\ssize {11000\atop 12210}$}   at 0 3
\put{$\ssize {01000\atop 01000}$}   at 0 1

\put{$\ssize {11100\atop 12210}$} at 1 4
\put{$\ssize {11000\atop 11110}$} at 1 2.8
\put{$\ssize {01000\atop 01100}$} at 1 2

\put{$\ssize {00000\atop 01110}$} at 2 5
\put{$\ssize {12100\atop 12210}$} at 2 3
\put{$\ssize {00000\atop 00100}$} at 2 1

\put{$\ssize {00000\atop 01111}$} at 3 6
\put{$\ssize {01000\atop 01110}$} at 3 4
\put{$\ssize {01100\atop 01100}$} at 3 2.8
\put{$\ssize {11100\atop 11210}$} at 3 2

\put{$\ssize {01000\atop 01111}$} at 4 5
\put{$\ssize {01100\atop 01210}$} at 4 3
\put{$\ssize {11100\atop 11110}$} at 4 1

\put{$\ssize {01100\atop 01211}$} at 5 4
\put{$\ssize {00000 \atop 00110}$} at 5 2.8
\put{$\ssize {01100 \atop 01110}$} at 5 2
\put{$\ssize {11110 \atop 11110}$} at 5 0

\put{$\ssize {00100 \atop 00100}$} at 6 5
\put{$\ssize {01100 \atop 01221}$} at  6 3
\put{$\ssize {01110 \atop 01110}$} at  6 1

\arr{4.3 0.7} {4.7 0.3} 
\arr{5.3 0.3} {5.7 0.7} 
\arr{0.3 1.3} {0.7 1.7} 
\arr{1.3 1.7} {1.7 1.3} 
\arr{2.3 1.3} {2.7 1.7} 
\arr{3.3 1.7} {3.7 1.3} 
\arr{4.3 1.3} {4.7 1.7} 
\arr{5.3 1.7} {5.7 1.3} 
\arr{0.3 2.7} {0.7 2.3} 
\arr{1.3 2.3} {1.7 2.7} 
\arr{2.3 2.7} {2.7 2.3} 
\arr{3.3 2.3} {3.7 2.7} 
\arr{4.3 2.7} {4.7 2.3} 
\arr{5.3 2.3} {5.7 2.7} 
\arr{0.3 3.3} {0.7 3.7} 
\arr{1.3 3.7} {1.7 3.3} 
\arr{2.3 3.3} {2.7 3.7} 
\arr{3.3 3.7} {3.7 3.3} 
\arr{4.3 3.3} {4.7 3.7} 
\arr{5.3 3.7} {5.7 3.3} 
\arr{0.3 4.7} {0.7 4.3} 
\arr{1.3 4.3} {1.7 4.7} 
\arr{2.3 4.7} {2.7 4.3} 
\arr{3.3 4.3} {3.7 4.7} 
\arr{4.3 4.7} {4.7 4.3} 
\arr{5.3 4.3} {5.7 4.7} 
\arr{2.3 5.3} {2.7 5.7} 
\arr{3.3 5.7} {3.7 5.3} 
\arr{0.3 2.94}{0.7 2.86}
\arr{1.3 2.86}{1.7 2.94}
\arr{2.3 2.94}{2.7 2.86}
\arr{3.3 2.86}{3.7 2.94}
\arr{4.3 2.94}{4.7 2.86}
\arr{5.3 2.86}{5.7 2.94}
\setdots<2pt>
\plot 0.7 1  1.3 1 /
\plot 2.7 1  3.3 1 /
\plot 0.7 5  1.3 5 /
\plot 4.7 5  5.3 5 /
\setsolid
\plot 0 1.3  0 2.7 /
\plot 0 3.3  0 4.7 /
\plot 6 1.3  6 2.7 /
\plot 6 3.3  6 4.7 /
\endpicture} 
\quad
\hbox{\beginpicture
\setcoordinatesystem units <.7cm,.7cm>
\put{} at 0 0
\put{} at 7 6
\multiput{$\smallsq2$} at -.1 4.8  -.1 5  -.1 5.2 /
\put{$\sssize\bullet$} at 0 5.2 
\multiput{$\smallsq2$} at -.2 2.7  -.2 2.9  -.2 3.1  -.2 3.3  0 2.9  0 3.1 /
\multiput{$\sssize\bullet$} at -.1 2.9  .1 2.9 /
\plot -.1 2.9  .1 2.9 /
\put{$\smallsq2$} at -.1 1
\put{$\sssize \bullet$} at 0 1

\multiput{$\smallsq2$} at .8 3.7  .8 3.9  .8 4.1  .8 4.3  1 3.9  1 4.1 /
\multiput{$\sssize\bullet$} at .9 4.1  1.1 4.1 /
\plot .9 4.1  1.1 4.1 /
\multiput{$\smallsq2$} at .9 2.5  .9 2.7  .9 2.9  .9 3.1 /
\put{$\sssize\bullet$} at 1 2.7
\multiput{$\smallsq2$} at .9 1.9  .9 2.1 /
\put{$\sssize\bullet$} at 1 1.9

\multiput{$\smallsq2$} at 1.9 4.8  1.9 5  1.9 5.2 /
\multiput{$\smallsq2$} at 1.8 2.7  1.8 2.9  1.8 3.1  1.8 3.3  2 2.9  2 3.1 /
\multiput{$\sssize\bullet$} at 1.9 3.1  2.1 3.1  2.1 2.9 /
\plot 1.9 3.1  2.1 3.1 /
\put{$\smallsq2$} at 1.9 1

\multiput{$\smallsq2$} at 2.9 5.7  2.9 5.9  2.9 6.1  2.9 6.3 /
\multiput{$\smallsq2$} at 2.9 3.8  2.9 4  2.9 4.2 /
\put{$\sssize\bullet$} at 3 3.8
\multiput{$\smallsq2$} at 2.9 2.7  2.9 2.9 /
\put{$\sssize\bullet$} at 3 2.9
\multiput{$\smallsq2$} at 2.8 1.7  2.8 1.9  2.8 2.1  2.8 2.3  3 2.1 /
\multiput{$\sssize\bullet$} at 2.9 2.1  3.1 2.1 /
\plot 2.9 2.1  3.1 2.1 /

\multiput{$\smallsq2$} at 3.9 4.7  3.9 4.9  3.9 5.1  3.9 5.3 /
\put{$\sssize\bullet$} at 4 4.7
\multiput{$\smallsq2$} at 3.8 2.8  3.8 3  3.8 3.2  4 3 /
\multiput{$\sssize\bullet$} at 3.9 3  4.1 3 /
\plot 3.9 3  4.1 3 /
\multiput{$\smallsq2$} at 3.9 .7  3.9 .9  3.9 1.1  3.9 1.3 /
\put{$\sssize\bullet$} at 4 1.1 

\multiput{$\smallsq2$} at 4.8 3.7  4.8 3.9  4.8 4.1  4.8 4.3  5 3.9 /
\multiput{$\sssize \bullet$} at 4.9 3.9  5.1 3.9 /
\plot 4.9 3.9  5.1 3.9 /
\multiput{$\smallsq2$} at 4.9 2.7  4.9 2.9 /
\multiput{$\smallsq2$} at 4.9 1.8  4.9 2  4.9 2.2 /
\put{$\sssize\bullet$} at 5 2
\multiput{$\smallsq2$} at 4.9 -.3  4.9 -.1  4.9 .1  4.9 .3 /
\put{$\sssize\bullet$} at 5 .3

\put{$\smallsq2$} at 5.9 5
\put{$\sssize \bullet$} at 6 5
\multiput{$\smallsq2$} at 5.8 2.7  5.8 2.9  5.8 3.1  5.8 3.3  6 2.9  6 3.1 /
\multiput{$\sssize \bullet$} at 5.9 2.9  6.1 2.9 /
\plot 5.9 2.9  6.1 2.9 /
\multiput{$\smallsq2$} at 5.9 .8  5.9 1  5.9 1.2 /
\put{$\sssize\bullet$} at 6 1.2

\arr{4.3 0.7} {4.7 0.3} 
\arr{5.3 0.3} {5.7 0.7} 
\arr{0.3 1.3} {0.7 1.7} 
\arr{1.3 1.7} {1.7 1.3} 
\arr{2.3 1.3} {2.7 1.7} 
\arr{3.3 1.7} {3.7 1.3} 
\arr{4.3 1.3} {4.7 1.7} 
\arr{5.3 1.7} {5.7 1.3} 
\arr{0.3 2.7} {0.7 2.3} 
\arr{1.3 2.3} {1.7 2.7} 
\arr{2.3 2.7} {2.7 2.3} 
\arr{3.3 2.3} {3.7 2.7} 
\arr{4.3 2.7} {4.7 2.3} 
\arr{5.3 2.3} {5.7 2.7} 
\arr{0.3 3.3} {0.7 3.7} 
\arr{1.3 3.7} {1.7 3.3} 
\arr{2.3 3.3} {2.7 3.7} 
\arr{3.3 3.7} {3.7 3.3} 
\arr{4.3 3.3} {4.7 3.7} 
\arr{5.3 3.7} {5.7 3.3} 
\arr{0.3 4.7} {0.7 4.3} 
\arr{1.3 4.3} {1.7 4.7} 
\arr{2.3 4.7} {2.7 4.3} 
\arr{3.3 4.3} {3.7 4.7} 
\arr{4.3 4.7} {4.7 4.3} 
\arr{5.3 4.3} {5.7 4.7} 
\arr{2.3 5.3} {2.7 5.7} 
\arr{3.3 5.7} {3.7 5.3} 
\arr{0.3 2.94}{0.7 2.86}
\arr{1.3 2.86}{1.7 2.94}
\arr{2.3 2.94}{2.7 2.86}
\arr{3.3 2.86}{3.7 2.94}
\arr{4.3 2.94}{4.7 2.86}
\arr{5.3 2.86}{5.7 2.94}
\setdots<2pt>
\plot 0.3 1  1.7 1 /
\plot 2.3 1  3.7 1 /
\plot 0.3 5  1.7 5 /
\plot 4.3 5  5.7 5 /
\setsolid
\plot 0 1.2  0 2.5 /
\plot 0 3.5  0 4.6 /
\plot 6 1.4  6 2.5 /
\plot 6 3.5  6 4.8 /
\setshadegrid span <1.5mm>
\vshade   0   1 5 <,z,,> 
          2   1 5  <z,z,,> 
          3   1 6  <z,z,,>
          4   1 5  <z,z,,> 
          5   0 5  <z,z,,> 
          6   1 5 /
\endpicture} 
$$
Here are the possible dimension pairs of the indecomposable triples in $\Cal S(4)$.
$$
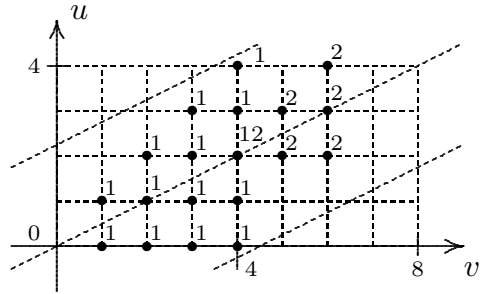
\hbox{\beginpicture
\setcoordinatesystem units <.5cm,.5cm>
\arr{0 -1}{0 5}
\arr{-1 0}{9 0}
\put{$u$} at  0.5 5.2
\put{$v$} at 9.2 -.5
\setdashes <.5mm> 
\setplotarea x from 0 to 8, y from 0 to 4 
\grid {8} {4}
\plot -1 -.5  9 4.5 /
\plot -1 1.75  4.5 4.5 /
\plot 3.5 -.5  9 2.25 /
\setsolid
\plot 4 -0.5  4 0.2 /
\put{$\ssize 4$} at 4.3 -.5
\plot 8 -0.2  8 0.2 /
\put{$\ssize 8$} at 8 -.5
\plot -0.2 4  0.2 4 /
\put{$\ssize 4$} at -.5 4 
\put{$\ssize 0$} at -.5 .3
\multiput{$\ssize \bullet$} at 1 0  2 0  3 0  4 0  1 1  2 1  3 1  4 1  
                             2 2  3 2  4 2  3 3  4 3  4 4  
                             4 2  5 2  6 2  5 3  6 3  6 4  /
\multiput{$\ssize 1$} at 1.2 .3   2.2 .3  3.2 .3   4.2 .3   1.2 1.3  2.2 1.4  3.2 1.3  
                            4.2 1.3  2.2 2.3  3.2 2.3  3.2 3.3  4.2 3.3  4.5 4.25 /
\multiput{$\ssize 2$} at 5.2 2.3  6.2 2.3  5.2 3.3  6.2 3.4  6.2 4.3 /
\multiput{$\ssize 12$} at 4.3 2.5 /
\endpicture}
$$

\medskip\noindent{\bf (6.5) The case $n = 5$.}
For the construction of the Auslander-Reiten quiver for $\Cal S(5)$
we use the methods of Section 1.  Here is the result.
$$
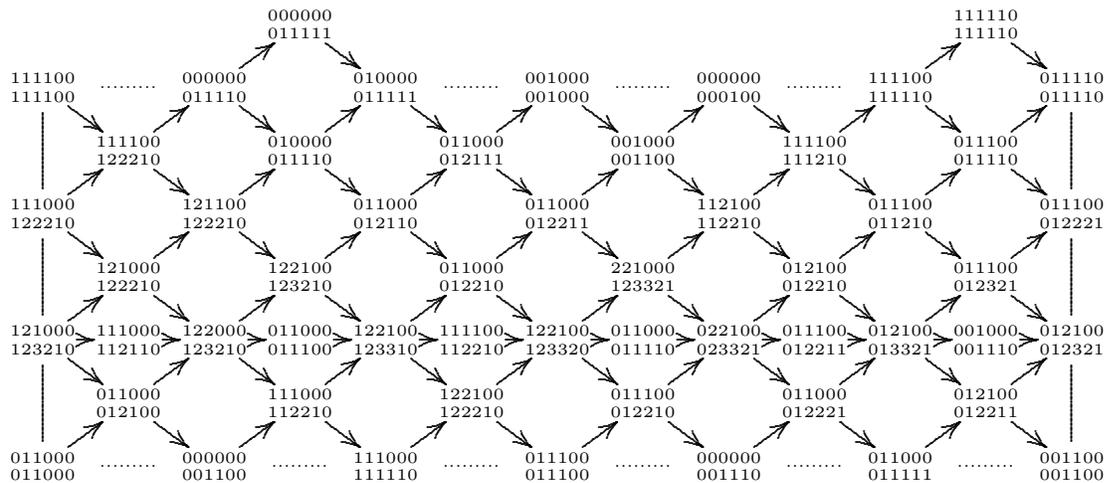
\hbox{\beginpicture
\setcoordinatesystem units <0.95cm,0.7cm>
\put{} at 0 0
\put{} at 12 7
\put{$\ssize {011000 \atop 011000}$}   at 0 0 
\put{$\ssize {121000 \atop 123210}$}   at 0 2
\put{$\ssize {111000 \atop 122210}$}   at 0 4
\put{$\ssize {111100 \atop 111100}$}   at 0 6
\put{$\ssize {011000 \atop 012100}$}   at 1 1
\put{$\ssize {111000 \atop 112110}$}   at 1 2
\put{$\ssize {121000 \atop 122210}$}   at 1 3
\put{$\ssize {111100 \atop 122210}$}   at 1 5
\put{$\ssize {000000 \atop 001100}$}   at 2 0
\put{$\ssize {122000 \atop 123210}$}   at 2 2
\put{$\ssize {121100 \atop 122210}$}   at 2 4
\put{$\ssize {000000 \atop 011110}$}   at 2 6
\put{$\ssize {111000 \atop 112210}$}   at 3 1
\put{$\ssize {011000 \atop 011100}$}   at 3 2
\put{$\ssize {122100 \atop 123210}$}   at 3 3
\put{$\ssize {010000 \atop 011110}$}   at 3 5
\put{$\ssize {000000 \atop 011111}$}   at 3 7
\put{$\ssize {111000 \atop 111110}$}   at 4 0
\put{$\ssize {122100 \atop 123310}$}   at 4 2
\put{$\ssize {011000 \atop 012110}$}   at 4 4
\put{$\ssize {010000 \atop 011111}$}   at 4 6
\put{$\ssize {122100 \atop 122210}$}   at 5 1
\put{$\ssize {111100 \atop 112210}$}   at 5 2 
\put{$\ssize {011000 \atop 012210}$}   at 5 3
\put{$\ssize {011000 \atop 012111}$}   at 5 5
\put{$\ssize {011100 \atop 011100}$}   at 6 0
\put{$\ssize {122100 \atop 123320}$}   at 6 2
\put{$\ssize {011000 \atop 012211}$}   at 6 4
\put{$\ssize {001000 \atop 001000}$}   at 6 6
\put{$\ssize {011100 \atop 012210}$}   at 7 1
\put{$\ssize {011000 \atop 011110}$}   at 7 2
\put{$\ssize {221000 \atop 123321}$}   at 7 3 
\put{$\ssize {001000 \atop 001100}$}   at 7 5
\put{$\ssize {000000 \atop 001110}$}   at 8 0
\put{$\ssize {022100 \atop 023321}$}   at 8 2
\put{$\ssize {112100 \atop 112210}$}   at 8 4
\put{$\ssize {000000 \atop 000100}$}   at 8 6
\put{$\ssize {011000 \atop 012221}$}   at 9 1
\put{$\ssize {011100 \atop 012211}$}   at 9 2
\put{$\ssize {012100 \atop 012210}$}   at 9 3
\put{$\ssize {111100 \atop 111210}$}   at 9 5 
\put{$\ssize {011000 \atop 011111}$}   at 10 0 
\put{$\ssize {012100 \atop 013321}$}   at 10 2
\put{$\ssize {011100 \atop 011210}$}   at 10 4
\put{$\ssize {111100 \atop 111110}$}   at 10 6
\put{$\ssize {012100 \atop 012211}$}   at 11 1
\put{$\ssize {001000 \atop 001110}$}   at 11 2
\put{$\ssize {011100 \atop 012321}$}   at 11 3
\put{$\ssize {011100 \atop 011110}$}   at 11 5 
\put{$\ssize {111110 \atop 111110}$}   at 11 7
\put{$\ssize {001100 \atop 001100}$}   at 12 0
\put{$\ssize {012100 \atop 012321}$}   at 12 2
\put{$\ssize {011100 \atop 012221}$}   at 12 4
\put{$\ssize {011110 \atop 011110}$}   at 12 6

\arr{0.3 0.3} {0.7 0.7}
\arr{2.3 0.3} {2.7 0.7}
\arr{4.3 0.3} {4.7 0.7}
\arr{6.3 0.3} {6.7 0.7}
\arr{8.3 0.3} {8.7 0.7}
\arr{10.3 0.3} {10.7 0.7}

\arr{1.3 0.7} {1.7 0.3}
\arr{3.3 0.7} {3.7 0.3}
\arr{5.3 0.7} {5.7 0.3}
\arr{7.3 0.7} {7.7 0.3}
\arr{9.3 0.7} {9.7 0.3}
\arr{11.3 0.7} {11.7 0.3}

\arr{0.3 1.7} {0.7 1.3} 
\arr{2.3 1.7} {2.7 1.3} 
\arr{4.3 1.7} {4.7 1.3} 
\arr{6.3 1.7} {6.7 1.3} 
\arr{8.3 1.7} {8.7 1.3} 
\arr{10.3 1.7} {10.7 1.3} 

\arr{1.3 1.3} {1.7 1.7} 
\arr{3.3 1.3} {3.7 1.7} 
\arr{5.3 1.3} {5.7 1.7} 
\arr{7.3 1.3} {7.7 1.7} 
\arr{9.3 1.3} {9.7 1.7} 
\arr{11.3 1.3} {11.7 1.7} 

\arr{0.4 2.0} {0.6 2.0}
\arr{1.4 2.0} {1.6 2.0}
\arr{2.4 2.0} {2.6 2.0}
\arr{3.4 2.0} {3.6 2.0}
\arr{4.4 2.0} {4.6 2.0}
\arr{5.4 2.0} {5.6 2.0}
\arr{6.4 2.0} {6.6 2.0}
\arr{7.4 2.0} {7.6 2.0}
\arr{8.4 2.0} {8.6 2.0}
\arr{9.4 2.0} {9.6 2.0}
\arr{10.4 2.0} {10.6 2.0}
\arr{11.4 2.0} {11.6 2.0}

\arr{0.3 2.3} {0.7 2.7} 
\arr{2.3 2.3} {2.7 2.7} 
\arr{4.3 2.3} {4.7 2.7} 
\arr{6.3 2.3} {6.7 2.7} 
\arr{8.3 2.3} {8.7 2.7} 
\arr{10.3 2.3} {10.7 2.7} 

\arr{1.3 2.7} {1.7 2.3} 
\arr{3.3 2.7} {3.7 2.3} 
\arr{5.3 2.7} {5.7 2.3} 
\arr{7.3 2.7} {7.7 2.3} 
\arr{9.3 2.7} {9.7 2.3} 
\arr{11.3 2.7} {11.7 2.3} 

\arr{0.3 3.7} {0.7 3.3} 
\arr{2.3 3.7} {2.7 3.3} 
\arr{4.3 3.7} {4.7 3.3} 
\arr{6.3 3.7} {6.7 3.3} 
\arr{8.3 3.7} {8.7 3.3} 
\arr{10.3 3.7} {10.7 3.3} 

\arr{1.3 3.3} {1.7 3.7} 
\arr{3.3 3.3} {3.7 3.7} 
\arr{5.3 3.3} {5.7 3.7} 
\arr{7.3 3.3} {7.7 3.7} 
\arr{9.3 3.3} {9.7 3.7} 
\arr{11.3 3.3} {11.7 3.7} 

\arr{0.3 4.3} {0.7 4.7} 
\arr{2.3 4.3} {2.7 4.7} 
\arr{4.3 4.3} {4.7 4.7} 
\arr{6.3 4.3} {6.7 4.7} 
\arr{8.3 4.3} {8.7 4.7} 
\arr{10.3 4.3} {10.7 4.7} 

\arr{1.3 4.7} {1.7 4.3} 
\arr{3.3 4.7} {3.7 4.3} 
\arr{5.3 4.7} {5.7 4.3} 
\arr{7.3 4.7} {7.7 4.3} 
\arr{9.3 4.7} {9.7 4.3} 
\arr{11.3 4.7} {11.7 4.3} 

\arr{0.3 5.7} {0.7 5.3} 
\arr{2.3 5.7} {2.7 5.3} 
\arr{4.3 5.7} {4.7 5.3} 
\arr{6.3 5.7} {6.7 5.3} 
\arr{8.3 5.7} {8.7 5.3} 
\arr{10.3 5.7} {10.7 5.3} 

\arr{1.3 5.3} {1.7 5.7} 
\arr{3.3 5.3} {3.7 5.7} 
\arr{5.3 5.3} {5.7 5.7} 
\arr{7.3 5.3} {7.7 5.7} 
\arr{9.3 5.3} {9.7 5.7} 
\arr{11.3 5.3} {11.7 5.7} 

\arr{2.3 6.3} {2.7 6.7}
\arr{3.3 6.7} {3.7 6.3}
\arr{10.3 6.3} {10.7 6.7}
\arr{11.3 6.7} {11.7 6.3}

\setdots<2pt>
\plot 0.7 0  1.3 0 /
\plot 2.7 0  3.3 0 /
\plot 4.7 0  5.3 0 /
\plot 6.7 0  7.3 0 /
\plot 8.7 0  9.3 0 /
\plot 10.7 0  11.3 0 /

\plot 0.7 6  1.3 6 /
\plot 4.7 6  5.3 6 /
\plot 6.7 6  7.3 6 /
\plot 8.7 6  9.3 6 /

\setsolid                       
\plot 0 0.4  0 1.6 /
\plot 0 2.4  0 3.6 /
\plot 0 4.4  0 5.6 /

\plot 12 0.4  12 1.6 /
\plot 12 2.4  12 3.6 /
\plot 12 4.4  12 5.6 /
\endpicture} 
$$

\def \ssq{$\smallsq2$}
\def \bul{$\sssize\bullet$}

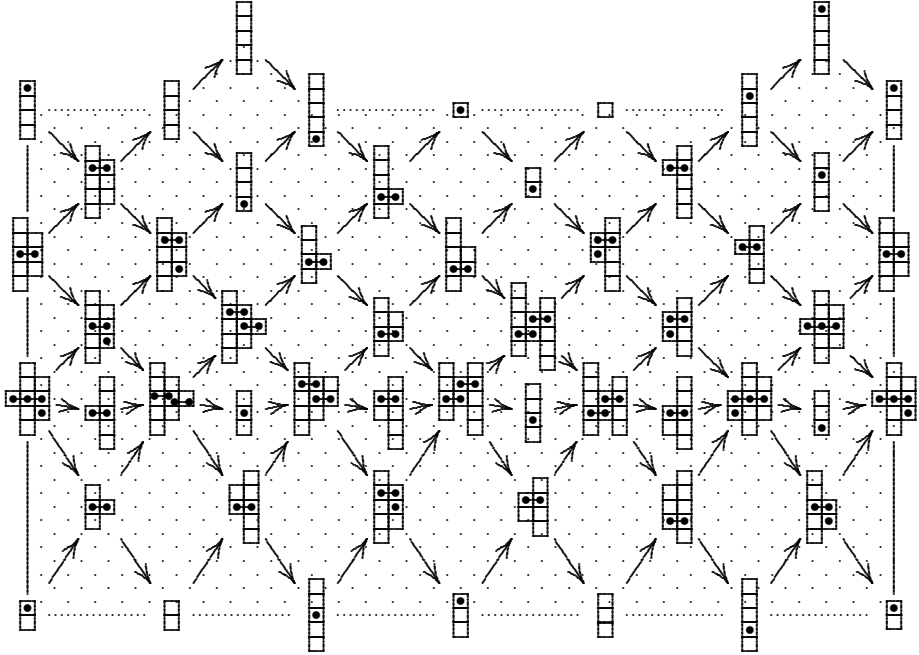
$$
\hbox{\beginpicture
\setcoordinatesystem units <.8cm,.8cm>
\multiput{\ssq} at -.1 -.1  -.1 .1 /
\put{\bul} at 0 .1 
\multiput{\ssq} at 1.9 -.1  1.9 .1 /  
\multiput{\ssq} at 3.9 -.4  3.9 -.2  3.9 0  3.9 .2  3.9 .4 /
\put{\bul} at 4 0
\multiput{\ssq} at 5.9 -.2  5.9 0  5.9 .2 /
\put{\bul} at 6 .2
\multiput{\ssq} at 7.9 -.2  7.9 0  7.9 .2 /
\multiput{\ssq} at 9.9 -.4  9.9 -.2  9.9 0  9.9 .2  9.9 .4 /
\put{\bul} at 10 -.2  
\multiput{\ssq} at 11.9 -.1  11.9 .1 /
\put{\bul} at 12 .1
\multiput{\ssq} at .8 1.3  .8 1.5  .8 1.7  1 1.5 
                3 1.1  3 1.3  3 1.5  3 1.7  3 1.9  2.8 1.5  2.8 1.7 
                4.8 1.1  4.8 1.3  4.8 1.5  4.8 1.7  4.8 1.9  5 1.3  5 1.5  
                                5 1.7
                7 1.2  7 1.4  7 1.6  7 1.8  6.8 1.4  6.8 1.6
                9 1.1  9 1.3  9 1.5  9 1.7  9 1.9  8.8 1.3  8.8 1.5  8.8 1.7
                10.8 1.1  10.8 1.3  10.8 1.5  10.8 1.7  10.8 1.9  
                                11 1.3  11 1.5 /
\multiput{\bul} at .9 1.5  1.1 1.5  2.9 1.5  3.1 1.5  4.9 1.7  5.1 1.7  5.1 1.5
              6.9 1.6  7.1 1.6  8.9 1.3  9.1 1.3  10.9 1.5  11.1 1.5  
                                11.1 1.3 /
\plot .9 1.5  1.1 1.5 /
\plot 2.9 1.5  3.1 1.5 /
\plot 4.9 1.7  5.1 1.7 /
\plot 6.9 1.6  7.1 1.6 /
\plot 8.9 1.3  9.1 1.3 /
\plot 10.9 1.5  11.1 1.5 /
%
%
\multiput{\ssq} at 
-.1 2.6  -.1 2.8  -.1 3  -.1 3.2  -.1 3.4  .1 2.8  .1 3  .1 3.2  -.3 3 
1 2.4  1 2.6  1 2.8  1 3  1 3.2  .8 2.8
1.7 2.6  1.7 2.8  1.7 3  1.7 3.2  1.7 3.4  1.9 2.8  1.9 3  1.9 3.2  2.1 3
2.9 2.6  2.9 2.8  2.9 3
3.7 2.6  3.7 2.8  3.7 3  3.7 3.2  3.7 3.4  3.9 2.8  3.9 3  3.9 3.2  
                4.1 3  4.1 3.2
5 2.4  5 2.6  5 2.8  5 3  5 3.2  4.8 2.8  4.8 3
5.7 2.8  5.7 3  5.7 3.2  5.7 3.4  6.1 2.6  6.1 2.8  6.1 3  6.1 3.2  6.1 3.4
                5.9 3  5.9 3.2
6.9 2.5  6.9 2.7  6.9 2.9  6.9 3.1
7.7 2.6  7.7 2.8  7.7 3  7.7 3.2  7.7 3.4  8.1 2.6  8.1 2.8  8.1 3  8.1 3.2
                7.9 2.8  7.9 3
9 2.4  9 2.6  9 2.8  9 3  9 3.2  8.8 2.6  8.8 2.8 
9.9 2.6  9.9 2.8  9.9 3  9.9 3.2  9.9 3.4  10.1 2.8  10.1 3  10.1 3.2  9.7 2.8
               9.7 3 
10.9 2.6  10.9 2.8  10.9 3
11.9 2.6  11.9 2.8  11.9 3  11.9 3.2  11.9 3.4  12.1 2.8  12.1 3  12.1 3.2
               11.7 3 /
\multiput{\bul} at -.2 3  0 3  .2 3  .2 2.8  .9 2.8  1.1 2.8  1.76 3.04  
        1.96 3.04
        2.04 2.96  2.24 2.96  3 2.8  3.8 3.2  4 3.2  4 3  4.2 3  4.9 3  5.1 3
        5.8 3  6 3  6 3.2  6.2 3.2  7 2.7  7.8 2.8  8 2.8  8 3  8.2 3
        8.9 2.8  9.1 2.8  9.8 2.8  9.8 3  10 3  10.2 3  11 2.6  11.8 3
        12 3  12.2 3  12.2 2.8 /
\plot -.2 3  .2 3 /
\plot .9 2.8  1.1 2.8 /
\plot 1.76 3.04  1.96 3.04 /
\plot 2.04 2.96  2.24 2.96 /
\plot 3.8 3.2  4 3.2 /
\plot 4 3  4.2 3 /
\plot 4.9 3  5.1 3 /
\plot 5.8 3  6 3 /
\plot 6 3.2  6.2 3.2 /
\plot 7.8 2.8  8 2.8 /
\plot 8 3  8.2 3 /
\plot 8.9 2.8  9.1 2.8 /
\plot 9.8 3  10.2 3 /
\plot 11.8 3  12.2 3 / 
\multiput{\ssq} at .8 3.6  .8 3.8  .8 4  .8 4.2  .8 4.4  1 3.8  1 4  1 4.2
        2.7 3.6  2.7 3.8  2.7 4  2.7 4.2  2.7 4.4  2.9 3.8  2.9 4  2.9 4.2
                                3.1 4
        4.8 3.7  4.8 3.9  4.8 4.1  4.8 4.3  5 3.9  5 4.1
        6.7 3.7  6.7 3.9  6.7 4.1  6.7 4.3  6.7 4.5  6.9 3.9  6.9 4.1  7.1 
                                3.5  7.1 3.7  7.1 3.9  7.1 4.1  7.1 4.3
        8.8 3.9  8.8 4.1  9 3.7  9 3.9  9 4.1  9 4.3
        10.7 4  10.9 3.6  10.9 3.8  10.9 4  10.9 4.2  10.9 4.4  11.1 3.8
                                11.1 4  11.1 4.2 /
\multiput{\bul} at .9 4  1.1 4  1.1 3.8  2.8 4.2  3 4.2  3 4  3.2 4  
        4.9 3.9  5.1 3.9  6.8 3.9  7 3.9  7 4.1  7.2 4.1  8.9 3.9  8.9 4.1 
        9.1 4.1  10.8 4  11 4  11.2 4 /
\plot .9 4  1.1 4 /
\plot 2.8 4.2  3 4.2 /
\plot 3 4  3.2 4 /
\plot 4.9 3.9  5.1 3.9 /
\plot 6.8 3.9 7 3.9 /
\plot 7 4.1  7.2 4.1 /
\plot 8.9 4.1  9.1 4.1 /
\plot 10.8 4  11.2 4 / 
\multiput{\ssq} at -.2 4.6  -.2 4.8  -.2 5  -.2 5.2  -.2 5.4  0 4.8  0 5  0 5.2
        1.8 4.6  1.8 4.8  1.8 5  1.8 5.2  1.8 5.4  2 4.8  2 5  2 5.2
        3.8 4.7  3.8 4.9  3.8 5.1  3.8 5.3  4 4.9
        5.8 4.6  5.8 4.8  5.8 5  5.8 5.2  5.8 5.4  6 4.8  6 5
        7.8 5  7.8 5.2  8 4.6  8 4.8  8 5  8 5.2  8 5.4
        9.8 5.1  10 4.7  10 4.9  10 5.1  10 5.3
        11.8 4.6  11.8 4.8  11.8 5  11.8 5.2  11.8 5.4  12 4.8  12 5  12 5.2 /
\multiput{\bul} at -.1 5  .1 5  1.9 5.2  2.1 5.2  2.1 4.8  3.9 4.9  4.1 4.9
        5.9 4.8  6.1 4.8  7.9 5  7.9 5.2  8.1 5.2  9.9 5.1  10.1 5.1  
        11.9 5  12.1 5 /
\plot -.1 5  .1 5  /
\plot 1.9 5.2  2.1 5.2 /
\plot 3.9 4.9  4.1 4.9 /
\plot 5.9 4.8  6.1 4.8 /
\plot 7.9 5.2  8.1 5.2 /
\plot 9.9 5.1  10.1 5.1 /
\plot 11.9 5  12.1 5 /
\multiput{\ssq} at .8 5.6  .8 5.8  .8 6  .8 6.2  .8 6.4  1 5.8  1 6  1 6.2
        2.9 5.7  2.9 5.9  2.9 6.1  2.9 6.3 
        4.8 5.6  4.8 5.8  4.8 6  4.8 6.2  4.8 6.4  5 5.8
        6.9 5.9  6.9 6.1
        8.8 6.2  9 5.6  9 5.8  9 6  9 6.2  9 6.4
        10.9 5.7  10.9 5.9  10.9 6.1  10.9 6.3 /
\multiput{\bul} at .9 6.2  1.1 6.2  3 5.7  4.9 5.8  5.1 5.8  7 5.9  8.9 6.2
        9.1 6.2  11 6.1 /
\plot .9 6.2  1.1 6.2 /
\plot 4.9 5.8  5.1 5.8 /
\plot 8.9 6.2  9.1 6.2 /
\multiput{\ssq} at -.1 6.7  -.1 6.9  -.1 7.1  -.1 7.3  1.9 6.7  1.9 6.9
        1.9 7.1  1.9 7.3
        2.9 7.6  2.9 7.8  2.9 8  2.9 8.2  2.9 8.4  
        3.9 6.6  3.9 6.8  3.9 7  3.9 7.2  3.9 7.4
        5.9 7
        7.9 7
        9.9 6.6  9.9 6.8  9.9 7  9.9 7.2  9.9 7.4
        10.9 7.6  10.9 7.8  10.9 8  10.9 8.2  10.9 8.4
        11.9 6.7  11.9 6.9  11.9 7.1  11.9 7.3 /
\multiput{\bul} at 0 7.3  4 6.6  6 7 10 7.2  11 8.4  12 7.3 /
\arr{0.3 0.45} {0.7 1.05} 
\arr{1.3 1.05} {1.7 0.45} 
\arr{2.3 0.45} {2.7 1.05} 
\arr{3.3 1.05} {3.7 0.45} 
\arr{4.3 0.45} {4.7 1.05} 
\arr{5.3 1.05} {5.7 0.45} 
\arr{6.3 0.45} {6.7 1.05} 
\arr{7.3 1.05} {7.7 0.45} 
\arr{8.3 0.45} {8.7 1.05} 
\arr{9.3 1.05} {9.7 0.45} 
\arr{10.3 0.45} {10.7 1.05} 
\arr{11.3 1.05} {11.7 0.45} 
\arr{0.3 2.55} {0.7 1.95} 
\arr{1.3 1.95} {1.6 2.4} 
\arr{2.3 2.55} {2.7 1.95} 
\arr{3.3 1.95} {3.6 2.4} 
\arr{4.3 2.55} {4.7 1.95} 
\arr{5.3 1.95} {5.7 2.55} 
\arr{6.4 2.4 } {6.7 1.95} 
\arr{7.3 1.95} {7.6 2.4 } 
\arr{8.4 2.4 } {8.7 1.95} 
\arr{9.3 1.95} {9.7 2.55} 
\arr{10.3 2.55} {10.7 1.95} 
\arr{11.3 1.95} {11.7 2.55} 
\arr{0.4 3.4} {0.7 3.7} 
\arr{1.3 3.7} {1.6 3.4} 
\arr{2.3 3.3} {2.6 3.6} 
\arr{3.3 3.7} {3.6 3.4} 
\arr{4.4 3.4} {4.7 3.7} 
\arr{5.3 3.7} {5.6 3.4} 
\arr{6.4 3.4} {6.6 3.6} 
\arr{7.4 3.6} {7.6 3.4} 
\arr{8.4 3.4} {8.7 3.7} 
\arr{9.3 3.7} {9.7 3.3} 
\arr{10.4 3.4} {10.7 3.7} 
\arr{11.4 3.6} {11.7 3.3} 
\arr{0.3 4.7} {0.7 4.3} 
\arr{1.3 4.3} {1.7 4.7} 
\arr{2.3 4.7} {2.6 4.4} 
\arr{3.3 4.3} {3.7 4.7} 
\arr{4.3 4.7} {4.7 4.3} 
\arr{5.3 4.3} {5.7 4.7} 
\arr{6.3 4.7} {6.6 4.4} 
\arr{7.4 4.4} {7.7 4.7} 
\arr{8.3 4.7} {8.7 4.3} 
\arr{9.3 4.3} {9.7 4.7} 
\arr{10.3 4.7} {10.7 4.3} 
\arr{11.4 4.4} {11.7 4.7} 
\arr{0.3 5.3} {0.7 5.7} 
\arr{1.3 5.7} {1.7 5.3} 
\arr{2.3 5.3} {2.7 5.7} 
\arr{3.3 5.7} {3.7 5.3} 
\arr{4.3 5.3} {4.7 5.7} 
\arr{5.3 5.7} {5.7 5.3} 
\arr{6.3 5.3} {6.7 5.7} 
\arr{7.3 5.7} {7.7 5.3} 
\arr{8.3 5.3} {8.7 5.7} 
\arr{9.3 5.7} {9.7 5.3} 
\arr{10.3 5.3} {10.7 5.7} 
\arr{11.3 5.7} {11.7 5.3} 
\arr{0.3 6.7} {0.7 6.3} 
\arr{1.3 6.3} {1.7 6.7} 
\arr{2.3 6.7} {2.7 6.3} 
\arr{3.3 6.3} {3.7 6.7} 
\arr{4.3 6.7} {4.7 6.3} 
\arr{5.3 6.3} {5.7 6.7} 
\arr{6.3 6.7} {6.7 6.3} 
\arr{7.3 6.3} {7.7 6.7} 
\arr{8.3 6.7} {8.7 6.3} 
\arr{9.3 6.3} {9.7 6.7} 
\arr{10.3 6.7} {10.7 6.3} 
\arr{11.3 6.3} {11.7 6.7} 

\arr{2.3 7.3} {2.7 7.7} 
\arr{3.3 7.7} {3.7 7.3} 
\arr{10.3 7.3} {10.7 7.7} 
\arr{11.3 7.7} {11.7 7.3} 

\arr{0.4 2.92}{0.7 2.86}
\arr{1.3 2.86}{1.6 2.92}
\arr{2.4 2.92}{2.7 2.86}
\arr{3.3 2.86}{3.6 2.92}
\arr{4.4 2.92}{4.7 2.86}
\arr{5.3 2.86}{5.6 2.92}
\arr{6.4 2.92}{6.7 2.86}
\arr{7.3 2.86}{7.6 2.92}
\arr{8.4 2.92}{8.7 2.86}
\arr{9.3 2.86}{9.6 2.92}
\arr{10.4 2.92}{10.7 2.86}
\arr{11.3 2.86}{11.6 2.92}
\setdots<2pt>
\plot 0.3 0  1.7 0 /
\plot 2.3 0  3.7 0 /
\plot 4.3 0  5.7 0 /
\plot 6.3 0  7.7 0 /
\plot 8.3 0  9.7 0 /
\plot 10.3 0  11.7 0 /
\plot 0.3 7  1.7 7 /
\plot 4.3 7  5.7 7 /
\plot 6.3 7  7.7 7 /
\plot 8.3 7  9.7 7 /
\setsolid
\plot 0 0.3  0 2.4 /
\plot 0 3.6  0 4.4 /
\plot 0 5.6  0 6.5 /
\plot 12 0.3  12 2.4 /
\plot 12 3.6  12 4.4 /
\plot 12 5.6  12 6.5 /
\setshadegrid span <1.5mm>
\vshade   0   0 7 <,z,,> 
          2   0 7  <z,z,,> 
          3   0 8  <z,z,,>
          4   0 7  <z,z,,> 
          10  0 7  <z,z,,> 
          11  0 8  <z,z,,> 
          12   0 7 /
\endpicture} 
$$

\bigskip
Let us indicate how $\Gamma(5)$ is constructed. 
We consider the algebra $\Psi$ given by
quiver and relations, as indicated.

$$
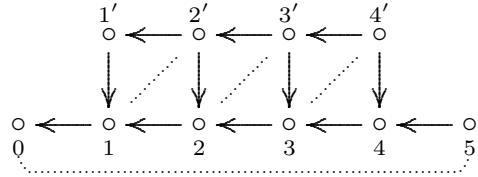
\hbox{\beginpicture
\setcoordinatesystem units <0.5cm,0.5cm>
\put{} at 0 -1.05
\put{} at 5 2
\put{$\circ$} at 0 0
\put{$\circ$} at 2 0
\put{$\circ$} at 4 0
\put{$\circ$} at 6 0
\put{$\circ$} at 8 0
\put{$\circ$} at 10 0
\put{$\circ$} at 2 2
\put{$\circ$} at 4 2
\put{$\circ$} at 6 2
\put{$\circ$} at 8 2
\arr{1.6 0}{0.4 0}
\arr{3.6 0}{2.4 0}
\arr{5.6 0}{4.4 0}
\arr{7.6 0}{6.4 0}
\arr{9.6 0}{8.4 0}

\arr{3.6 2}{2.4 2}
\arr{5.6 2}{4.4 2}
\arr{7.6 2}{6.4 2}

\arr{2 1.6}{2 0.4}
\arr{4 1.6}{4 0.4}
\arr{6 1.6}{6 0.4}
\arr{8 1.6}{8 0.4}


\put{$\ssize 1'$} at  2 2.5
\put{$\ssize 2'$} at  4 2.5
\put{$\ssize 3'$} at  6 2.5
\put{$\ssize 4'$} at  8 2.5

\put{$\ssize 0$} at 0 -.5
\put{$\ssize 1$} at 2 -.5
\put{$\ssize 2$} at 4 -.5
\put{$\ssize 3$} at 6 -.5
\put{$\ssize 4$} at 8 -.5
\put{$\ssize 5$} at 10 -.5

\setdots<2pt>
\plot 3.5 1.5  2.5 0.5 /
\plot 5.5 1.5  4.5 0.5 /
\plot 7.5 1.5  6.5 0.5 /
\plot 0 -0.7  0.1 -0.9  0.2 -1  0.3 -1.05  9.7 -1.05  9.8 -1  9.9 -0.9  10 -0.7 /
\endpicture}
$$

\medskip
The algebra $\Psi$ has a preprojective component; here is a sketch of the 
left hand part.

$$
\hbox{\beginpicture
\setcoordinatesystem units <0.45cm,0.45cm>
\put{} at 0 0
\put{} at 24 10
\put{$\bullet$}  at  0 3
\put{$\ssize P(0)$}  at -1 3
\put{$\bullet$}  at  1 4
\put{$\ssize P(1)$}  at  0 4
\put{$\bullet$}  at  2 5
\put{$\ssize P(2)$}  at  1 5
\put{$\bullet$}  at  3 6
\put{$\ssize P(3)$}  at  2 6
\put{$\bullet$}  at  4 7
\put{$\ssize P(4)$}  at  3 7
\put{$\bullet$}  at  11 8
\put{$\ssize P(3')$} at 11.5 8.5
\put{$\bullet$}  at  14 9
\put{$\ssize P(5)$}  at 14 9.5
\put{$\bullet$}  at  16 9
\put{$\ssize S(1')$} at 15.75 9.5
\put{$\bullet$}  at  23 10
\put{$\ssize P(4')=I(0)$} at 23 10.5
\put{$\bullet$}  at  2 4
\put{$\ssize P(1')$} at  2.5 3.5
\put{$\bullet$}  at  5 2
\put{$\ssize P(2')$} at  5 1.5

\put{$\bullet$}  at  17 2
\put{$\ssize P(2')[1]$} at 17 1.5
\put{$\bullet$}  at  24 9
\put{$\ssize P(3')[1]$} at 24.8 9.5

\put{$\circ$} at 7 2
\put{$\circ$} at 9 2
\put{$\circ$} at 11 2
\put{$\circ$} at 13 2
\put{$\circ$} at 15 2
\put{$\circ$} at 19 2
\put{$\circ$} at 21 2
\put{$\circ$} at 23 2

\put{$\circ$} at 2 3
\put{$\circ$} at 4 3
\put{$\circ$} at 6 3
\put{$\circ$} at 8 3
\put{$\circ$} at 10 3
\put{$\circ$} at 12 3
\put{$\circ$} at 14 3
\put{$\circ$} at 16 3
\put{$\circ$} at 18 3
\put{$\circ$} at 20 3
\put{$\circ$} at 22 3

\put{$\circ$} at 3 4
\put{$\circ$} at 4 4
\put{$\circ$} at 5 4
\put{$\circ$} at 6 4
\put{$\circ$} at 7 4
\put{$\circ$} at 8 4
\put{$\circ$} at 9 4
\put{$\circ$} at 10 4
\put{$\circ$} at 11 4
\put{$\circ$} at 12 4
\put{$\circ$} at 13 4
\put{$\circ$} at 14 4
\put{$\circ$} at 15 4
\put{$\circ$} at 16 4
\put{$\circ$} at 17 4
\put{$\circ$} at 18 4
\put{$\bullet$} at 19 4
\put{$\ssize A_4$} at 19 3.5
\put{$\circ$} at 20 4
\put{$\circ$} at 21 4
\put{$\circ$} at 22 4
\put{$\circ$} at 23 4

\put{$\circ$} at 4 5
\put{$\circ$} at 6 5
\put{$\circ$} at 8 5
\put{$\circ$} at 10 5
\put{$\circ$} at 12 5
\put{$\circ$} at 14 5
\put{$\circ$} at 16 5
\put{$\bullet$} at 18 5
\put{$\ssize A_3$} at 17.3 5
\put{$\circ$} at 20 5
\put{$\circ$} at 22 5

\put{$\circ$} at 5 6
\put{$\circ$} at 7 6
\put{$\circ$} at 9 6
\put{$\circ$} at 11 6
\put{$\circ$} at 13 6
\put{$\circ$} at 15 6
\put{$\bullet$} at 17 6
\put{$\ssize A_2$} at 16.3 6
\put{$\circ$} at 19 6
\put{$\bullet$} at 21 6
\put{$\ssize C_3$} at 21.7 6
\put{$\circ$} at 23 6

\put{$\circ$} at 6 7
\put{$\circ$} at 8 7
\put{$\circ$} at 10 7
\put{$\circ$} at 12 7
\put{$\circ$} at 14 7
\put{$\bullet$} at 16 7
\put{$\ssize A_1$} at 15.3 7
\put{$\circ$} at 18 7
\put{$\bullet$} at 20 7
\put{$\ssize C_2$} at 21.2 7
\put{$\circ$} at 22 7

\put{$\circ$} at 13 8
\put{$\bullet$} at 15 8
\put{$\ssize A_0$} at 14.3 8
\put{$\circ$} at 17 8
\put{$\bullet$} at 19 8
\put{$\ssize C_1$} at 19.7 8
\put{$\circ$} at 21 8
\put{$\circ$} at 23 8

\put{$\circ$} at 20 9
\put{$\bullet$} at 18 9
\put{$\ssize C_0$} at 18 9.5
\put{$\circ$} at 22 9

\setsolid 
\plot 5.3 2.3 5.7 2.7 /
\plot 7.3 2.3 7.7 2.7 /
\plot 9.3 2.3 9.7 2.7 /
\plot 11.3 2.3 11.7 2.7 /
\plot 13.3 2.3 13.7 2.7 /
\plot 15.3 2.3 15.7 2.7 /
\plot 17.3 2.3 17.7 2.7 /
\plot 19.3 2.3 19.7 2.7 /
\plot 21.3 2.3 21.7 2.7 /
\plot 23.3 2.3 23.7 2.7 /

\plot 4.3 2.7 4.7 2.3 /
\plot 6.3 2.7 6.7 2.3 /
\plot 8.3 2.7 8.7 2.3 /
\plot 10.3 2.7 10.7 2.3 /
\plot 12.3 2.7 12.7 2.3 /
\plot 14.3 2.7 14.7 2.3 /
\plot 16.3 2.7 16.7 2.3 /
\plot 18.3 2.7 18.7 2.3 /
\plot 20.3 2.7 20.7 2.3 /
\plot 22.3 2.7 22.7 2.3 /

\plot 0.3 3.3 0.7 3.7 /
\plot 4.3 3.3 4.7 3.7 /
\plot 6.3 3.3 6.7 3.7 /
\plot 8.3 3.3 8.7 3.7 /
\plot 10.3 3.3 10.7 3.7 /
\plot 12.3 3.3 12.7 3.7 /
\plot 14.3 3.3 14.7 3.7 /
\plot 16.3 3.3 16.7 3.7 /
\plot 18.3 3.3 18.7 3.7 /
\plot 20.3 3.3 20.7 3.7 /
\plot 22.3 3.3 22.7 3.7 /
 
\plot 1.3 3.7 1.7 3.3 /
\plot 3.3 3.7 3.7 3.3 /
\plot 5.3 3.7 5.7 3.3 /
\plot 7.3 3.7 7.7 3.3 /
\plot 9.3 3.7 9.7 3.3 /
\plot 11.3 3.7 11.7 3.3 /
\plot 13.3 3.7 13.7 3.3 /
\plot 15.3 3.7 15.7 3.3 /
\plot 17.3 3.7 17.7 3.3 /
\plot 19.3 3.7 19.7 3.3 /
\plot 21.3 3.7 21.7 3.3 /
\plot 23.3 3.7 23.7 3.3 /

\plot 1.3 4 1.6 4 /
\plot 2.4 4 2.7 4 /
\plot 3.3 4 3.7 4 /
\plot 4.3 4 4.7 4 /
\plot 5.3 4 5.7 4 /
\plot 6.3 4 6.7 4 /
\plot 7.3 4 7.7 4 /
\plot 8.3 4 8.7 4 /
\plot 9.3 4 9.7 4 /
\plot 10.3 4 10.7 4 /
\plot 11.3 4 11.7 4 /
\plot 12.3 4 12.7 4 /
\plot 13.3 4 13.7 4 /
\plot 14.3 4 14.7 4 /
\plot 15.3 4 15.7 4 /
\plot 16.3 4 16.7 4 /
\plot 17.3 4 17.7 4 /
\plot 18.3 4 18.7 4 /
\plot 19.3 4 19.7 4 /
\plot 20.3 4 20.7 4 /
\plot 21.3 4 21.7 4 /
\plot 22.3 4 22.7 4 /
\plot 23.3 4 23.7 4 /

\plot 1.3 4.3 1.7 4.7 /
\plot 3.3 4.3 3.7 4.7 /
\plot 5.3 4.3 5.7 4.7 /
\plot 7.3 4.3 7.7 4.7 /
\plot 9.3 4.3 9.7 4.7 /
\plot 11.3 4.3 11.7 4.7 /
\plot 13.3 4.3 13.7 4.7 /
\plot 15.3 4.3 15.7 4.7 /
\plot 17.3 4.3 17.7 4.7 /
\plot 19.3 4.3 19.7 4.7 /
\plot 21.3 4.3 21.7 4.7 /
\plot 23.3 4.3 23.7 4.7 /

\plot 2.3 4.7 2.7 4.3 /
\plot 4.3 4.7 4.7 4.3 /
\plot 6.3 4.7 6.7 4.3 /
\plot 8.3 4.7 8.7 4.3 /
\plot 10.3 4.7 10.7 4.3 /
\plot 12.3 4.7 12.7 4.3 /
\plot 14.3 4.7 14.7 4.3 /
\plot 16.3 4.7 16.7 4.3 /
\plot 18.3 4.7 18.7 4.3 /
\plot 20.3 4.7 20.7 4.3 /
\plot 22.3 4.7 22.7 4.3 /

\plot 2.3 5.3 2.7 5.7 /
\plot 4.3 5.3 4.7 5.7 /
\plot 6.3 5.3 6.7 5.7 /
\plot 8.3 5.3 8.7 5.7 /
\plot 10.3 5.3 10.7 5.7 /
\plot 12.3 5.3 12.7 5.7 /
\plot 14.3 5.3 14.7 5.7 /
\plot 16.3 5.3 16.7 5.7 /
\plot 18.3 5.3 18.7 5.7 /
\plot 20.3 5.3 20.7 5.7 /
\put{$\ssize a$} at 20.9 5.5
\plot 22.3 5.3 22.7 5.7 /

\plot 3.3 5.7 3.7 5.3 /
\plot 5.3 5.7 5.7 5.3 /
\plot 7.3 5.7 7.7 5.3 /
\plot 9.3 5.7 9.7 5.3 / 
\plot 11.3 5.7 11.7 5.3 /
\plot 13.3 5.7 13.7 5.3 /
\plot 15.3 5.7 15.7 5.3 /
\plot 17.3 5.7 17.7 5.3 /
\plot 19.3 5.7 19.7 5.3 /
\plot 21.3 5.7 21.7 5.3 /
\plot 23.3 5.7 23.7 5.3 /

\plot 3.3 6.3 3.7 6.7 /
\plot 5.3 6.3 5.7 6.7 /
\plot 7.3 6.3 7.7 6.7 /
\plot 9.3 6.3 9.7 6.7 /
\plot 11.3 6.3 11.7 6.7 /
\plot 13.3 6.3 13.7 6.7 /
\plot 15.3 6.3 15.7 6.7 /
\plot 17.3 6.3 17.7 6.7 /
\plot 19.3 6.3 19.7 6.7 /
\put{$\ssize a$} at 19.8 6.6
\plot 21.3 6.3 21.7 6.7 /
\plot 23.3 6.3 23.7 6.7 /

\plot 4.3 6.7 4.7 6.3 /
\plot 6.3 6.7 6.7 6.3 /
\plot 8.3 6.7 8.7 6.3 /
\plot 10.3 6.7 10.7 6.3 /
\plot 12.3 6.7 12.7 6.3 /
\plot 14.3 6.7 14.7 6.3 /
\plot 16.3 6.7 16.7 6.3 /
\plot 18.3 6.7 18.7 6.3 /
\plot 20.3 6.7 20.7 6.3 /
\plot 22.3 6.7 22.7 6.3 /

\plot 10.3 7.3 10.7 7.7 /
\plot 12.3 7.3 12.7 7.7 /
\plot 14.3 7.3 14.7 7.7 /
\plot 16.3 7.3 16.7 7.7 /
\plot 18.3 7.3 18.7 7.7 /
\put{$\ssize a$} at 18.9 7.5
\plot 20.3 7.3 20.7 7.7 /
\plot 22.3 7.3 22.7 7.7 /

\plot 11.3 7.7 11.7 7.3 /
\plot 13.3 7.7 13.7 7.3 /
\plot 15.3 7.7 15.7 7.3 /
\plot 17.3 7.7 17.7 7.3 /
\plot 19.3 7.7 19.7 7.3 /
\plot 21.3 7.7 21.7 7.3 /
\plot 23.3 7.7 23.7 7.3 /

\plot 13.3 8.3 13.7 8.7 /
\plot 15.3 8.3 15.7 8.7 /
\plot 17.3 8.3 17.7 8.7 /
\put{$\ssize a$} at 17.9 8.5
\plot 19.3 8.3 19.7 8.7 /
\plot 21.3 8.3 21.7 8.7 /
\plot 23.3 8.3 23.7 8.7 /

\plot 14.3 8.7 14.7 8.3 /
\plot 16.3 8.7 16.7 8.3 /
\plot 18.3 8.7 18.7 8.3 /
\plot 20.3 8.7 20.7 8.3 /
\plot 22.3 8.7 22.7 8.3 /
\plot 24.3 8.7 24.7 8.3 /

\plot 22.3 9.3 22.7 9.7 /

\plot 23.3 9.7 23.7 9.3 /

\setdots<2pt>
\plot 0.4 3   1.6 3 /
\plot 2.4 3   3.6 3 /
\plot 5.4 2   6.6 2 /
\plot 7.4 2   8.6 2 /
\plot 9.4 2   10.6 2 /
\plot 11.4 2  12.6 2 /
\plot 13.4 2  14.6 2 /
\plot 15.4 2  16.2 2 /
\plot 17.8 2  18.6 2 /
\plot 19.4 2  20.6 2 /
\plot 21.4 2  22.6 2 /
\plot 23.4 2  24.6 2 /

\plot  4.4 7   5.6 7 /
\plot  6.4 7   7.6 7 /
\plot  8.4 7   9.6 7 /
\plot 11.4 8  12.6 8 /
\plot 14.4 9  15.6 9 /
\plot 16.4 9  17.6 9 /
\plot 18.4 9  19.6 9 /
\plot 20.4 9  21.6 9 /
\plot 24.4 9  25.6 9 /
  
\setshadegrid span <0.8mm>
\vshade 14.5 9.5 9.5 <,z,,> 15.5 8.5 10.5 <z,,,> 21.5 2.5 4.5 
<z,,,> 22.5 1.5 4.5 <z,,,> 23 1 5 <z,,,> 24 2 6 /

\plot 2.6 1 11.6 10 /
\plot 13.6 0 23.6 10 /

\betweenarrows{$\Cal D'$} from 2.6 1.0 to 14.6 1.0
\betweenarrows{$\Cal C$} from -1.0 0.3 to 13.9 0.3 

\endpicture}
$$

In $\Cal D'$, the set of indecomposable $\Psi$-modules 
located between the two diagonals in the
diagram, each module $M$ for which one of the vertical 
arrows $\beta_i$ is not represented by a monomorphism lies 
on the sectional path starting at the simple module $S(1')$. 
By deleting this path in $\Cal D'$ 
we obtain $\Cal D=\Cal D'\cap\Cal S(\widetilde 5)$.
It follows from the three claims below that the Auslander-Reiten quiver
is as pictured.

\medskip\noindent
{\bf Claim 1.} {\it
The set $\Cal D$ is the fundamental domain in $\Cal S(\widetilde 5)$ for the shift.
}

\smallskip\noindent
{\it Proof: }\/
If $M\in\ind\Cal S(\widetilde 5)$ 
satisfies that $M_z=0=M_{z'}$ for $z<0$ and $M_0\neq 0$
then either $\beta_0$ is not an epimorphism and there is a path (even a non-zero map)
$M[1]\to P(5)$ in the category $\Cal C$, as seen
in the proof of Claim~1 in (1.4)
or else $\beta_0$ is an isomorphism and there is a non-zero map 
$M\to I(0)=P(4')$ in the category $\Cal C$.
As a consequence, either $M[1]$ or $M$ occurs in $\Cal C$.
Since all modules in $\Cal C\backslash \Cal D'$ have a translate in $\Cal D$
with respect to the shift, $\Cal D$ is a full set of 
representatives for the orbits under the shift
of the indecomposable representations in $\Cal S(\widetilde 5)$.  \qed

\medskip\noindent
{\bf Claim 2.} {\it
Each source map in $\mod \Psi$ starting at a $\Psi$-module in $\Cal D$ is a source map 
in $\Cal S(\widetilde 5)$.
}

\smallskip\noindent
{\it Proof: }\/
Let $X$ be an object in $\Cal D$ and $f:X\to Y$ a source map in $\mod \Psi$.
We show that any non-isomorphism $t:X\to T$ for $T\in\ind \Cal S(\widetilde 5)$
factors over $f$.
Let $j$ be such that $T[j]$ is in $\Cal D$.
If $j>0$ then $X[j]$ and $T[j]$ are both modules over an algebra $\Psi'$ which is
obtained from $\Psi$ by performing one point extensions at radicals of 
projective modules, as in (1.2).  We obtain:
$\Hom_{\Cal S}(X,T)=\Hom_{\Psi'}(X[j],T[j])=0$.
Similarly, if $j\leq 0$ then $X$ and $T$ are both modules over a suitable 
extension algebra $\Psi'$ of $\Psi$; in this case $t$ factors over $f$ since 
$f$ is a source map in the category $\mod \Psi'$. \qed

\medskip \noindent
{\bf Claim 3.} {\it
Let $X$ be an object in $\Cal D$.  Either the Auslander-Reiten
sequence in $\mod \Psi$ starting at $X$ is in $\Cal S(\widetilde 5)$,
or $X$ is the projective-injective module $P(5)$,
or else $X$ is one of the modules $A_i$, $0\leq i\leq 3$, in the above diagram,
and the relative Auslander-Reiten sequence starting at $X$ has the 
form $0\to A_i\to A_{i+1}\oplus C_{i-1}\to C_i\to 0$ where $C_{-1}=0$.
}

\smallskip\noindent
{\it Proof: }\/
See Lemma (1.3.1) for the source map of $P(5)$; 
the four maps labelled ``a'' in the diagram
are left $\Cal S(\widetilde 5)$-approximations,
so each Auslander-Reiten sequence starting at one of the $A_i$ is as claimed.
\qed

We conclude this section with the picture of the dimension pairs, and the corresponding
Krull-Remak-Schmidt multiplicities of the total space, of 
the indecomposable triples in $\Cal S(5)$.
$$
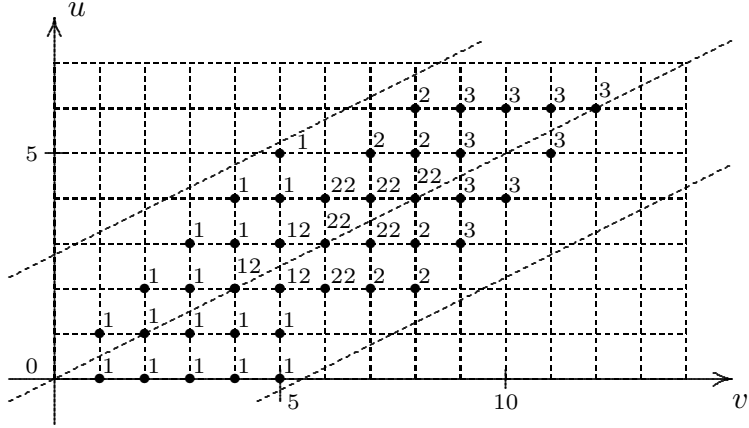
\hbox{\beginpicture
\setcoordinatesystem units <.5cm,.5cm>
\arr{0 -1}{0 8}
\arr{-1 0}{15 0}
\put{$u$} at  0.5 8.2
\put{$v$} at 15.2 -.5
\setdashes <.5mm> 
\setplotarea x from 0 to 14, y from 0 to 7 
\grid {14} {7}
\plot -1 -.5  15 7.5 /
\plot -1 2.25  9.5 7.5  /
\plot 4.5 -.5  15 4.75 /
\setsolid
\plot 5 -0.5  5 0.2 /
\put{$\ssize 5$} at 5.3 -.5
\plot 10 -0.2  10 0.2 /
\put{$\ssize 10$} at 10 -.5
\plot -0.2 5  0.2 5 /
\put{$\ssize 5$} at -.5 5 
\put{$\ssize 0$} at -.5 .3 
\multiput{$\ssize \bullet$} at 1 0  2 0  3 0  4 0  5 0  
                             1 1  2 1  3 1  4 1  5 1
                             2 2  3 2  4 2  5 2  3 3  4 3  5 3  4 4  5 4  5 5  
                             4 2  5 2  6 2  5 3  6 3  6 4  
                             7 2  8 2  7 3  8 3  7 4  8 4  7 5  8 5  8 6 
                             9 3  9 4  9 5  9 6  10 4  10 6  11 5  11 6  12 6 /
\multiput{$\ssize 1$} at 1.2 .3   2.2 .3  3.2 .3   4.2 .3   5.2 .3 
                      1.2 1.3  2.2 1.4  3.2 1.3  4.2 1.3  5.2 1.3  
                    2.2 2.3  3.2 2.3   3.2 3.3  4.2 3.3  4.2 4.3  5.2 4.3  5.5 5.25 /
\multiput{$\ssize 2$} at 7.2 2.3  8.2 2.3  8.2 3.3  
                        7.2 5.3  8.2 5.3  8.2 6.3 /
\multiput{$\ssize 12$} at 5.4 2.3  5.4 3.3  4.3 2.5 /
\multiput{$\ssize 22$} at 6.4 2.3  6.3 3.5  7.4 3.3  6.4 4.3  7.4 4.3  8.3 4.5 /
\multiput{$\ssize 3$} at 9.2 3.3  9.2 4.3  9.2 5.3  9.2 6.3  10.2 4.3  10.2 6.3
                         11.2 5.3  11.2 6.3  12.2 6.4 /
\endpicture}
$$

\medskip\noindent{\bf (6.6) A formula for the number of indecomposables.}
The number $s(n)$ of indecomposable objects in $\Cal S(n)$ is
predicted by the formula $s(n)=2+2(n-1)\frac 6{6-n}$.
If the stable part of $\Gamma(n)$ is 
of type $\Bbb Z\Delta/\varphi$, then the Auslander-Reiten quiver of 
$\Cal S(\widetilde n)$ has stable type $\Bbb Z\Delta$
and we read off from the sequence of boundary modules that within each
interval of six $\tau_{\Cal S}$-translations in $\Bbb Z\Delta$ 
there will occur $6-n$ translates of each boundary module. 
This accounts for the factor $\frac6{6-n}$
in the formula.
In fact for $n<6$ and for every indecomposable non-projective module $M$ in the covering
category $\Cal S(\widetilde n)$, 
the quotient $q_n$ of the $\tau_{\Cal S}$-period, up to the shift, and the
number of $\tau_{\Cal S}$-orbits containing one of the $M[i]$
is equal to $\frac 6{6-n}$ (and independent of $M$).
By checking each case, we see that the Dynkin diagram $\Delta$ 
has $2(n-1)$ points
--- the other factor in the formula. The product of the two factors is the number
$2(n-1)\cdot\frac6{6-n}$ of indecomposables in the stable part of the Auslander-Reiten
quiver for $\Cal S(n)$. 
Two more indecomposables are not represented in
the stable part, namely the projective-injective ones.

\smallskip
In the table below we list some numerical data about the Auslander-Reiten
quivers of the categories $\Cal S(n)$.
By $\sigma$ and $\rho$ we denote automorphisms
of the underlying Dynkin diagram of order 2 and 3, respectively.

\smallskip
\centerline{
\vbox{\offinterlineskip
\halign{\strut#&#\hfil&\vrule#& \quad\hfil#\hfil&\;\hfil#\hfil&\;\hfil#\hfil&
          \;\hfil#\hfil&\;\hfil#\hfil&\;\hfil#\hfil\cr
& $n$  &height14pt& $1$ & $2$ & $3$ & $4$ & $5$ & $6$ \cr
\noalign{\hrule}
& $2(n-1)$  &height14pt& $0$ & $2$ & $4$ & $6$ & $8$ & $10$ \cr
& stable part of $\Gamma(\widetilde n)$ 
        &height14pt& $-$ & $\Bbb ZA_2$ & $\Bbb ZD_4$ & $\Bbb ZE_6$ & $\Bbb ZE_8$ 
  &  ${\beginpicture
\setcoordinatesystem units <0.5cm,0.5cm>
\put{$\ssize E_8$} at 0 0.3
        \circulararc 300 degrees from 0.4 0.3 center at -.05 0.35
        \endpicture}$ \cr
& $\frac 6{6-n}$  &height14pt& $\frac 65$ & $\frac 32$ & $2$ & $3$ & $6$ & $\infty$ \cr
& $q_n$           &height14pt& $-$ & $\frac 32$ & $2$ & $3$ & $6$ & $-$ \cr
& stable part of $\Gamma(n)$ &height14pt& $-$ &
        $\Bbb ZA_2/\tau^{\frac32}\sigma$ & 
        $\Bbb ZD_4/\tau^2\rho$ & 
        $\Bbb ZE_6/\tau^3\sigma$ & 
        $\Bbb ZE_8/\tau^6$ &  ${\beginpicture
\setcoordinatesystem units <0.5cm,0.5cm>
\put{$\ssize E_8$} at 0 0.3
        \circulararc 300 degrees from 0.4 0.3 center at -.05 0.35
        \endpicture}$ \cr
& $2+2(n-1)\frac6{6-n}$ &height14pt&$2$ & $5$ & $10$ & $20$ & $50$ & $\infty$ 
 \cr
& \# indec.\ in $\Gamma(n)$ &height14pt& $2$ & $5$ & $10$ & $20$ & $50$ & $\infty$ \cr
}}}

\medskip
We would like to point out that the formula 
$s(n)=2+2(n-1)\frac 6{6-n}$ predicts the number of indecomposables
correctly for all values $0\leq n\leq 6$.  
In the case where $n=6$ (and where $s(6)$ is predicted correctly as $\infty$),
even the first factor $2(n-1)$ in the formula 
can be interpreted as the number of points of the underlying
diagram.  In fact, we have seen that the category $\Cal S(6)$ arises as a 
coextension of a tame concealed algebra of type $\widetilde E_8$ and in
this sense is associated with $8+1+1=10=2(6-1)$ points.

        \bigskip\bigskip
{\bf 7. Three Remarks about Wild Submodule Categories.}
        \medskip
In the cases where $n\geq 7$, the categories $\Cal S(\widetilde n)$ 
and $\Cal S(n)$
have wild representation type.  In fact, it is shown in [10] that for any 
commutative local uniserial ring $\Lambda$ of Loewy length 7
the category $\Cal S(\Lambda)$ is controlled $k$-wild where $k$ is the radical
factor field of $\Lambda$. 

        \smallskip 
In this section we use the boundary modules and our knowledge about the 
periodicity of the translation to describe the shape of the connected components of the 
Auslander-Reiten quiver for $\Cal S(\widetilde n)$ and $\Cal S(n)$ 
where $n\geq 7$.
Moreover, we construct an indecomposable object in 
$\Cal S(\widetilde n)^{\lfd}$ which does not have finite support.
Thus for $n\geq 7$ the category $\Cal S(\widetilde n)$ 
is not locally support finite.

        \smallskip
We have seen that for $\Cal S(6)$, the region in $\Bbb Z\times \Bbb Z$ containing
the possible dimension pairs of indecomposable objects is contained in a stripe.
We show in (7.3) that this is no longer the case in $\Cal S(n)$ when $n\geq 7$.

\medskip\noindent{\bf (7.1) The connected components of the
Auslander-Reiten quiver.}

\medskip\noindent
{\bf (7.1.1) Proposition.} {\it Let $n\geq 7$.
\item{\rm (a)} Each component of the Auslander-Reiten quiver 
   $\Gamma(\widetilde n)$ for $\Cal S(\widetilde n)$
   has stable part of type $\Bbb ZA_\infty$. All components in 
        $\Gamma(\widetilde n)$ are stable
   with the exception of $n-6$ components which contain boundary modules.
\item{\rm (b)} The stable part of each component of the Auslander-Reiten quiver 
   $\Gamma(n)$ for
   $\Cal S(n)$ is a tube of rank 1, 2, 3, or 6.
   All components are stable with the exception of one tube of rank 6 which
   contains the boundary modules.

}

\smallskip\noindent
{\it Proof: }\/
(b) Let $\Cal C$ be a connected component of $\Gamma(n)$. The length function
on $\Cal C$ gives rise to a subadditive function on the stable part of 
$\Cal C$ which is an unbounded function since $\Cal C$ is not a finite
component. Under the Auslander-Reiten translation, every 
indecomposable module
in $\Gamma(n)$ has a period which is a divisor of 6 ([11], Corollary~6.5),
and hence $\Cal C$ contains periodic modules.
Thus [7], Theorem~2, implies that the diagram of the stable part of 
$\Cal C$
is $\Bbb ZA_\infty/G$ where $G$ is 
an admissible group of automorphisms of $\Bbb ZA_\infty$.
More precisely, the stable
part of $\Cal C$ is a tube of rank a divisor of 6.
Since the boundary modules have $\tau_{\Cal S}$-period 6 whenever $n\geq 3$,
we obtain that 
both projective-injective indecomposables occur in one tube of rank 6.

(a) A connected component $\widetilde{\Cal C}$ of 
$\Gamma(\widetilde n)$ cannot contain any 
periodic modules if $n\geq 7$, since for every object $M$ in such a component
there is an isomorphism $\tau^6_{\Cal S}M\cong M[n-6]$.
Since the component $\Cal C$ of the Auslander-Reiten quiver for $\Cal S(n)$
corresponding to $\widetilde{\Cal C}$ under the covering functor 
has stable part a tube of 
rank a divisor of 6, the category $\widetilde{\Cal C}$ must have stable part 
$\Bbb ZA_\infty$.  In $\Cal S(\widetilde n)$, there are $n-6$ 
$\tau_{\Cal S}$-orbits of boundary modules, they
give rise to the $n-6$ non-stable components in $\Gamma(\widetilde n)$.
\qed

\medskip\noindent{\bf (7.2) Objects in $\Cal S(\widetilde n)$ with large
support.}

\medskip\noindent
{\bf (7.2.1) Proposition.} {\it
For $n\geq 7$, the category $\Cal S(\widetilde n)$ 
is not locally support finite.
}

\smallskip\noindent
{\it Proof: }\/ We construct a finite dimensional indecomposable
representation $N$ and a locally finite dimensional 
indecomposable representation $M$
such that there is a short exact sequence
$$0\longrightarrow M[1]\longrightarrow M\longrightarrow N\longrightarrow 0.$$
We present our construction of $M$ and $N$ for the case $n=7$ and the 
component of $\Gamma(\widetilde 7)$ containing
the boundary modules (but this construction works for any $n\geq 7$ and 
any other component, too).
$$
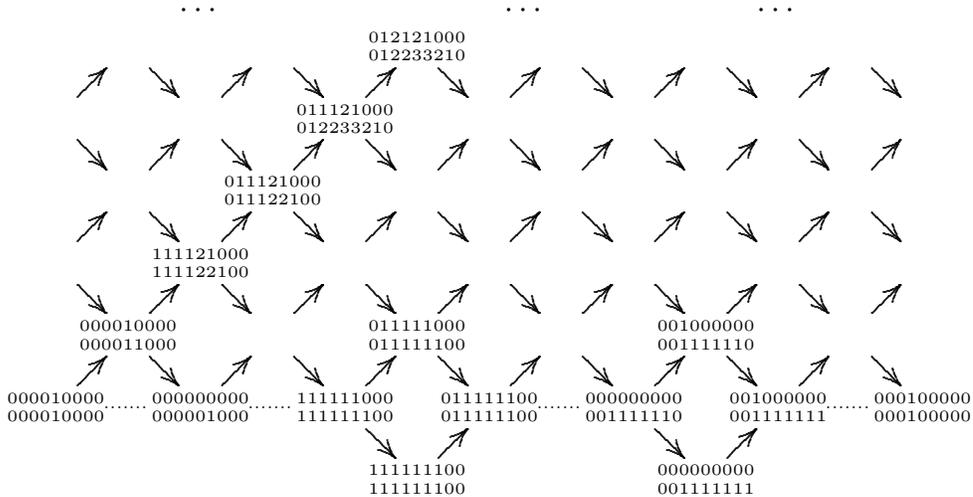
\hbox{\beginpicture
\setcoordinatesystem units <0.8cm,0.8cm>
\put{} at 0 0
\put{} at 0 5

\put{$\ssize {000010000 \atop 000010000}$} at 0 0
\put{$\ssize {000010000 \atop 000011000}$} at 1 1
\put{$\ssize {111121000 \atop 111122100}$} at 2 2
\put{$\ssize {011121000 \atop 011122100}$} at 3 3
\put{$\ssize {011121000 \atop 012233210}$} at 4 4
\put{$\ssize {012121000 \atop 012233210}$} at 5 5
\put{$\ssize {000000000 \atop 000001000}$} at 2 0
\put{$\ssize {111111000 \atop 111111100}$} at 4 0
\put{$\ssize {011111000 \atop 011111100}$} at 5 1
\put{$\ssize {111111100 \atop 111111100}$} at 5 -1
\put{$\ssize {011111100 \atop 011111100}$} at 6 0
\put{$\ssize {000000000 \atop 001111110}$} at 8 0
\put{$\ssize {000000000 \atop 001111111}$} at 9 -1
\put{$\ssize {001000000 \atop 001111110}$} at 9 1
\put{$\ssize {001000000 \atop 001111111}$} at 10 0
\put{$\ssize {000100000 \atop 000100000}$} at 12 0

\arr{4.3 -0.3} {4.7 -0.7} 
\arr{5.3 -0.7} {5.7 -0.3} 
\arr{8.3 -0.3} {8.7 -0.7} 
\arr{9.3 -0.7} {9.7 -0.3} 
\arr{0.3 0.3} {0.7 0.7} 
\arr{1.3 0.7} {1.7 0.3} 
\arr{2.3 0.3} {2.7 0.7} 
\arr{3.3 0.7} {3.7 0.3} 
\arr{4.3 0.3} {4.7 0.7} 
\arr{5.3 0.7} {5.7 0.3} 
\arr{6.3 0.3} {6.7 0.7} 
\arr{7.3 0.7} {7.7 0.3} 
\arr{8.3 0.3} {8.7 0.7} 
\arr{9.3 0.7} {9.7 0.3} 
\arr{10.3 0.3} {10.7 0.7} 
\arr{11.3 0.7} {11.7 0.3} 
\arr{0.3 1.7} {0.7 1.3} 
\arr{1.3 1.3} {1.7 1.7} 
\arr{2.3 1.7} {2.7 1.3} 
\arr{3.3 1.3} {3.7 1.7} 
\arr{4.3 1.7} {4.7 1.3} 
\arr{5.3 1.3} {5.7 1.7} 
\arr{6.3 1.7} {6.7 1.3} 
\arr{7.3 1.3} {7.7 1.7} 
\arr{8.3 1.7} {8.7 1.3} 
\arr{9.3 1.3} {9.7 1.7} 
\arr{10.3 1.7} {10.7 1.3}
\arr{11.3 1.3} {11.7 1.7} 
\arr{0.3 2.3} {0.7 2.7} 
\arr{1.3 2.7} {1.7 2.3} 
\arr{2.3 2.3} {2.7 2.7} 
\arr{3.3 2.7} {3.7 2.3} 
\arr{4.3 2.3} {4.7 2.7} 
\arr{5.3 2.7} {5.7 2.3} 
\arr{6.3 2.3} {6.7 2.7} 
\arr{7.3 2.7} {7.7 2.3} 
\arr{8.3 2.3} {8.7 2.7} 
\arr{9.3 2.7} {9.7 2.3} 
\arr{10.3 2.3} {10.7 2.7} 
\arr{11.3 2.7} {11.7 2.3} 
\arr{0.3 3.7} {0.7 3.3} 
\arr{1.3 3.3} {1.7 3.7} 
\arr{2.3 3.7} {2.7 3.3} 
\arr{3.3 3.3} {3.7 3.7} 
\arr{4.3 3.7} {4.7 3.3} 
\arr{5.3 3.3} {5.7 3.7} 
\arr{6.3 3.7} {6.7 3.3} 
\arr{7.3 3.3} {7.7 3.7} 
\arr{8.3 3.7} {8.7 3.3} 
\arr{9.3 3.3} {9.7 3.7} 
\arr{10.3 3.7} {10.7 3.3}
\arr{11.3 3.3} {11.7 3.7} 
\arr{0.3 4.3} {0.7 4.7} 
\arr{1.3 4.7} {1.7 4.3} 
\arr{2.3 4.3} {2.7 4.7} 
\arr{3.3 4.7} {3.7 4.3} 
\arr{4.3 4.3} {4.7 4.7} 
\arr{5.3 4.7} {5.7 4.3} 
\arr{6.3 4.3} {6.7 4.7} 
\arr{7.3 4.7} {7.7 4.3} 
\arr{8.3 4.3} {8.7 4.7} 
\arr{9.3 4.7} {9.7 4.3} 
\arr{10.3 4.3} {10.7 4.7} 
\arr{11.3 4.7} {11.7 4.3} 
\setdots<2pt>
\plot 0.7 0  1.3 0 /
\plot 2.7 0  3.3 0 /
\plot 6.7 0  7.3 0 /
\plot 10.7 0  11.3 0 /
\put{$\cdots$} at 2 5.5 
\put{$\cdots$} at 6.5 5.5
\put{$\cdots$} at 10 5.5 
\endpicture}
$$
We start with the module $N$ with dimension vector
$\ssize {012121000 \atop 012233210}$ and with a sectional path
$K(4) \to N$ formed by irreducible maps.
Note that $\tau^{6}N$ is just the module $N[1]$, and 
the Auslander-Reiten component exhibits a non-split extension
$0 \to N[1] \to N^2[1] \to N \to 0.$  Here we sketch a part of the
component:
$$
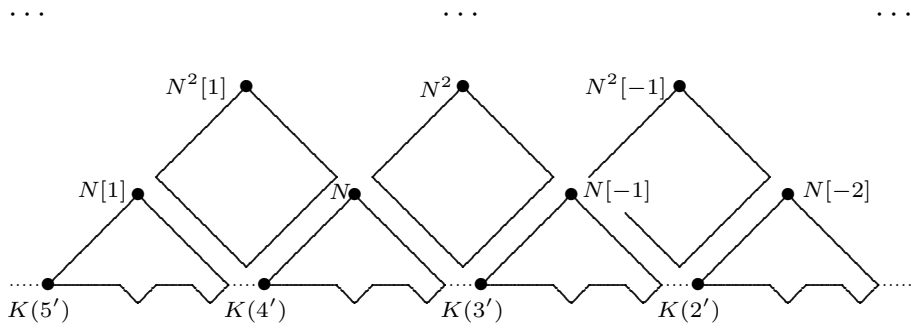
\hbox{\beginpicture
\setcoordinatesystem units <0.2cm,0.2cm>
\put{} at -13 0
\put{} at 35 15
\plot 0 0  4 0  5 -1  6 0  8 0  9 -1  10 0  5 5  0 0 /
\plot 12 0  16 0  17 -1  18 0  20 0  21 -1  22 0  17 5  12 0 /
\plot 24 0  28 0  29 -1  30 0  32 0  33 -1  34 0  29 5  24 0 /
\plot -12 0  -8 0  -7 -1  -6 0  -4 0  -3 -1  -2 0  -7 5  -12 0 /
\plot 6 6  11 1  16 6  11 11  6 6 /
\plot -6 6  -1 1  4 6  -1 11  -6 6 /
\plot 20 4  23 1  28 6  23 11  18 6 /
\put{$\bullet$} at 0 0 
\put{$\bullet$} at 5 5
\put{$\bullet$} at 11 11 
\put{$\bullet$} at 12 0
\put{$\bullet$} at -12 0
\put{$\bullet$} at -7 5
\put{$\bullet$} at 24 0
\put{$\bullet$} at 29 5
\put{$\bullet$} at 17 5
\put{$\bullet$} at -1 11
\put{$\bullet$} at 23 11
 
\put{$\ssize K(4')$} at -0.5 -1.4
\put{$\ssize K(3')$} at 11.5 -1.4
\put{$\ssize K(5')$} at -12.5 -1.4
\put{$\ssize K(2')$} at 23.5 -1.4
\put{$\ssize N$} at 4.2 5.2
\put{$\ssize N[1]$} at -9 5.2
\put{$\ssize N[-1]$} at 19.5 5.2
\put{$\ssize N[-2]$} at 31.7 5.2
\put{$\ssize N^2$} at 9.5 11
\put{$\ssize N^2[1]$} at -3.8 11
\put{$\ssize N^2[-1]$} at 20 11
\setdots<2pt>
\plot -14 0  -12 0 /
\plot -2 0  0 0 /
\plot 10 0 12 0 /
\plot 22 0  24 0 /
\plot 34 0  36 0 /
\put{$\cdots$} at 11 15
\put{$\cdots$} at -13 15
\put{$\cdots$} at 35 15
\endpicture}
$$
Similarly we obtain for each $m\in \Bbb N$ an indecomposable module
$N^{m+1}[m]$ and a  non-split extension
$$0\longrightarrow N[m]\longrightarrow N^{m+1}[m]\lto{\phi_m}N^m[m-1]
\longrightarrow 0$$
and define $(M,\psi_m:M\to N^m[m-1])$ to be the limit of the inverse system
$(\phi_i)_{i\in\Bbb N}$.
The module $M$ is locally finite dimensional with dimension vector
$\smallmatrix 
0\, 1\, 3\, 4\, 6\,\phantom{0}7\,\phantom{0}7\,\phantom{0}7\,\phantom{0}7\, \cdots\cr
0\, 1\, 3\, 5\, 8\,11\,13\,14\,14\, \cdots \endsmallmatrix$, 
and each map $\psi_m$ is the canonical map modulo $M[m]$. 
Moreover, the mapping property of the inverse limit yields
a map $M[1]\to M$ which has cokernel $N$. 

We claim that $M$ is indecomposable in $\Cal S(\widetilde 7)^{\lfd}$. 
Let $M=U\oplus V$ be a direct sum decomposition in $\Cal S(\widetilde 7)^{\lfd}$. We assume 
that $U$ and $V$ are both non-zero and fix $\ell \in\Bbb N$ such that
$U_\ell\neq 0$ and $V_\ell\neq 0$. 
Since $\psi_\ell(M)_\ell=M_\ell$ we obtain a proper decomposition
$\psi_\ell(M)=\psi_\ell(U)\oplus\psi_\ell(V)$, in contradiction to the
indecomposability of $\psi_\ell(M)=N^\ell[\ell-1]$.      \qed

        \medskip\noindent
{\bf (7.3) Dimension pairs in $\Cal S(n)$ for $n\geq 7$.}

\medskip\noindent
{\bf (7.3.1) Proposition.} {\it
For $n\geq 7$, the dimension pairs of indecomposable triples 
in the category $\Cal S(n)$ are not contained in any stripe.
}

        \smallskip\noindent
{\it Proof:} We specify two families of indecomposable triples which 
have dimension pairs that lie on two lines which have different slope.
More precisely, for each $f>0$ we construct two modules $\bar M^+$ and 
$\bar M^-$ in $\Cal S(7)$ such that $\bdim \bar M^+=f\cdot(13,6)$ 
and $\bdim \bar M^-=f\cdot(13,7)$, 
so $\bdim \bar M^+$ and $\bdim \bar M^-$ have vertical distance $f$.
This shows that there is no stripe which contains all dimension pairs.

        \smallskip
Given $f>0$, let $M$ be an $f$-fold selfextension of a module on the mouth
of one of the homogeneous tubes in $\Cal U$ in $\Cal S(\widetilde 6)$,
say corresponding to parameter $c\in k\backslash\{0,1\}$,  
so $M$ has dimension vector $\bdim M=f\cdot \bh^\infty$. 
Let $M^+$ be the representation in $\Cal S(\widetilde 7)$ obtained from $M$
by replacing $M_7=0$ by $M^+_7=M_6$ and by taking for $M^+_{\alpha_7}$ the
identity map on $M_6$. 
Similarly, let $M^-$ be obtained from $M$ by replacing $M_0=0=M_{0'}$
by $M^-_0=M_1=M_{1'}=M^-_{0'}$ and by taking as maps 
$M^-_{\beta_0}=M^-_{\alpha_1}=M^-_{\alpha_1'}=1_{M_1}$.
Then $M^+$ and $M^-$ have dimension vectors
$$\bdim M^+=f\cdot {\ssize {01221\phantom{000} \atop 01233211}}\quad\text{and}\quad
  \bdim M^-=f\cdot {\ssize {11221\phantom{000}\atop 11233210}}, $$
respectively.  
Both modules are indecomposable since a decomposition would give rise to a 
decomposition of $M$.

        \smallskip
For $\bar M^+$ and $\bar M^-$ take the image of $M^+$ and $M^-$ under the 
covering functor, this yields indecomposable objects by Proposition (2.1.1).
In case $f=2$, they can be pictured as follows 
(for the notation see Section 2.3):

{
\def\sixfourtwo{\multiput{} at 0 -2  0 7 /
        \multiput{\sq} at 0 0  0 1  0 2  0 3  0 4  0 5 
                1 1  1 2  1 3  1 4    2 2  2 3 /
        \multiput{$\bullet$} at .5 2  1.5 2 / 
        \plot .5 2  1.5 2 / }
\def\sixfourtwoplus{\multiput{} at 0 -2  0 7 /
        \multiput{\sq} at 0 0  0 1  0 2  0 3  0 4  0 5  0 6
                1 1  1 2  1 3  1 4    2 2  2 3 /
        \multiput{$\bullet$} at .5 2  1.5 2 / 
        \plot .5 2  1.5 2 / }
\def\sixfourtwominus{\multiput{} at 0 -2  0 7 /
        \multiput{\sq} at 0 -1  0 0  0 1  0 2  0 3  0 4  0 5  
                1 1  1 2  1 3  1 4    2 2  2 3 /
        \multiput{$\bullet$} at .5 2  1.5 2 / 
        \plot .5 2  1.5 2 / }
\def\lambdaone{\multiput{} at 0 -2  0 7 /
        \multiput{$\bullet$} at 2.5 3 /
        \put{$\sssize c$} at .5 3
        \put{$\sssize c'$} at 1.5 3.1 
        \plot .8 3  1.1 3 /
        \plot 1.8 3  2.5 3 / }
\def\lambdaonelow{\multiput{} at 0 -2  0 7 /
        \multiput{$\bullet$} at 2.3 2.8 /
        \put{$\sssize c$} at .5 2.8
        \put{$\sssize c'$} at 1.5 2.9
        \plot .8 2.8  1.1 2.8 /
        \plot 1.8 2.8  2.3 2.8 / }
\def\lambdatwohigh{\multiput{} at 0 -2  0 7 /
        \multiput{$\bullet$} at -.3 3.15  2.3 3.15 /
        \put{$\sssize c$} at .5 3.15
        \put{$\sssize c'$} at 1.5 3.15 
        \plot -.3 3.15  .2 3.15 /
        \plot .8 3.15  1.1 3.15 /
        \plot 1.8 3.15  2.3 3.15 / }
\def\lambdatwolow{\multiput{} at 0 -2  0 7 /
        \multiput{$\bullet$} at -.3 2.8  2.3 2.8 /
        \put{$\sssize c$} at .5 2.8
        \put{$\sssize c'$} at 1.5 2.9
        \plot -.3 2.8  .2 2.8 /
        \plot .8 2.8  1.1 2.8 /
        \plot 1.8 2.8  2.3 2.8 / }

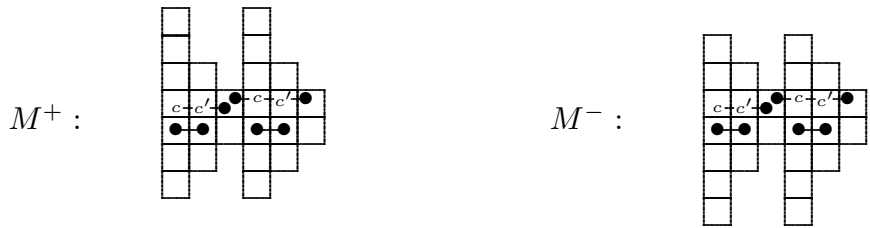
\def\scale#1{\beginpicture\setcoordinatesystem units <3mm,3mm>#1
        \endpicture}
%
%
$$\hbox{\beginpicture\setcoordinatesystem units <1cm,1cm>
        \put{$M^+:$} at 1 0
        \put{\scale{\multiput{\sixfourtwoplus} at 0 0  3 0 /
                    \put{\lambdaonelow} at 0 0 
                    \put{\lambdatwohigh} at 3 0 }} at 2.75 0 
        \put{$M^-:$} at 7 0
        \put{\scale{\multiput{\sixfourtwominus} at 0 0  3 0 /
                    \put{\lambdaonelow} at 0 0 
                    \put{\lambdatwohigh} at 3 0 }} at 8.75 0 
        \endpicture}$$
}
The dimension pairs of $\bar M^+$ and $\bar M^-$ are as stated above. \qed

        \bigskip\bigskip
{\bf Appendix. Operators which are not nilpotent.} 
        \medskip
We conclude with a remark on invariant subspaces of operators which
are not necessarily nilpotent. In this section we deal with the more general case 
where the linear operator does not act nilpotently on the vector space.  
We assume that the minimal polynomial of the
operator is a product of linear factors and demonstrate that the classification
of the invariant subspaces can be reduced to the nilpotent case.

\smallskip Let $(V,U,T)$ be a triple consisting
of a finite dimensional $k$-vector space $V$, a subspace $U$ of $V$
and a linear map $T\:V\to V$ which maps $U$ into $U$.

\smallskip 
Suppose that the minimal polynomial $\min_T(x)$ for $T$ factors as
$$\min_T(x)=p_1^{n_1}(x)\cdots p_s^{n_s}(x)$$
where the $p_i$ are pairwise distinct monic irreducible polynomials.
Then the vector space $V$ has primary decomposition
$$V=V_{p_1}\oplus\cdots\oplus V_{p_s}\T{where the}
V_{p_i}\;=\;\big\{v\in V\;\big|\; p_i^{n_i}(T)(v)=0\big\}$$
are invariant subspaces of $V$ with minimal polynomials
$\min_{T|{V_{p_i}}} (x)=p_i^{n_i}(x)$.

\smallskip
Since $U$ is an invariant subspace of $V$, the minimal polynomial 
for $T|U$ is a divisor of the minimal polynomial for $T$ and hence
$U$ has a corresponding primary decomposition
$U=U_{p_1}\oplus\cdots\oplus U_{p_s}$. Thus,
$$(V,U,T) \;=\; (V_{p_1}, U_{p_1}, T|{V_{p_1}}) \;\oplus\;\cdots
        \;\oplus \;(V_{p_s}, U_{p_s}, T|{V_{p_s}})$$
is a decomposition where each operator $T|V_{p_i}$ has a minimal
polynomial which is a power of an irreducible polynomial. 

\smallskip
If $p_i(x)$ is a linear polynomial, then the operator
$S_i=p_i(T)$ is a linear operator which acts nilpotently on $V_{p_i}$ 
and which fixes the subspace $U_{p_i}$.

\smallskip 
In conclusion, we have reduced the decomposition problem for 
the triple $(V,U,T)$
where $T$ is an arbitrary $k$-linear operator
whose minimal polynomial is a product of linear factors,
 to the corresponding problem
for the triples $(V_{p_i},U_{p_i},S_i)$ where $S_i$ is a nilpotent
linear operator.

        \bigskip\bigskip
\frenchspacing
\noindent
{\bf References}        {\baselineskip=9pt \rmk
\parindent=1.5truecm
                                         \medskip\smallskip\noindent
 \item{[1]} 
        {\itk M.\ Auslander, I\. Reiten, S.~O.~Smal\o,}
        Representation Theory of Artin Algebras,
        Cambridge University Press (1995).
\item{[2]} 
        {\itk G.~Birkhoff,}
        Subgroups of abelian groups,
        Proc.~Lond.~Math.~Soc., II.~Ser.~{\bfk 38} (1934), 385--401.
\item{[3]} 
        {\itk K.~Bongartz and P.~Gabriel,} 
        Covering spaces in representation theory,
        Inventiones Math.\ {\bfk 65} (1982), 331--378.
\item{[4]} 
        {\itk N.~Bourbaki,} Groupes et alg\`ebres de Lie, Chapitres 4,5 et 6,
        Masson, Paris, 1981.
\item{[5]} 
        {\itk P.\ Dowbor, H.~Lenzing, A.~Skowronski,} Galois coverings
        by algebras of locally support-finite categories, in:
        Representation theory I, Proc.\ Conf.\ Ottawa (1984), 91--93.
\item{[6]} 
        {\itk P.~Gabriel,}
        The universal cover of a representation-finite algebra.
        in: Representations of Algebras, 
        Springer LNM {\bfk 903} (1981), 68--105.
\item{[7]} 
        {\itk D.~Happel, U.~Preiser, and C.~M.~Ringel,}
        Vinberg's characterization of Dynkin diagrams using subadditive
        functions with application to {\itk D}{\rmk Tr}-periodic modules. 
        in: Representation Theory II, Proc.\ Conf.\ Ottawa 1979, 
        Springer LNM {\bfk 832} (1980), 280--294.
\item{[8]} 
        {\itk F.~Richman and E.~A.~Walker,}
        Subgroups of p${}^5$-bounded groups,
        in: Abelian groups and modules, Trends Math., Birkh\"auser, Basel, (1999),
        55--73.
\item{[9]} 
        {\itk C.~M.~Ringel,} Tame Algebras and Integral Quadratic Forms,
        Springer LNM {\bfk 1099} (1984).
\item{[10]} 
        {\itk C.~M.~Ringel and M.~Schmidmeier,}
        Submodule categories of wild representation type,
        J.~Pure and Applied Alg.~{\bfk 205} (2006), 412-422.
\item{[11]} 
        {\itk C.~M.~Ringel and M.~Schmidmeier,}
        The Auslander-Reiten translation in
        submodule categories, to appear in: Transactions of the AMS. 
\item{[12]} 
        {\itk C.~M.~Ringel and M.~Schmidmeier,}
        Invariant subspaces of nilpotent linear operators. II. In preparation.
\item{[13]} 
        {\itk M.~Schmidmeier,} Bounded submodules of modules,
         J.~Pure and Applied Alg.~{\bfk 203} (2005), 45--82.
\item{[14]} 
        {\itk D.~Simson,} Chain categories of modules and 
        subprojective representations 
        of posets over uniserial algebras,
        Rocky Mountain J.~Math.~{\bfk 32} (2002), 1627--1650.  

\par}

\bigskip\noindent
{\rmk Claus Michael Ringel,
Fakult\"at f\"ur Mathematik, Universit\"at Bielefeld,
\par\noindent  POBox 100\,131, \ D-33\,501 Bielefeld 
\par\noindent {\ttk ringel\@math.uni-bielefeld.de}
        \medskip\noindent
Markus Schmidmeier, 
Department of Mathematical Sciences, Florida Atlantic University,
\par\noindent Boca Raton, Florida 33431-0991
\par\noindent {\ttk markus\@math.fau.edu}
}
\bye